# INTRODUCTION TO LINEAR BIALGEBRA

W. B. Vasantha Kandasamy

Florentin Smarandache

K. Ilanthenral

# INTRODUCTION
# TO LINEAR BIALGEBRA

**W. B. Vasantha Kandasamy**
Department of Mathematics
Indian Institute of Technology, Madras
Chennai – 600036, India
e-mail: vasantha@iitm.ac.in
web: http://mat.iitm.ac.in/~wbv

**Florentin Smarandache**
Department of Mathematics
University of New Mexico
Gallup, NM 87301, USA
e-mail: smarand@gallup.unm.edu

**K. Ilanthenral**
Editor, Maths Tiger, Quarterly Journal
Flat No.11, Mayura Park,
16, Kazhikundram Main Road, Tharamani,
Chennai – 600 113, India
e-mail: ilanthenral@gmail.com

# CONTENTS

## Preface

The algebraic structure, linear algebra happens to be one of the subjects which yields itself to applications to several fields like coding or communication theory, Markov chains, representation of groups and graphs, Leontief economic models and so on. This book has for the first time, introduced a new algebraic structure called linear bialgebra, which is also a very powerful algebraic tool that can yield itself to applications.

With the recent introduction of bimatrices (2005) we have ventured in this book to introduce new concepts like linear bialgebra and Smarandache neutrosophic linear bialgebra and also give the applications of these algebraic structures.

It is important to mention here it is a matter of simple exercise to extend these to linear n-algebra for any n greater than 2; for n = 2 we get the linear bialgebra.

This book has five chapters. In the first chapter we just introduce some basic notions of linear algebra and S-linear algebra and their applications. Chapter two introduces some new algebraic bistructures. In chapter three we introduce the notion of linear bialgebra and discuss several interesting properties about them. Also, application of linear bialgebra to bicodes is given. A remarkable part of our research in this book is the introduction of the notion of birepresentation of bigroups.

The fourth chapter introduces several neutrosophic algebraic structures since they help in defining the new concept of neutrosophic linear bialgebra, neutrosophic bivector spaces, Smarandache neutrosophic linear bialgebra and Smarandache neutrosophic bivector spaces. Their probable applications to real-world models are discussed. We have aimed to make this book engrossing and illustrative and supplemented it with nearly 150 examples. The final chapter gives 114 problems which will be a boon for the reader to understand and develop the subject.

The main purpose of this book is to familiarize the reader with the applications of linear bialgebra to real-world problems.

Finally, we express our heart-felt thanks to Dr.K.Kandasamy whose assistance and encouragement in every manner made this book possible.

W.B.VASANTHA KANDASAMY
FLORENTIN SMARANDACHE
K. ILANTHENRAL

## Chapter One

# INTRODUCTION TO LINEAR ALGEBRA AND S-LINEAR ALGEBRA

In this chapter we just give a brief introduction to linear algebra and S-linear algebra and its applications. This chapter has three sections. In section one; we just recall the basic definition of linear algebra and some of the important theorems. In section two we give the definition of S-linear algebra and some of its basic properties. Section three gives a few applications of linear algebra and S-linear algebra.

## 1.1 Basic properties of linear algebra

In this section we give the definition of linear algebra and just state the few important theorems like Cayley Hamilton theorem, Cyclic Decomposition Theorem, Generalized Cayley Hamilton Theorem and give some properties about linear algebra.

**DEFINITION 1.1.1:** *A vector space or a linear space consists of the following:*

*i.* *a field F of scalars.*
*ii.* *a set V of objects called vectors.*

*iii.* *a rule (or operation) called vector addition; which associates with each pair of vectors $\alpha, \beta \in V$; $\alpha + \beta$ in V, called the sum of $\alpha$ and $\beta$ in such a way that*

*a.* *addition is commutative $\alpha + \beta = \beta + \alpha$.*

*b.* *addition is associative $\alpha + (\beta + \gamma) = (\alpha + \beta) + \gamma$.*

*c.* *there is a unique vector 0 in V, called the zero vector, such that*

$$\alpha + 0 = \alpha$$

*for all $\alpha$ in V.*

*d.* *for each vector $\alpha$ in V there is a unique vector $-\alpha$ in V such that*

$$\alpha + (-\alpha) = 0.$$

*e.* *a rule (or operation), called scalar multiplication, which associates with each scalar c in F and a vector $\alpha$ in V a vector $c \bullet \alpha$ in V, called the product of c and $\alpha$, in such a way that*

1. $1 \bullet \alpha = \alpha$ *for every* $\alpha$ *in V.*
2. $(c_1 \bullet c_2) \bullet \alpha = c_1 \bullet (c_2 \bullet \alpha)$.
3. $c \bullet (\alpha + \beta) = c \bullet \alpha + c \bullet \beta$.
4. $(c_1 + c_2) \bullet \alpha = c_1 \bullet \alpha + c_2 \bullet \alpha$ .

*for $\alpha, \beta \in V$ and $c, c_1 \in F$.*

It is important to note as the definition states that a vector space is a composite object consisting of a field, a set of 'vectors' and two operations with certain special properties. V is a linear algebra if V has a multiplicative closed binary operation '.' which is associative; i.e., if $v_1, v_2 \in V$, $v_1.v_2 \in V$. The same set of vectors may be part of a number of distinct vectors.

We simply by default of notation just say V a vector space over the field F and call elements of V as vectors only as matter of convenience for the vectors in V may not bear much resemblance to any pre-assigned concept of vector, which the reader has.

**THEOREM (CAYLEY HAMILTON):** *Let T be a linear operator on a finite dimensional vector space V. If f is the characteristic polynomial for T, then f(T) = 0, in other words the minimal polynomial divides the characteristic polynomial for T.*

**THEOREM: (CYCLIC DECOMPOSITION THEOREM):** *Let T be a linear operator on a finite dimensional vector space V and let $W_0$ be a proper T-admissible subspace of V. There exist non-zero vectors $\alpha_1, \ldots, \alpha_r$ in V with respective T-annihilators $p_1, \ldots, p_r$ such that*

i. $V = W_0 \oplus Z(\alpha_1; T) \oplus \ldots \oplus Z(\alpha_r; T)$.
ii. $p_t$ *divides* $p_{t-1}$, $t = 2, \ldots, r$.

*Further more the integer r and the annihilators $p_1, \ldots, p_r$ are uniquely determined by*
*(i) and (ii) and the fact that $\alpha_t$ is 0.*

**THEOREM (GENERALIZED CAYLEY HAMILTON THEOREM):** *Let T be a linear operator on a finite dimensional vector space V. Let p and f be the minimal and characteristic polynomials for T, respectively*

i. *p divides f.*
ii. *p and f have the same prime factors except the multiplicities.*
iii. *If* $p = f_1^{\alpha_1} \ldots f_t^{\alpha_t}$ *is the prime factorization of p, then* $f = f_1^{d_1} f_2^{d_2} \ldots f_t^{d_t}$ *where $d_i$ is the nullity of $f_i(T)^{\alpha}$ divided by the degree of $f_i$.*

The following results are direct and left for the reader to prove.

Here we take vector spaces only over reals i.e., real numbers. We are not interested in the study of these properties in case of complex fields. Here we recall the concepts of linear functionals, adjoint, unitary operators and normal operators.

**DEFINITION 1.1.2:** *Let F be a field of reals and V be a vector space over F. An inner product on V is a function which assigns to each ordered pair of vectors $\alpha$ , $\beta$ in V a scalar $\langle \alpha / \beta \rangle$ in F in such a way that for all $\alpha$, $\beta$, $\gamma$ in V and for all scalars c.*

i. $\langle \alpha + \beta \mid \gamma \rangle = \langle \alpha \mid \gamma \rangle + \langle \beta \mid \gamma \rangle$.
ii. $\langle c\, \alpha \mid \beta \rangle = c\langle \alpha \mid \beta \rangle$.
iii. $\langle \beta \mid \alpha \rangle = \langle \alpha \mid \beta \rangle$.
iv. $\langle \alpha \mid \alpha \rangle > 0$ *if* $\alpha \neq 0$.
v. $\langle \alpha \mid c\beta + \gamma \rangle = c\langle \alpha \mid \beta \rangle + \langle \alpha \mid \gamma \rangle$.

*Let $Q^n$ or $F^n$ be a n dimensional vector space over Q or F respectively for $\alpha, \beta \in Q^n$ or $F^n$ where*

$$\alpha = \langle \alpha_1, \alpha_2, \ldots, \alpha_n \rangle \text{ and}$$
$$\beta = \langle \beta_1, \beta_2, \ldots, \beta_n \rangle$$
$$\langle \alpha \mid \beta \rangle = \sum_j \alpha_j \beta_j .$$

*Note:* We denote the positive square root of $\langle \alpha \mid \alpha \rangle$ by $||\alpha||$ and $||\alpha||$ is called the norm of $\alpha$ with respect to the inner product $\langle\ \rangle$.

We just recall the notion of quadratic form.

The quadratic form determined by the inner product is the function that assigns to each vector $\alpha$ the scalar $||\alpha||^2$.

Thus we call an inner product space is a real vector space together with a specified inner product in that space. A finite dimensional real inner product space is often called a Euclidean space.

The following result is straight forward and hence the proof is left for the reader.

**Result 1.1.1:** If V is an inner product space, then for any vectors $\alpha$, $\beta$ in V and any scalar c.

i. $||c\alpha|| = |c|\ ||\alpha||$.
ii. $||\alpha|| > 0$ for $\alpha \neq 0$.
iii. $|\langle \alpha \mid \beta \rangle| \leq ||\alpha||\ ||\beta||$.
iv. $||\alpha + \beta|| \leq ||\alpha|| + ||\beta||$.

Let $\alpha$ and $\beta$ be vectors in an inner product space V. Then $\alpha$ is orthogonal to $\beta$ if $\langle \alpha \mid \beta \rangle = 0$ since this implies $\beta$ is orthogonal to $\alpha$, we often simply say that $\alpha$ and $\beta$ are orthogonal. If S is a set of vectors in V, S is called an orthogonal set provided all pair of distinct vectors in S are orthogonal. An orthogonal set S is an orthonormal set if it satisfies the additional property $||\alpha|| = 1$ for every $\alpha$ in S.

**Result 1.1.2:** An orthogonal set of non-zero vectors is linearly independent.

**Result 1.1.3:** If $\alpha$ and $\beta$ is a linear combination of an orthogonal sequence of non-zero vectors $\alpha_1, \ldots, \alpha_m$ then $\beta$ is the particular linear combinations

$$\beta = \sum_{t=1}^{m} \frac{\langle \beta \mid \alpha_t \rangle}{|| \alpha_t ||^2} \alpha_t \ .$$

**Result 1.1.4:** Let V be an inner product space and let $\beta_1, \ldots, \beta_n$ be any independent vectors in V. Then one may construct orthogonal vectors $\alpha_1, \ldots, \alpha_n$ in V such that for each t = 1, 2, …, n the set $\{\alpha_1, \ldots, \alpha_t\}$ is a basis for the subspace spanned by $\beta_1, \ldots, \beta_t$.

This result is known as the Gram-Schmidt orthgonalization process.

**Result 1.1.5:** Every finite dimensional inner product space has an orthogonal basis.

One of the nice applications is the concept of a best approximation. A best approximation to β by vector in W is a vector α in W such that

$$||\beta - \alpha|| \leq ||\beta - \gamma||$$

for every vector γ in W.

The following is an important concept relating to the best approximation.

**THEOREM 1.1.1:** *Let W be a subspace of an inner product space V and let β be a vector in V.*

i. *The vector α in W is a best approximation to β by vectors in W if and only if β – α is orthogonal to every vector in W.*

ii. *If a best approximation to β by vectors in W exists it is unique.*

iii. *If W is finite dimensional and {$\alpha_1, \ldots, \alpha_t$} is any orthogonal basis for W, then the vector*

$$\alpha = \sum_t \frac{(\beta \mid \alpha_t)\alpha_t}{||\alpha_t||^2}$$

*is the unique best approximation to β by vectors in W.*

Let V be an inner product space and S any set of vectors in V. The orthogonal complement of S is that set $S^{\perp}$ of all vectors in V which are orthogonal to every vector in S.

Whenever the vector α exists it is called the orthogonal projection of β on W. If every vector in V has orthogonal projection on W, the mapping that assigns to each vector in V its orthogonal projection on W is called the orthogonal projection of V on W.

**Result 1.1.6:** Let V be an inner product space, W is a finite dimensional subspace and E the orthogonal projection of V on W.

Then the mapping

$$\beta \to \beta - E\beta$$

is the orthogonal projection of V on $W^{\perp}$.

**Result 1.1.7:** Let W be a finite dimensional subspace of an inner product space V and let E be the orthogonal projection of V on W. Then E is an idempotent linear transformation of V onto W, $W^{\perp}$ is the null space of E and $V = W \oplus W^{\perp}$. Further $1 - E$ is the orthogonal projection of V on $W^{\perp}$. It is an idempotent linear transformation of V onto $W^{\perp}$ with null space W.

**Result 1.1.8:** Let $\{\alpha_1, \ldots, \alpha_t\}$ be an orthogonal set of non-zero vectors in an inner product space V.

If $\beta$ is any vector in V, then

$$\sum_t \frac{|(\beta, \alpha_t)|^2}{\|\alpha_t\|^2} \leq \|\beta\|^2$$

and equality holds if and only if

$$\beta = \sum_t \frac{(\beta \mid \alpha_t)}{\|\alpha_t\|^2} \alpha_t .$$

Now we prove the existence of adjoint of a linear operator T on V, this being a linear operator $T^*$ such that $(T\alpha \mid \beta) = (\alpha \mid T^*\beta)$ for all $\alpha$ and $\beta$ in V.

We just recall some of the essential results in this direction.

**Result 1.1.9:** Let V be a finite dimensional inner product space and f a linear functional on V. Then there exists a unique vector $\beta$ in V such that $f(\alpha) = (\alpha \mid \beta)$ for all $\alpha$ in V.

**Result 1.1.10:** For any linear operator T on a finite dimensional inner product space V there exists a unique linear operator $T^*$ on V such that

$$(T\alpha \mid \beta) = (\alpha \mid T^*\beta)$$

for all $\alpha$, $\beta$ in V.

**Result 1.1.11:** Let V be a finite dimensional inner product space and let $B = \{\alpha_1, \ldots, \alpha_n\}$ be an (ordered) orthonormal basis for V. Let T be a linear operator on V and let A be the matrix of T in the ordered basis B. Then

$$A_{ij} = (T\alpha_j \mid \alpha_i).$$

Now we define adjoint of T on V.

**DEFINITION 1.1.3:** *Let T be a linear operator on an inner product space V. Then we say that T has an adjoint on V if there exists a linear operator $T^*$ on V such that*

$$(T\alpha \mid \beta) = (\alpha \mid T^*\beta)$$

*for all $\alpha$, $\beta$ in V.*

It is important to note that the adjoint of T depends not only on T but on the inner product as well.

The nature of $T^*$ is depicted by the following result.

**THEOREM 1.1.2:** *Let V be a finite dimensional inner product real vector space. If T and U are linear operators on V and c is a scalar*

*i.* $(T + U)^* = T^* + U^*$.

*ii.* $(cT)^* = cT^*$.
*iii.* $(TU)^* = U^*T^*$.
*iv.* $(T^*)^* = T$.

*A linear operator T such that* $T = T^*$ *is called self adjoint or Hermitian.*

Results relating the orthogonal basis is left for the reader to explore.

*Let V be a finite dimensional inner product space and T a linear operator on V. We say that T is normal if it commutes with its adjoint i.e.* $TT^* = T^*T$.

## 1.2 Introduction to S-linear algebra

In this section we first recall the definition of Smarandache R-module and Smarandache k-vectorial space. Then we give different types of Smarandache linear algebra and Smarandache vector space.

Further we define Smarandache vector spaces over the finite rings which are analogous to vector spaces defined over the prime field $Z_p$. Throughout this section $Z_n$ will denote the ring of integers modulo n, n a composite number $Z_n[x]$ will denote the polynomial ring in the variable x with coefficients from $Z_n$.

**DEFINITION [27, 40]:** *The Smarandache-R-Module (S-R-module) is defined to be an R-module (A, +, ×) such that a proper subset of A is a S-algebra (with respect with the same induced operations and another '×' operation internal on A), where R is a commutative unitary Smarandache ring (S-ring) and S its proper subset that is a field. By a proper subset we understand a set included in A, different from the empty set, from the unit element if any and from A.*

**DEFINITION [27, 40]:** *The Smarandache k-vectorial space (S-k-vectorial space) is defined to be a k-vectorial space,*

*(A, +, •) such that a proper subset of A is a k-algebra (with respect with the same induced operations and another '×' operation internal on A) where k is a commutative field. By a proper subset we understand a set included in A different from the empty set from the unit element if any and from A. This S-k-vectorial space will be known as type I, S-k-vectorial space.*

Now we proceed on to define the notion of Smarandache k-vectorial subspace.

**DEFINITION 1.2.1:** *Let A be a k-vectorial space. A proper subset X of A is said to be a Smarandache k-vectorial subspace (S-k-vectorial subspace) of A if X itself is a Smarandache k-vectorial space.*

**THEOREM 1.2.1:** *Let A be a k-vectorial space. If A has a S-k-vectorial subspace then A is a S-k-vectorial space.*

*Proof:* Direct by the very definitions.

Now we proceed on to define the concept of Smarandache basis for a k-vectorial space.

**DEFINITION 1.2.2:** *Let V be a finite dimensional vector space over a field k. Let B = $\{v_1, v_2, \ldots, v_n\}$ be a basis of V. We say B is a Smarandache basis (S-basis) of V if B has a proper subset say A, $A \subset B$ and $A \neq \phi$ and $A \neq B$ such that A generates a subspace which is a linear algebra over k that is if W is the subspace generated by A then W must be a k-algebra with the same operations of V.*

**THEOREM 1.2.2:** *Let V be a vector space over the field k. If B is a S-basis then B is a basis of V.*

*Proof:* Straightforward by the very definitions.

The concept of S-basis leads to the notion of Smarandache strong basis which is not present in the usual vector spaces.

**DEFINITION 1.2.3:** *Let V be a finite dimensional vector space over a field k. Let B = {$x_1$, …, $x_n$} be a basis of V. If every proper subset of B generates a linear algebra over k then we call B a Smarandache strong basis (S-strong basis) for V.*

Now having defined the notion of S-basis and S-strong basis we proceed on to define the concept of Smarandache dimension.

**DEFINITION 1.2.4:** *Let L be any vector space over the field k. We say L is a Smarandache finite dimensional vector space (S-finite dimensional vector space) of k if every S-basis has only finite number of elements in it. It is interesting to note that if L is a finite dimensional vector space then L is a S-finite dimensional space provided L has a S-basis.*

It can also happen that L need not be a finite dimensional space still L can be a S-finite dimensional space.

**THEOREM 1.2.3:** *Let V be a vector space over the field k. If A = {$x_1$, …, $x_n$} is a S-strong basis of V then A is a S-basis of V.*

*Proof:* Direct by definitions, hence left for the reader as an exercise.

**THEOREM 1.2.4:** *Let V be a vector space over the field k. If A = {$x_1$, …, $x_n$ } is a S-basis of V. A need not in general be a S-strong basis of V.*

*Proof:* By an example. Let V = Q [x] be the set of all polynomials of degree less than or equal to 10. V is a vector space over Q.

Clearly A = {1, x, $x^2$, …, $x^{10}$ } is a basis of V. In fact A is a S-basis of V for take B = {1, $x^2$, $x^4$, $x^6$, $x^8$, $x^{10}$}. Clearly B generates a linear algebra. But all subsets of A do not form a S-basis of V, so A is not a S-strong basis of V but only a S-basis of V.

We will define Smarandache eigen values and Smarandache eigen vectors of a vector space.

**DEFINITION 1.2.5:** *Let V be a vector space over the field F and let T be a linear operator from V to V. T is said to be a Smarandache linear operator (S-linear operator) on V if V has a S-basis, which is mapped by T onto another Smarandache basis of V.*

**DEFINITION 1.2.6:** *Let T be a S-linear operator defined on the space V. A characteristic value c in F of T is said to be a Smarandache characteristic value (S-characteristic value) of T if the characteristic vector of T associated with c generate a subspace, which is a linear algebra that is the characteristic space, associated with c is a linear algebra. So the eigen vector associated with the S-characteristic values will be called as Smarandache eigen vectors (S-eigen vectors) or Smarandache characteristic vectors (S-characteristic vectors).*

Thus this is the first time such Smarandache notions are given about S-basis, S-characteristic values and S-characteristic vectors. For more about these please refer [43, 46].

Now we proceed on to define type II Smarandache k-vector spaces.

**DEFINITION 1.2.7:** *Let R be a S-ring. V be a module over R. We say V is a Smarandache vector space of type II (S-vector space of type II) if V is a vector space over a proper subset k of R where k is a field.*

We have no means to interrelate type I and type II Smarandache vector spaces.

However in case of S-vector spaces of type II we define a stronger version.

**DEFINITION 1.2.8:** *Let R be a S-ring, M a R-module. If M is a vector space over every proper subset k of R which is a field; then we call M a Smarandache strong vector space of type II (S-strong vector space of type II).*

**THEOREM 1.2.5:** *Every S-strong vector space of type II is a S-vector space of type II.*

*Proof:* Direct by the very definition.

***Example 1.2.1:*** Let $Z_{12}$ [x] be a module over the S-ring $Z_{12}$. $Z_{12}$ [x] is a S-strong vector space of type II.

***Example 1.2.2:*** Let $M_{2\times2} = \{(a_{ij}) \mid a_{ij} \in Z_6\}$ be the set of all $2 \times 2$ matrices with entries from $Z_6$. $M_{2\times2}$ is a module over $Z_6$ and $M_{2\times2}$ is a S-strong vector space of type II.

***Example 1.2.3:*** Let $M_{3\times5} = \{(a_{ij}) \mid a_{ij} \in Z_6\}$ be a module over $Z_6$. Clearly $M_{3\times5}$ is a S-strong vector space of type II over $Z_6$.

Now we proceed on to define Smarandache linear algebra of type II.

**DEFINITION 1.2.9:** *Let R be any S-ring. M a R-module. M is said to be a Smarandache linear algebra of type II (S-linear algebra of type II) if M is a linear algebra over a proper subset k in R where k is a field.*

**THEOREM 1.2.6:** *All S-linear algebra of type II is a S-vector space of type II and not conversely.*

*Proof:* Let M be an R-module over a S-ring R. Suppose M is a S-linear algebra II over $k \subset R$ (k a field contained in R) then by the very definition M is a S-vector space II.

To prove converse we have show that if M is a S-vector space II over $k \subset R$ (R a S-ring and k a field in R) then M in general need not be a S-linear algebra II over k contained in R. Now by example 1.2.3 we see the collection $M_{3\times5}$ is a S-vector space II over the field k {0, 2, 4} contained in $Z_6$. But clearly $M_{3\times5}$ is not a S-linear algebra II over $\{0, 2, 4\} \subset Z_6$.

We proceed on to define Smarandache subspace II and Smarandache subalgebra II.

**DEFINITION 1.2.10:** *Let M be an R-module over a S-ring R. If a proper subset P of M is such that P is a S-vector space of type II over a proper subset k of R where k is a field then we call P a Smarandache subspace II (S-subspace II) of M relative to P.*

It is important to note that even if M is a R-module over a S-ring R, and M has a S-subspace II still M need not be a S-vector space of type II.

On similar lines we will define the notion of Smarandache subalgebra II.

**DEFINITION 1.2.11:** *Let M be an R-module over a S-ring R. If M has a proper subset P such that P is a Smarandache linear algebra II (S-linear algebra II) over a proper subset k in R where k is a field then we say P is a S-linear subalgebra II over R.*

Here also it is pertinent to mention that if M is a R-module having a subset P that is a S-linear subalgebra II then M need not in general be a S-linear algebra II. It has become important to mention that in all algebraic structure, S if it has a proper substructure P that is Smarandache then S itself is a Smarandache algebraic structure. But we see in case of R-Modules M over the S-ring R if M has a S-subspace or S-subalgebra over a proper subset k of R where k is a field still M in general need not be a S-vector space or a S-linear algebra over k; $k \subset R$.

Now we will illustrate this by the following examples.

***Example 1.2.4:*** Let $M = R[x] \times R[x]$ be a direct product of polynomial rings, over the ring $R \times R$. Clearly $M = R[x] \times R[x]$ is a S-vector space over the field $k = R \times \{0\}$.

It is still more interesting to see that M is a S-vector space over k = {0} × Q, Q the field of rationals. Further M is a S-strong vector space as M is a vector space over every proper subset of R × R which is a field.

Still it is important to note that M = R [x] × R [x] is a S-strong linear algebra. We see Q[x] × Q[x] = P ⊂ M is a S-subspace over $k_1$ = Q × {0} and {0} × Q but P is not a S-subspace over $k_2$ = R × {0} or {0} × R.

Now we will proceed on to define Smarandache vector spaces and Smarandache linear algebras of type III.

**DEFINITION 1.2.12:** *Let M be a any non empty set which is a group under '+'. R any S-ring. M in general need not be a module over R but a part of it is related to a section of R. We say M is a Smarandache vector space of type III (S-vector space III) if M has a non-empty proper subset P which is a group under '+', and R has a proper subset k such that k is a field and P is a vector space over k.*

Thus this S-vector space III links or relates and gets a nice algebraic structure in a very revolutionary way.

We illustrate this by an example.

***Example 1.2.5:*** Consider M = Q [x] × Z [x]. Clearly M is an additively closed group. Take R = Q × Q; R is a S-ring. Now P = Q [x] × {0} is a vector space over k = Q × {0}. Thus we see M is a Smarandache vector space of type III.

So this definition will help in practical problems where analysis is to take place in such set up.

Now we can define Smarandache linear algebra of type III in an analogous way.

**DEFINITION 1.2.13:** *Suppose M is a S-vector space III over the S-ring R. We call M a Smarandache linear algebra of type III (S-linear algebra of type III) if P ⊂ M which is a vector space over k ⊂ R (k a field) is a linear algebra.*

Thus we have the following naturally extended theorem.

**THEOREM 1.2.7:** *Let M be a S-linear algebra III for $P \subset M$ over R related to the subfield $k \subset R$. Then clearly P is a S-vector space III.*

*Proof:* Straightforward by the very definitions.

To prove that all S-vector space III in general is not a S-linear algebra III we illustrate by an example.

***Example 1.2.6:*** Let $M = P_1 \cup P_2$ where $P_1 = M_{3\times3} = \{(a_{ij}) \mid a_{ij} \in Q\}$ and $P_2 = M_{2\times2} = \{(a_{ij}) \mid a_{ij} \in Z\}$ and R be the field of reals. Now take the proper subset $P = P_1$, $P_1$ is a S-vector space III over $Q \underset{\neq}{\subset} R$. Clearly $P_1$ is not a S-linear algebra III over Q.

Now we proceed on to define S-subspace III and S-linear algebra III.

**DEFINITION 1.2.14:** *Let M be an additive Aeolian group, R any S-ring. $P \subset M$ be a S-vector space III over a field $k \subset R$. We say a proper subset $T \subset P$ to be a Smarandache vector subspace III (S-vector subspace III) or S-subspace III if T itself is a vector space over k.*

*If a S-vector space III has no proper S-subspaces III relative to a field $k \subset R$ then we call M a Smarandache simple vector space III (S-simple vector space III).*

On similar lines one defines Smarandache sublinear algebra III and S-simple linear algebra III.

Yet a new notion called Smarandache super vector spaces are introduced for the first time.

**DEFINITION 1.2.15:** *Let R be S-ring. V a module over R. We say V is a Smarandache super vector space (S-super vector space) if V is a S-k-vector space over a proper set k, $k \subset R$ such that k is a field.*

**THEOREM 1.2.8:** *All S-super spaces are S-k-vector spaces over the field k, k contained in R.*

*Proof:* Straightforward.

Almost all results derived in case of S-vector spaces type II can also be derived for S-super vector spaces.

Further for given V, a R-module of a S-ring R we can have several S-super vector spaces.

Now we just give the definition of Smarandache super linear algebra.

**DEFINITION 1.2.16:** *Let R be a S-ring. V a R module over R. Suppose V is a S-super vector space over the field k, $k \subset R$. we say V is a S-super linear algebra if for all a, $b \in V$ we have $a \bullet b \in V$ where '$\bullet$' is a closed associative binary operation on V.*

Almost all results in case of S-linear algebras can be easily extended and studied in case of S-super linear algebras.

**DEFINITION 1.2.17:** *Let V be an additive abelian group, $Z_n$ be a S-ring (n a composite number). V is said to be a Smarandache vector space over $Z_n$ (S-vector space over $Z_n$) if for some proper subset T of $Z_n$ where T is a field satisfying the following conditions:*

i. *$vt$, $tv \in V$ for all $v \in V$ and $t \in T$.*
ii. *$t(v_1 + v_2) = tv_1 + tv_2$ for all $v_1$ $v_2 \in V$ and $t \in T$.*
iii. *$(t_1 + t_2)v = t_1v + t_2v$ for all $v \in V$ and $t_1$, $t_2 \in T$.*
iv. *$t_1(t_2 u) = (t_1 t_2) u$ for all $t_1, t_2 \in T$ and $u \in V$.*
v. *if e is the identity element of the field T then $ve = ev = v$ for all $v \in V$.*

*In addition to all these if we have an multiplicative operation on V such that $u \bullet v_1 \in V$ for all $u \bullet v_1 \in V$ then we call V a Smarandache linear algebra (S-linear algebra) defined over finite S-rings.*

It is a matter of routine to check that if V is a S-linear algebra then obviously V is a S-vector space. We however will illustrate by an example that all S-vector spaces in general need not be S-linear algebras.

***Example 1.2.7*:** Let $Z_6 = \{0, 1, 2, 3, 4, 5\}$ be a S-ring (ring of integers modulo 6). Let $V = M_{2\times 3} = \{(a_{ij}) \mid a_{ij} \in Z_6\}$.

Clearly V is a S-vector space over $T = \{0, 3\}$. But V is not a S-linear algebra. Clearly V is a S-vector space over $T_1 = \{0, 2, 4\}$. The unit being 4 as $4^2 \equiv 4 \pmod 6$.

***Example 1.2.8*:** Let $Z_{12} = \{0, 1, 2, \ldots, 10, 11\}$ be the S-ring of characteristic two. Consider the polynomial ring $Z_{12}[x]$ . Clearly $Z_{12}[x]$ is a S-vector space over the field $k = \{0, 4, 8\}$ where 4 is the identity element of k and k is isomorphic with the prime field $Z_3$.

***Example 1.2.9*:** Let $Z_{18} = \{0, 1, 2, \ldots, 17\}$ be the S-ring. $M_{2\times 2} = \{(a_{ij}) \mid a_{ij} \in Z_{18}\}$ $M_{2\times 2}$ is a finite S-vector space over the field $k = \{0, 9\}$. What is the basis for such space?

Here we see $M_{2\times 2}$ has basis

$$\begin{bmatrix} 1 & 0 \\ 0 & 0 \end{bmatrix}, \begin{bmatrix} 0 & 1 \\ 0 & 0 \end{bmatrix}, \begin{bmatrix} 0 & 0 \\ 0 & 1 \end{bmatrix} \text{and} \begin{bmatrix} 0 & 0 \\ 1 & 0 \end{bmatrix}.$$

Clearly $M_{2x2}$ is not a vector space over $Z_{18}$ as $Z_{18}$ is only a ring.

Now we proceed on to characterize those finite S-vector spaces, which has only one field over which the space is defined. We call such S-vector spaces as Smarandache unital vector spaces. The S-vector space $M_{2\times 2}$ defined over $Z_{18}$ is a S-unital vector space. When the S-vector space has more than one S-vector space defined over more than one field we call the S-vector space as Smarandache multi vector space (S-multi vector space).

For consider the vector space $Z_6[x]$. $Z_6[x]$ is the polynomial ring in the indeterminate x with coefficients from $Z_6$. Clearly $Z_6[x]$ is a S-vector space over . k = {0, 3}; k is a field isomorphic with $Z_2$ and $Z_6[x]$ is also a S-vector space over $k_1$ = {0, 2, 4} a field isomorphic to $Z_3$ . Thus $Z_6[x]$ is called S-multi vector space.

Now we have to define Smarandache linear operators and Smarandache linear transformations. We also define for these finite S-vector spaces the notion of Smarandache eigen values and Smarandache eigen vectors and its related notions.

Throughout this section we will be considering the S-vector spaces only over finite rings of modulo integers $Z_n$ (n always a positive composite number).

**DEFINITION 1.2.18:** *Let U and V be a S-vector spaces over the finite ring $Z_n$. i.e. U and V are S-vector space over a finite field P in $Z_n$. That is $P \subset Z_n$ and P is a finite field. A Smarandache linear transformation (S-linear transformation) T of U to V is a map given by $T(c\,\alpha + \beta) = c\,T(\alpha) + T(\beta)$ for all $\alpha, \beta \in U$ and $c \in P$. Clearly we do not demand c to be from $Z_n$ or the S-vector spaces U and V to be even compatible with the multiplication of scalars from $Z_n$.*

***Example 1.2.10:*** Let $Z_{15}^{8}[x]$ and $M_{3\times 3} = \{(a_{ij}) \mid a_{ij} \in Z_{15}\}$ be two S-vector spaces defined over the finite S-ring. Clearly both $Z_{15}^{8}[x]$ and $M_{3\times 3}$ are S-vector spaces over P = {0, 5, 10} a field isomorphic to $Z_3$ where 10 serves as the unit element of P. $Z_{15}^{8}[x]$ is a additive group of polynomials of degree less than or equal to 8 and $M_{3\times 3}$ is the additive group of matrices.

Define T: $Z^{8}_{15}[x] \to M_{3x3}$ by

$$T(p_0 + p_1x + \ldots + p_8x^8) = \begin{bmatrix} p_0 & p_1 & p_2 \\ p_3 & p_4 & p_5 \\ p_6 & p_7 & p_8 \end{bmatrix}.$$

Thus T is a S-linear transformation. Both the S-vector spaces are of dimension 9.

Now we see the groups $Z_{15}^{8}[x]$ and $M_{3\times3}$ are also S-vector spaces over $P_1 = \{0, 3, 6, 9, 12\}$, this is a finite field isomorphic with $Z_5$, 6 acts as the identity element.

Thus we see we can have for S-vector spaces more than one field over which they are vector spaces.

Thus we can have a S-vector spaces defined over finite ring, we can have more than one base field. Still they continue to have all properties.

***Example 1.2.11:*** Let $M_{3\times3} = \{(a_{ij}) \mid a_{ij} \in \{0, 3, 6, 9, 12\} \subset Z_{15}\}$. $M_{3\times3}$ is a S-vector space over the S-ring $Z_{15}$. i.e.$M_{3\times3}$ is a S-vector space over $P = \{0, 3, 6, 9, 12\}$ where P is the prime field isomorphic to $Z_{15}$.

***Example 1.2.12:*** $V = Z_{12} \times Z_{12} \times Z_{12}$ is a S-vector space over the field, $P = \{0, 4, 8\} \subset Z_{12}$.

**DEFINITION 1.2.19:** *Let $Z_n$ be a finite ring of integers. V be a S-vector space over the finite field P, $P \subset Z_n$. We call V a Smarandache linear algebra (S-linear algebra) over a finite field P if in V we have an operation '$\bullet$' such that for all a, b $\in V$, $a \bullet b \in V$.*

It is important to mention here that all S-linear algebras over a finite field is a S-vector space over the finite field. But however every S-vector space over a finite field, in general need not be a S-linear algebra over a finite field k.

We illustrate this by an example.

***Example 1.2.13:*** Let $M_{7\times3} = \{(a_{ij}) \mid a_{ij} \in Z_{18}\}$ i.e. the set of all $7 \times 3$ matrices. $M_{7\times3}$ is an additive group. Clearly $M_{7\times3}$ is a S-vector space over the finite field, $P = \{0, 9\} \subset Z_{18}$. It is easily verified that $M_{7\times3}$ is not a S-linear algebra.

Now we proceed on to define on the collection of S-linear transformation of two S-vector spaces relative to the same field P in $Z_n$. We denote the collection of all S-linear transformation from two S-vector spaces U and V relative to the field $P \subset Z_n$ by $SL_P$ (U, V). Let V be a S-vector space defined over the finite field P, $P \subset Z_n$. A map $T_P$ from V to V is said to be a Smarandache linear operator (S-linear operator) of V if $T_P(c\alpha + \beta) = c\, T_P(\alpha) + T_P(\beta)$ for all $\alpha, \beta \in V$ and $c \in P$. Let $SL_P(V, V)$ denote the collections of all S-linear operators from V to V.

**DEFINITION 1.2.20:** *Let V be a S-vector space over a finite field $P \subset Z_n$. Let T be a S-linear operator on V. A Smarandache characteristic value (S-characteristic value) of T is a scalar c in P (P a finite field of the finite S-ring $Z_n$) such that there is a non-zero vector $\alpha$ in V with $T\alpha = c\alpha$. If c is a S-characteristic value of T, then*

i. *Any $\alpha$ such that $T\alpha = c\alpha$ is called a S-characteristic vector of T associated with the S-characteristic value c.*

ii. *The collection of all $\alpha$ such that $T\alpha = c\alpha$ is called the S-characteristic space associated with c.*

Almost all results studied and developed in the case of S-vector spaces can be easily defined and analyzed, in case of S-vector spaces over finite fields, P in $Z_n$.

Thus in case of S-vector spaces defined over $Z_n$ the ring of finite integers we can have for a vector space V defined over $Z_n$ we can have several S-vector spaces relative to $P_i \subset Z_n$, $P_i$ subfield of $Z_n$ Each $P_i$ will make a S-linear operator to give distinct S-characteristic values and S-characteristic vectors. In some cases we may not be even in a position to have all characteristic roots to be present in the same field $P_i$ such situations do not occur in our usual vector spaces they are only possible in case of Smarandache structures.

Thus a S-characteristic equation, which may not be reducible over one of the fields, $P_i$ may become reducible over some other field $P_j$. This is the quality of S-vector spaces over finite rings $Z_n$.

Study of projections $E_i$, primary decomposition theorem in case of S-finite vector spaces will yield several interesting results. So for a given vector space V over the finite ring $Z_n$ V be S-vector spaces over the fields $P_1, \ldots, P_m$, where $P_i \subset Z_n$, are fields in $Z_n$ and V happen to be S-vector space over each of these $P_i$ then we can have several decomposition of V each of which will depend on the fields $P_i$. Such mixed study of a single vector space over several fields is impossible except for the Smarandache imposition.

Now we can define inner product not the usual inner product but inner product dependent on each field which we have defined earlier. Using the new pseudo inner product once again we will have the modified form of spectral theorem. That is, the Smarandache spectral theorem which we will be describing in the next paragraph for which we need the notion of Smarandache new pseudo inner product on V.

Let V be an additive abelian group. $Z_n$ be a ring of integers modulo n, n a finite composite number. Suppose V is a S-vector space over the finite fields $P_1, \ldots, P_t$ in $Z_n$ where each $P_i$ is a proper subset of $Z_n$ which is a field and V happens to be a vector space over each of these $P_i$. Let $\langle \, , \rangle_{P_i}$ be an inner product defined on V relative to each $P_i$. Then $\langle \, , \rangle_{P_i}$ is called the Smarandache new pseudo inner product on V relative to $P_i$.

Now we just define when is a Smarandache linear operator T, Smarandache self-adjoint. We say T is Smarandache self adjoint (S- self adjoint) if $T = T^*$.

***Example 1.2.14:*** Let V = $Z_6^2[x]$ be a S-vector space over the finite field, P = {0, 2, 4}, {1, x, $x^2$} is a S-basis of V,

$$A = \begin{bmatrix} 4 & 0 & 0 \\ 0 & 2 & 2 \\ 0 & 2 & 2 \end{bmatrix}$$

be the matrix associated with a linear operator T.

$$\lambda - AI = \begin{bmatrix} \lambda - 4 & 0 & 0 \\ 0 & \lambda - 2 & 4 \\ 0 & 4 & \lambda - 2 \end{bmatrix}$$

$$\begin{aligned} &= (\lambda - 4)\,[(\lambda - 2)\,(\lambda - 2) - 4] \\ &= (\lambda - 4)\,[(\lambda - 2)^2] - 4\,(\lambda - 4) = 0 \\ &= \lambda^3 - 2\lambda^2 + 4\lambda = 0 \end{aligned}$$

$\lambda = 0, 4, 4$ are the S-characteristic values. The S-characteristic vector for $\lambda = 4$ are

$$V_1 = (0, 4, 4)$$
$$V_2 = (4, 4, 4).$$

For $\lambda = 0$ the characteristic vector is (0, 2, 4). So

$$A = \begin{bmatrix} 4 & 0 & 0 \\ 0 & 2 & 2 \\ 0 & 2 & 2 \end{bmatrix} = A^*.$$

Thus T is S-self adjoint operator.

$W_1$ is the S-subspace generated by {(0, 4, 4), (4, 4, 4)}. $W_2$ is the S-subspace generated by {(0, 2, 4)}.

$$\begin{aligned} V &= W_1 + W_2. \\ T &= c_1E_1 + c_2 E_2. \\ c_1 &= 4. \\ c_2 &= 0. \end{aligned}$$

**THEOREM (SMARANDACHE SPECTRAL THEOREM FOR S-VECTOR SPACES OVER FINITE RINGS $Z_N$):** *Let $T_i$ be a Smarandache self adjoint operator on the S-finite* *dimensional pseudo inner product space V = $Z_n[x]$, over each of the finite fields $P_1, P_2, \ldots, P_t$ contained in $Z_n$.*

*Let $c_1, c_2, \ldots, c_k$ be the distinct S-characteristic values of $T_i$ . Let $W_i$ be the S-characteristic space associated with $c_i$ and $E_i$ the orthogonal projection of V on $W_i$, then $W_i$ is orthogonal to $W_j$, $i \neq j$; V is the direct sum of $W_1, \ldots, W_k$ and $T_i = c_1 E_1 + \ldots + c_k E_k$ (we have several such decompositions depending on the number of finite fields in $Z_n$ over which V is defined ).*

*Proof:* Direct as in case of S-vector spaces.

Further results in this direction can be developed as in case of other S-vector spaces.

One major difference here is that V can be S-vector space over several finite fields each finite field will reflect its property.

## 1.3. Some application of S-linear algebra

In this section we just indicate how the study of Markov chains can be defined and studied as Smarandache Markov chains; as in the opinion of the author such study in certain real world problems happens to be better than the existing ones. Further we deal with a Smarandache analogue of Leontief economic models.

### 1.3.1 Smarandache Markov chains

Suppose a physical or a mathematical system is such that at any moment it can occupy one of a finite number of states. When we view them as stochastic process or Markov chains we make a assumption that the system moves with time from one state to another, so that a schedule of observation times and keep the states of the system at these times. But when we tackle real world problems say even for simplicity the emotions of a persons it need not fall under the category of sad, happy, angry, …, many a times the emotions of a person may be very unpredictable depending largely on the situation, such study cannot fall under Markov chains for such states cannot be included and given as next observation, this largely affects the very transition matrix P = $[p_{ij}]$ with nonnegative entries for which each of the column sums are one and all of whose entries are positive. This has relevance as even the policy makers are humans and their view is ultimate and this rules the situation. So to over come this problem when we have indecisive situations we give negative values so that our transition matrix column sums do not add to one and all entries may not be positive.

Thus we call the new transition matrix, which is a square matrix which can have negative entries also falling in the interval [–1, 1] and whose column sums can also be less than 1 as the Smarandache transition matrix (S-transition matrix).

Further the Smarandache probability vector (S-probability vector) will be a column vector, which can take entries from [–1, 1] whose sum can lie in the interval [–1, 1]. The Smarandache probability vectors $x^{(n)}$ for n = 0, 1, 2 , … , are said to be the Smarandache state vectors (S-state vectors) of a Smarandache Markov process (S-Markov process). Clearly if P is a S-transition matrix of a S-Markov process and $x^{(n)}$ is the S-state vectors at the $n^{th}$ observation then $x^{(n+1)} \neq p\, x^{(n)}$ in general.

The interested reader is requested to develop results in this direction.

### 1.3.2 Smarandache Leontief economic models

Matrix theory has been very successful in describing the interrelations between prices, outputs and demands in an economic model. Here we just discuss some simple models based on the ideals of the Nobel-laureate Massily Leontief. Two types of models discussed are the closed or input-output model and the open or production model each of which assumes some economic parameter which describe the inter relations between the industries in the economy under considerations. Using matrix theory we evaluate certain parameters.

The basic equations of the input-output model are the following:

$$\begin{bmatrix} a_{11} & a_{12} & \cdots & a_{1n} \\ a_{21} & a_{22} & \cdots & a_{2n} \\ \vdots & & & \vdots \\ a_{n1} & a_{n2} & \cdots & a_{nn} \end{bmatrix} \begin{bmatrix} p_1 \\ p_2 \\ \vdots \\ p_n \end{bmatrix} = \begin{bmatrix} p_1 \\ p_2 \\ \vdots \\ p_n \end{bmatrix}$$

each column sum of the coefficient matrix is one

i. $p_i \geq 0$, $i = 1, 2, \ldots, n$.
ii. $a_{ij} \geq 0$, $i, j = 1, 2, \ldots, n$.
iii. $a_{ij} + a_{2j} + \ldots + a_{nj} = 1$

for $j = 1, 2, \ldots, n$.

$$p = \begin{bmatrix} p_1 \\ p_2 \\ \vdots \\ p_n \end{bmatrix}$$

are the price vector. $A = (a_{ij})$ is called the input-output matrix

$$Ap = p \text{ that is, } (I - A)\, p = 0.$$

Thus A is an exchange matrix, then Ap = p always has a nontrivial solution p whose entries are nonnegative. Let A be an exchange matrix such that for some positive integer m, all of the entries of $A^m$ are positive. Then there is exactly only one linearly independent solution of $(I - A)\, p = 0$ and it may be chosen such that all of its entries are positive in Leontief open production model.

In contrast with the closed model in which the outputs of k industries are distributed only among themselves, the open model attempts to satisfy an outside demand for the outputs. Portions of these outputs may still be distributed among the industries themselves to keep them operating, but there is to be some excess some net production with which to satisfy the outside demand. In some closed model, the outputs of the industries were fixed and our objective was to determine the prices for these outputs so that the equilibrium condition that expenditures equal incomes was satisfied.

$x_i$ = monetary value of the total output of the $i^{th}$ industry.

$d_i$ = monetary value of the output of the $i^{th}$ industry needed to satisfy the outside demand.

$\sigma_{ij}$ = monetary value of the output of the $i^{th}$ industry needed by the $j^{th}$ industry to produce one unit of monetary value of its own output.

With these qualities we define the production vector.

$$x = \begin{bmatrix} x_1 \\ x_2 \\ \vdots \\ x_k \end{bmatrix}$$

the demand vector

$$d = \begin{bmatrix} d_1 \\ d_2 \\ \vdots \\ d_k \end{bmatrix}$$

and the consumption matrix,

$$C = \begin{bmatrix} \sigma_{11} & \sigma_{12} & \cdots & \sigma_{1k} \\ \sigma_{21} & \sigma_{22} & \cdots & \sigma_{2k} \\ \vdots & & & \\ \sigma_{k1} & \sigma_{k2} & \cdots & \sigma_{kk} \end{bmatrix}.$$

By their nature we have

$$x \geq 0, d \geq 0 \text{ and } C \geq 0.$$

From the definition of $\sigma_{ij}$ and $x_j$ it can be seen that the quantity

$$\sigma_{i1} x_1 + \sigma_{i2} x_2 + \ldots + \sigma_{ik} x_k$$

is the value of the output of the $i^{th}$ industry needed by all k industries to produce a total output specified by the production vector x.

Since this quantity is simply the $i^{th}$ entry of the column vector Cx, we can further say that the $i^{th}$ entry of the column vector x – Cx is the value of the excess output of the $i^{th}$ industry available to satisfy the outside demand. The value of the outside demand for the output of the $i^{th}$ industry is the $i^{th}$ entry of the demand vector d; consequently; we are led to the following equation:

$$x - Cx = d \text{ or}$$
$$(I - C)\, x = d$$

for the demand to be exactly met without any surpluses or shortages. Thus, given C and d, our objective is to find a production vector $x \geq 0$ which satisfies the equation $(I - C)x = d$.

A consumption matrix C is said to be productive if $(1 - C)^{-1}$ exists and $(1 - C)^{-1} \geq 0$.

A consumption matrix C is productive if and only if there is some production vector $x \geq 0$ such that $x > Cx$.

A consumption matrix is productive if each of its row sums is less than one. A consumption matrix is productive if each of its column sums is less than one.

Now we will formulate the Smarandache analogue for this, at the outset we will justify why we need an analogue for those two models.

Clearly, in the Leontief closed Input – Output model, $p_i$ = price charged by the $i^{th}$ industry for its total output in reality need not be always a positive quantity for due to competition to capture the market the price may be fixed at a loss or the demand for that product might have fallen down so badly so that the industry may try to charge very less than its real value just to market it.

Similarly $a_{ij} \geq 0$ may not be always be true. Thus in the Smarandache Leontief closed (Input – Output) model (S-Leontief closed (Input-Output) model) we do not demand $p_i \geq 0$, $p_i$ can be negative; also in the matrix $A = (a_{ij})$,

$$a_{1j} + a_{2j} + \ldots + a_{kj} \neq 1$$

so that we permit $a_{ij}$'s to be both positive and negative, the only adjustment will be we may not have $(I - A)\, p = 0$, to have only one linearly independent solution, we may have more than one and we will have to choose only the best solution.

As in this complicated real world problems we may not have in practicality such nice situation. So we work only for the best solution.

On similar lines we formulate the Smarandache Leontief open model (S-Leontief open model) by permitting that $x \geq 0$ , $d \geq 0$ and $C \geq 0$ will be allowed to take $x \leq 0$ or $d \leq 0$ and or $C \leq 0$ . For in the opinion of the author we may not in reality have the monetary total output to be always a positive quality for all industries and similar arguments for $d_i$'s and $C_{ij}$'s.

When we permit negative values the corresponding production vector will be redefined as Smarandache production vector (S-production vector) the demand vector as Smarandache demand vector (S-demand vector) and the consumption matrix as the Smarandache consumption matrix (S-consumption matrix). So when we work out under these assumptions we may have different sets of conditions

We say productive if $(1 – C)^{-1} \geq 0$, and non-productive or not up to satisfaction if $(1 – C)^{-1} < 0$.

The reader is expected to construct real models by taking data's from several industries. Thus one can develop several other properties in case of different models.

## Chapter Two

# INTRODUCTORY CONCEPTS ON BASIC BISTRUCTURES AND S-BISTRUCTURES

In this chapter we recall some of the basic concepts on bistructures and S-bistructures used in this book to make the book a self contained one. This chapter has two sections. In section one, we give the basic definition in section 2 some basic notions of S-bistructure.

## 2.1 Basic concepts of bigroups and bivector spaces

This section is devoted to the recollection of bigroups, sub-bigroups and we illustrate it with examples. [40] was the first one to introduce the notion of bigroups in the year 1994. As there is no book on bigroups we give all algebraic aspects of it.

**DEFINITION [40]:** *A set (G, +, •) with two binary operation '+' and '•' is called a bigroup if there exist two proper subsets $G_1$ and $G_2$ of G such that*

i. $G = G_1 \cup G_2$.
ii. $(G_1, +)$ *is a group.*
iii. $(G_2, \bullet)$ *is a group.*

*A subset H ($\neq \phi$ ) of a bigroup (G, + , •) is called a sub-bigroup, if H itself is a bigroup under '+' and '•' operations defined on G.*

**THEOREM 2.1.1:** *Let (G, +, •) be a bigroup. The subset H $\neq \phi$ of a bigroup G is a sub-bigroup, then (H, + ) and (H, •) in general are not groups.*

*Proof*: Given (G, +, •) is a bigroup and H $\neq \phi$ of G is a sub-bigroup of G. To show (H, +) and (H, •) are not groups.

We give an example to prove this. Consider the bigroup G = { …, –2, –1, 0 1, 2, …} ∪ {i, –j} under the operations '+' and '•'. G = $G_1 \cup G_2$ where ($G_1$, •) = {–1, 1, i, – i} under product and $G_2$ = { …, –2, –1, 0 1, 2, …} under '+' are groups.

Take H = {–1, 0, 1}. H is a proper subset of G and H = $H_1 \cup H_2$ where $H_1$ = {0} and $H_2$ = {–1, 1}; $H_1$ is a group under '+' and $H_2$ is a group under product i.e. multiplication.

Thus H is a sub-bigroup of G but H is not a group under '+' or '•'.

Now we get a characterization theorem about sub-bigroup in what follows:

**THEOREM [46]:** *Let (G, +, •) be a bigroup. Then the subset H ($\neq \phi$ ) of G is a sub-bigroup of G if and only if there exists two proper subsets $G_1$, $G_2$ of G such that*

i. *G = $G_1 \cup G_2$ where ($G_1$, +) and ($G_2$, •) are groups.*
ii. *($H \cap G_1$, +) is subgroup of ($G_1$, +).*
iii. *($H \cap G_2$, •) is a subgroup of ($G_2$, •).*

*Proof*: Let H ($\neq \phi$) be a sub-bigroup of G, then (H, +, •) is a bigroup. Therefore there exists two proper subsets $H_1$, $H_2$ of H such that

i. H = $H_1 \cup H_2$.

ii. $(H_1, +)$ is a group.
iii. $(H_2, \bullet)$ is a group.

Now we choose $H_1$ as $H \cap G_1$ then we see that $H_1$ is a subset of $G_1$ and by (ii) $(H_1, +)$ is itself a group. Hence ($H_1 = H \cap (G_1, +)$) is a subgroup of $(G_1, +)$. Similarly ($H_2 = H \cap G_2, \bullet$) is a subgroup of $(G_2, \bullet)$. Conversely let (i), (ii) and (iii) of the statement of theorem be true. To prove $(H, +, \bullet)$ is a bigroup it is enough to prove $(H \cap G_1) \cup (H \cap G_2) = H$.

Now, $(H \cap G_1) \cup (H \cap G_2) = [(H \cap G_1) \cup H] \cap [(H \cap G_1) \cup G_2]$

$= [(H \cup H) \cap (G_1 \cup H)] \cap [(H \cup G_2) \cap (G_1 \cup G_2)]$
$= [H \cap (G_1 \cup H)] \cap [(H \cup G_2) \cap G]$
$= H \cap (H \cap G_2)$ (since $H \subseteq G_1 \cup H$ and $H \cup G_2 \subseteq G$)
$= H$ (since $H \subset H \cup G_2$).

Hence the theorem is proved.

It is important to note that in the above theorem just proved the condition (i) can be removed. We include this condition only for clarification or simplicity of understanding.

Another natural question would be can we have at least some of the classical theorems and some more classical concepts to be true in them.

**DEFINITION 2.1.1:** *Let $(G, +, \bullet)$ be a bigroup where $G = G_1 \cup G_2$; bigroup G is said to be commutative if both $(G_1, +)$ and $(G_2, \bullet)$ are commutative.*

***Example 2.1.1***: Let $G = G_1 \cup G_2$ where $G_1 = Q \setminus \{0\}$ with usual multiplication and $G_2 = \langle g \mid g^2 = I \rangle$ be a cyclic group of order two. Clearly G is a commutative bigroup.

We say the order of the bigroup $G = G_1 \cup G_2$ is finite if the number of elements in them is finite; otherwise we say the bigroup G to be of infinite order.

The bigroup given in example 2.1.1 is a bigroup of infinite order.

***Example 2.1.2*:** Let G = $G_1 \cup G_2$ where $G_1 = Z_{10}$ group under addition modulo 10 and $G_2 = S_3$ the symmetric group of degree3. Clearly G is a non-commutative bigroup of finite order. $|G| = 16$.

***Example 2.1.3*:** Let G = $G_1 \cup G_2$ be a bigroup with $G_1$ = {set of all n × n matrices under '+' over the field of reals} and $G_2$ = { set of all n × n matrices A with $|A| \neq \{0\}$ with entries from Q}, ($G_1$, +) and ($G_2$, ×) are groups and G = $G_1 \cup G_2$ is a non-commutative bigroup of infinite order.

In this section we introduce the concept of bivector spaces and S-bivector spaces. The study of bivector spaces started only in 1999 [106]. Here we recall these definitions and extend it to the Smarandache bivector spaces.

**DEFINITION 2.1.2:** *Let V = $V_1 \cup V_2$ where $V_1$ and $V_2$ are two proper subsets of V and $V_1$ and $V_2$ are vector spaces over the same field F that is V is a bigroup, then we say V is a bivector space over the field F.*

*If one of $V_1$ or $V_2$ is of infinite dimension then so is V. If $V_1$ and $V_2$ are of finite dimension so is V; to be more precise if $V_1$ is of dimension n and $V_2$ is of dimension m then we define dimension of the bivector space V = $V_1 \cup V_2$ to be of dimension m + n. Thus there exists only m + n elements which are linearly independent and has the capacity to generate V = $V_1 \cup V_2$.*

The important fact is that same dimensional bivector spaces are in general not isomorphic.

***Example 2.1.4:*** Let V = $V_1 \cup V_2$ where $V_1$ and $V_2$ are vector spaces of dimension 4 and 5 respectively defined over rationals where $V_1 = \{(a_{ij}) / a_{ij} \in Q\}$, collection of all 2 × 2 matrices with entries from Q. $V_2$ = {Polynomials of degree less than or equal to 4}.

Clearly V is a finite dimensional bivector space of dimension 9. In order to avoid confusion we can follow the following convention whenever essential. If $v \in V = V_1 \cup$ $V_2$ then $v \in V_1$ or $v \in V_2$ if $v \in V_1$ then v has a representation of the form $(x_1, x_2, x_3, x_4, 0, 0, 0, 0, 0)$ where $(x_1, x_2, x_3, x_4) \in V_1$ if $v \in V_2$ then $v = (0, 0, 0, 0, y_1, y_2, y_3, y_4, y_5)$ where $(y_1, y_2, y_3, y_4, y_5) \in V_2$.

Thus we follow the notation.

***Notation***: Let $V = V_1 \cup V_2$ be the bivector space over the field F with dimension of V to be m + n where dimension of $V_1$ is m and that of $V_2$ is n. If $v \in V = V_1 \cup V_2$, then $v \in V_1$ or $v \in V_2$ if $v \in V_1$ then $v = (x_1, x_2, \ldots, x_m, 0, 0, \ldots, 0)$ if $v \in V_2$ then $v = (0, 0, \ldots, 0, y_1, y_2, \ldots, y_n)$.

We never add elements of $V_1$ and $V_2$. We keep them separately as no operation may be possible among them. For in example 2.1.4 we had $V_1$ to be the set of all 2 × 2 matrices with entries from Q where as $V_2$ is the collection of all polynomials of degree less than or equal to 4. So no relation among elements of $V_1$ and $V_2$ is possible. Thus we also show that two bivector spaces of same dimension need not be isomorphic by the following example:

***Example 2.1.5:*** Let $V = V_1 \cup V_2$ and $W = W_1 \cup W_2$ be any two bivector spaces over the field F. Let V be of dimension 8 where $V_1$ is a vector space of dimension 2, say $V_1 = F \times F$ and $V_2$ is a vector space of dimension say 6 all polynomials of degree less than or equal to 5 with coefficients from F. W be a bivector space of dimension 8 where $W_1$ is a vector space of dimension 3 i.e. $W_1$ = {all polynomials of degree less than or equal to 2} with coefficients from F and $W_2 = F \times F \times F \times F \times F$ a vector space of dimension 5 over F.

Thus any vector in V is of the form $(x_1, x_2, 0, 0, 0, \ldots, 0)$ or $(0, 0, y_1, y_2, \ldots, y_6)$ and any vector in W is of the form $(x_1, x_2, x_3, 0, \ldots, 0)$ or $(0, 0, 0, y_1, y_2, \ldots, y_5)$. Hence no isomorphism can be sought between V and W in this set up.

This is one of the marked differences between the vector spaces and bivector spaces. Thus we have the following theorems, the proof of which is left for the reader to prove.

**THEOREM 2.1.2:** *Bivector spaces of same dimension defined over same fields need not in general be isomorphic.*

**THEOREM 2.1.3:** *Let $V = V_1 \cup V_2$ and $W = W_1 \cup W_2$ be any two bivector spaces of same dimension over the same field F. Then V and W are isomorphic as bivector spaces if and only if the vector space $V_1$ is isomorphic to $W_1$ and the vector space $V_2$ is isomorphic to $W_2$, that is dimension of $V_1$ is equal to dimension $W_1$ and the dimension of $V_2$ is equal to dimension $W_2$.*

**THEOREM 2.1.4:** *Let $V = V_1 \cup V_2$ be a bivector space over the field F. W any non empty set of V. $W = W_1 \cup W_2$ is a sub-bivector space of V if and only if $W \cap V_1 = W_1$ and $W \cap V_2 = W_2$ are sub spaces of $V_1$ and $V_2$ respectively.*

**DEFINITION 2.1.3:** *Let $V = V_1 \cup V_2$ and $W = W_1 \cup W_2$ be two bivector spaces defined over the field F of dimensions $p = m + n$ and $q = m_1 + n_1$ respectively.*

*We say the map $T: V \to W$ is a bilinear transformation (transformation bilinear) of the bivector spaces if $T = T_1 \cup T_2$ where $T_1 : V_1 \to W_1$ and $T_2 : V_2 \to W_2$ are linear transformations from vector spaces $V_1$ to $W_1$ and $V_2$ to $W_2$ respectively satisfying the following three rules:*

i. *$T_1$ is always a linear transformation of vector spaces whose first co ordinates are non-zero and $T_2$ is a linear transformation of the vector space whose last co ordinates are non zero.*

ii. *$T = T_1 \cup T_2$ '$\cup$' is just only a notational convenience.*

iii. *$T(v) = T_1(v)$ if $v \in V_1$ and $T(v) = T_2(v)$ if $v \in V_2$.*

Yet another marked difference between bivector spaces and vector spaces are the associated matrix of an operator of bivector spaces which has $m_1 + n_1$ rows and $m + n$ columns where dimension of V is m + n and dimension of W is $m_1 + n_1$ and T is a linear transformation from V to W. If A is the associated matrix of T then.

$$A = \begin{bmatrix} B_{m_1 \times m} & O_{n_1 \times m} \\ O_{m_1 \times n} & C_{n_1 \times n} \end{bmatrix}$$

where A is a $(m_1 + n_1) \times (m + n)$ matrix with $m_1 + n_1$ rows and m + n columns. $B_{m_1 \times m}$ is the associated matrix of $T_1 : V_1 \to W_1$ and $C_{n_1 \times n}$ is the associated matrix of $T_2 : V_2 \to W_2$ and $O_{n_1 \times m}$ and $O_{m_1 \times n}$ are non zero matrices.

Without loss of generality we can also represent the associated matrix of T by the bimatrix $B_{m_1 \times m} \cup C_{n_1 \times n}$.

***Example 2.1.6:*** Let $V = V_1 \cup V_2$ and $W = W_1 \cup W_2$ be two bivector spaces of dimension 7 and 5 respectively defined over the field F with dimension of $V_1 = 2$, dimension of $V_2 = 5$, dimension of $W_1 = 3$ and dimension of $W_2 = 2$. T be a linear transformation of bivector spaces V and W. The associated matrix of $T = T_1 \cup T_2$ where $T_1 : V_1 \to W_1$ and $T_2 : V_2 \to W_2$ given by

$$A = \begin{bmatrix} 1 & -1 & 2 & 0 & 0 & 0 & 0 & 0 \\ -1 & 3 & 0 & 0 & 0 & 0 & 0 & 0 \\ 0 & 0 & 0 & 2 & 0 & 1 & 0 & 0 \\ 0 & 0 & 0 & 3 & 3 & -1 & 2 & 1 \\ 0 & 0 & 0 & 1 & 0 & 1 & 1 & 2 \end{bmatrix}$$

where the matrix associated with $T_1$ is given by

$$A_1 = \begin{bmatrix} 1 & -1 & 2 \\ -1 & 3 & 0 \end{bmatrix}$$

and that of $T_2$ is given by

$$A_2 = \begin{bmatrix} 2 & 0 & 1 & 0 & 0 \\ 3 & 3 & -1 & 0 & 1 \\ 1 & 0 & 1 & 1 & 2 \end{bmatrix}$$

Thus $A = A_1 \cup A_2$ is the bimatrix representation for T for it saves both time and space. When we give the bimatrix representation. We can also call the linear operator T as linear bioperator and denote it as $T = T_1 \cup T_2$.

We call $T : V \to W$ a linear operator of both the bivector spaces if both V and W are of same dimension. So the matrix A associated with the linear operator T of the bivector spaces will be a square matrix. Further we demand that the spaces V and W to be only isomorphic bivector spaces. If we want to define eigen bivalues and eigen bivectors associated with T.

The eigen bivector values associated with are the eigen values associated with $T_1$ and $T_2$ separately. Similarly the eigen bivectors are that of the eigen vectors associated with $T_1$ and $T_2$ individually. Thus even if the dimension of the bivector spaces V and W are equal still we may not have eigen bivalues and eigen bivectors associated with them.

***Example 2.1.7:*** Let T be a linear operator of the bivector spaces – V and W. $T = T_1 \cup T_2$ where $T_1 : V_1 \to W_1$ dim $V_1$ = dim $W_1$ = 3 and $T_2 : V_2 \to W_2$ where dim $V_2$ = dim $W_2$ = 4. The associated matrix of T is

$$A = \begin{bmatrix} 2 & 0 & -1 & 0 & 0 & 0 & 0 \\ 0 & 1 & 0 & 0 & 0 & 0 & 0 \\ -1 & 0 & 3 & 0 & 0 & 0 & 0 \\ 0 & 0 & 0 & 2 & -1 & 0 & 6 \\ 0 & 0 & 0 & -1 & 0 & 2 & 1 \\ 0 & 0 & 0 & 0 & 2 & -1 & 0 \\ 0 & 0 & 0 & 6 & 1 & 0 & 3 \end{bmatrix}$$

The eigen bivalues and eigen bivectors can be calculated.

**DEFINITION 2.1.4:** *Let T be a linear operator on a bivector space V. We say that T is diagonalizable if $T_1$ and $T_2$ are diagonalizable where $T = T_1 \cup T_2$.*

*The concept of symmetric operator is also obtained in the same way, we say the linear operator $T = T_1 \cup T_2$ on the bivector space $V = V_1 \cup V_2$ is symmetric if both $T_1$ and $T_2$ are symmetric.*

**DEFINITION 2.1.5:** *Let $V = V_1 \cup V_2$. be a bivector space over the field F. We say $\langle\, ,\rangle$ is an inner product on V if $\langle\, ,\rangle = \langle\, ,\rangle_1 \cup \langle\, ,\rangle_2$ where $\langle\, ,\rangle_1$ and $\langle\, ,\rangle_2$ are inner products on the vector spaces $V_1$ and $V_2$ respectively.*

*Note that in $\langle\, ,\rangle = \langle\, ,\rangle_1 \cup \langle\, ,\rangle_2$ the '$\cup$' is just a conventional notation by default.*

**DEFINITION 2.1.6:** *Let $V = V_1 \cup V_2$ be a bivector space on which is defined an inner product $\langle\, ,\rangle$. If $T = T_1 \cup T_2$ is a linear operator on the bivector spaces V we say $T^*$ is an adjoint of T if $\langle T\alpha / \beta \rangle = \langle \alpha / T^* \beta \rangle$ for all $\alpha, \beta \in V$ where $T^* = T_1^* \cup T_2^*$ are $T_1^*$ is the adjoint of $T_1$ and $T_2^*$ is the adjoint of $T_2$.*

The notion of normal and unitary operators on the bivector spaces are defined in an analogous way. T is a unitary operator on the bivector space $V = V_1 \cup V_2$ if and only if $T_1$ and $T_2$ are unitary operators on the vector space $V_1$ and $V_2$ .

Similarly T is a normal operator on the bivector space if and only if $T_1$ and $T_2$ are normal operators on $V_1$ and $V_2$ respectively. We can extend all the notions on bivector spaces $V = V_1 \cup V_2$ once those properties are true on $V_1$ and $V_2$.

The primary decomposition theorem and spectral theorem are also true is case of bivector spaces. The only problem with bivector spaces is that even if the dimension of bivector spaces are the same and defined over the same field still they are not isomorphic in general.

Now we are interested in the collection of all linear transformation of the bivector spaces $V = V_1 \cup V_2$ to $W = W_1 \cup W_2$ where V and W are bivector spaces over the same field F.

We denote the collection of linear transformation by B-$Hom_F$ (V, W).

**THEOREM 2.1.5:** *Let V and W be any two bivector spaces defined over F. Then B-$Hom_F$ (V, W) is a bivector space over F.*

*Proof:* Given $V = V_1 \cup V_2$ and $W = W_1 \cup W_2$ be two bivector spaces defined over the field F. B-$Hom_F$ (V, W) = $\{T_1 : V_1 \to W_1\} \cup \{T_2 : V_2 \to W_2\} = Hom_F (V_1, W_1) \cup Hom_F (V_2, W_2)$. So clearly B- $Hom_F$ (V,W) is a bivector space as $Hom_F (V_1, W_1)$ and $Hom_F (V_2, W_2)$ are vector spaces over F.

**THEOREM 2.1.6:** *Let $V = V_1 \cup V_2$ and $W = W_1 \cup W_2$ be two bivector spaces defined over F of dimension m + n and $m_1 + n_1$ respectively. Then B-$Hom_F$ (V,W) is of dimension $mm_1 + nn_1$.*

*Proof*: Obvious by the associated matrices of T.

## 2.2 Introduction of S-bigroups and S-bivector spaces

In this section we introduce the concept of Smarandache bigroups. Bigroups were defined and studied in the year 1994 [40]. But till date Smarandache bigroups have not been defined. Here we define Smarandache bigroups and try to obtain several of the classical results enjoyed by groups. Further the study of Smarandache bigroups will throw several interesting features about bigroups in general.

**DEFINITION 2.2.1:** *Let (G, •, ∗) be a non-empty set such that $G = G_1 \cup G_2$ where $G_1$ and $G_2$ are proper subsets of G. (G, •, ∗) is called a Smarandache bigroup (S-bigroup) if the following conditions are true.*

- *i.* *(G, •) is a group.*
- *ii.* *(G, ∗) is a S-semigroup.*

***Example 2.2.1***: Let $G = \{g^2, g^4, g^6, g^8, g^{10}, g^{12} = 1\} \cup S(3)$ where S(3) is the symmetric semigroup. Clearly $G = G_1 \cup G_2$ where $G_1 = \{1, g^2, g^4, g^6, g^8, g^{10}\}$, group under '×' and $G_2 = S(3)$; S-semigroup under composition of mappings. G is a S-bigroup.

**THEOREM 2.2.1:** *Let G be a S-bigroup, then G need not be a bigroup.*

*Proof*: By an example consider the bigroup given in example 2.2.1. G is not a bigroup only a S-bigroup.

***Example 2.2.2***: Let $G = Z_{20} \cup S_5$; $Z_{20}$ is a S-semigroup under multiplication modulo 20 and $S_3$ is the symmetric group of degree 3. Clearly G is a S-bigroup; further G is not a bigroup.

**DEFINITION 2.2.2:** *Let $G = G_1 \cup G_2$ be a S-bigroup, a proper subset $P \subset G$ is said to be a Smarandache sub-bigroup of G if $P = P_1 \cup P_2$ where $P_1 \subset G_1$ and $P_2 \subset G_2$ and $P_1$ is a group or a S-semigroup under the operations of $G_1$ and $P_2$ is a group or S-semigroup under the operations of $G_2$ i.e. either $P_1$ or $P_2$ is a S-semigroup i.e. one of $P_1$ or $P_2$ is a group, or in short P is a S-bigroup under the operation of $G_1$ and $G_2$.*

**THEOREM 2.2.2:** *Let $G = G_1 \cup G_2$ be a S-bigroup. Then G has a proper subset H such that H is a bigroup.*

*Proof*: Given $G = G_1 \cup G_2$ is a S-bigroup. $H = H_1 \cup H_2$ where if we assume $G_1$ is a group then $H_1$ is a subgroup of $G_1$ and if we have assumed $G_2$ is a S-semigroup then $H_2$ is proper subset of $G_2$ and $H_2$ which is a subgroup of $G_2$.

Thus $H = H_1 \cup H_2$ is a bigroup.

**COROLLARY:** *If $G = G_1 \cup G_2$ a S-bigroup then G contains a bigroup.*

Study of S-bigroups is very new as in literature we do not have the concept of Smarandache groups we have only the concept of S-semigroups.

**DEFINITION 2.2.3:** *Let $G = G_1 \cup G_2$ be a S-bigroup we say G is a Smarandache commutative bigroup (S-commutative bigroup) if $G_1$ is a commutative group and every proper subset S of $G_2$ which is a group is a commutative group.*

*If both $G_1$ and $G_2$ happens to be commutative trivially G becomes a S-commutative bigroup.*

***Example 2.2.3:*** Let $G = G \cup S(3)$ where $G = \langle g \mid g^9 = 1\rangle$ Clearly G is a S-bigroup. In fact G is not a S-commutative bigroup.

***Example 2.2.4:*** Let $G = G \cup S(4)$ where $G = \langle g \mid g^{27} = 1\rangle$ and S(4) the symmetric semigroup which is a S-semigroup this S-bigroup is also non-commutative.

**DEFINITION 2.2.4:** *Let $G = G_1 \cup G_2$ be a S-bigroup where $G_1$ is a group and $G_2$ a S-semigroup we say G is a S-weekly commutative bigroup if the S-semigroup $G_2$ has at least one proper subset which is a commutative group.*

**DEFINITION 2.2.5:** *Let $A = A_1 \cup A_2$ be a k-bivector space. A proper subset X of A is said to be a Smarandache k-bivectorial space (S-k-bivectorial space) if X is a biset and*

*$X = X_1 \cup X_2 \subset A_1 \cup A_2$ where each $X_i \subset A_i$ is S-k-vectorial space.*

**DEFINITION 2.2.6:** *Let A be a k-vectorial bispace. A proper sub-biset X of A is said to be a Smarandache-k-vectorial bi-subspace (S-k-vectorial bi-subspace) of A if X itself is a S-k-vectorial subspace.*

**DEFINITION 2.2.7:** *Let V be a finite dimensional bivector space over a field K. Let $B = B_1 \cup B_2 = \{(x_1, \ldots, x_k, 0 \ldots 0)\} \cup \{(0,0, \ldots, 0, y_1 \ldots y_n)\}$ be a basis of V. We say B is a Smarandache basis (S-basis) of V if B has a proper subset A, $A \subset B$ and $A \neq \phi$, $A \neq B$ such that A generates a bisubspace which is bilinear algebra over K; that is W is the sub-bispace generated by A then W must be a k-bi-algebra with the same operations of V.*

**THEOREM 2.2.3:** *Let A be a k-bivectorial space. If A has a S-k-vectorial sub-bispace then A is a S-k-vectorial bispace.*

*Proof:* Straightforward by the very definition.

**THEOREM 2.2.4:** *Let V be a bivector space over the field K. If B is a S-basis of V then B is a basis of V.*

*Proof:* Left for the reader to verify.

**DEFINITION 2.2.8:** *Let V be a finite dimensional bivector space over a field K. Let $B = \{\upsilon_1, \ldots , \upsilon_n\}$ be a basis of V. If every proper subset of B generates a bilinear algebra over K then we call B a Smarandache strong basis (S-strong basis) for V.*

*Let V be any bivector space over the field K. We say L is a Smarandache finite dimensional bivector space (S-finite dimensional bivector space) over K if every S-basis has only finite number of elements in it.*

All results proved for bivector spaces can be extended by taking the bivector space $V = V_1 \cup V_2$ both $V_1$ and $V_2$ to be S-vector space. Once we have $V = V_1 \cup V_2$ to be a S-bivector space i.e. $V_1$ and $V_2$ are S-vector spaces, we see all properties studied for bivector spaces are easily extendable in case of S-bivector spaces with appropriate modifications.

Further the notion of Smarandache-k-linear algebra can be defined with appropriate modifications.

**Chapter Three**

# LINEAR BIALGEBRA, S-LINEAR BIALGEBRA AND THEIR PROPERTIES

This chapter introduces the notions of linear bialgebra and S-linear bialgebra and their properties and also give some applications. This chapter has nine section. In the first section we introduce the notion of basic properties of linear bialgebra which is well illustrated by several examples. Section two introduces the concept of linear bitransformation and linear bioperators. Section three introduces bivector spaces over finite fields. The concept of representation of finite bigroups is given in section four. Section five gives the application of bimatrix to bigraphs. Jordan biforms are introduced in section six. For the first time the application of bivector spaces to bicodes in given in section seven. The eighth section gives the best biapproximation and applies it to bicodes to find the closest sent message. The final section indicates about Markov bichains / biprocesss.

## 3.1 Basic Properties of Linear Bialgebra

In this section we for the first time introduce the notion of linear bialgebra, prove several interesting results and illustrate them also with example.

**DEFINITION 3.1.1:** *Let $V = V_1 \cup V_2$ be a bigroup. If $V_1$ and $V_2$ are linear algebras over the same field F then we say V is a linear bialgebra over the field F.*

*If both $V_1$ and $V_2$ are of infinite dimension vector spaces over F then we say V is an infinite dimensional linear bialgebra over F. Even if one of $V_1$ or $V_2$ is infinite dimension then we say V is an infinite dimensional linear bialgebra. If both $V_1$ and $V_2$ are finite dimensional linear algebra over F then we say $V = V_1 \cup V_2$ is a finite dimensional linear bialgebra.*

***Examples 3.1.1:*** Let $V = V_1 \cup V_2$ where $V_1$ = {set of all n × n matrices with entries from Q} and $V_2$ be the polynomial ring Q [x]. $V = V_1 \cup V_2$ is a linear bialgebra over Q and the linear bialgebra is an infinite dimensional linear bialgebra.

***Example 3.1.2:*** Let $V = V_1 \cup V_2$ where $V_1 = Q \times Q \times Q$ abelian group under '+', $V_2$ = {set of all 3 × 3 matrices with entries from Q} then $V = V_1 \cup V_2$ is a bigroup. Clearly V is a linear bialgebra over Q. Further dimension of V is 12 V is a 12 dimensional linear bialgebra over Q.

The standard basis is {(0 1 0), (1 0 0), (0 0 1)} ∪

$$\left\{\begin{array}{l}\begin{pmatrix}1&0&0\\0&0&0\\0&0&0\end{pmatrix},\begin{pmatrix}0&1&0\\0&0&0\\0&0&0\end{pmatrix},\begin{pmatrix}0&0&1\\0&0&0\\0&0&0\end{pmatrix},\begin{pmatrix}0&0&0\\1&0&0\\0&0&0\end{pmatrix},\begin{pmatrix}0&0&0\\0&1&0\\0&0&0\end{pmatrix},\\ \begin{pmatrix}0&0&0\\0&0&1\\0&0&0\end{pmatrix},\begin{pmatrix}0&0&0\\0&0&0\\1&0&0\end{pmatrix},\begin{pmatrix}0&0&0\\0&0&0\\0&1&0\end{pmatrix},\begin{pmatrix}0&0&0\\0&0&0\\0&0&1\end{pmatrix}\end{array}\right\}$$

***Example 3.1.3:*** Let $V = V_1 \cup V_2$ where $V_1$ is a collection of all 8 × 8 matrices over Q and $V_2$ = {the collection of all 3 × 2 matrices over Q}. Clearly V is a bivector space of dimension 70 and is not a linear bialgebra.

From this example it is evident that there exists bivector spaces which are not linear bialgebras.

Now if in a bivector space $V = V_1 \cup V_2$ one of $V_1$ or $V_2$ is a linear algebra then we call V as a semi linear bialgebra. So the vector space given in example 3.1.3 is a semi linear bialgebra.

We have the following interesting theorem.

**THEOREM 3.1.1:** *Every linear bialgebra is a semi linear bialgebra. But a semi linear bialgebra in general need not be a linear bialgebra.*

*Proof:* The fact that every linear bialgebra is a semi linear bialgebra is clear from the definition of linear bialgebra and semi linear bialgebra.

To prove that a semi linear bialgebra need not in general be a linear bialgebra. We consider an example. Let $V = V_1 \cup V_2$ where $V_1 = Q \times Q$ and $V_2$ the collection of all $3 \times 2$ matrices with entries from Q clearly $V = Q \times Q$ is a linear algebra of dimension 2 over Q and $V_2$ is not a linear algebra but only a vector space of dimension 6, as in $V_2$ we cannot define matrix multiplication. Thus $V = V_1 \cup V_2$ is not a linear bialgebra but only a semi linear bialgebra.

Now we have another interesting result.

**THEOREM 3.1.2:** *Every semi linear bialgebra over Q is a bivector space over Q but a bivector space in general is not a semi linear bialgebra.*

*Proof:* Every semi linear bialgebra over Q is clearly by the very definition a bivector space over Q. But to show a bivector space over Q in general is not a semi linear bialgebra we give an example. Let $V = V_1 \cup V_2$ where $V_1$ = {the set of all $2 \times 5$ matrices with entries from Q} and $V_2$ = {all polynomials of degree less than or equal to 5 with entries from Q}. Clearly both $V_1$ and $V_2$ are only vector spaces over Q and none of them are linear algebra. Hence $V = V_1 \cup V_2$ is only a bivector space and not a semi linear bialgebra over Q.

Hence the claim.

Now we define some more types of linear bialgebra. Let $V = V_1 \cup V_2$ be a bigroup. Suppose $V_1$ is a vector space over Q $\left(\sqrt{2}\right)$ and $V_2$ is a vector space over Q $\left(\sqrt{3}\right)$ (Q $\left(\sqrt{2}\right)$ and Q $\left(\sqrt{3}\right)$ are fields).

Then V is said to be a strong bivector space over the bifield Q $\left(\sqrt{3}\right) \cup$ Q $\left(\sqrt{2}\right)$. Similarly if $V = V_1 \cup V_2$ be a bigroup and if V is a linear algebra over F and $V_2$ is a linear algebra over K, K ≠ F K ∩ F ≠ F or K. i.e. if K ∪ F is a bifield then we say V is a strong linear bialgebra over the bifield.

Thus now we systematically give the definitions of strong bivector space and strong linear bialgebra.

**DEFINITION 3.1.2:** *Let $V = V_1 \cup V_2$ be a bigroup. $F = F_1 \cup F_2$ be a bifield. If $V_1$ is a vector space over $F_1$ and $V_2$ is a vector space over $F_2$ then $V = V_1 \cup V_2$ is called the strong bivector space over the bifield $F = F_1 \cup F_2$. If $V = V_1 \cup V_2$ is a bigroup and if $F = F_1 \cup F_2$ is a bifield. If $V_1$ is a linear algebra over $F_1$ and $V_2$ is a linear algebra over $F_2$. Then we say $V = V_1 \cup V_2$ is a strong linear bialgebra over the field $F = F_1 \cup F_2$.*

***Example 3.1.4:*** Consider the bigroup $V = V_1 \cup V_2$ where $V_1 = Q\left(\sqrt{2}\right) \times Q\left(\sqrt{2}\right)$ and $V_2$ = {set of all 3 × 3 matrices with entries from $Q\left(\sqrt{3}\right)$}.

Let $F = Q\left(\sqrt{2}\right) \cup Q\left(\sqrt{3}\right)$, Clearly F is a bifield. $V_1$ is a linear algebra over $Q\left(\sqrt{2}\right)$ and $V_2$ is a linear algebra over $Q\left(\sqrt{3}\right)$. So $V = V_1 \cup V_2$ is a strong linear bialgebra over the bifield $F = Q\left(\sqrt{2}\right) \cup Q\left(\sqrt{3}\right)$.

It is interesting to note that we do have the notion of weak linear bialgebra and weak bivector space.

***Example 3.1.5:*** Let $V = V_1 \cup V_2$ where $V_1 = Q \times Q \times Q \times Q$ be the group under '+' and $V_2$ = The set of all polynomials over the field $Q\left(\sqrt{2}\right)$. Now $V_1$ is a linear algebra over Q and $V_2$ is a linear algebra over $Q\left(\sqrt{2}\right)$. Clearly $Q \cup Q\left(\sqrt{2}\right)$ is a not a bifield as $Q \subseteq Q\left(\sqrt{2}\right)$. Thus $V = V_1 \cup V_2$ is a linear bialgebra over Q we call $V = V_1 \cup V_2$ to be a weak linear bialgebra over $Q \cup Q\left(\sqrt{2}\right)$. For $Q \times Q \times Q \times Q = V_1$ is not a linear algebra over $Q\left(\sqrt{2}\right)$. It is a linear algebra only over Q. Based on this now we give the definition of weak linear bialgebra over $F = F_1 \cup F_2$. where $F_1$ and $F_2$ are fields and $F_1 \cup F_2$ is not a bifield.

**DEFINITION 3.1.3:** *Let $V = V_1 \cup V_2$ be a bigroup. Let $F = F_1 \cup F_2$. Clearly F is not a bifield (For $F_1 \underset{\neq}{\subset} F_2$ or $F_2 \underset{\neq}{\subset} F_1$) but $F_1$ and $F_2$ are fields. If $V_1$ is a linear algebra over $F_1$ and $V_2$ is a linear algebra over $F_2$, then we call $V_1 \cup V_2$ a weak linear bialgebra over $F_1 \cup F_2$. One of $V_1$ or $V_2$ is not a linear algebra over $F_2$ or $F_1$ respectively.*

On similar lines we can define weak bivector space.

***Example 3.1.6:*** Let $V = V_1 \cup V_2$ be a bigroup. Let $F = Q \cup Q\left(\sqrt{2},\sqrt{3}\right)$ be a union of two fields, for $Q \subset Q\left(\sqrt{2},\sqrt{3}\right)$ so $Q \cup Q\left(\sqrt{2},\sqrt{3}\right) = Q\left(\sqrt{2},\sqrt{3}\right)$ so is a field. This cannot be always claimed, for instance if $F = Q\left(\sqrt{2}\right) \cup Q\left(\sqrt{3}\right)$ is not a field only a bifield. Let $V_1 = Q \times Q \times Q$ and $V_2 = Q\left(\sqrt{2},\sqrt{3}\right)[x]$. $V_1$ is a vector space, in fact a linear algebra over Q but $V_1$ is not a vector space over the field $Q\left(\sqrt{2},\sqrt{3}\right)$. $V_2$ is a vector space or linear algebra over Q or $Q\left(\sqrt{2},\sqrt{3}\right)$. Since $Q \subset Q\left(\sqrt{2},\sqrt{3}\right)$, we see V is a weak bivector space over $Q \cup Q\left(\sqrt{2},\sqrt{3}\right)$ infact V is a weak linear bialgebra over $Q \cup Q\left(\sqrt{2},\sqrt{3}\right)$.

***Example 3.1.7:*** Let $V = V_1 \cup V_2$ be a bigroup. Let $F = Q\left(\sqrt{2},\sqrt{3}\right) \cup Q\left(\sqrt{5},\sqrt{7},\sqrt{11}\right)$ be a bifield. Suppose $V_1 = Q\left(\sqrt{2},\sqrt{3}\right)[x]$ be a linear algebra over $Q\left(\sqrt{2},\sqrt{3}\right)$ and $V_2 = Q\left(\sqrt{5},\sqrt{7},\sqrt{11}\right) \times Q\left(\sqrt{5},\sqrt{7},\sqrt{11}\right)$ be a linear algebra over $Q\left(\sqrt{5},\sqrt{7},\sqrt{11}\right)$. Clearly $V = V_1 \cup V_2$ is a strong linear bialgebra over the bifield $Q\left(\sqrt{2},\sqrt{3}\right) \cup Q\left(\sqrt{5},\sqrt{7},\sqrt{11}\right)$.

Now we have the following results.

**THEOREM 3.1.3:** *Let $V = V_1 \cup V_2$ be a bigroup and V be a strong linear bialgebra over the bifield $F = F_1 \cup F_2$. V is not a linear bialgebra over F.*

*Proof:* Now we analyze the definition of strong linear bialgebra and the linear bialgebra. Clearly the strong linear bialgebra has no relation with the linear bialgebra or a linear bialgebra has no relation with strong linear bialgebra for linear bialgebra is defined over a field where as the strong linear bialgebra is defined over a bifield, hence no relation can ever be derived. In the similar means one cannot derive any form of relation between the weak linear bialgebra and linear bialgebra.

All the three notions, weak linear bialgebra, linear bialgebra and strong linear bialgebra for a weak linear bialgebra is defined over union of fields $F = F_1 \cup F_2$ where $F_1 \subset F_2$ or $F_2 \subset F_1$; $F_1$ and $F_2$ are fields; linear bialgebras are defined over the same field where as the strong linear bialgebras are defined over bifields. Thus these three concepts are not fully related.

It is important to mention here that analogous to weak linear bialgebra we can define weak bivector space and analogous to strong linear bialgebra we have the notion of strong bivector spaces.

***Example 3.1.8:*** Let $V = V_1 \cup V_2$ where $V_1$ = {set of all linear transformation of a n dimensional vector space over Q to a m dimensional vector space W over Q} and $V_2$ = {All polynomials of degree ≤ 6 with coefficients from R}.

Clearly $V = V_1 \cup V_2$ is a bigroup. V is a weak bivector space over $Q \cup R$.

***Example 3.1.9:*** Let $V = V_1 \cup V_2$ be a bigroup. $V_1$ = {set of all polynomials of degree less than or equal to 7 over $Q\left(\sqrt{2}\right)$} and $V_2$ = {set of all 5 × 2 matrices with entries from $Q\left(\sqrt{3},\sqrt{7}\right)$}. $V = V_1 \cup V_2$ is a strong bivector space over the bifield $F = Q\left(\sqrt{2}\right) \cup Q\left(\sqrt{3},\sqrt{7}\right)$. Clearly $V = V_1 \cup V_2$ is not a strong linear bialgebra over F.

Now we proceed on to define linear subbialgebra and subbivector space.

**DEFINITION 3.1.4:** *Let $V = V_1 \cup V_2$ be a bigroup. Suppose V is a linear bialgebra over F. A non empty proper subset W of V is said to be a linear subbialgebra of V over F if*

(1) *$W = W_1 \cup W_2$ is a subbigroup of $V = V_1 \cup V_2$.*
(2) *$W_1$ is a linear subalgebra over F.*
(3) *$W_2$ is a linear subalgebra over F.*

***Example 3.1.10:*** Let $V = V_1 \cup V_2$ where $V_1 = Q\times Q\times Q\times Q$ and $V_2$ ={set of all 4 × 4 matrices with entries from Q}. $V = V_1 \cup V_2$ is a bigroup under '+' V is a linear bialgebra over Q.

Now consider $W = W_1 \cup W_2$ where $W_1 = Q \times \{0\} \times Q \times \{0\}$ and $W_2$ = {collection of all upper triangular matrices with entries from Q}. $W = W_1 \cup W_2$ is a subbigroup of $V = V_1 \cup V_2$. Clearly W is a linear subbialgebra of V over Q.

**DEFINITION 3.1.5:** *Let $V = V_1 \cup V_2$ be a bigroup. $F = F_1 \cup F_2$ be a bifield. Let V be a strong linear bialgebra over F. A non empty subset $W = W_1 \cup W_2$ is said to be a strong linear subbialgebra of V over F if*

*(1) $W = W_1 \cup W_2$ is a subbigroup of $V = V_1 \cup V_2$.*
*(2) $W_1$ is a linear algebra over $F_1$ and*
*(3) $W_2$ is a linear algebra over $F_2$.*

***Example 3.1.11:*** Let $V = V_1 \cup V_2$ where $V_1$ = {Set of all 2 × 2 matrices with entries from Q $\left(\sqrt{2}\right)$} be a group under matrix addition and $V_2$ = {collection of all polynomials in the variable x over Q $\left(\sqrt{5},\sqrt{3}\right)$}. $V_2$ under polynomial addition is a group. Thus $V = V_1 \cup V_2$ is a strong linear bialgebra over the bifield F = Q $\left(\sqrt{2}\right) \cup$ Q $\left(\sqrt{3},\sqrt{5}\right)$.

Consider $W = W_1 \cup W_2$ where $W_1 = \left\{\begin{pmatrix} 0 & 0 \\ 0 & a \end{pmatrix} \middle| a \in Q\left(\sqrt{2}\right)\right\}$ and $W_2$ = {all polynomials of even degree i.e. p (x) = $\sum_{i=1}^{n} p_i x^{2i} n \in Q(\sqrt{3},\sqrt{5})$. Clearly, $W = W_1 \cup W_2$ is a subbigroup of $V = V_1 \cup V_2$. Clearly W is a strong linear subbialgebra of V over the bifield F = $Q\left(\sqrt{2}\right) \cup Q\left(\sqrt{5},\sqrt{3}\right)$.

**DEFINITION 3.1.6:** *Let $V = V_1 \cup V_2$ be a bigroup. Let V be a weak linear bialgebra over the union of fields $F_1 \cup F_2$ (i.e. $F_1$ is a subfield of $F_2$ or $F_2$ is a subfield of $F_1$) where $V_1$ is a linear algebra over $F_1$ and $V_2$ is a linear algebra over $F_2$. A non empty subset $W = W_1 \cup W_2$ is a weak linear sub bialgebra of V over $F_1 \cup F_2$ if*

*(1) W is a subbigroup of V.*
*(2) $W_1$ is a linear subalgebra of $V_1$.*
*(3) $W_2$ is a linear subalgebra of $V_2$.*

***Example 3.1.12:*** Let V = $V_1 \cup V_2$ be a bigroup. Let F = Q $\cup$ Q $\left(\sqrt{3},\sqrt{7}\right)$ be the union of fields. $V_1$ be a linear algebra over Q and $V_2$ be a linear algebra over Q $\left(\sqrt{3},\sqrt{7}\right)$; where $V_1$ = {all polynomials in the variable x with coefficient from Q} and $V_2$ = {3 × 3 matrices with entries from $Q\left(\sqrt{2},\sqrt{3}\right)$}. V = $V_1 \cup V_2$ is a weak linear bialgebra over Q $\cup$ $Q\left(\sqrt{3},\sqrt{7}\right)$. Consider W = $W_1 \cup W_2$, $W_1$ = {all polynomial of only even degree in the variable x with coefficients from Q} and

$$W_2 = \left\{ \begin{pmatrix} a & b & 0 \\ e & d & 0 \\ 0 & 0 & 0 \end{pmatrix} \middle| a,b,c,d \in Q\left(\sqrt{3},\sqrt{7}\right) \right\}.$$

Clearly W = $W_1 \cup W_2$ is a subbigroup of V = $V_1 \cup V_2$. Further W = $W_1 \cup W_2$ is a weak linear subbialgebra of V over F = Q $\cup$ Q $\left(\sqrt{2},\sqrt{3}\right)$.

It is very interesting to note that at times a weak linear bialgebra can have a subbigroup which happens to be just a linear subbialgebra over $F_1$ (if $F_1 \subset F_2$ i.e., $F_1$ is a subfield of $F_2$) and not a weak linear subbialgebra. We define such structures as special weak linear bialgebra.

**DEFINITION 3.1.7:** *Let V = $V_1 \cup V_2$ be a bigroup. F = $F_1 \cup F_2$ (where either $F_1$ is a subfield of $F_2$ or $F_2$ is a subfield of $F_1$) be union of fields and V is a weak linear bialgebra over F. Suppose W = $W_1 \cup W_2$ is a subbigroup of V = $V_1 \cup V_2$ and if W is a linear bialgebra only over $F_1$ (or $F_2$) (which* *ever is the subfield of the other). Then we call $V = V_1 \cup V_2$ a special weak linear bialgebra.*

***Example 3.1.13:*** Let $V = V_1 \cup V_2$ be a bigroup where $V_1 = Q \times Q \times Q$ and $V_2$ = {2 × 2 matrices with entries from $Q\left(\sqrt{2}\right)$}. $V = V_1 \cup V_2$ is a weak linear bialgebra over $Q \cup Q\left(\sqrt{2}\right)$.

Consider $W = W_1 \cup W_2$ where $W_1 = Q \times \{0\} \times Q$ and $W_2$ = {set of all 2 × 2 matrices over Q}. Then $W = W_1 \cup W_2$ is a linear bialgebra over Q. Clearly $W_2$ is note linear algebra over $Q\left(\sqrt{2}\right)$. We call $V = V_1 \cup V_2$ a special weak linear algebra over $Q \cup Q\left(\sqrt{2}\right)$.

**DEFINITION 3.1.8:** *Let $V = V_1 \cup V_2$ be a bigroup, $F = F_1 \cup F_2$ be a union fields (i.e., $F_1 = F$ or $F_2 = F$), V be the weak linear bialgebra over F. We say V is a strong special weak linear bialgebra if V has a proper subset W such that $W = W_1 \cup W_2$ is a subbigroup of V and there exists subfields of $F_1$ and $F_2$ say $K_1$ and $K_2$ respectively such that $K = K_1 \cup K_2$ is a bifield and W is a strong linear bialgebra over K.*

*Note:* Every weak linear bialgebra need not in general be a strong special weak linear bialgebra. We illustrate both by examples.

***Example 3.1.14:*** Let $V = V_1 \cup V_2$ where $V_1 = Q \times Q$ and $V_2$ = {Collection of all 5 × 5 matrices with entries from $Q\left(\sqrt{2}\right)$}.V is a weak linear bialgebra over $F= Q \cup Q\left(\sqrt{2}\right)$. Clearly V is never a strong special weak linear bialgebra.

Thus all weak linear bialgebras need not in general be a strong special weak linear bialgebras.

Thus all weak linear bialgebras need not in general be a strong special weak linear bialgebras. We will prove a theorem to this effect.

***Example 3.1.15:*** Let $V = V_1 \cup V_2$ be a bigroup. $V_1 = Q\left(\sqrt{2},\sqrt{3}\right)$ [x] be the polynomial ring over $Q\left(\sqrt{2},\sqrt{3}\right) = F_1$. Clearly $V_1$ is a abelian group under addition. Let $Q\left(\sqrt{2},\sqrt{3},\sqrt{7},\sqrt{11}\right) = F_2$ be the field. Let $V_2$ = {the set of all $3 \times 3$ matrices with entries from $F_2$}. $V_2$ is an abelian group under matrix addition. $V = V_1 \cup V_2$ is a weak linear bialgebra over $F_1 \cup F_2$.

Now consider the subbigroup, $W = W_1 \cup W_2$ where $W_1$ = {Set of all polynomials in x with coefficients from the field $Q\left(\sqrt{2}\right)$} and $W_2$ = {set of all $3 \times 3$ matrices with entries from $Q\left(\sqrt{3},\sqrt{7}\right)$}. $F = Q\left(\sqrt{2}\right) \cup Q\left(\sqrt{3},\sqrt{7}\right)$ is a bifield and the subbigroup $W = W_1 \cup W_2$ is a strong linear bialgebra over the bifield $F = Q\left(\sqrt{2}\right) \cup Q\left(\sqrt{3},\sqrt{7}\right)$. Thus V has a subset W such that W is a strong linear bialgebra so V is a strong special weak linear bialgebra.

Now we give condition for a weak linear bialgebra to be strong special weak linear bialgebra.

**THEOREM 3.1.4:** *Let $V = V_1 \cup V_2$ be a bigroup which has a proper subbigroup. Suppose V is a weak linear bialgebra over $F = F_1 \cup F_2$. V is not a strong special weak linear bialgebra if and only if one of $F_1$ or $F_2$ is a prime field.*

*Proof:* Let $V = V_1 \cup V_2$ be a bigroup and V be a weak linear bialgebra over the field $F = F_1 \cup F_2$. Suppose we assume $F_1$ is a prime field then clearly $F_2$ has $F_1$ to be its subfield. So $F = F_1 \cup F_2$ has no subset say $K = K_1 \cup K_2$ such that K is a bifield. Since F has no subset $K = K_1 \cup K_2$ such that K is a bifield, we see for no subbigroup W can be a strong linear bialgebra.

Conversely if $F = F_1 \cup F_2$ and if one of $F_1$ or $F_2$ is a prime field then there does not exist $K \subset F$, $K = K_1 \cup K_2$

such that K is a bifield, then V cannot have a subbigroup which is a strong linear bialgebra.

Now we give the conditions for the weak linear bialgebra to be a strong special weak linear bialgebra.

(1) $V = V_1 \cup V_2$ should have proper subbigroups.
(2) In $F = F_1 \cup F_2$ ($F = F_1$ or $F = F_2$) the subfield must not be a prime field and the extension field must contain some other non prime subfield other than the already mentioned subfield i.e. V should have subset which is a bifield.
(3) V should have a subbigroup which is a strong linear bialgebra over the bifield.

**THEOREM 3.1.5:** *Let $V = V_1 \cup V_2$ be a bigroup. $F = F_1 \cup F_2$ be the union of field. V be a weak linear bialgebra over F. If V has no proper subbigroup then*

*(1) V is not a strong special weak linear bialgebra over F.*
*(2) V is not a special weak linear bialgebra.*

The proof of the above theorem is left as an exercise for the reader.

## 3.2 Linear Bitransformation and Linear Bioperators

This section introduces the notion of bitransformation and bioperators. The notion of bidiagonlizable is also introduced and studied. We define inner biproduct on a bivector space leading to analogous results as in case of vector spaces.

The homomorphism T of a linear bialgebra $V = V_1 \cup V_2$ to another linear bialgebra $U = U_1 \cup U_2$ both U and V are defined over the same field F is defined to be a bitransformation of V to U if T satisfies the following conditions i.e., $T : V \to U$ is such that (1) $T = T_1 \cup T_2$ where '$\cup$' is just a symbol a convenience of notation, where $T_1 : V_1 \to U_1$ is a linear transformation and $T_2 : V_2 \to U_2$ is a linear transformation.

*Note:* If $\text{Hom}_{F_B}^{B}(V_1 \cup V_2, U_1 \cup U_2) = \text{Hom}_{F_B}^{B}(U, V)$ denotes the collection of all linear Bitransformation, it can be easily verified that $\text{Hom}_{F_B}^{B}(V, U)$ is a bivector space over F.

For we have with every transformation associated a m × n matrix and with every m × n matrix we have a transformation so we have with every bitransformation a bimatrix and every bimatrix represents a bitransformation so it can be easily proved $\text{Hom}_{F_B}^{B}(V, U)$ is a bivector space over F. To this end first we give some examples.

***Example 3.2.1:*** Let $V = V_1 \cup V_2$ where $V_1$ = {The collection of all m × n matrices with entries from Q} and $V_2$ = {The collection of all p × q matrices with entries from Q}. Clearly $V = V_1 \cup V_2$ is a bivector space over Q as $V_1$ and $V_2$ are vector spaces over Q.

***Example 3.2.2:*** Let $T : V \to U$ be a linear bitransformation of the linear bialgebra $V = V_1 \cup V_2$ where $V_1 = Q \times Q$ and $V_2 = Q_3 [ x ]$ = {all polynomials of degree less than or equal to three}. $U = U_1 \cup U_2$ where $U_1 = Q \times Q \times Q$ and $U_2$ = {all polynomials of degree less than or equal to five} = $Q_5 [ x ]$. U and V are bivector spaces over Q.

$$T = T_1 \cup T_2 : V \to U$$

is defined by $T_1 : V_1 \to U_1$ and $T_2 : V_2 \to U_2$ by

$$T_1 (x, y) = (x + y, x - y, y)$$
$$T_2 (\nu_1, \nu_2, \nu_3, \nu_4) = (\nu_1, \nu_2, \nu_3, \nu_4, \nu_1 + \nu_2 + \nu_3).$$

The related matrix with

$$T_1 \text{ is } \begin{bmatrix} 1 & 1 \\ 1 & -1 \\ 0 & 1 \end{bmatrix}, \text{ i.e., } T_1(x, y) = (x + y, x - y, y).$$

The matrix related with

$$T_2 = \begin{bmatrix} 1 & 0 & 0 & 0 \\ 0 & 1 & 0 & 0 \\ 0 & 0 & 1 & 0 \\ 0 & 0 & 0 & 1 \\ 1 & 1 & 1 & 0 \end{bmatrix}$$

$$T_2 (v_1, v_2, v_3, v_4) = (v_1, v_2, v_3, v_4, v_1 + v_2 + v_3).$$

Thus

$$T = \begin{bmatrix} 1 & 1 \\ 1 & -1 \\ 0 & 1 \end{bmatrix} \cup \begin{bmatrix} 1 & 0 & 0 & 0 \\ 0 & 1 & 0 & 0 \\ 0 & 0 & 1 & 0 \\ 0 & 0 & 0 & 1 \\ 1 & 1 & 1 & 0 \end{bmatrix}.$$

Thus to every linear bitransformation is associated with a bimatrix and to each bimatrix we have a linear bitransformation associated with it. Thus the collection of all bimatrices of the form T over Q is a bivector space over Q.

***Example 3.2.3:*** Let $V = V_1 \cup V_2$, where $V_1$ = {set of all 7 × 3 matrices with entries from Q} and $V_2$ = {set of all 4 × 4 matrices with entries from Q} Clearly $V = V_1 \cup V_2$ is a bivector space over Q. Dimension of V is 37 = 21 + 16.

***Example 3.2.4:*** Let $V = V_1 \cup V_2$, where $V_1$ = {set of all 6 × 6 matrices with entries from Q} and $V_2$ = {set of all 2 × 2 matrices with entries from R}. Clearly $V = V_1 \cup V_2$ is a linear bialgebra over Q. It is to be noted V is an infinite dimensional bivector space over Q. For $V_1$ is of dimension 36 and $V_2$ is of infinite dimension over Q. So V is an infinite dimensional linear bialgebra over Q.

***Example 3.2.5:*** Let $V = V_1 \cup V_2$ where $V_1$ = {set of all 6 × 6 matrices with entries from R} and $V_2$ = {set of all 2 × 2 matrices with entries from R}. $V = V_1 \cup V_2$ is a linear bialgebra over R.

Dimension of V as a linear bialgebra over R is finite, and dimension is 40 = 36 + 4 over R. If we assume $V = V_1 \cup V_2$ to be a linear bialgebra over Q then V is a linear bialgebra of infinite dimension over Q.

Thus we can analogously derive or prove the following theorems.

**THEOREM 3.2.1:** *Let $V = V_1 \cup V_2$ and $U = U_1 \cup U_2$ be finite dimensional linear bialgebras over the same field F. To every linear bitransformation $T: V \to U$ we have a bimatrix $M = M_1 \cup M_2$ with entries from F. $M_1$ a $n \times m$ matrix and $M_2$ a $p \times q$ matrix is associated. Conversely given any bimatrix $P = P_1 \cup P_2$ where we have $P_1$ to be a $n \times m$ matrix and $P_2$ to be a $p \times q$ matrix with entries from F we can always associate a linear bitransformation $T: V \to U$ for a suitable basis of V and U.*

The proof is left as an exercise for the reader as the working is straight forward.

**THEOREM 3.2.2:** *Let V and U be finite dimensional linear bialgebras over the field F. Let $Hom_F^B(V,U)$ be the collection of all linear bitransformations of V to U; then $Hom_F^B(V,U)$ is a bivector space over F which is finite dimensional.*

The proof of this theorem is left as an exercise to the reader.

Now we proceed on to define the notion of linear bioperator of linear bialgebras.

**DEFINITION 3.2.1:** *Let $V = V_1 \cup V_2$ be a linear bialgebra of finite dimension over the field F. A linear bitransformation $T : V \to V$ is said to be a linear bioperator of V.*

***Example 3.2.6:*** Let $V = V_1 \cup V_2$ where $V_1 = Q \times Q \times Q$ and $V_2 = Q^6 [ x ]$ = {all polynomials of degree less than or equal to 6 with entries from Q}. Clearly $V = V_1 \cup V_2$ is a finite dimensional linear bialgebra over Q. We see dimension of V is (3 + 7) i.e. 10.

Now let $T : V \to V$ be a linear bioperator defined by

$$T = (T_1 \cup T_2) : V \to V$$

$$T_1 (x, y, z) = (x + y, x - z, x + y + z)$$

and

$$T_2 (x_1, x_2, x_3, x_4, x_5, x_6, x_7)$$

$$= (x_1 + x_2, x_2 + x_3, x_3 + x_4, x_4 + x_5, x_5 + x_6, x_6 + x_7, x_7 + x_1).$$

$T = T_1 \cup T_2$ is a linear bioperator, the bimatrix associated with T, is given by $M = M_1 \cup M_2$ where

$$M_1 = \begin{bmatrix} 1 & 1 & 0 \\ 1 & 0 & -1 \\ 1 & 1 & 1 \end{bmatrix}$$

and

$$M_2 = \begin{bmatrix} 1 & 1 & 0 & 0 & 0 & 0 & 0 \\ 0 & 1 & 1 & 0 & 0 & 0 & 0 \\ 0 & 0 & 1 & 1 & 0 & 0 & 0 \\ 0 & 0 & 0 & 1 & 1 & 0 & 0 \\ 0 & 0 & 0 & 0 & 1 & 1 & 0 \\ 0 & 0 & 0 & 0 & 0 & 1 & 1 \\ 1 & 0 & 0 & 0 & 0 & 0 & 1 \end{bmatrix}.$$

Now we make the following observations:

1. We have intentionally defined the transformations from bivector space V to bivector space U as linear bitransformations as the linearity is preserved and if we say bilinear it would become mixed up with the bilinear forms and multilinear forms of the linear algebra / vector spaces.

2. When $V = V_1 \cup V_2$ is a finite dimensional linear bialgebra of dimension say (m, n) i.e. m + n then the linear bialgebra of all linear bioperator $\text{Hom}_F^B(V,V)$ is a $m \times m + n \times n$ dimensional linear bialgebra over F i.e. $\text{Hom}_F^B(V,V) \cong$ {all $m \times m$ matrices $\cup$ all $n \times n$ matrices for both entries are taken from F} so dim $\text{Hom}_F^B(V,V) = m^2 + n^2$.

***Example 3.2.7:*** Let $V = V_1 \cup V_2$ where $V_1 = Q \times Q \times Q \times Q$ and $V_2$ = {collection of all polynomials of degree less than or equal to 2 with coefficients from Q}. $V = V_1 \cup V_2$ is a linear bialgebra of dimension 4 + 3 = 7.
Now

$$\text{Hom}_Q^B(V,V) \cong$$

$$\left\{ \begin{pmatrix} a_{11} & \cdots & a_{14} \\ a_{21} & & a_{24} \\ \vdots & & \\ a_{41} & & a_{44} \end{pmatrix} \cup \begin{pmatrix} b_{11} & b_{12} & b_{13} \\ b_{21} & b_{22} & b_{23} \\ b_{31} & b_{32} & b_{33} \end{pmatrix} \middle| \, a_{ij}\, b_{ij}, \in Q \begin{array}{c} 1 \le i, j \le 4 \\ 1 \le i_1, j_1 \le 3 \end{array} \right\}$$

Clearly dimension of $\text{Hom}_Q^B(V,V)$ is $4^2 + 3^2$. $\text{Hom}_Q^B(V,V)$ is a linear bialgebra of dimension $4^2 + 3^2$.

It is still interesting and important to note that in case of linear operators of a vector space / linear algebra of dimension n; we have dimension of $\text{Hom}_F(V,V)$ is $n^2$. But in case of linear bialgebra / bivector space of dimension (m, n) i.e. m + n the dimension of the space of linear bioperators is $m^2 + n^2$ and not $(m + n)^2$. This is an important difference and thus it is not the dimension square as in case of vector spaces. One may wonder how addition and composition of maps take place in $\text{Hom}_F^B(V,V)$. If $T = T_1 \cup T_2$ and $S = S_1 \cup S_2$ be linear bioperator on V to V where $V = V_1 \cup V_2$. Then $T + S = (T_1 \cup T_2) + (S_1 \cup S_2) = (T_1 + T_2) \cup (S_1 + S_2)$

We know $T_1, S_1 : V_1 \to V_1$ and $T_2, S_2 : V_2 \to V_2$ so $T_1 + S_1$ and $T_2 + S_2$ are well defined linear operators from $V_1 \to V_1$ and $V_2 \to V_2$ respectively so $T + S : V \to V$ is a linear bioperator from V to V.

Now the composition of linear bioperator are defined by, if $T = T_1 \cup T_2$ and $S = S_1 \cup S_2$ be any two linear bioperators on V; Then $T \text{ o } S = (T_1 \cup T_2) \text{ o } (S_1 \cup S_2) = (T_1 \text{ o } S_1) \cup (T_2 \text{ o } S_2)$. We know $T_1, S_1 : V_1 \to V_1$ and $T_2, S_2 : V_2 \to V_2$ are linear operators so $T_1 \text{ o } S_1$ and $T_2 \text{ o } S_2$ are linear operators on $V_1$ and $V_2$ respectively.

Hence $T \text{ o } S = (T_1 \cup T_2) \text{ o } (S_1 \cup S_2) = (T_1 \text{ o } S_1) \cup (T_2 \text{ o } S_2)$ is a linear bioperator from V to V. Now what about scalar multiplication,

$$
\begin{aligned}
CT &= C\,(T_1 \cup T_2) \\
&= CT_1 \cup CT_2.
\end{aligned}
$$

$$
\begin{aligned}
(T + S)\,(\nu) & \\
&= (T + S)\,(\nu_1 \cup \nu_2),\ \nu \in V = V_1 \cup V_2 \\
&= (T_1 \cup T_2 + S_1 \cup S_2)\,(\nu_1 \cup \nu_2) \\
&= ((T_1 + S_1) \cup (T_2 + S_2))\,(\nu_1 \cup \nu_2) \\
&= (T_1 + S_1)\,(\nu_1) \cup (T_2 + S_2)\,(\nu_2) \\
&= (T_1\,(\nu_1) + S_1\,(\nu_1)) \cup (\,T_2\,(\nu_2) + S_2\,(\nu_2)) \\
&= \nu_1^1 \cup \nu_2^1\ .
\end{aligned}
$$

$(\because T_1\,(\nu_1) \text{ and } S_1\,(\nu_1) \in V_1 \text{ and } T_2\,(\nu_2) \text{ and } S_2\,(\nu_2) \in V_2)$

$$
\begin{aligned}
(T + S)\nu &= T\,(\nu) + S\,(\nu) \\
&= \nu^1.
\end{aligned}
$$

Thus it is easily verified that as in case of linear operators, linear bioperators also form a linear bialgebra over the same underlying field.

We define a linear bitransformation T : V → U (V and U linear bialgebras / bivector spaces over the same field F) to be invertible if there exists a linear bitransformation S from U into V such that S o T is the identity linear bioperator on V.

(1) This is possible if and only if T is 1 – 1 i.e. $T = T_1 \cup T_2$ then $T_i$ is 1 – 1 for i = 1, 2 that is $T(\nu) = T(\upsilon)$ implies $T(\nu_1 \cup \nu_2) = T(\upsilon_1 \cup \upsilon_2)$ i.e. $(T_1 \cup T_2)(\nu_1 \cup \nu_2) = (T_1 \cup T_2)(\upsilon_1 \cup \upsilon_2)$ i.e. $T_1(\nu_1) \cup T_2(\nu_2) = T_1(\upsilon_1) \cup T_2(\upsilon_2)$ i.e. if and only if $\nu_1 = \upsilon_1$ and $\nu_2 = \upsilon_2$.

(2) T is onto that is range of T is all of U i.e. range of $T_1$ is all of $U_1$ and range of $T_2$ is all of $U_2$. We call a linear bioperator T to be non-singular if $T(\nu) = 0$ implies $\nu = 0$ i.e. $\nu_1 = 0$ and $\nu_2 = 0$ $(\nu = \nu_1 \cup \nu_2)$ i.e. $T(\nu) = 0$ implies $T_1(\nu_1) = 0$ and $T_2(\nu_2) = 0$.

Now we proceed on to define the null space of a linear bitransformation $T = T_1 \cup T_2 : V \to U$. Null space of T = (null space of $T_1$) ∪ (null space of $T_2$) where $T_1: V_1 \to U_1$, a linear transformation of the linear algebra / vector space $V_1$ to $U_1$ and $T_2: V_2 \to U_2$ is a linear transformation of the linear algebra / vector space. Clearly null bispace of T is a subbispace or bisubspace or linear subbialgebra of $V = V_1 \cup V_2$.

If the linear bitransformation T is non singular than the null bispace of T is $\{0\} \cup \{0\}$. It can be proved analogous to linear operators.

**THEOREM 3.2.3:** *Let $V = V_1 \cup V_2$ and $U = U_1 \cup U_2$ be finite dimensional bivector spaces over the field F such that dim $V_i$ = dim $U_i$, i = 1, 2. If T is a linear bitransformation from V into U the following are equivalent:*

*(i) T is invertible*
*(ii) T is non singular*

(iii) *T is onto, that is range of T is* $U = U_1 \cup U_2$.

The proof is left as an exercise for the reader.

Now when we say for any linear bialgebra / bivector spaces U and V over F the dimension are identical we mean only dim $V_i$ = dim $U_i$ for i = 1, 2 where $U = U_1 \cup U_2$ and $V = V_1 \cup V_2$.

Thus it is very important to note that in case of finite dimensional linear bialgebras or bivector spaces over the same field F, if dim V = m + n and dim U = p + q even if m + n = p + q but p ≠ m or n and q ≠ m or n then we don't say the dimension of U and V are identical but they are same dimensional. So like in linear algebra the number of elements in the basis set of $U_1$ and $U_2$ equal to number of elements in the basis set of $V_1$ and $V_2$ does not imply the linear bialgebras or bivector spaces are identical but we say they are of same dimension.

***Example 3.2.8:*** Let $V = V_1 \cup V_2$ where $V_1$ = {set of all 2 × 2 matrices with entries from Q} and $V_2$ = { Q × Q × Q} be a bivector space over Q. Clearly V is a finite dimensional vector space over Q and dimension of V = 4 + 3 = 7. Consider $U = U_1 \cup U_2$ where $U_1$ = {Q × Q} and $U_2$ = {set of all polynomials of degree less than or equal to 4}. $U = U_1 \cup U_2$ is a finite dimensional bivector space over Q and dim U = 2 + 5 = 7. But we don't say though dim U = dim V. The dimensions are same but they are not bivector spaces of identical dimension.

***Example 3.2.9:*** Let $V = V_1 \cup V_2$ where $V_1$ = {set of all 2 × 3 matrices with entries from Q} and $V_2$ = {Q}. Clearly $V = V_1 \cup V_2$ is a finite dimensional bivector space and dim V = dim $V_1$ + dim $V_2$ = 6 + 1 = 7. Let $U = U_1 \cup U_2$ where $U_1$ = {Q × Q × Q} and $U_2$ = {set of all 2 × 2 matrices with entries from Q}. dim U = dim V = 7, but they have same dimension but not identical dimension over Q.

Now we proceed on to give an example of a identical dimension bivector space.

***Example 3.2.10:*** Let $V = V_1 \cup V_2$ where $V_1$ = {All polynomials of degree less than or equal to 3 with coefficient from Q} and $V_2 = \{Q \times Q\}$. $V = V_1 \cup V_2$ is a finite dimensional bivector space over Q; dimension of $(V = V_1 \cup V_2) = 4 + 2 = 6$. Consider the bivector space $U = U_1 \cup U_2$ where $U_1$ = {2 × 2 matrices with entries from Q} and $U_2 = \{ \begin{pmatrix} x \\ y \end{pmatrix} \mid x, y \in Q\}$ be the column vector. $U = U_1 \cup U_2$ is a finite dimensional bivector space of dimension 4 + 2 = 6. The bivector spaces are infact identical.

*Note:* The important theorem in finite dimensional vector spaces over same field F; i.e. all vector spaces of same dimension are isomorphic is not true in case of finite dimensional bivector spaces of same dimension. Only bivector spaces of identical dimensions are isomorphic defined over the same field. This is the marked difference between vector spaces / linear algebras and bivector spaces / linear bialgebras.

Now we define two types of subspaces.

**DEFINITION 3.2.2:** *Let $V = V_1 \cup V_2$ be a bivector space over F. Let $W = W_1 \cup W_2$ be a bisubspace of V. If both $W_1$ and $W_2$ are proper subsets of V then we say W is a bisubspace or subbispace of V. If one of $W_1$ or $W_2$ is not a proper subset of $V_1$ or $V_2$ (respectively) then we call W a semi subbispace of V.*

***Example 3.2.11:*** Let $V = V_1 \cup V_2$ be a bivector space over Q, where $V_1$ = {set of all 3 × 3 matrices with entries from Q} and $V_2 = Q \times Q \times Q \times Q$. Clearly $V = V_1 \cup V_2$ is a bivector space over Q. Let $W = W_1 \cup W_2$ where $W_1$ = {3 × 3 upper triangular matrices with entries from Q} and $W_2 = Q \times \{0\} \times Q \times \{0\}$. Clearly $W = W_1 \cup W_2$ is a subbispace of V.

Let $W^1 = W_1^1 \cup W_2^1$ where $W_1^1 = V_1$ and $W_2 = Q \times \{0\} \times \{0\} \times \{0\}$. $W^1$ is called the semi bisubspace of $V = V_1 \cup$ $V_2$. As our main motivation is only study of linear bialgebra and not much of bivector spaces we just give an outline of some of the definitions, its properties and related results. Let T : V $\rightarrow$ V be a linear bioperator on $V = V_1 \cup V_2$ a bivector space over F.

Now we define the characteristic bivalue of T to be a scalar $C = C_1 \cup C_2$ in F, such that there is a non zero vector $\alpha = \alpha_1 \cup \alpha_2$ in $V_1 \cup V_2$ with $T\alpha = C\,\alpha$ i.e. $T(\alpha_1 \cup \alpha_2) = (C_1 \cup C_2)(\alpha_1 \cup \alpha_2)$ i.e. $T_1\,\alpha_1 \cup T_2\,\alpha_2 = C_1\,\alpha_1 \cup C_2\,\alpha_2$. $C = C_1 \cup C_2$ ($C_1$, $C_2$ in F) is called the characteristic bivalue of $T = T_1 \cup T_2$ then:

1. Any $\alpha = \alpha_1 \cup \alpha_2$ such that $T\alpha = C\alpha$ is called a characteristic bivector of T associated with the characteristic bivalue $C = C_1 \cup C_2$.
2. The collection of all $\alpha = \alpha_1 \cup \alpha_2$ such that $T\alpha = C\alpha$ is called the characteristic bispace associated with $C = C_1 \cup C_2$.Characteristic bivalues are also called as characteristic biroots, latent biroots, eigen bivalues, proper bivalues or spectral bivalues. If T is any linear bioperator and $C = C_1 \cup C_2$ is any biscalar the set of bivectors $\alpha = \alpha_1 \cup \alpha_2$ in $V = V_1 \cup V_2$ such that $T\alpha = C\alpha$ is a subbispace of V.

It can be easily proved by any interested reader that the bisubspace is the null bispace of the linear bioperator $T - CI = (T_1 - C_1\, I^1) \cup (T_2 - C_2\, I^2)$ we call $C = C_1 \cup C_2$ is a characteristics bivalue of $T = T_1 \cup T_2$, if this subbispace is different from the zero subbispace i.e. if $(T - CI) = (T_1 \cup T_2 - C_1 \cup C_2\,(I^1 \cup I^2)) = (T_1 - C_1\, I^1) \cup (T_2 - C_2\, I^2)$ fails to be one to one.

It is left as an exercise for the reader to prove the following theorem.

**THEOREM 3.2.4:** *Let $T = T_1 \cup T_2$ be a linear bioperator or a finite dimensional bivector space $V = V_1 \cup V_2$, let $C = C_1 \cup C_2$ be a scalar. The following are equivalent*

(i) $C = C_1 \cup C_2$ is a characteristic bivalue of $T = T_1 \cup T_2$.

(ii) the bioperator $(T - CI) = (T_1 \cup T_2 - (C_1 \cup C_2)(I^1 \cup I^2) = (T_1 - C_1 I^1) \cup (T_2 - C_2 I^2)$ is singular.

(iii) $\det (T - CI) = 0$ i.e. $\det (T_1 \cup T_2 - (C_1 \cup C_2)(I^1 \cup I^2)) = \det (T_1 - C_1 I^1) \cup \det (T_2 - C_2 I^2) = 0 \cup 0$.

Now we proceed on to discuss it in case of bimatrix as every bimatrix and a linear bioperator have a one-to-one relation. So now we can give the characteristic bivalue of a square bimatrix and a mixed square bimatrix.

**DEFINITION 3.2.3:** *Let $A = A_1 \cup A_2$ be a square $n \times n$ bimatrix over a field F, a characteristic bivalue of A in F is a scalar $C = C_1 \cup C_2$ such that the matrix $(A - CI) = (A_1 - C_1 I^1) \cup (A_2 - C_2 I^2)$ is singular.*

**DEFINITION 3.2.4:** *Let V be a bivector space over the field F and let $T = T_1 \cup T_2$ be a linear bioperator on $V = V_1 \cup V_2$. A characteristic bivalue of $T = T_1 \cup T_2$ is a scalar $C = C_1 \cup C_2$ in F such that there is a non zero bivector $\alpha = \alpha_1 \cup \alpha_2$ in V with $T\alpha = c\alpha$ i.e.*

$$T(\alpha_1 \cup \alpha_2) = (T_1 \cup T_2)(\alpha_1 \cup \alpha_2)$$

$$T_1 \alpha_1 \cup T_2 \alpha_1 = C_1 \alpha_1 \cup C_2 \alpha_2.$$

*If $C = C_1 \cup C_2$ is a characteristic bivalue of the bioperator $T = T_1 \cup T_2$ then*

*(a) any $\alpha = \alpha_1 \cup \alpha_2$ such that $T\alpha = C\alpha$ i.e. $T_1 \alpha_1 \cup T_2 \alpha_2 = C_1 \alpha_1 \cup C_2 \alpha_2$ is called the characteristic bivector of $T = T_1 \cup T_2$ associated with the characteristic bivalue $C = C_1 \cup C_2$.*

*(b) The collection of all $\alpha = \alpha_1 \cup \alpha_2$ in Vsuch that $T\alpha = T_1 \alpha_1 \cup T_2 \alpha_2 = C_1 \alpha_1 \cup C_2 \alpha_2$ is called the characteristic bispace associated with the scalar $C = C_1 \cup C_2$.*

(Note a bimatrix $A = A_1 \cup A_2$ is singular (or non invertible) if and only if both $A_1$ and $A_2$ and singular (non invertible), $C = C_1 \cup C_2$ is the characteristic bivalue of $A = A_1 \cup A_2$ if and only if det $(A - CI) = 0$ i.e. det $(A_1 = C_1I^1) \cup$ det $(A_2 = C_2I^2) = 0 \cup 0$. $(I = I^1 \cup I^2)$ or equivalently if and only if det $(CI - A) = 0$ i.e. det $(C_1 I^1 - A^1) \cup$ det $(C_2 I^2 - A_2) = 0 \cup 0$, we form the bimatrix $(xI - A) = (x_1 I^1 - A_1) \cup (x_2 I^2 - A_2)$ with bipolynomial entries, consider the bipolynomial $f(x_1, x_2) = \det (x_1 I^1 - A_1) \cup \det (x_2 I^2 - A_2)$. Clearly the characteristic bivalues of $A = A_1 \cup A_2$ in F are just the scalars, $C = C_1 \cup C_2$ in F such that $f(C) = f(C_1, C_2) = 0$ for this reason f is called the characteristic bipolynomial of $A = A_1 \cup A_2$. It is important to note that $f(C) = f(C_1) \cup f_2(C_2)$ is a monic bipolynomial which has degree exactly (n, n).

We illustrate this by the following example.

***Example 3.2.12:*** Let $T = T_1 \cup T_2$ be a linear bioperator on the bivector space $V = (Q \times Q \times Q) \cup$ {all polynomials of degree less than or equal to two} which is represented in the standard ordered basis by the bimatrix A.

$$A = A_1 \cup A_2 = \begin{bmatrix} 0 & 1 & 0 \\ 2 & -2 & 2 \\ 2 & -3 & 2 \end{bmatrix} \cup \begin{bmatrix} 3 & 1 & -1 \\ 2 & 2 & -1 \\ 2 & 2 & 0 \end{bmatrix}.$$

The characteristic bipolynomial for $T = T_1 \cup T_2$ (or for $A = A_1 \cup A_2$) is det $(xI - A) = \det (x_1 I^1 - A_1) \cup \det (x_1 I^2 - A_2)$

$$= \begin{vmatrix} x_1 & -1 & 0 \\ -2 & x_1+2 & -2 \\ -2 & 3 & x_1-2 \end{vmatrix} \cup \begin{vmatrix} x_2-3 & -1 & 1 \\ -2 & x_2-2 & 1 \\ -2 & -2 & x_2 \end{vmatrix}$$

$$= x_1^3 \cup x_2^3 - 5x_2^2 + 8x_2 - 4 .$$

The characteristic bivalues of $A = A_1 \cup A_2$ are $(0, 0, 0) \cup (1, 2, 2)$.

We give another example.

***Example 3.2.13:*** Let $T = T_1 \cup T_2$ be a linear bioperator on the bivector space $V = Q \times Q \cup \{ \begin{bmatrix} a \\ b \end{bmatrix}$ - a column vector with entries from Q} over Q. Let the bioperator be represented by the bimatrix $A = A_1 \cup A_2$ in the standard basis.

$$A = \begin{bmatrix} 0 & -1 \\ 1 & 0 \end{bmatrix} \cup \begin{bmatrix} 1 & -1 \\ 2 & 2 \end{bmatrix}.$$

The characteristic bipolynomial for T (or for $A = A_1 \cup A_2$) is $\det (xI - A) = \det (x_1 I^1 - A_1) \cup \det (x_2 I^2 - A_2)$

$$= \begin{vmatrix} x_1 & 1 \\ -1 & x_1 \end{vmatrix} \cup \begin{vmatrix} x_2 - 1 & 1 \\ -2 & x_2 - 2 \end{vmatrix}$$

$$= \left(x_1^2 + 1\right) \cup (x_2^2 - 3x_2 + 4).$$

Since this bipolynomial has no real roots $T = T_1 \cup T_2$ has no characteristic bivalues. Thus we have illustrated how characteristic bivalues are determined in case of bivector spaces / bimatrices.

Now we illustrate this by an example how the characteristic bivalues look in case of mixed square bimatrices.

***Example 3.2.14:*** Let T be a linear bioperator on the bivector space $V = (V_1 \cup V_2) = (Q \times Q) \cup$ (set of all polynomials of degree less than or equal to 2) over Q.

The mixed square bimatrix associated with T relative to the standard basis is given by $A = A_1 \cup A_2$.

$$= \begin{bmatrix} 0 & -1 \\ 1 & 0 \end{bmatrix} \cup \begin{bmatrix} 3 & 1 & -1 \\ 2 & 2 & -1 \\ 2 & 2 & 0 \end{bmatrix}.$$

The related characteristic bipolynomial for the bimatrix A is

$$\begin{vmatrix} x_1 & 1 \\ -1 & x_1 \end{vmatrix} \cup \begin{vmatrix} x_2 - 3 & -1 & 1 \\ -2 & x_2 - 2 & 1 \\ -2 & -2 & x_2 \end{vmatrix}$$

$$= (x_1^2 + 1) \cup (x_2 - 1)\,(x_2 - 2)^2 .$$

Clearly only one set of characteristic bivalues does not exist for the polynomial $x_1^2 + 1$ has no real roots. The other set of characteristic bivalues are 1, 2 and 2.

Thus we see in case of characteristic bipolynomials, three possibilities occur.

(i) The characteristic bivalues does not exist for both $T_1$ and $T_2$.
(ii) Characteristic bivalues exist for one of $T_1$ or $T_2$.
(iii) Characteristic bivalues exists for both $T_1$ and $T_2$.

In case of (ii) we give a new definition.

**DEFINITION 3.2.5:** *Let $T = T_1 \cup T_2$ be a linear bioperator on $V = V_1 \cup V_2$. If the characteristic bivalues exists only for $T_1$ or $T_2$ then we say the linear bioperator has semi characteristic bivalues.*

Now we proceed on to define the notion of bidiagonalizable and semi bidiagonalizable in case of linear bioperators.

**DEFINITION 3.2.6:** *Let $T = T_1 \cup T_2$ be a linear bioperator on $V = V_1 \cup V_2$, V a finite dimensional bivector space. We say $T = T_1 \cup T_2$ is bidiagonalizable if there exists a basis* *for $V = V_1 \cup V_2$ each vector of which is the characteristic bivector of $T = T_1 \cup T_2$.*

We will illustrate the method of finding the characteristic bivectors for a linear bioperator.

***Example 3.2.15:*** Let $T = T_1 \cup T_2$ be a linear bioperator on the bivector space $V = V_1 \cup V_2$. Let the bimatrix associated with T be given by

$$\begin{bmatrix} 5 & -6 & -6 \\ -1 & 4 & 2 \\ 3 & -6 & -4 \end{bmatrix} \cup \begin{bmatrix} -1 & 0 & 0 \\ 2 & 1 & 0 \\ 0 & 1 & 4 \end{bmatrix}.$$

The characteristic bipolynomial associated with T is

$$\begin{vmatrix} x_1 - 5 & 6 & 6 \\ 1 & x_1 - 4 & 2 \\ -3 & 6 & x_1 + 4 \end{vmatrix} \cup \begin{vmatrix} x_2 + 1 & 0 & 0 \\ -2 & x_2 - 1 & 0 \\ 0 & -1 & x_2 - 4 \end{vmatrix}$$

$= \{(x_1 - 2)^2 (x_1 - 1)\} \cup \{(x_2 + 1) (x_2 - 1) (x_2 - 4)\}$.

The associated characteristic bivalues of the bioperator is $\{2, 2, 1\} \cup \{1, -1, 4\}$.

Now we calculate the related characteristic bivectors related to $\{2, 1\}$. Consider $[A_1 - 2I] \cup [A_2 - I]$.

Now we find the characteristic bivector associated with the bivalue $\{2, 1\}$.

$$\begin{bmatrix} 3 & -6 & -6 \\ -1 & 2 & 2 \\ 3 & -6 & -6 \end{bmatrix} \begin{bmatrix} \alpha_1^1 \\ \alpha_2^1 \\ \alpha_3^1 \end{bmatrix} \cup \begin{bmatrix} -2 & 0 & 0 \\ 2 & 0 & 0 \\ 0 & 1 & 3 \end{bmatrix} \begin{bmatrix} \alpha_1^2 \\ \alpha_2^2 \\ \alpha_3^2 \end{bmatrix}$$

This gives a pair of bivectors for (2, 1) viz {(2, 1 0) ∪ (0, –3, 1)} and (2 0 1) ∪ (0 –3 1).

Now the bimatrix related with the bivalues for {2, –1} is:

$$\begin{bmatrix} 3 & -6 & -6 \\ -1 & 2 & 2 \\ 3 & -6 & -6 \end{bmatrix} \begin{bmatrix} \alpha_1^1 \\ \alpha_2^1 \\ \alpha_3^1 \end{bmatrix} \cup \begin{bmatrix} 0 & 0 & 0 \\ 2 & 2 & 0 \\ 0 & 1 & 5 \end{bmatrix} \begin{bmatrix} \alpha_1^2 \\ \alpha_2^2 \\ \alpha_3^2 \end{bmatrix}$$

The associated bivector is {(2, 1, 0) (2 0 1)} ∪ {(5, –5, 1)}. Now for the bivalue (1, 4), we have the following bimatrix.

$$\begin{bmatrix} 4 & -6 & -6 \\ -1 & 3 & 2 \\ 3 & -6 & -5 \end{bmatrix} \begin{bmatrix} \alpha_1^1 \\ \alpha_2^1 \\ \alpha_3^1 \end{bmatrix} \cup \begin{bmatrix} -5 & 0 & 0 \\ 2 & -3 & 0 \\ 0 & 1 & 0 \end{bmatrix}$$

The related bivector is given by {(3, –1, 3) ∪ (0 0 1)}.
Thus we have the bivector {(3, –1, 3), (2, 1, 0), (2, 0, 1)} ∪ {(0 –3, 1), (0 0 1), (5, –5, 1)} associated with the linear bioperator T.

Now we find whether the bimatrix $A = A_1 \cup A_2$ associated with T is bidiagonalizable. Since the characteristic bivectors associated with the linear bioperator are linearly independent; we see the given linear bioperator is bidiagonalizable.

For

$$P_1 = \begin{bmatrix} 3 & 2 & 2 \\ -1 & 1 & 0 \\ 3 & 0 & 1 \end{bmatrix}$$

and

$$P_2 = \begin{bmatrix} 0 & 5 & 0 \\ -3 & -5 & 0 \\ 1 & 1 & 1 \end{bmatrix}$$

are such that

$$P_1 A_1 P_1^{-1} \cup P_2 A_2 P_2^{-1} = \begin{bmatrix} 1 & 0 & 0 \\ 0 & 2 & 0 \\ 0 & 0 & 2 \end{bmatrix} \cup \begin{bmatrix} -1 & 0 & 0 \\ 0 & 1 & 0 \\ 0 & 0 & 4 \end{bmatrix}$$

Thus if $P = P_1 \cup P_2$, then

$$\begin{aligned} P A P^{-1} &= (P_1 \cup P_2)(A_1 \cup A_2)(P_1 \cup P_2)^{-1} \\ &= P_1 A_1 P_1^{-1} \cup P_2 A_2 P_2^{-1} \\ &= D_1 \cup D_2. \end{aligned}$$

Now analogous to symmetric bilinear forms we define bisymmetric bilinear biforms.

**DEFINITION 3.2.7:** *Let $f = f_1 \cup f_2$ be a bilinear biform on the bivector space $V = V_1 \cup V_2$. We say that f is bisymmetric bilinear biform if $f(\alpha, \beta) = f(\beta\ \alpha)$ for all vectors $\alpha, \beta \in V_1 \cup V_2$ with $\alpha = \alpha_1 \cup \alpha_2$, $\beta = \beta_1 \cup \beta_2$ : $(f_1 \cup f_2)(\alpha_1 \cup \alpha_2, \beta_1 \cup \beta_2) = (f_1 \cup f_2)(\beta_2 \cup \beta_1, \alpha_1 \cup \alpha_2)$.*

$$f_1(\alpha_1, \beta_1) \cup f_2(\alpha_2, \beta_2) = f_1(\beta_1, \alpha_1) \cup f_2(\beta_2, \alpha_2).$$

*If $V = V_1 \cup V_2$ is finite dimensional, the bilinear biform $f = f_1 \cup f_2$ is bisymmetric if and only if its bimatrix $A = A_1 \cup A_2$ in some ordered basis is bisymmetric.*

*[i.e., $A^t = A_1^t \cup A_2^t = A = A1 \cup A2$ (A1t = A1, At2= A2)]*

$$\begin{aligned} f(X\ Y) &= (f_1 \cup f_2)(X_1 \cup X_2, Y_1 \cup Y_2) \\ &= f_1(X_1, Y_1) \cup f_2(X_2\ Y_2) \\ &= X_1^t A_1 Y_1 \cup X_2^t A_2 Y_2. \end{aligned}$$

*If $f = f_1 \cup f_2$ is a bisymmetric bilinear biform, the biquadratic biform associated with $f = f_1 \cup f_2$ is the function $q = q_1 \cup q_2$ from $V_1 \cup V_2$ onto $F \cup F$ defined by $q(\alpha) = (q_1 \cup q_2)(\alpha_1 \cup \alpha_2)$.*

*If $f_i(\alpha_i\ \alpha_i) > 0$ if $\alpha_I \neq 0$ $i = 1, 2$. Then the bilinear biform is called bipositive definite. Results as in case of bilinear forms can also be derived in case of bilinear biforms.*

*Throughout this section the bivector space is defined over the field F. A bilinear biform* $f = f_1 \cup f_2$ *on* $V = V_1 \cup V_2$ *is called skew bisymmetric if*

$$f(\alpha, \beta) = -f(\beta, \alpha) \text{ i.e. } (f_1 \cup f_2)(\alpha_1 \cup \alpha_2, \beta_1 \cup \beta_2)$$
$$= -(f_1 \cup f_2)(\alpha_1 \cup \alpha_2, \beta_1 \cup \beta_2)$$
$$\text{i.e. } f_1(\alpha_1, \beta_1) \cup f_2(\alpha_2, \beta_2)$$
$$= -f_1(\beta_1, \alpha_1) \cup -f_2(\beta_2, \alpha_2)$$

*for all* $\alpha = \alpha_1 \cup \alpha_2$, $\beta = \beta_1 \cup \beta_2$ *in* $V = V_1 \cup V_2$.

*Suppose* $f = f_1 \cup f_2$ *is any bilinear biform on* $V = V_1 \cup V_2$.

*If we let*

$$g(\alpha,\beta) = \frac{1}{2}\left[ f(\alpha,\beta) + f(\beta\alpha)\right]$$

$$h(\alpha,\beta) = \frac{1}{2}\left[ f(\alpha,\beta) - f(\beta\alpha)\right]$$

$g = g_1 \cup g_2$ *and* $h = h_1 \cup h_2$, $\alpha = \alpha_1 \cup \alpha_2$, $\beta = \beta_1 \cup \beta_2$ *in* $V = V_1 \cup V_2$ *that* $g = g_1 \cup g_2$ *is a bisymmetric bilinear biform on* $V = V_1 \cup V_2$ *and* $h = h_1 \cup h_2$ *is a skew bisymmetric bilinear biform on* $V = V_1 \cup V_2$.

*Clearly if* $f = g + h$ *i.e.* $f = f_1 \cup f_2 = (g_1 \cup g_2) + (h_1 \cup h_2)$
$= (g_1 + h_1) \cup (g_2 + h_2)$.

Further more, this expression for $f = f_1 \cup f_2$ as the sum of a bisymmetric and a skew- bisymmetric biform is unique as in case of symmetric and skew symmetric form. Thus the space. $L^B(V, V, F) = L^B(V_1 \cup V_2, V_1 \cup V_2, F)$ is the direct sum of the subbispaces of bisymmetric biforms and skew bisymmetric biforms.

Now we can give the analogous results on bimatrices. If $A = A_1 \cup A_2$ is a bimatrix, A is skew bisymmetric if

$$A^t = -A$$

i.e. $(A_1 \cup A_2)^t = -(A_1 \cup A_2)$

i.e. $A_1^t \cup A_2^t = -(A_1 \cup A_2)$.

All analogous results can be derived as in case of symmetric bilinear forms. Now we proceed onto define the notion of an inner biproduct on a bivector space $V = V_1 \cup V_2$ over F.

**DEFINITION 3.2.8:** *Let F be a field of real numbers or complex numbers and $V = V_1 \cup V_2$ be a bivector space over F. An inner biproduct on V is a bifunction $( / ) = ( / )_1 \cup ( / )_2$ which assigns to each ordered pair of bivectors $\alpha = \alpha_1 \cup \alpha_2$, $\beta = \beta_1 \cup \beta_2$ in $V = V_1 \cup V_2$ $\beta_i, \alpha_i \in V_i$ ( i = 1, 2 ) a pair of scalars $(\alpha / \beta) = (\alpha_1 / \beta_1)_1 \cup (\alpha_2 / \beta_2)_2$ in F in such a way that for all $\alpha, \beta, \gamma$ in $V = V_1 \cup V_2$ and all scalars $C = C_1 \cup C_2$.*

(i) $(\alpha + \beta / \gamma) = (\alpha_1 \cup \alpha_2 + \beta_1 \cup \beta_2 \mid (\gamma_1 \cup \gamma_2))$
$= \{(\alpha_1 + \beta_1) \cup (\alpha_2 + \beta_2)\} \mid (\gamma_1 \cup \gamma_2))$
$= (\alpha / \gamma) + (\beta / \gamma)$
$= (\alpha_1 \cup \alpha_2 \mid \gamma_1 \cup \gamma_2) + (\beta_1 \cup \beta_2 \mid \gamma_1 \cup \gamma_2)$
$= (\alpha_1 + \beta_1 \mid \gamma_1)_1 \cup (\alpha_2 + \beta_2 \mid \gamma_2)_2$
$= (\alpha_1 \mid \gamma_1)_1 \cup (\alpha_2 \mid \gamma_2)_2 + (\beta_1 \mid \gamma_1)_1 \cup (\beta_2 \mid \gamma_2)_2$
$= (\alpha_1 \mid \gamma_1)_1 + (\beta_1 \mid \gamma_1)_1 \cup (\alpha_2 \mid \gamma_2)_2 + (\beta_2 \mid \gamma_2)_2$
$= (\alpha_1 + \beta_1 \mid \gamma_1)_1 \cup (\alpha_2 + \beta_2 \mid \gamma_2)_2$.

(ii) $(C\alpha \mid \beta) = C(\alpha \mid \beta)$
i.e. $((C_1 \cup C_2)(\alpha_1 \cup \alpha_2) \mid (\beta_1 \cup \beta_2))$
$= (C_1 \alpha_1 \mid \beta_1)_1 \cup (C_2 \alpha_2 \mid \beta_2)_2$
$= C_1 (\alpha_1 \mid \beta_1)_1 \cup C_2 (\alpha_2 \mid \beta_2)_2$
$= (C_1 \cup C_2)(\alpha_1 \cup \alpha_2 \mid (\beta_1 \cup \beta_2)$
$= C(\alpha \mid \beta)$.

(iii) $(\beta \mid \alpha) = \overline{(\alpha \mid \beta)}$ (bar-denoting the complex conjugation)
(d) $(\alpha \mid \alpha) > 0$ if $\alpha \neq 0$.
i.e. $(\alpha_1 \cup \alpha_2 \mid \alpha_1 \cup \alpha_2) = (\alpha \mid \alpha)$
$= (\alpha_1 \mid \alpha_1)_1 \cup (\alpha_2 \cup \alpha_2)_2$
$(\alpha_i \mid \alpha_i) > 0$ if $\alpha_i \neq 0$  i = 1, 2.

***Example 3.2.16:*** Let $V = V_1 \cup V_2$ be a bivector space over the field Q. $V_1 = Q \times Q$ and $V_2$ = {all polynomials of degree less than three with coefficients from Q}. Let ( / ) be a bifunction i.e.

$$( / ) = ( / )_1 \cup ( / )_2.$$

For $\alpha = \alpha_1 \cup \alpha_2$, and $\beta = \beta_1 \cup \beta_2$ $\alpha_1, \beta_1 \in V_1$ and $\alpha_2, \beta_2 \in V_2$ where $\alpha_1 = (x_1, x_2)$ and $\beta_1 = (y_1 , y_2)$ define $(\alpha \mid \beta) = (\alpha_1 \mid \beta_1)_1 \cup (\alpha_2 \mid \beta_2)_2$

$$= (x_1 y_1 + x_2 y_2) \cup \int_0^1 \alpha_2 \beta_2 \, dt ,$$

clearly $(\alpha \mid \beta)$ is a inner biproduct on $V = V_1 \cup V_2$.

**DEFINITION 3.2.9:** *An inner biproduct space is a real or complex bivector space together with a specified inner biproduct defined on that space. A finite dimensional real inner biproduct space is called a Euclidean bispace. A complex inner product space is often referred to as a unitary bispace.*

All results in case of inner product space can be derived with simple modifications in case of inner product bispace. The binorm of $\alpha = \alpha_1 \cup \alpha_2$ is defined to be

$$(\alpha \mid \alpha) = (\alpha_1 \mid \alpha_1)_1 \cup (\alpha_2 \mid \alpha_2)_2 .$$
$$= \|\alpha_1\|_1 \cup \|\alpha_2\|_2 .$$
$$\|\alpha\| = \|\alpha_1\|_1 \cup \|\alpha_2\|_2$$

is called the binorm of $\alpha$. Let $\alpha = \alpha_1 \cup \alpha_2$, $\beta = \beta_1 \cup \beta_2$ be bivectors in an inner biproduct space $V = V_1 \cup V_2$. Then $\alpha$ is biorthogonal to $\beta$ if $(\alpha \mid \beta) = 0$ i.e. $(\alpha \mid \beta) = (\alpha_1 \mid \beta_1)_1 \cup (\alpha_2 \mid \beta_2)_2 = 0 \cup 0$.

So $(\alpha \mid \beta) = 0$ implies $\beta$ is biorthogonal to $\alpha$ and $\alpha$ and $\beta$ are biorthogonal.

If $(\alpha \mid \beta) \neq 0$ but $(\alpha_1 \mid \beta_1)_1$ or $(\alpha_2 \mid \beta_2 )_2$ is zero in $(\alpha \mid \beta) = (\alpha_1 \mid \beta_1)_1 \cup (\alpha_2 \mid \beta_2)_2$ then we say $\alpha$ and $\beta$ are semi biorthogonal.

Let $S = S_1 \cup S_2$ be a set of vectors in $V = V_1 \cup V_2$. We say S is a biorthogonal set provided all pairs of distinct vectors in $S_1$ and all pairs distinct vectors in $S_2$ are orthogonal. An biorthonormal set is an biorthogonal set S with an additional property that

$$\| \alpha \| = \| \alpha_1 \|_1 \cup \| \alpha_2 \|_2 = 1 \cup 1 \text{ i.e.} \| \alpha_i \|_i = 1$$

for i = 1, 2 for every $\alpha_i$ in $S_i$ , i = 1, 2.

Clearly zero vector is biorthogonal with every bivector in V. It is an additional property that $\|\alpha\| = \|\alpha_1\|_1 \cup \|\alpha_2\|_2$

$$= 1 \cup 1 \text{ i.e. } \|\alpha_i\| = 1$$

for i = 1, 2 for every $\alpha_i$ in $S_i$ , i = 1, 2.

Clearly zero vector is biorthogonal with every bivector in V. It is left as an exercise for the reader to prove.

An biorthogonal set of non zero bivectors is linearly independent. Note when we say a set of vectors $\nu_1, \ldots, \nu_k$ in $V = V_1 \cup V_2$ is linearly independent, we mean if

$\nu_i = \nu_i^1 \cup \nu_i^2$, i = 1, 2, …, k

then $\left(\nu_1^1, \ldots, \nu_k^1\right)$ and $\left(\nu_1^2, \ldots, \nu_k^2\right)$ form a linearly independent set of $V_1$ and $V_2$ respectively. Now we give a result analogous to Gram-Schmidt orthogonalization process.

We say a set of bivectors $\left(\nu_1, \ldots, \nu_k\right)$ in $V = V_1 \cup V_2$ is semilinearly independent if $\nu_i = \nu_i^1 \cup \nu_i^2$, i = 1, 2, … k, and $\left(\nu_1^1, \ldots, \nu_k^1\right)$ or $\left(\nu_1^2, \ldots, \nu_k^2\right)$ alone is a linearly independent set of $V_1$ or $V_2$ respectively i.e. one set is linearly dependent and another set is linearly independent.

***Example 3.2.17:*** Let $V = V_1 \cup V_2$ be a bivector space over Q. Let $V_1 = \{Q \times Q \times Q\}$ and

$$V_2 = \left\{ \begin{pmatrix} a & b \\ c & d \end{pmatrix} \middle| a, b, c, d, \in Q \right\}.$$

Clearly V is a bivector space over Q.

Take the set $(\nu_1, \nu_2, \nu_3, \nu_4)$ in $V = V_1 \cup V_2$ where

$$\nu_1 = (1\ 0\ 0) \cup \begin{pmatrix} 0 & 0 \\ 0 & 1 \end{pmatrix}$$

$$\nu_2 = (0\ 1\ 1) \cup \begin{pmatrix} 1 & 0 \\ 0 & 0 \end{pmatrix}$$

$$\nu_3 = (1\ 1\ 1) \cup \begin{pmatrix} 0 & 1 \\ 0 & 0 \end{pmatrix}$$

$$\nu_4 = (0\ 1\ 0) \cup \begin{pmatrix} 0 & 0 \\ 1 & 0 \end{pmatrix}.$$

Clearly the set {(1 0 0), (0 1 1), (1 1 1), (0 1 0)} is a linearly dependent set where as

$$\left\{\begin{pmatrix} 0 & 0 \\ 0 & 1 \end{pmatrix}, \begin{pmatrix} 1 & 0 \\ 0 & 0 \end{pmatrix}, \begin{pmatrix} 0 & 1 \\ 0 & 0 \end{pmatrix}, \begin{pmatrix} 0 & 0 \\ 1 & 0 \end{pmatrix}\right\}$$

is a linearly independent set.

Thus $\{\nu_1, \nu_2, \nu_3, \nu_4\}$ is a semilinearly independent set of V and clearly $\{\nu_1, \nu_2, \nu_3, \nu_4\}$ is not a linearly independent set of V.

***Example 3.2.18:*** Let $V = V_1 \cup V_2$ where $V_1 = \{Q \times Q \times Q \times Q\}$ and $V_2$ = {all polynomials of degree less than or equal to 5 with coefficients from Q}. Clearly $V = V_1 \cup V_2$ is a bivector space over Q.

Consider the set of vectors $\{\nu_1, \nu_2, \nu_3, \nu_4\}$ with

$$\begin{aligned} \nu_1 &= (1\ 0\ 0\ 0) \cup \{1\} \\ \nu_2 &= (0\ 1\ 0\ 0) \cup \{x\} \\ \nu_3 &= (0\ 0\ 1\ 0) \cup \{x^3\} \\ \nu_4 &= (0\ 0\ 0\ 1) \cup \{x^5\}. \end{aligned}$$

Clearly $\{v_1, v_2, v_3, v_4\}$ is a linearly independent set of V over Q.

Now we have the following interesting theorem.

**THEOREM 3.2.5:** *Let $V = V_1 \cup V_2$ be a finite dimensional bivector space. Let dim V = (m, n) i.e. dim $V_1$ = m and dim $V_2$ = n. Any set of non zero bivectors $\{v_1, \ldots, v_k\}$ of V (with K > m and K > n) where $v_i = v_i^1 \cup v_i^2$ (i = 1, 2, …, K) is always a linearly dependent set of V.*

*Proof:* Follows from the fact the set of elements $\left(v_1^1,\ldots,v_k^1\right)$ in $V_1$ and $\left(v_1^2,\ldots,v_k^2\right)$ in $V_2$ are a linearly dependent sets of $V_1$ and $V_2$ respectively.

Now we give a theorem for semi-linearly dependent set of the bivector space $V = V_1 \cup V_2$.

**THEOREM 3.2.6:** *Let $V = V_1 \cup V_2$ be a bivector space over the field F. Let dim V = (m , n) i.e. dim $V_1$ = m and dim $V_2$ = n. A set of non zero bivectors $(v_1, \ldots, v_k)$ with (K > m and K < n 'or' K > n and K < m) is a linearly independent set if and only if $\left(v_1^2, \ldots, v_k^2\right)$ forms an independent set of $V_2$ (or $\left(v_1^1, \ldots, v_k^1\right)$ forms an independent set of $V_1$).*

*Proof:* Follows as in case of the theorem 3.2.5.

Suppose $(v_1, \ldots, v_m) \in V = V_1 \cup V_2$ be a linearly independent biorthogonal set of a inner biproduct bivector space and $\beta = \beta_1 \cup \beta_2$ is a linear combination of $(v_1, \ldots, v_m)$ then β is a particular linear combination of the form

$$\beta = \sum_{K=1}^{m} \frac{(\beta \mid \alpha_K)\alpha_K}{\|\alpha_k\|^2}$$

$$\alpha_K = \alpha_K^1 \cup \alpha_K^2 \text{ and } \|\alpha_K\|^2 = \|\alpha_K\|_1^2 \cup \|\alpha_K\|_2^2$$

$$\beta = \beta_1 \cup \beta_2 = \sum_{K=1}^{m} \frac{(\beta \mid \alpha_K)\ \alpha_K}{\|\alpha_K\|_2^2}.$$

Now we proceed on to give the modified form of Gram Schmidt orthogonalization process which we choose to call as Gram Schmidt biorthogonalization process.

**THEOREM 3.2.7:** *Let $V = V_1 \cup V_2$ be an inner biproduct space and let $\beta_1,\ldots, \beta_K$ be any independent bivectors in $V = V_1 \cup V_2$. Then one may construct biorthogonal bivectors $\alpha_1,\ldots,\alpha_K$ in $V = V_1 \cup V_2$ such that the set $\{\alpha_1,\ldots,\alpha_K\}$ is a basis for the subbispace spanned by $\beta_1,\ldots, \beta_K$.*

*Proof:* The bivectors $\alpha_1, \ldots, \alpha_K$ will be obtained by means of a new construction known as the Gram-Schmidt biorthogonalization process.
Let

$$\alpha_i = \alpha_i^1 \cup \alpha_i^2; \quad i = 1, 2, \ldots, K$$
$$\beta_i = \beta_i^1 \cup \beta_i^2, \quad i = 1, 2, \ldots, K$$

First let $\alpha_1 = \beta_1$ i.e. $\alpha_1^1 = \beta_1^1$ and $\alpha_2^2 = \beta_2^2$.
The other bivectors are given inductively as follows: Suppose $\alpha_1, \ldots, \alpha_m$ $(1 \leq m < K)$ have been chosen so that for every K, $\{\alpha_1, \ldots, \alpha_K\}$, $1 \leq K < m$ is an biorthogonal basis for the subbispace of V spanned by $\beta_1,\ldots, \beta_K$.
To construct the next vector

$$\alpha_{m+1} = \alpha_{m+1}^1 \cup \alpha_{m+1}^2.$$

Let

$$\alpha_{m+1} = \beta_{m+1} - \sum_{K=1}^{m} \frac{(\beta_{m+1} \mid \alpha_K)\alpha_K}{\|\alpha_K\|^2}$$

that is

$$\alpha_{m+1}^1 \cup \alpha_{m+1}^2 = \left(\beta_{m+1}^1 \cup \beta_{m+1}^2\right) -$$

$$\sum_{K=1}^{m}\frac{\left(\beta_{m+1}^{1}|\alpha_{K}^{1}\right)\alpha_{K}^{1}}{\left\|\alpha_{K}^{1}\right\|_{1}^{2}} \cup \sum_{K=1}^{m}\frac{\left(\beta_{m+1}^{2}|\alpha_{K}^{2}\right)\alpha_{K}^{2}}{\left\|\alpha_{K}^{2}\right\|_{2}^{2}}$$

$$= \left(\beta_{m+1}^{1}-\sum_{K=1}^{m}\frac{\left(\beta_{m+1}^{1}|\alpha_{K}^{1}\right)\alpha_{K}^{1}}{\left\|\alpha_{K}^{1}\right\|_{1}^{2}}\right) \cup \left(\beta_{m+1}^{2}-\sum_{K=1}^{m}\frac{\left(\beta_{m+1}^{2}|\alpha_{K}^{2}\right)\alpha_{K}^{2}}{\left\|\alpha_{K}^{2}\right\|_{2}^{2}}\right).$$

Then $\alpha_{m+1} \neq 0$. For otherwise $\beta_{m+1}$ is a linear combination of $\alpha_1,\ldots,\alpha_m$ and hence a linear combination of $\beta_1,\ldots, \beta_m$. Further more if $1 \leq j \leq m$, then

$$\left(\alpha_{m+1}|\alpha_j\right)$$

$$= \quad \left(\beta_{m+1}|\alpha_j\right)-\sum_{K=1}^{m}\frac{\left(\beta_{k+1}|\alpha_k\right)\left(\alpha_k|\alpha_j\right)}{\left\|\alpha_k\right\|^2}$$

$$\left(\alpha_{m+1}^{1}|\alpha_{j}^{1}\right) \cup \left(\alpha_{m+1}^{2}|\alpha_{j}^{2}\right)$$

$$= \quad \left(\beta_{m+1}^{1}|\alpha_{j}^{1}\right) \cup \left(\beta_{m+1}^{2}|\alpha_{j}^{2}\right) \; -$$

$$\sum_{K=1}^{m}\frac{\left(\beta_{K+1}^{1}|\alpha_{K}^{1}\right)\left(\alpha_{K}^{1}\,|\,\alpha_{j}^{1}\right)}{\left\|\alpha_{K}^{1}\right\|_{1}^{2}} \cup \sum_{K=1}^{m}\frac{\left(\beta_{K+1}^{2}|\alpha_{K}^{2}\right)\left(\alpha_{K}^{2}\,|\,\alpha_{j}^{2}\right)}{\left\|\alpha_{K}^{2}\right\|_{2}^{2}}$$

$$= \left(\beta_{m+1}^{1}|\alpha_{j}^{1}\right) \cup \left(\beta_{m+1}^{2}|\alpha_{j}^{2}\right) \; - \left(\beta_{m+1}^{1}|\alpha_{j}^{1}\right) \cup \left(\beta_{m+1}^{2}|\alpha_{j}^{2}\right)$$

$$= 0 \cup 0.$$

Therefore $(\alpha_1,\ldots, \alpha_{m+1})$ is a biorthogonal set consisting of m + 1 non zero bivectors in the sub bispace spanned by $\beta_1, \ldots, \beta_{m+1}$, it is a basis for this bisubspace. Now the bivectors $\alpha_1,\ldots, \alpha_m$ can be constructed one after the other in accordance with the following rule $\alpha_1 = \beta_1$ .

$$\alpha_1^1 \cup \alpha_1^2 = \beta_1^1 \cup \beta_1^2$$

$$\alpha_2 = \alpha_2^1 \cup \alpha_2^2 = \beta_2 - \frac{\left(\beta_2|\alpha_1\right)\alpha_1}{\left\|\alpha_1\right\|^2}$$

$$= \quad (\beta_2^1 \cup \beta_2^2) - \left( \frac{(\beta_2^1 | \alpha_1^1)\alpha_1^1}{\|\alpha_1^1\|_1^2} \cup \frac{(\beta_2^2 | \alpha_2^2)\alpha_1^2}{\|\alpha_1^2\|_2^2} \right)$$

$$= \quad \left( \beta_2^1 - \frac{(\beta_2^1 | \alpha_1^1)\alpha_1^1}{\|\alpha_1^1\|_1^2} \right) \cup \left( \beta_2^2 - \frac{(\beta_2^2 | \alpha_1^2)\alpha_1^2}{\|\alpha_1^2\|_2^2} \right).$$

$$\alpha_3 \quad = \quad \beta_3 - \frac{(\beta_3 | \alpha_1)\alpha_1}{\|\alpha_1\|^2} - \frac{(\beta_3 | \alpha_2)\alpha_2}{\|\alpha_2\|^2}$$

$$\alpha_3^1 \cup \alpha_3^2 = \quad (\beta_3^1 \cup \beta_3^2) - \left( \frac{(\beta_3^1 | \alpha_1^1)\alpha_1^1}{\|\alpha_1^1\|_1^2} \cup \frac{(\beta_2^2 | \alpha_1^2)\alpha_1^2}{\|\alpha_1^2\|_2^2} \right)$$

$$- \left( \frac{(\beta_3^1 | \alpha_2^1)\alpha_2^1}{\|\alpha_2^1\|_1^2} \cup \frac{(\beta_3^2 | \alpha_2^2)\alpha_2^2}{\|\alpha_2^2\|_2^2} \right)$$

$$= \quad \left( \beta_3^1 - \frac{(\beta_3^1 | \alpha_1^1)\alpha_1^1}{\|\alpha_1^1\|_1^2} - \frac{(\beta_3^1 | \alpha_2^1)\alpha_2^1}{\|\alpha_2^1\|_1^2} \right) \cup$$

$$\left( \beta_3^2 - \frac{(\beta_2^2 | \alpha_1^2)\alpha_1^2}{\|\alpha_1^2\|_2^2} - \frac{(\beta_3^2 | \alpha_2^2)\alpha_2^2}{\|\alpha_2^2\|_2^2} \right)$$

and so on.

Thus we give an example how we construct biorthogonal basis.

***Example 3.2.19:*** Consider the bivector space $V = V_1 \cup V_2$ over Q; where $V_1 = Q \times Q \times Q$ and $V_2$ = {the set of all polynomials of degree less than or equal to 2 with coefficients from Q}. $V = V_1 \cup V_2$ is a bivector space over Q. Consider the set of bivectors which forms an linearly independent set given by $\{\nu_1, \nu_2, \nu_3\}$ where

$$\nu_1 = (3\ 0\ 4) \cup \{ 1 \}$$
$$\nu_2 = (-1, 0, 7) \cup \{ x \}$$
$$\nu_3 = (2\ 9\ 11) \cup \{x^2\}.$$

Let V be endowed with the inner biproduct ( / ) = $( / )_1 \cup ( / )_2$ where $( / )_1$ is the standard inner product and $( / )_2$

$$\int_0^1 f(x)\, g(x)\, dx\ .$$

Now applying the new Gram Schmidt biorthogonalization process we obtain the following bivector.

$$\alpha_1 = (3\ 0\ 4) \cup \{ 1 \}$$

$$\alpha_2 = \left\{(-1,0,7) - \frac{((-107)|(304))}{25}(3,0,4)\right\} \cup \left\{x - \int_0^1 \frac{1.x.dx}{1}\right\}$$

$$= (-4, 0, 3) \cup (x - ½ ).$$

$$\alpha_3 = \left\{(2,9,11) - \frac{(2,9,11)|(304)}{25}(3,0,4) - \frac{((2,9,11)|(-4,0,3))}{25}(-4,0,3) \cup \left\{x^2 - \left(\int_0^1 \frac{x^2.1}{1}\right)^{dx^1} - \left(\int_0^1 \frac{x^2\left(x - \frac{1}{2}\right)dx}{1/12}\right)^{x-\frac{1}{2}}\right\}\right.$$

$$= (0, 9, 0) \cup \left(x^2 - x + \tfrac{1}{6}\right)$$

Thus {(3 0 4), (-4, 0, 3), (0, 9, 0)} $\cup \left\{1, x - \tfrac{1}{2}, x^2 - x + \tfrac{1}{6}\right\}$ forms a biorthogonal bivectors of $V = V_1 \cup V_2$ over Q.

Now we proceed on to define the notion of best biapproximation to the bivector β by bivectors in the subbivector space W of a bivector α in W.

**DEFINITION 3.2.10:** *Let V = $V_1 \cup V_2$ be a bivector space over the field F. Let on V be defined a inner biproduct* $\| \; \| = \| \; \|_1 \cup \| \; \|_2$. *Let W be a subbivector space of V = $V_1 \cup V_2$ and let $\beta = \beta_1 \cup \beta_2$ be an arbitrary bivector in V. A best biapproximation to $\beta = \beta_1 \cup \beta_2$ by bivectors in W is a bivector $\alpha = \alpha_1 \cup \alpha_2$ in W such that*

$$\|\beta - \alpha\| \leqq \|\beta - \gamma\|$$

*i.e.* $\|\beta_1 - \alpha_1\|_1 \cup \|\beta_2 - \alpha_2\|_2 \leq \|\beta_1 - \gamma_1\|_1 \cup \|\beta_2 - \gamma_2\|_2$

*for every bivector $\gamma = \gamma_1 \cup \gamma_2$ in W = $W_1 \cup W_2$*

*(i)* *If a best biapproximation to $\beta = \beta_1 \cup \beta_2$ by bivectors in W = $W_1 \cup W_2$ exists, it is unique.*

*(ii)* *If W = $W_1 \cup W_2$ is finite dimensional and $\{\alpha_1,\ldots,\alpha_k\}$ is any biorthonormal basis for W = $W_1 \cup W_2$ then the vector*

$$\alpha = \sum \frac{(\beta \,|\, \alpha_k)\alpha_k}{\|\alpha_k\|^2}$$

$$= \alpha_1 \cup \alpha_2 = \left( \sum_k \frac{(\beta_1 \,|\, \alpha_k^1)\alpha_k^1}{\|\alpha_k^1\|_1^2} \right) \cup \left( \sum_k \frac{(\beta_2 \,|\, \alpha_k^2)\alpha_k^2}{\|\alpha_k^2\|_2^2} \right)$$

*is the unique best biapproximation to $\beta = \beta_1 \cup \beta_2$ by bivectors in W = $W_1 \cup W_2$.*

Now are proceed on to define the notion of biorthogonal bicomplement of an inner biproduct bivector space V = $V_1 \cup V_2$.

**DEFINITION 3.2.11:** *Let V = $V_1 \cup V_2$ be a inner biproduct bivector space and S = $S_1 \cup S_2$ any set of bivectors in V = $V_1 \cup V_2$. The biorthogonal bicomplement of S = $S_1 \cup S_2$ is* *the set* $S^{\perp} = S_1^{\perp} \cup S_2^{\perp}$ *of all bivectors in* $V = V_1 \cup V_2$ *which are biorthogonal to every bivector in* $S = S_1 \cup S_2$*. The biorthogonal bicomplement of the space* $V = V_1 \cup V_2$ *is the zero bisubspace* $\{0\} \cup \{0\}$ *and conversely* $\{0\}^{\perp} = \{0\}^{\perp} \cup \{0\}^{\perp}$ $= V_1 \cup V_2$.

Now we proceed on to define the notion of idempotent bioperator on a linear bivector space.

**DEFINITION 3.2.12:** *Let* $V = V_1 \cup V_2$ *be a linear bivector space over the field F. A linear bioperator* $E : V \to V$ *is said to be an idempotent bioperator, if* $E = E_1 \cup E_2$ *and each of* $E_1$ *and* $E_2$ *are idempotent operators on* $V_1$ *and* $V_2$ *respectively i.e.,* $E^2 = E$ *as* $= E_1^2 = E_1$ *and* $E_2^2 = E_2$ *i.e.* $E = E_1 \cup E_2$ *and* $E^2 = E$.

We now proceed on to define the notion of biadjoints analogous to adjoints. Let V be a bivector space over F. T be a linear bioperator in V i.e. $T = T_1 \cup T_2 : V_1 \cup V_2 \to V_1 \cup V_2$. Let $V_1 \cup V_2$ be a inner biproduct bivector space with an inner biproduct $( \mid ) = ( \mid )_1 \cup ( \mid )_2$ defined on $V = V_1 \cup V_2$. The biadjoint of the linear bioperator $T^*$ on V is defined such that $(T\alpha \mid \beta) = (\alpha \mid T^* \beta)$ where $T^* = T_1^* \cup T_2^*$ and $\alpha = \alpha_1 \cup \alpha_2$ and $\beta = \beta_1 \cup \beta_2$ with

$$
\begin{aligned}
(T\alpha \mid \beta) &= ((T_1 \cup T_2)(\alpha_1 \cup \alpha_2) \mid (\beta_1 \cup \beta_2)) \\
&= (T_1 \alpha_1 \mid \beta_1) \cup (T_2 \alpha_2 \mid \beta_2) \\
&= ((\alpha_1 \cup \alpha_2) \mid (T_1^* \cup T_2^*)(\beta_1 \cup \beta_2) \\
&= (\alpha_1 \mid T_1^* \beta_1) \cup (\alpha_1 \mid T_2^* \beta_2). \\
&= (\alpha \mid T^* \beta)
\end{aligned}
$$

$T^* = T_1^* \cup T_2^*$ (where $T = T_1 \cup T_2$) is called the biadjoint of T.

A linear bioperator T such that $T = T^*$ is called the self biadjoint on the inner biproduct bivector space $V = V_1 \cup V_2$.

Now we proceed on to define the notion of unitary bioperators on a bivector space.

**DEFINITION 3.2.13:** *Let $A = A_1 \cup A_2$ be a $n \times n$ complex bimatrix. (A bimatrix A is said to be complex if it takes entries from the complex field). A is called biunitary if $A^* A = I$ i.e. $A_1^* A_1 \cup A_2^* A_2 = I^1_{n \times n} \cup I^2_{n \times n}$. A real or complex $n \times n$ bimatrix $A = A_1 \cup A_2$ is said to be biorthogonal if $A^t A = I$ i.e. $A_1^t A_1 \cup A_2^t A_2 = I^1_{n \times n} \cup I^2_{n \times n}$. A complex $n \times n$ bimatrix $A = A_1 \cup A_2$ is called binormal if $AA^* = A^* A$ i.e. $A_1 A_1^* \cup A_2 A_2^* = A_1^* A_1 \cup A_2^* A_2$.*

*Now even if we have a complex mixed square bimatrix $A = A_1 \cup A_2$ where $A_1$ is $n \times n$ matrix and $A_2$ is a $m \times m$ matrix still the notion of biunitary and binormal can be defined without any difficulty.*

*Thus if V is a finite dimensional inner biproduct space and $T = T_1 \cup T_2$, a linear bioperator on $V = V_1 \cup V_2$ we say that T is binormal if $T = T_1 \cup T_2$ commutes with its adjoint i.e. $TT^* = T^*T$ and T is biunitary if $TT^* = T^*T = I$. i.e. $(T_1 T_1^* \cup T_2 T_2^*) = (T_1^* T_1 \cup T_2^* T_2) = I^1 \cup I^2$.*

Several other routine relations can be derived for biunitary and binormal bioperators as in case of unitary and normal operators.

A linear bioperator T on a finite dimensional inner biproduct space $V = V_1 \cup V_2$ is non negative if $T = T_1 \cup T_2 = T^* = T_1^* \cup T_2^*$ and $(T \alpha \mid \alpha) = (T \alpha_1 \mid \alpha_1) \cup (T \alpha_2 \mid \alpha_2) \geq 0 \cup 0$ for all $\alpha = \alpha_1 \cup \alpha_2$ in $V = V_1 \cup V_2$. A positive linear bioperator is such that $T = T^*$ and $(T \alpha \mid \alpha) = (T \alpha_1 \mid \alpha_1) \cup (T \alpha_2 \mid \alpha_2) > 0 \cup 0$ for all $\alpha \neq 0$.

Now we proceed on to define biunitary bitransformation. Let $V = V_1 \cup V_2$ and $V^1 = V_1^1 \cup V_2^1$ be inner biproduct spaces over the same field. A linear bitransformation

$$U : V \to V^1$$

i.e. $U = U_1 \cup U_2 : V_1 \cup V_2 \rightarrow V_1^1 \cup V_2^1$ is called a biunitary bitransformation if it maps V onto $V^1$ and preserves inner biproduct. If $T = T_1 \cup T_2$ is a linear bioperator on $V = V_1 \cup V_2$ and $T^1 = T_1^1 \cup T_2^1$ a linear bioperator on $V^1 = V_1^1 \cup V_2^1$ then $T = T_1 \cup T_2$ is unitarily equivalent to $T^1 = T_1^1 \cup T_2^1$ if there exists a biunitary bitransformation U of V onto $V^1$ such that $UTU^{-1} = T^1$. i.e., $(U_1\ T_1\ U_1^{-1}) \cup (U_2\ T_2 U_2^{-1}) = T_1^1 \cup T_2^1$.

Most of the interrelated results can be derived in case of bioperators as in case of operators. Now we study linear bialgebra in an entirely different way, i.e., over finite fields.

## 3.3 Bivector spaces over finite fields

We introduce the notion of bivector spaces over finite fields $Z_p$. This study forces us to define pseudo inner product on these spaces.

Throughout this section $Z_p$ will denote the prime field of characteristic p, p a prime i.e. $Z_p = \{0, 1, 2, \ldots, p - 1\}$. Clearly $Z_p$ is a finite field having just p elements including 0. $Z_p \times Z_p \times Z_p \times \ldots \times Z_p = Z_p^m$ denotes the direct product of $Z_p$ taken m times m > 1 and $Z_p^m$ is a group under addition modulo p.

***Example 3.3.1:*** Let $V = V_1 \cup V_2$ be a bivector space over $Z_p$ where $V_1 = Z_p \times Z_p \times Z_p$ and $V_2 = Z_p [x]$. V is a bivector space over $Z_p$.

***Example 3.3.2:*** Let $V = V_1 \cup V_2$ where $V_1 = Z_3 \times Z_3 \times Z_3$ and $V_2 = Z_3[x]$. V is a bivector space over $Z_3$ and dimension of V is infinite.

**DEFINITION 3.3.1:** *Let $Z_{p_1} \cup Z_{p_2}$ ($p_1 \neq p_2$, $p_1$, $p_2$ be primes) be a bifield. $V = V_1 \cup V_2$ be a bigroup say*

*$V_1 = Z_{p_1} \times Z_{p_1}$ and $V_2$ = $Z_{p_2} \times Z_{p_2} \times Z_{p_2}$ we call V the bivector space over the bifield $Z_{p_1} \cup Z_{p_2}$. This bivector space is called the strong bivector space.*

It is interesting to note that we in general cannot derive any form of relation between bivector spaces.

***Example 3.3.3:*** Let F = $Z_5 \cup Z_7$ be the bifield, V = $(Z_5 \times Z_5 \times Z_5) \cup (Z_7 \times Z_7)$ be the bigroup. V is a strong bivector space over F = $Z_5 \cup Z_7$.

Clearly {(1 0 0), (0 1 0) (0 0 1)} $\cup$ {(0, 1), (1, 0)} is a basis of V over the bifield $Z_5 \cup Z_7$.

As in case of bivector spaces over field of characteristic zero, we cannot define inner biproduct in case of bivector spaces over field of characteristic p. So in case of bivector spaces over finite field or even in case of strong bivector spaces over finite bifields we have to define a new inner biproduct which we choose to call as pseudo inner biproduct. The main reason for defining *pseudo inner biproduct* is we may not have

$(\alpha / \alpha)(\alpha_1 / \alpha_1)_1 \cup (\alpha_2 / \alpha_2)_2 > 0$ where $\alpha = \alpha_1 \cup \alpha_2 \neq 0 \cup 0$.

For in the case of finite fields it may happen $(\alpha / \alpha) = ( / )_1 \cup ( / )_2 = 0 \cup 0$ even without $\alpha$ being equal to zero, but all other conditions will be satisfied .

***Example 3.3.4:*** Let V = $V_1 \cup V_2$ be a bivector space over $Z_{11}$ where $V_1 = Z_{11} \times Z_{11}$ and $V_2$ = {all polynomials of degree less than or equal to 3 with coefficients from $Z_{11}$}. Clearly V is bivector space over $Z_{11}$.

Define inner biproduct $( / ) = ( / )_1 \cup ( / )_2$ where $( / )_1$ is the inner product $2x_1 y_1 + 4 x_2 y_2$ and $(p (x) / q(x))_2 = p_0 q_0 + \ldots + p_4 q_4$, where $x = (x_1, x_2)$ and $y = (y_1, y_2)$ and

$$p (x) = (p_0, \ldots, p_4)$$
$$q (x) = (q_0, \ldots, q_4).$$

Suppose consider $\alpha = (2, 3) \cup (p(x)) = 2 + 2x + 3x^2 + 4x^3)$
$= \alpha_1 \cup \alpha_2$

$$
\begin{aligned}
(\alpha / \alpha) &= (\alpha_1 / \alpha_1)_1 \cup (\alpha_2 / \alpha_2)_2 \\
&= (2.2.2 + 4.3.3) \cup (2.2 + 2.2 + 3.3 + 4.4) \\
&\qquad\qquad \text{(addition modulo 11)} \\
&= 44 \text{ (mod 11)} \cup 33 \text{ (mod 11)} \\
&= 0 \cup 0 \text{ (mod 11).}
\end{aligned}
$$

But $\alpha = \alpha_1 \cup \alpha_2 \neq 0$ with $(\alpha / \alpha) = 0$. Thus we may have an inner biproduct such that $(\alpha / \alpha) = 0$ without $\alpha = 0$. Thus this inner biproduct is called as the pseudo inner biproduct on V. Determination of several results for pseudo inner biproducts is left as an exercise to the reader.

## 3.4 Representation of finite bigroup

Now we proceed on to a very interesting concept called Representation of finite bigroups. This is the first time such notion is defined. We have extensively studied bigroups but have not ventured to study their representation. Here we bring in the representation of finite bigroups and it's inter-relation with bivector spaces and its properties. Throughout this section $G = G_1 \cup G_2$ be a bigroup and $V = V_1 \cup V_2$ be a bivector space.

**DEFINITION 3.4.1:** *A birepresentation of $G = G_1 \cup G_2$ on $V = V_1 \cup V_2$ is a mapping $\rho = \rho_1^1 \cup \rho_2^2$ from G to invertible linear bi transformation on $V = V_1 \cup V_2$ such that $\rho_{xy} = \rho_x \circ \rho_y$ for all $x, y \in G_1 \cup G_2$ here $x = x_1 \cup x_2$ and $y = y_1 \cup y_2$. So*

$$\rho_{xy} = \rho^1_{x_1 y_1} \cup \rho^2_{x_2 y_2} = \left(\rho^1_{x_1} \cup \rho^2_{x_2}\right) \circ \left(\rho^1_{y_1} \cup \rho^2_{y_2}\right).$$

*Here we use $\rho_x = \rho^1_{x_1} \cup \rho^2_{x_2}$ to describe the invertible linear bitransformation on $V = V_1 \cup V_2$ associated to $x = x_1 \cup x_2$ in $G = G_1 \cup G_2$. So we can write*

$$\rho_x(\upsilon) = \rho^1_{x_1}(\upsilon_1) \cup \rho^2_{x_2}(\upsilon_2)$$

*for the image of a bivector $\upsilon = \upsilon_1 \cup \upsilon_2$ in V under $\rho_x = \rho^1_{x_1} \cup \rho^2_{x_2}$. Thus $\rho_e = I = \rho^1_{e_1} \cup \rho^2_{e_2}$ where $e_1$ is the identity element of $G_1$ and $e_2$ is the identity element of $G_2$ and $I = I_1 \cup I_2$ denotes the identity bitransformation on $V = V_1 \cup V_2$ and*

$$\rho_{x^{-1}} = \rho^1_{x_1^{-1}} \cup \rho^2_{x_2^{-1}}$$

$$= \left(\rho^1_{x_1}\right)^{-1} \cup \left(\rho^2_{x_2}\right)^{-1}$$

*for all $x = x_1 \cup x_2$ in $G = G_1 \cup G_2$.*

*In other words a birepresentation of the bigroup $G = G_1 \cup G_2$ on $V = V_1 \cup V_2$ is a bihomomorphism from $G = G_1 \cup G_2$ into GL $(V = V_1 \cup V_2)$ = GL $(V_1) \cup$ GL $(V_2)$. The bi dimension (m, n) of V (i.e. dim of $V_1 = m$ and dim $V_2 = n$) is called the bidegree of the birepresentation.*

Basic examples of birepresentations are left biregular birepresentation and right biregular birepresentation over a field K. This is defined as follows:

We take $V = V_1 \cup V_2$ to be a bivector space of bifunction on the bigroup $G = G_1 \cup G_2$ with values in K. For the left biregular birepresentation we define;

$$L_x = L^1_{x_1} \cup L^2_{x_2} : V \to V$$

i.e. $L^1_{x_1} : V_1 \to V_1$ and $L^2_{x_2} : V_2 \to V_2$ where $x = x_1 \cup x_2$ for each x in $G = G_1 \cup G$ by

$$L_x (f) (z) = \left(L^1_{x_1} \cup L^2_{x_2}\right) (f_1 \cup f_2) (z)$$

for each function f(z) in V (i.e. $f_1 (z_1)$ in $V_1$ and $f_2 (z_2)$ in $V_2$) is equal to

$$f(x_z^{-1}) = f_1\left(x_1^{-1} z_1\right) \cup f_2\left(x_2^{-1} z_2\right).$$

For right biregular birepresentation we define

$$R_x = R^1_{x_1} \cup R^2_{x_2} : V \to V$$

for each $x = x_1 \cup x_2$ in $G_1 \cup G_2$.

By

$$R_x (f)(z) = f(zx)$$

i.e.

$$\left(R^1_{x_1} \cup R^2_{x_2}\right)\left(f_1 \cup f_2\right)\left(z_1 \cup z_2\right)$$

$$= \quad R^1_{x_1}(f)(z_1) \cup R^2_{x_2}(f_2)(z_2)$$

$$= \quad f_1(z_1 x_1) \cup f_2(z_2 x_2)$$

for each bifunction

$$f(z) = f_1(z_1) \cup f_2(z_2)$$

in $V = V_1 \cup V_2$.

Thus if x and y are elements of $G = G_1 \cup G_2$ then

$$\left(L_x \circ L_y\right) f(z) = \left(L^1_{x_1} \cup L^2_{x_2}\right) \circ \left(L^1_{y_1} \cup L^2_{y_2}\right)(f_1(z_1) \cup f_2(z_2))$$

by

$$R_x (f)(z) = f(zx)$$

i.e. $\left(R^1_{x_1} \cup R^2_{x_2}\right)\left(f_1 \cup f_2\right)\left(z_1 \cup z_2\right)$

$$= \quad R^1_{x_1}(f)(z_1) \cup R^2_{x_2}(f_2)(z_2)$$

$$= \quad f_1\left(z_1 x_1\right) \cup f_2\left(z_2 x_2\right)$$

for each bifunction

$$f(z) = f_1(z_1) \cup f_2(z_2)$$

in $V = V_1 \cup V_2$.

Thus if x and y are elements of $G = G_1 \cup G_2$ then

$$\left(L_x \circ L_y\right) f(z) = \left(L^1_{x_1} \cup L^2_{x_2}\right) \circ \left(L^1_{y_1} \cup L^2_{y_2}\right)\left(f_1(z_1) \cup f_2(z_2)\right)$$

$$= \quad \left(L_{x_1} \cup L_{x_2}\right)\left[\left(L^1_{y_1}(f_1)(z_1)\right) \cup \left(L^2_{y_2}(f_2)(z_2)\right)\right]$$

$$= \quad L^1_{x_1} L^1_{y_1}(f_1)(z_1) \cup L^2_{x_2} L^2_{y_2}(f_2)(z_2)$$

$$= \quad \left(L^1_{y_1} \cup L^2_{y_2}\right)\left(f_1 \cup f_2\right)(x^{-1}z)$$

$$= L^1_{y_1}(f_1)\left(x_1^{-1} z_1\right) \cup L^2_{y_2}(f_2)\left(x_2^{-1} z_2\right)$$
$$= f_1\left(y_1^{-1} x_1^{-1} z_1\right) \cup f_2\left(x\ y_2^{-1} x_2^{-1} z_2\right)$$
$$= f_1\left[\left(x_1 y_1\right)^{-1} z_1\right] \cup f_2\left[\left(x_2 y_2\right)^{-1} z_2\right]$$
$$= L^1_{x_1 y_1}(f_1)(z_1) \cup L^2_{x_2 y_2}(f_2)(z_2)$$
$$= L^1_{y_1}(f_1)(x_1^{-1} z_1) \cup L^2_{y_2}(f_2)(x_2^{-1} z_2)$$
$$= f_1\left(y_1^{-1} x_1^{-1} z_1\right) \cup f_2\left(y_2^{-1} x_2^{-1} z_2\right)$$
$$= f_1\left[\left(x_1 y_1\right)^{-1} z_1\right] \cup f_2\left[\left(x_2 y_2\right)^{-1} z_2\right]$$
$$= L^1_{x_1 y_1}(f_1)(z_1) \cup L^2_{x_2 y_2}(f_2)(z_2)$$

and

$(R_x \text{ o } R_y)(f)(z)$

$$= \left(R^1_{x_1} \cup R^2_{x_2}\right) \circ \left(R^1_{y_1} \cup R^2_{y_2}\right)\ \left(f_1 \cup f_2\right)\left(z_1 \cup z_2\right)$$
$$= R^1_{x_1}\left(R^1_{y_1}(f_1)\right)(z_1) \cup R^2_{x_1}\left(R^2_{y_2}(f_2)\right)(z_2)$$
$$= R^1_{y_1}(f_1)(z_1 x_1) \cup R^2_{y_2}(f_2)(z_2 x_2)$$
$$= f_1(z_1\ x_1\ y_1) \cup f_2(z_2\ x_2\ y_2)$$
$$= R^1_{x_1 y_1}(f_1)(z_1) \cup R^2_{x_2 y_2}(f_2)(z_2).$$

Another description of these birepresentation which can be convenient is the following. For each $w_1 \cup w_2$ in $G_1 \cup G_2 = G$ define bifunction

$$\phi_w(z) = \phi^1_{w_1}(z_1) \cup \phi^2_{w_2}(z_2)$$

on $G = G_1 \cup G_2$ by $\phi^1_{w_1}(z_1) = 1$ where $z_1 = w_1$ , $\phi^2_{w_2}(z_2) = 1$ when $z_2 = w_2$ $\phi^1_{w_1}(z_1) = 0$ when $z_1 \neq w_1$, $\phi^2_{w_2}(z_2) = 0$ when $z_2 \neq w_2$.

This function $\phi_w = \phi^1_{w_1} \cup \phi^2_{w_2}$ for $w = w_1 \cup w_2$ in $G = G_1 \cup G_2$ form a basis for the space of bifunction on the bigroup $G = G_1 \cup G_2$.

One can check

$$L_x\ (\phi_w) = \phi_{xw}$$

$$\text{i.e. } L^1_{x_1}\left(\phi^1_{w_1}\right) \cup L^2_{x_2}\left(\phi^2_{w_2}\right) = \phi^1_{x_1 w_1} \cup \phi^2_{x_2 w_2}$$

and

$$R_x\left(\phi_w\right) = \phi_{w\,x}$$

$$\text{i.e. } R^1_{x_1}\left(\phi^1_{w_1}\right) \cup R^2_{x_2}\left(\phi^2_{w_2}\right) = \phi^1_{w_1 x_1} \cup \phi^2_{w_2 x_2}$$

for all $x = x_1 \cup x_2$ in $G = G_1 \cup G_2$.

Observe that

$$L_x \text{ o } R_y = R_y \text{ o } L_x$$

$$\text{i.e. } \left(L^1_{x_1} \circ R^1_{y_1}\right) \cup \left(L^1_{x_2} \circ R^2_{y_2}\right) = \left(R^1_{y_1} \circ L^1_{x_1}\right) \cup \left(R^2_{y_2} \circ L^2_{x_2}\right)$$

for all x, y in G.

More generally suppose that we have a bihomomorphism from the bigroup $G = G_1 \cup G_2$ to the bipermutation on a non empty finite set $E = E_1 \cup E_2$ (i.e. bipermutation on a non empty set $E = E_1 \cup E_2$ we mean the permutations on a non empty finite set $E_1 \cup$ permutation on a finite set $E_2$). That is suppose for each $x = x_1 \cup x_2$ in $G = G_1 \cup G_2$ we have a bipermutation $\pi_x = \pi^1_x \cup \pi^2_x$ on $E_1 \cup E_2$ i.e. $\pi^1_{x_1}$ a one to one mapping from $E_1$ on to $E_1$ and $\pi^2_{x_2}$ a one to one mapping from $E_2$ on to $E_2$ such that

$$\begin{aligned}
\pi_x \circ \pi_y & \\
&= \left(\pi^1_{x_1} \circ \pi^1_{y_1}\right) \cup \left(\pi^2_{x_2} \circ \pi^2_{y_2}\right) \\
&= \pi^1_{x_1 y_1} \cup \pi^2_{x_2 y_2}.
\end{aligned}$$

As usual $\pi e = \pi^1_{e_1} \cup \pi^2_{e_2}$ is the identity bimapping of $E_1$ and $E_2$ (respectively) and

$$\pi_{x^{-1}} = \pi^1_{x_1^{-1}} \cup \pi^2_{x_2^{-1}}$$

is the inverse bimapping of

$$\pi_x = \pi^1_{x_1} \cup \pi^2_{x_2}$$

on $E = E_1 \cup E_2$. Let $V = V_1 \cup V_2$ be a bivector space of K-valued functions on $E_1 \cup E_2$.

Then we get a birepresentation of $G = G_1 \cup G_2$ on $V = V_1 \cup V_2$ by associating to each $x_1$ in $G_1 \cup x_2$ in $G_2$ the linear bimapping $\pi_x = \pi^1_{x_1} \cup \pi^2_{x_2} : V_1 \cup V_2 \to V_1 \cup V_2$ defined by

$$\pi^1_{x_1}(f_1)(a_1) \cup \pi^2_{x_2}(f_2)(a_2)$$
$$= f_1\left(\pi^1_{x_1} I^1(a_1)\right) \cup f_2\left(\pi^2_{x_2} I^2(a_2)\right)$$

for every bifunction $f(a) = f^1(a_1) \cup f^2(a_2)$ in $V = V_1 \cup V_2$.

Thus is called bipermutation birepresentation corresponding to the bihomomorphism $x \mapsto \pi_x$ i.e.$(x_1 \cup x_2) \mapsto \pi^1_{x_1} \cup \pi^2_{x_2}$ from $G = G_1 \cup G_2$ to bipermutations on $E = E_1 \cup E_2$.

It is indeed a birepresentation because for each $x = x_1 \cup x_2$ and $y = y_1 \cup y_2$ in $G = G_1 \cup G_2$ and each bifunction $f(a) = f\left(a_1^1\right) \cup f\left(a_2^2\right)$ in $V = V_1 \cup V_2$,

we have

$$\left(\pi_x \circ T_y\right)(f)(a)$$
$$= \left[\left(\pi^1_{x_1} \cup \pi^2_{x_2}\right) \circ \left(\pi^1_{y_1} \cup \pi^2_{y_2}\right)\right] (f_1 \cup f_2)(a_1 \cup a_2).$$
$$= \pi_x(\pi_y(f)(a))$$
$$= \left(\pi^1_{x_1} \cup \pi^2_{x_2}\right)\left(\pi^1_{y_1} \cup \pi^2_{y_2}\right)\left(f_1 \cup f_2\right)\left(a_1 \cup a_2\right)$$
$$= f(\pi_y(f)\left(\pi_x I(a)\right))$$
$$= \left(\pi^1_{y_1} \cup \pi^2_{y_2}\right)\left(f_1 \cup f_2\right)\left(\pi^1_{x_1} \cup \pi^2_{x_2}\right)\left(I^1(a_1) \cup I^2(a_2)\right)$$
$$= f\left(\pi_y I\left(\pi_x I(a)\right)\right)$$
$$= ((f_1 \cup f_2)\left(\pi^1_{y_1} \cup \pi^2_{y_2}\right)\left(I^1 \cup I^2\right)\left(\pi^1_{x_1} \cup \pi^2_{x_2}\right) I^1(a_1) \cup I^2(a_2))$$

$$= \quad f(\pi_{(xy)} I(a))$$
$$= \quad f_1\left(\pi^1_{(x_1y_1)} I^1(a_1)\right) \cup f_2\left(\pi^2_{(x_2y_2)} I^2(a_2)\right).$$

For more refer linear algebra [48].

Let G be a finite bigroup. V be a bivector space over a symmetric field K and let $\rho = \rho^1 \cup \rho^2$ be a birepresentation of $G = G_1 \cup G_2$ on $V = V_1 \cup V_2$. If $\langle \, , \rangle$ is a inner biproduct on V then $\langle \, , \rangle = \langle \, , \rangle_1 \cup \langle \, , \rangle_2$ is said to be invariant under the birepresentation $\rho = \rho^1 \cup \rho^2$ or simply $\rho$ - biinvariant if every $\rho_x = \rho^1_{x_1} \cup \rho^2_{x_2} : V \to V$, $x = x_1 \cup x_2$ in $G = G_1 \cup G_2$ preserves inner biproduct that is if

$$\left\langle \rho_x(\nu), \rho_x(w)\right\rangle = \left\langle \nu, w\right\rangle$$

i.e. $\left\langle \rho^1_{x_1}(\nu_1), \rho^1_{x_1}(w_1)\right\rangle_1 \cup \left\langle \rho^2_{x_2}(\nu_2), \rho^2_{x_2}(w_2)\right\rangle_2$

$$= \quad \langle \nu_1, w_1\rangle_1 \cup \langle \nu_2, w_2\rangle_2$$

for all $x = x_1 \cup x_2$ in $G = G_1 \cup G_2$ and $\nu = \nu_1 \cup \nu_2$ and $w = w_1 \cup w_2$ in $V = V_1 \cup V_2$. If $\langle \, , \rangle_o$ is any inner biproduct on $V = V_1 \cup V_2$ then we can obtain an invariant biproduct

$$\langle \, , \rangle = \langle \, , \rangle_1 \cup \langle \, , \rangle_2$$

from o by setting

$$\langle \nu, w\rangle = \langle \nu_1, w_1\rangle_1 \cup \langle \nu_2, w_2\rangle_2 =$$
$$\left(\sum_{y_1 \in G_1} \left\langle \rho^1_{y_1}(\nu_1), \rho^1_{y_1}(w_1)\right\rangle_{o_1}\right) \cup \left(\sum_{y_2 \in G_2} \left\langle \rho^2_{y_2}(\nu_2), \rho^2_{y_2}(w_2)\right\rangle_{o_2}\right).$$

It is easy to check that this does not define a inner biproduct on $V = V_1 \cup V_2$ which is invariant under the birepresentation $\rho = \rho^1 \cup \rho^2$.

One can thus obtain several important relations as in case of representation of finite groups.

### 3.5 Applications of Bimatrix to Bigraphs

In this section we recall some of the basic applications of linear bialgebra. For the very concept of linear bialgebra is itself new we have only very few applications. Interested researcher can develop several applications depending on his creativity. For example, consider the matrix representation existing between members of a set introduced, we have studied the theory of directed graphs to mathematically model the types of sets and relation. On similar lines we can introduce the directed bigraphs to mathematically model the types of sets and relation.

***Example 3.5.1:*** The bigraph given in the figures 3.5.1 and 3.5.2 has the following corresponding vertex bimatrix.

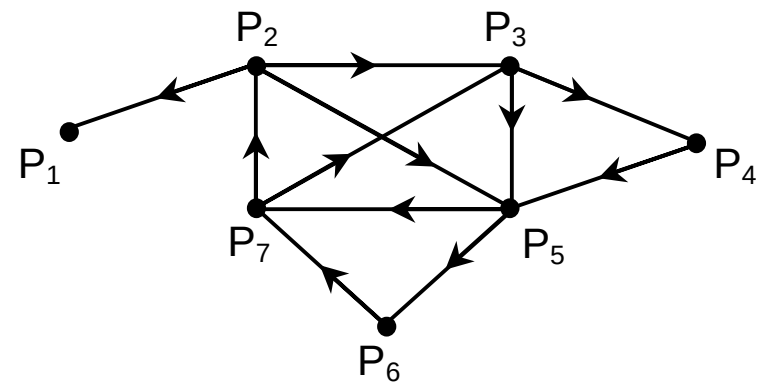

**G₁**

FIGURE: 3.5.1

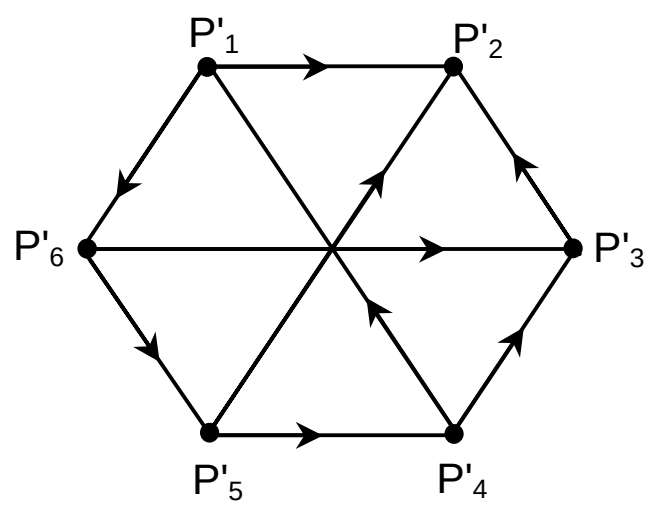

**G₂**

FIGURE: 3.5.2

$G = G_1 \cup G_2$ is the bigraph the related vertex bimatrix is a mixed square bimatrix given by $M = M_1 \cup M_2$. (Here $M_1$ is

the vertex matrix associated with graph the $G_1$ and $M_2$ is the vertex matrix associated with graph the $G_2$).

$$M = \begin{bmatrix} 0 & 0 & 0 & 0 & 0 & 0 & 0 \\ 1 & 0 & 1 & 0 & 1 & 0 & 0 \\ 0 & 0 & 0 & 1 & 1 & 0 & 0 \\ 0 & 0 & 0 & 0 & 1 & 0 & 0 \\ 0 & 0 & 0 & 0 & 0 & 1 & 1 \\ 0 & 0 & 0 & 0 & 0 & 0 & 1 \\ 0 & 1 & 1 & 0 & 0 & 0 & 0 \end{bmatrix} \cup \begin{bmatrix} 0 & 0 & 0 & 0 & 0 & 1 \\ 1 & 0 & 0 & 0 & 0 & 0 \\ 0 & 1 & 0 & 0 & 0 & 0 \\ 1 & 0 & 1 & 0 & 0 & 0 \\ 0 & 1 & 0 & 1 & 0 & 0 \\ 0 & 0 & 1 & 0 & 1 & 0 \end{bmatrix}$$

We give another example of the bigraph and its related vertex bimatrix.

***Example 3.5.2:*** $G = G_1 \cup G_2$ is the bigraph given by the following figures.

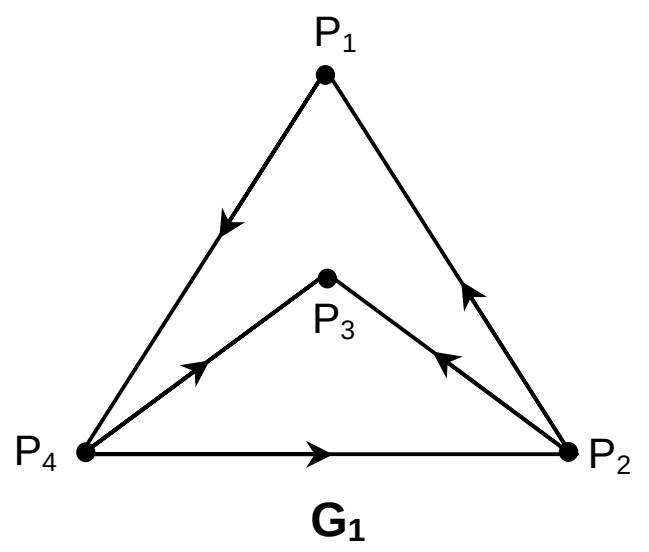

FIGURE: 3.5.3

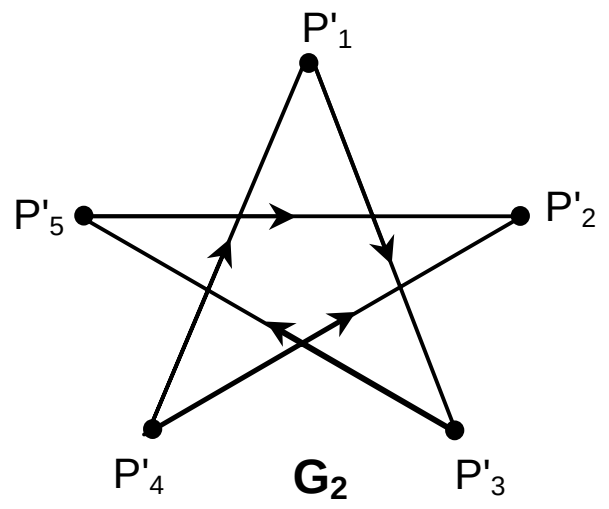

FIGURE: 3.5.4

The related / corresponding vertex bimatrix $M = M_1 \cup M_2$ is

$$\begin{bmatrix} 0 & 0 & 0 & 1 \\ 1 & 0 & 1 & 0 \\ 0 & 0 & 0 & 0 \\ 0 & 1 & 1 & 0 \end{bmatrix} \cup \begin{bmatrix} 0 & 0 & 1 & 0 & 0 \\ 0 & 0 & 0 & 1 & 0 \\ 0 & 0 & 0 & 0 & 1 \\ 1 & 0 & 0 & 0 & 0 \\ 0 & 1 & 0 & 0 & 0 \end{bmatrix}.$$

We get yet another example of bigraph which is connected by a vertex $V(G_1) = \left(p_1^1, p_2^1, p_3^1, p_4^1\right)$; $V(G_2) = (p_1, p_2, p_3, p_4, p_5)$ the vertex $p_1 = p_3^1$.

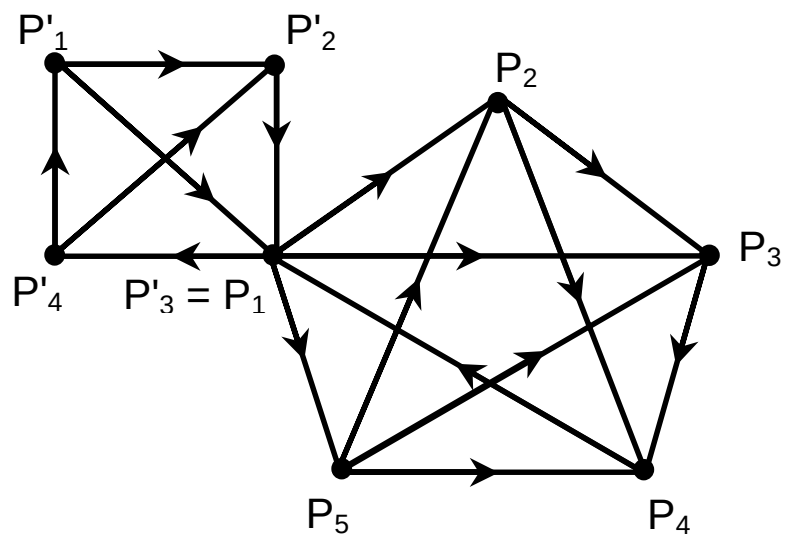

FIGURE: 3.5.5

Now we get the related vertex bimatrix $M = M_1 \cup M_2$

$$\begin{bmatrix} 0 & 1 & 0 & 0 \\ 0 & 0 & 1 & 0 \\ 0 & 0 & 0 & 1 \\ 1 & 0 & 0 & 0 \end{bmatrix} \cup \begin{bmatrix} 0 & 1 & 1 & 0 & 0 \\ 0 & 0 & 1 & 1 & 0 \\ 0 & 0 & 0 & 1 & 1 \\ 1 & 0 & 0 & 0 & 1 \\ 1 & 1 & 0 & 0 & 0 \end{bmatrix}.$$

Study in this direction is very interesting and innovative. For several properties like dominance of a directed graph, power of a bivertex of dominance etc can only be analysed using bimatrices.

Further bimatrices play a role in case of (labelled) bigraphs for the bigraphs can be given its adjacency bimatrix representation which can be used to analyze the number of walk and so on.

***Example 3.5.3:*** Let $G = G_1 \cup G_2$ be the bigraph given by the following figure:

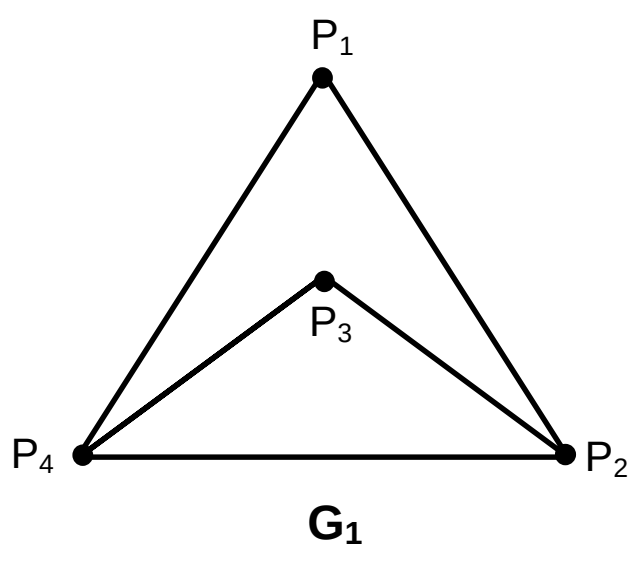

FIGURE: 3.5.6

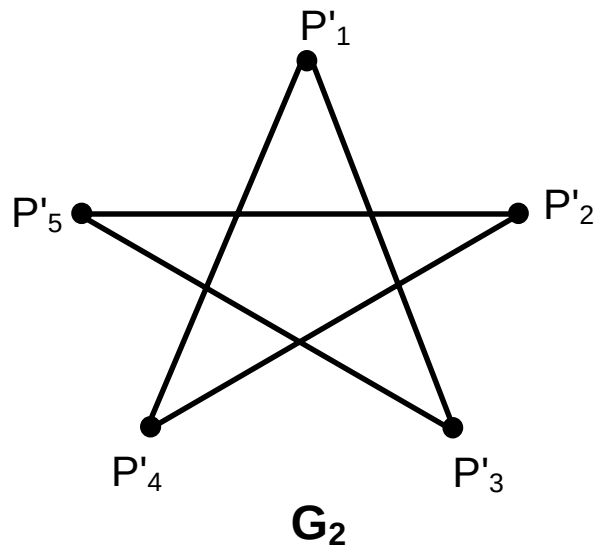

FIGURE: 3.5.7

The related adjacency bimatrix is given by

$$\begin{bmatrix} 0 & 1 & 0 & 1 \\ 1 & 0 & 1 & 1 \\ 0 & 1 & 0 & 1 \\ 1 & 1 & 1 & 0 \end{bmatrix} \cup \begin{bmatrix} 0 & 0 & 1 & 1 & 0 \\ 0 & 0 & 0 & 1 & 1 \\ 1 & 0 & 0 & 0 & 1 \\ 1 & 1 & 0 & 0 & 0 \\ 0 & 1 & 1 & 0 & 0 \end{bmatrix}.$$

Clearly one can see the vertex bimatrix of a directed bigraph is different from the adjacency bimatrix of a bigraph.

The applications of bigraphs have been elaborately discussed in [ ]. We just give one more example of the adjacency bimatrix before we proceed on to give other applications of linear bialgebra.

***Example 3.5.4:***

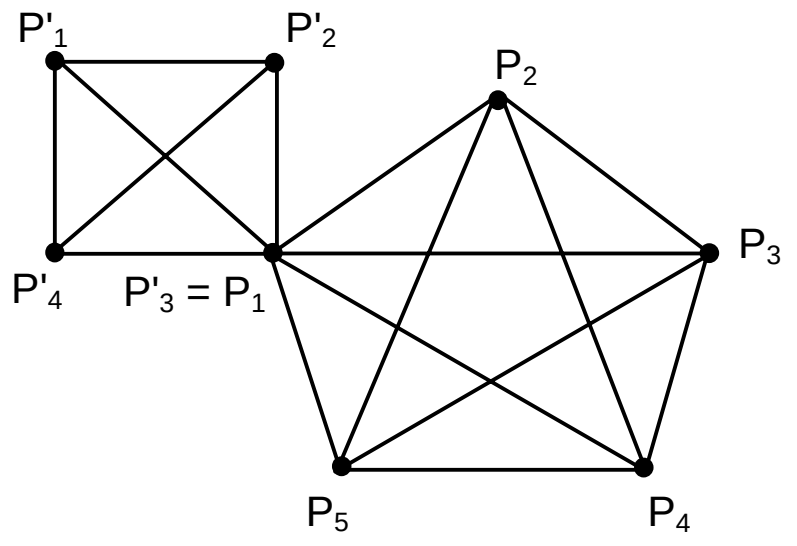

FIGURE: 3.5.8

The adjacency bimatrix of the above bigraph is given by M = $M_1 \cup M_2$.

$$\begin{bmatrix} 0 & 1 & 1 & 1 \\ 1 & 0 & 1 & 1 \\ 1 & 1 & 0 & 1 \\ 1 & 1 & 1 & 0 \end{bmatrix} \cup \begin{bmatrix} 0 & 1 & 1 & 1 & 1 \\ 1 & 0 & 1 & 1 & 1 \\ 1 & 1 & 0 & 1 & 1 \\ 1 & 1 & 1 & 0 & 1 \\ 1 & 1 & 1 & 1 & 0 \end{bmatrix} = M = M_1 \cup M_2.$$

The bigraph G whose lines and cycles are labelled can be given a cycle bimatrix representation.
For example consider the bigraph given by $G = G_1 \cup G_2$.

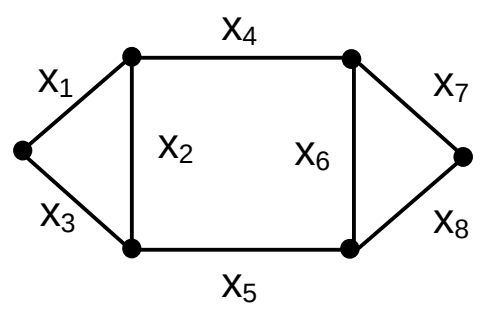

FIGURE: 3.5.9

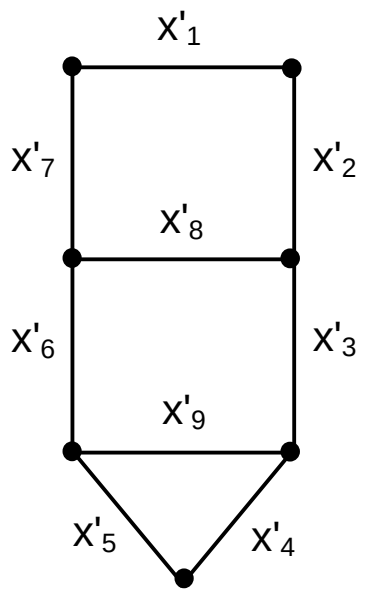

FIGURE: 3.5.10

The bicycles of $G_1 \cup G_2$ is (the cycles of $G_1$) $\cup$ (cycles of $G_2$). Cycle of $G_1$ are as follows:

$$
\begin{aligned}
Z_1 &= \{x_1, x_2, x_3\} \\
Z_2 &= \{x_2, x_4, x_5, x_6\} \\
Z_3 &= \{x_6, x_7, x_8\} \\
Z_4 &= \{x_1, x_3, x_4, x_5, x_6\} \\
Z_5 &= \{x_2, x_4, x_5, x_7, x_8\} \\
Z_6 &= \{x_1, x_3, x_4, x_5, x_7, x_8\}.
\end{aligned}
$$

Now cycles of $G_2$ are

$$
\begin{aligned}
Z_1^1 &= \{x_1^1\ x_2^1\ x_8^1\ x_7^1\} \\
Z_2^1 &= \{x_3^1\ x_6^1\ x_8^1\ x_9^1\} \\
Z_3^1 &= \{x_5^1\ x_9^1\ x_4^1\} \\
Z_4^1 &= \{x_1^1\ x_2^1\ x_3^1\ x_9^1\, x_6^1\ x_7^1\} \\
Z_5^1 &= \{x_1^1\ x_2^1\ x_3^1\ x_4^1\, x_5^1\ x_6^1\, x_7^1\} \\
Z_6^1 &= \{x_8^1\ x_3^1\ x_4^1\ x_5^1\, x_6^1\}.
\end{aligned}
$$

The bicycle of the bigraph $C = C_1 \cup C_2$ where $C_1 = \{z_1, z_2, z_3, z_4, z_5, z_6\}$ and $C_2 = \left\{z_1^1, z_2^1, z_3^1, z_4^1, z_5^1, z_6^1\right\}$.

Now we describe the bicycle bimatrix of the bigraph $G = G_1 \cup G_2$.

$$\begin{array}{c} \quad x_1 \;\; x_2 \;\; x_3 \;\; x_4 \;\; x_5 \;\; x_6 \;\; x_7 \;\; x_8 \\ \begin{bmatrix} 1 & 1 & 1 & 0 & 0 & 0 & 0 & 0 \\ 0 & 1 & 0 & 1 & 1 & 1 & 0 & 0 \\ 0 & 0 & 0 & 0 & 0 & 1 & 1 & 1 \\ 1 & 0 & 1 & 1 & 1 & 1 & 0 & 0 \\ 0 & 1 & 0 & 1 & 1 & 0 & 1 & 1 \\ 1 & 0 & 1 & 1 & 1 & 0 & 1 & 1 \end{bmatrix} \begin{array}{c} z_1 \\ z_2 \\ z_3 \\ z_4 \\ z_5 \\ z_6 \end{array} \end{array}$$

$$\cup \quad \begin{array}{c} \quad x_1^1 \;\; x_2^1 \;\; x_3^1 \;\; x_4^1 \;\; x_5^1 \;\; x_6^1 \;\; x_7^1 \;\; x_8^1 \;\; x_9^1 \\ \begin{bmatrix} 1 & 1 & 0 & 0 & 0 & 0 & 1 & 1 & 0 \\ 0 & 0 & 1 & 0 & 0 & 1 & 0 & 1 & 1 \\ 0 & 0 & 0 & 1 & 1 & 0 & 0 & 0 & 1 \\ 1 & 1 & 1 & 0 & 0 & 1 & 1 & 0 & 1 \\ 1 & 1 & 1 & 1 & 1 & 1 & 1 & 0 & 0 \\ 0 & 0 & 1 & 1 & 1 & 1 & 0 & 1 & 0 \end{bmatrix} \begin{array}{c} z_1 \\ z_2 \\ z_3 \\ z_4 \\ z_5 \\ z_6 \end{array} \end{array}$$

Several properties of the bigraph can be analysed using the bicycle bimatrix.

### 3.6 Jordan biform

Let us now proceed to define the notion of Jordan biform. The notion of Jordan form in each equivalence class of matrices under similarity and so it has been proved that to each endomorphism of a finite dimensional Vector space a unique matrix is attached. This work was developed by Camille Jordan, a late nineteenth century mathematician.

Let A be a n × n matrix defined over the field K with characteristic polynomial p (x). Suppose the characteristic polynomial

$$p(x) = (x - \alpha_1)^{d_1} \ldots (x - \alpha_k)^{d_k}$$

where $\alpha_i \in K$, $1 \leq i \leq K$ and $d_1 + \ldots + d_k = n$. We know by spectral theorem there are subspaces $V_1, \ldots, V_k$ of $K^n$ such that minimum polynomial for A on $V_i$ is $(x - \alpha_i)^{d_i}$ for some $1 \leq d_i^1 \leq m_i, 1 \leq i \leq K$ and

$$V \cong K^n = V_1 \oplus \ldots \oplus V_k = V \text{ and } \dim V_i = r_i \text{ with } r_1 + \ldots + r_k = n.$$

So now we define the notion of Jordan biform for bivector spaces defined over a field. Let $V = V_1 \cup V_2$ be a bivector space of dimension (m, n) over K i.e. $\dim V_1 = m$ and $\dim V_2 = n$. Let $A = A_1 \cup A_2$ be a bimatrix defined over K where $A_1$ is a $m \times m$ matrix and $A_2$ is a $n \times n$ matrix. Let $p^B(x) = p_1(x) \cup p_2(x)$ be the characteristic bipolynomial associated with $A = A_1 \cup A_2$.

Suppose $\alpha_i = \alpha_i^1 \cup \alpha_i^2$ be the eigen bivalues associated with $p^B(x)$. Let us assume the characteristic bipolynomial has the form i.e.

$$p^B(x) = (x - \alpha_1^1)^{d_1^1} \cdots (x - \alpha_{K_1}^1)^{d_{k_1}^1} \cup (x - \alpha_1^2)^{d_1^2} \cdots (x - \alpha_{K_2}^2)^{d_{k_2}^2}$$

where $(\alpha_1^1, \ldots, \alpha_{k_1}^1)$ are the distinct eigen values associated with $A_1$ and $(\alpha_1^2, \ldots, \alpha_{k_2}^2)$ are the distinct eigen values associated with $A_2$ such that $d_1^1 + \ldots + d_{k_1}^1 = m$ and $d_{k_1}^2 + \ldots + d_{k_2}^2 = n$.

There are subbispaces $(V_1^1, \ldots, V_{k_1}^1) \cup (V_1^2, \ldots, V_{k_2}^2)$ such that

$$V = V_1 \cup V_2 = (V_1^1 \oplus \ldots \oplus V_k^1) \cup (V_1^2 \oplus \ldots \oplus V_{k_2}^2)$$

and the minimal bipolynomial for A on

$$V_i = V_i^1 \cup V_i^2 \text{ is } (x - \alpha_i)^{t_i} = (x - \alpha_i^1)^{t_i^1} \cup (x - \alpha_i^2)^{t_i^2}$$

for some $t_i = t_i^1 \cup t_i^2$ satisfying

$$1 \leq t_i^1 \leq d_i^1$$

and

$$1 \leq t_i^2 \leq d_i^2.$$

So using similar results as in case of Jordan forms we can define for any linear bioperator T on $V = V_1 \cup V_2$ the $\phi$ - stable (i.e., $\phi_1$ stable $\cup$ $\phi_2$ stable-subspaces) subbispaces given by the bispectral theorem and get the related bimatrix of T in the form,

$$A = A_1 \cup A_2 = \begin{pmatrix} c_1^1 & 0 & \cdots & 0 \\ 0 & c_2^1 & & 0 \\ & & & \\ & \cdots & & c_{k_1}^1 \end{pmatrix} \cup \begin{pmatrix} c_1^2 & 0 & \cdots & 0 \\ 0 & c_2^2 & & 0 \\ & & & \\ & \cdots & & c_{k_2}^2 \end{pmatrix}.$$

where for each i = 1, 2,…, $k_1$ and j = 1, 2, …, $k_2$, there exist a bibasis $B = B_1 \cup B_2$ of $V = V_1 \cup V_2$ where each $c_t = c_{t_1}^1 \cup c_{t_2}^2$ is a $\left(r_{t_1}^1 \times r_{t_1}^1\right), \left(r_{t_2}^2 \times r_{t_2}^2\right)$ bimatrix and all other entries are zero. Also $T = T_1 \cup T_2$ restricted to $c_t = c_t^1 \cup c_t^2$ relative to the bibasis $B = B_1 \cup B_2$.

***Example 3.6.1:*** Consider a bivector space $V = V_1 \cup V_2$ over $\mathbb{C}$, where $V_1$ is of dimension 3 and $V_2$ is of dimension 4. Let $T = T_1 \cup T_2$ be a linear bioperator on V with related mixed square bimatrix,

$$A = \begin{bmatrix} 2 & 0 & 0 \\ 1 & 2 & 0 \\ 0 & 0 & -1 \end{bmatrix} \cup \begin{bmatrix} 2 & 0 & 0 & 0 \\ 1 & 2 & 0 & 0 \\ 0 & 0 & 2 & 0 \\ 0 & 0 & 1 & 2 \end{bmatrix}.$$

A has the following characteristic bipolynomial

$$\begin{aligned} p(x_1, x_2) &= p(x) \\ &= p_1(x_1) \cup p_2(x_2) \\ &= (x_1 - 2)^2 (x_1 + 1) \cup (x_2 - 2)^4. \end{aligned}$$

The minimal bipolynomial

$$\begin{aligned} q(x_1, x_2) &= q(x) \\ &= q_1(x_1) \cup q_2(x_2) \\ &= (x_1 - 2)(x_1 + 1) \cup (x_2 - 2)^2. \end{aligned}$$

The bimatrix $A = A_1 \cup A_2$ is the Jordan biform.

Now we proceed on to define the Jordan block bimatrix.

**DEFINITION 3.6.1:** *The $(n \times n)$ square bimatrix $A = A_1 \cup A_2$ is called the Jordan block bimatrix, if A is of the form*

$$\begin{bmatrix} e_1 & 1 & 0 & 0 & \dots & 0 \\ 0 & e_1 & 1 & 0 & \dots & \vdots \\ 0 & & \dots & & \dots & 1 \\ 0 & & \dots & & \dots & e_1 \end{bmatrix} \cup \begin{bmatrix} e_2 & 1 & 0 & 0 & \dots & 0 \\ 0 & e_2 & 1 & 0 & \dots & \vdots \\ 0 & & \dots & & \dots & 1 \\ 0 & & \dots & & \dots & e_2 \end{bmatrix}.$$

*That the bimatrix has the scalar $e = e_1 \cup e_2$ at each main diagonal entry and 1 at each entry on the super diagonal and zero at all other entries.*

*A mixed square bimatrix $A = A_1 \cup A_2$ is called the Jordan block bimatrix of A and is of the form*

$$A = \begin{bmatrix} e_1 & 1 & 0 & 0 & \dots & 0 \\ 0 & e_1 & 1 & 0 & \dots & \vdots \\ 0 & & \dots & & \dots & 1 \\ 0 & & \dots & & \dots & e_1 \end{bmatrix}_{m \times m} \cup \begin{bmatrix} e_2 & 1 & 0 & 0 & \dots & 0 \\ 0 & e_2 & 1 & 0 & \dots & \vdots \\ 0 & & \dots & & \dots & 1 \\ 0 & & \dots & & \dots & e_2 \end{bmatrix}_{n \times n}$$

We just give an example of a Jordan block bimatrix.

***Example 3.6.2:*** Let $A = A_1 \cup A_2$ .

$$A = \begin{bmatrix} 3 & 1 & 0 \\ 0 & 3 & 1 \\ 0 & 0 & 3 \end{bmatrix} \cup \begin{bmatrix} 2 & 1 & 0 & 0 \\ 0 & 2 & 1 & 0 \\ 0 & 0 & 2 & 1 \\ 0 & 0 & 0 & 2 \end{bmatrix}$$

A is the mixed square Jordan block bimatrix.

*Note:* It is sometimes also defined as the elementary Jordan bimatrix. It can be a super diagonal or sub diagonal. For even a square bimatrix of the form $A = A_1 \cup A_2$ where

$$A = \begin{bmatrix} c_1 & 0 & & \cdots & 0 \\ 1 & c_1 & 0 & & 0 \\ 0 & 1 & c_1 & 0 & \vdots \\ 0 & & 1 & c_1 & 0 \\ 0 & \cdots & & 1 & c_1 \end{bmatrix} \cup \begin{bmatrix} c_2 & 0 & & \cdots & 0 \\ 1 & c_2 & 0 & & 0 \\ 0 & 1 & c_2 & 0 & \vdots \\ 0 & & 1 & c_2 & 0 \\ 0 & \cdots & 0 & 1 & c_2 \end{bmatrix}$$

$$= \quad J(c_1) \cup J(c_2).$$

Thus this is also a Jordan bimatrix with characteristic bivalue $c = c_1 \cup c_2$. Now if the characteristic bivalues take $c_1, \ldots, c_k$, with $c_i = c_i^1 \cup c_i^2$ , $1 \leq i \leq K$ then the Jordan biform

$$A = A_1 \cup A_2$$

$$= \begin{bmatrix} J(c_1^1) & \cdots & 0 \\ \vdots & J(c_1^1) & \vdots \\ 0 & \cdots & J(c_k^1) \end{bmatrix} \cup \begin{bmatrix} J(c_1^2) & \cdots & 0 \\ \vdots & J(c_2^2) & \vdots \\ 0 & \cdots & J(c_k^2) \end{bmatrix}.$$

***Example 3.6.3:*** Let $A = A_1 \cup A_2$ be a bimatrix where $A_1$ is a $8 \times 8$ matrix and $A_2$ is a $6 \times 6$ matrix. The Jordan biform for A is given by

$$\begin{bmatrix} 3 & 0 & 0 & 0 & 0 & 0 & 0 & 0 \\ 1 & 3 & 0 & 0 & 0 & 0 & 0 & 0 \\ 0 & 1 & 3 & 0 & 0 & 0 & 0 & 0 \\ 0 & 0 & 0 & 2 & 0 & 0 & 0 & 0 \\ 0 & 0 & 0 & 1 & 2 & 0 & 0 & 0 \\ 0 & 0 & 0 & 0 & 0 & -1 & 0 & 0 \\ 0 & 0 & 0 & 0 & 0 & 1 & -1 & 0 \\ 0 & 0 & 0 & 0 & 0 & 0 & 1 & -1 \end{bmatrix} \cup \begin{bmatrix} 4 & 0 & 0 & 0 & 0 & 0 \\ 1 & 4 & 0 & 0 & 0 & 0 \\ 0 & 1 & 4 & 0 & 0 & 0 \\ 0 & 0 & 1 & 4 & 0 & 0 \\ 0 & 0 & 0 & 0 & 3 & 0 \\ 0 & 0 & 0 & 0 & 1 & 3 \end{bmatrix}$$

$$= \begin{bmatrix} J(3) & 0 & 0 & 0 & 0 & 0 \\ 0 & 0 & J(2) & 0 & 0 & 0 \\ 0 & 0 & 0 & 0 & 0 & J(-1) \end{bmatrix} \cup \begin{bmatrix} J(4) & 0 \\ 0 & J(3) \end{bmatrix}$$

All related results of Jordan form can be derived and proved in case of Jordan biform.

## 3.7 Application of bivector spaces to bicodes

Now for the first time we give the results of bimatrices in the field of coding theory. Very recently we have defined the notion of bicodes here we give the applications of linear bicodes and finite bimatrix bigroups. [51]

A bicode is a pair of sets of markers or labels built up from a finite "alphabet" which can be attached to some or all of the entities of a system of study. A mathematical bicode is a bicode with some mathematical structure. In many cases as in case of codes the mathematical structure of a bicodes refers to the fact that it forms a bivector space over finite field and in this case the bicode is linear. Let $Z_q(m)$ =

$Z_{q_1}(m_2) \cup Z_{q_2}(m_2)$. Denote the m-dimensional linear bicode in the bivector space $V = V_1 \cup V_2$ where $V_1$ is a $m_1$ dimensional vector space over $Z_p$ where $q_1$ is a power of the prime p and $V_2$ is a $m_2$ – dimensional vector space over $Z_p$ where $q_2$ is a power of the prime p, then $Z_q$ (m) = $Z_{q_1}(m_1) \cup Z_{q_2}(m_2)$ is a bicode over $Z_p$.

Thus a bicode can also be defined as a 'union' of two codes $C_1$ and $C_2$ where union is just the symbol.

For example we can have a bicode over $Z_2$.

***Example 3.7.1:*** Let $C = C_1 \cup C_2$ be a bicode over $Z_2$ given by

$C = C_1 \cup C_2$
= {( 0 0 0 0 0 0), (0 1 1 0 1 1 ), (1 1 0 1 1 0), (0 0 1 1 1 0), (1 0 0 0 1 1), (1 1 1 0 0 0), (0 1 0 1 0 1), (1 0 1 1 0 1)} $\cup$ {(0 0 0 0), (1 1 1 0), (1 0 0 1), (0 1 1 1), (0 1 0 1), (0 0 1 0), (1 1 0 0), (1 0 1 1)}

over $Z_2$.

These codes are subbispaces of the bivector space over $Z_2$ of dimension (6, 4). Now the bicodes are generated by bimatrices and we have parity check bimatrices to play a role in finding the check symbols. Thus we see we have the applications of linear bialgebra / bivector spaces in the study of bicodes.

$$C (n_1 \cup n_2, k_1, k_2) = C_1 (n_1, k_1) \cup C_2 (n_2, k_2)$$

is a linear bicode if and only if both $C_1$ $(n_1, k_1)$ and $C_2$ $(n_2, k_2)$ are linear codes of length $n_1$ and $n_2$ with $k_1$ and $k_2$ message symbols respectively with entries from the same field $Z_q$. The check symbols can be obtained from the $k_1$ and $k_2$ messages in such a way that the bicode words $x = x^1 \cup x^2$ satisfy the system of linear biequations.

i.e. $Hx^T = (0)$

i.e. $(H_1 \cup H_2) (x^1 \cup x^2)^T = (0) \cup (0)$

i.e. $H_1 (x^1)^T \cup H_2 (x^2)^T = (0) \cup (0)$

where $H = H_1 \cup H_2$ is a given mixed bimatrix of order ($n_1 - k_1 \times n_1$, $n_2 - k_2 \times n_2$) with entries from the same field $Z_q$, $H = H_1 \cup H_2$ is the parity check bimatrix of the bicode C.

The standard form of H is $(A_1, I_{n_1-k_1}) \cup (A_2, I_{n_2-k_2})$ with $A_1$ a $n_1 - k_1 \times k_1$ matrix and $A_2$ a $n_2 - k_2 \times k_2$ matrix. $I_{n_1-k_1}$ and $I_{n_2-k_2}$ are $n_1 - k_1 \times n_1 - k_1$ and $n_2 - k_2 \times n_2 - k_2$ identity matrices respectively. The bimatrix H is called the parity check bimatrix of the bicode $C = C_1 \cup C_2$. C is also called the linear $(n, k) = (n_1 \cup n_2, k_1 \cup k_2) = (n_1, k_1) \cup (n_2, k_2)$ bicode.

***Example 3.7.2:*** Let $C(n, k) = C(6 \cup 7, 3 \cup 4) = C_1(6, 3) \cup C_2(7, 4)$ be a bicode got by the parity check bimatrix, $H = H_1 \cup H_2$ where

$$H = \begin{pmatrix} 0 & 1 & 1 & 1 & 0 & 0 \\ 1 & 0 & 1 & 0 & 1 & 0 \\ 1 & 1 & 0 & 0 & 0 & 1 \end{pmatrix} \cup \begin{pmatrix} 1 & 1 & 1 & 0 & 1 & 0 & 0 \\ 0 & 1 & 1 & 1 & 0 & 1 & 0 \\ 1 & 1 & 0 & 1 & 0 & 0 & 1 \end{pmatrix}.$$

The bicodes are obtained by solving the equations.

$$Hx^T = H_1 x_1^T \cup H_2 x_2^T = (0) \cup (0).$$

There are $2^3 \cup 2^4$ bicode words given by

{(0 0 0 0 0 0), (0 1 1 0 1 1), (1 1 0 1 1 0), (0 0 1 1 1 1 0), (1 0 0 0 1 1), (1 1 1 0 0 0), (0 1 0 1 0 1), (1 0 1 1 0 1)} ∪ {(0 0 0 0 0 0 0), (1 0 0 0 1 0 1), (0 1 0 0 1 1 1), (0 0 1 0 1 1 0), (0 0 0 1 0 1 1), (1 1 0 0 0 1 0), (1 0 1 0 0 1 1), (1 0 0 1 1 1 1), (0 1 1 0 0 0 1), (0 1 0 1 1 0 0), (0 0 1 1 1 0 1), (1 1 1 0 1 0 0), (1 1 0 1 0 0 1), (1 0 1 1 0 0 0), (0 1 1 1 0 1 0), (1 1 1 1 1 1 1)}.

Clearly this is a bicode over $Z_2$. Now the main advantage of a bicode is that at a time two codes of same length or of different length are sent simultaneously and the resultant can be got.

As in case of codes if $x = x_1 \cup x_2$ is a sent message and $y = y_1 \cup y_2$ is a received message using the parity check bimatrix $H = H_1 \cup H_2$ one can easily verify whether the received message is a correct one or not.

For if we consider $Hx^T = H_1 x_1^T \cup H_2 x_2^T$ then $Hx^T = (0) \cup (0)$ for $x = x_1 \cup x_2$ is the bicode word which was sent. Let $y = y_1 \cup y_2$ be the received word. Consider $Hy^T = H_1 y_1^T \cup H_2 y_2^T$ if $Hy^T = (0) \cup (0)$ then the received bicode word is a correct bicode if $Hy^T \neq (0)$, then we say there is error in the received bicode word. So we define by syndrome $S^B(y) = Hy^T$ for any bicode $y = y_1 \cup y_2$.

Thus the bisyndrome

$$S^B(y) = S_1(y_1) \cup S_2(y_2) = H_1 y_1^T \cup H_2 y_2^T .$$

If the bisyndrome $S^B(y) = (0)$ then we say y is the bicode word. If $S^B(y) \neq (0)$ then the word y is not a bicode word. This is the fastest and the very simple means to check whether the received bicode word is a bicode word or not.

Now we proceed on to define the notion of bigenerator bimatrix or just the generator bimatrix of a bicode $C(n, k) = C_1(n_1, k_1) \cup C_2(n_2, k_2)$.

**DEFINITION 3.7.1:** *The generator bimatrix* $G = (I_k, -A^T) = G_1 \cup G_2 = \left(I_{k_1}^1 - A_1^T\right) \cup \left(I_{k_2}^2 - A_2^T\right)$ *is called the canonical generator bimatrix or canonical basic bimatrix or encoding bimatrix of a linear bicode,* $C(n, k) = C(n_1 \cup n_2, k_1 \cup k_2)$ *with parity check bimatrix;* $H = H_1 \cup H_2 = (A_1, I_{n_1-k_1}^1) \cup (A_2, I_{n_2-k_2}^2)$. *We have* $GH^T = (0)$ *i.e.* $G_1 H_1^T \cup G_2 H_2^T = (0) \cup (0)$.

We now illustrate by an example how a generator bimatrix of a bicode functions.

***Example 3.7.3:*** Consider $C(n, k) = C_1(n_1, k_1) \cup C_2(n_2, k_2)$ a bicode over $Z_2$ where $C_1(n_1, k_1)$ is a (7, 4) code and $C_2(n_2, k_2)$ is a (9, 3) code given by the generator bimatrix

$G = G_1 \cup G_2 =$

$$\begin{pmatrix} 1 & 0 & 0 & 0 & 1 & 0 & 1 \\ 0 & 1 & 0 & 0 & 1 & 1 & 1 \\ 0 & 0 & 1 & 0 & 1 & 1 & 0 \\ 0 & 0 & 0 & 1 & 0 & 1 & 1 \end{pmatrix} \cup \begin{pmatrix} 1 & 0 & 0 & 1 & 0 & 0 & 1 & 0 & 0 \\ 0 & 1 & 0 & 0 & 1 & 0 & 0 & 1 & 0 \\ 0 & 0 & 1 & 0 & 0 & 1 & 0 & 0 & 1 \end{pmatrix}.$$

One can obtain the bicode words by using the rule $x_1 = aG$ where a is the message symbol of the bicode. $x_1 \cup x_2 = a^1G_1 \cup a^2G_2$ with $a = a^1 \cup a^2$.
where

$$a^1 = a_1^1 a_2^1 a_3^1$$

and

$$a^2 = a_1^2\ a_2^2\ a_3^2\ a_4^2\ a_5^2\ a_6^2 .$$

We give another example in which we calculate the bicode words.

***Example 3.7.4:*** Consider a $C(n, k) = C_1(6, 3) \cup C_2(4, 2)$ bicode over $Z_2$.

The $2^3 \cup 2^2$ code words $x^1$ and $x^2$ of the binary bicode can be found using the generator bimatrix $G = G_1 \cup G_2$.

$$G = \begin{pmatrix} 1 & 0 & 0 & 0 & 1 & 1 \\ 0 & 1 & 0 & 1 & 0 & 1 \\ 0 & 0 & 1 & 1 & 1 & 0 \end{pmatrix} \cup \begin{pmatrix} 1 & 0 & 1 & 1 \\ 0 & 1 & 0 & 1 \end{pmatrix}.$$

$x = aG$

$x^1 \cup x^2 = a^1G_1 \cup a^2G_2$ where $a^1 = a_1^1\ a_2^1\ a_3^1$ and $a^2 = a_1^2\, a_2^2$ where the message symbols are

$$\begin{Bmatrix} 000, & 001, & 100, & 010, \\ 110 & 011 & 101 & 111 \end{Bmatrix} \cup \begin{Bmatrix} 01 & 10 \\ 00 & 11 \end{Bmatrix}.$$

Thus $x = x^1 \cup x^2 = a^1G_1 \cup a^2G_2$.
We get the bicodes as follows:

$$\begin{Bmatrix} 000000 & 011011 & 110110 & 001110 \\ 100011 & 111000 & 010101 & 101101 \end{Bmatrix} \cup \begin{Bmatrix} 0000, & 1010 \\ 0101, & 1110 \end{Bmatrix}$$

Now we proceed on to define the notion of repetition bicode and parity check bicode.

**DEFINITION 3.7.2:** *If each bicode word of a bicode consists of only one message symbol $a = a_1 \cup a_2 \in F_2$ and the $(n_1 - 1) \cup (n_2 - 1)$ check symbols and $x_2^2 = .... = x_{n_1}^1 = a_1$ and $x_2^2 = x_3^2 = ... = x_{n_2}^2 = a_2$, $a_1$ repeated $n_1 - 1$ times and $a_2$ is repeated $n_2 - 1$ times. We obtain the binary $(n, 1) = (n_1, 1) \cup (n_2, 1)$ repetition bicode with parity check bimatrix, $H = H_1 \cup H_2$*

$$\text{i.e., } H = \begin{pmatrix} 110 & \cdots & 0 \\ 101 & \cdots & 0 \\ \vdots & & \vdots \\ 100 & & 1 \end{pmatrix}_{n_1-1 \times n_1} \cup \begin{pmatrix} 110 & \ldots & 0 \\ 101 & & 0 \\ \vdots & \ldots & \vdots \\ 100 & & 1 \end{pmatrix}_{n_2-1 \times n_2}$$

There are only 4 bicode words in a repetition bicode.
$\{(1\ 1\ 1\ 1\ldots1), (0\ 0\ 0\ldots0)\} \cup \{(1\ 1\ 1\ 1\ 1\ 1\ldots1), (0\ 0\ \ldots\ 0)\}$

***Example 3.7.5:*** Consider the Repetition bicode $(n, k) = (5, 1) \cup (4, 1)$. The parity check bimatrix is given by

$$\begin{pmatrix} 1 & 1 & 0 & 0 & 0 \\ 1 & 0 & 1 & 0 & 0 \\ 1 & 0 & 0 & 1 & 0 \\ 1 & 0 & 0 & 0 & 1 \end{pmatrix}_{4 \times 5} \cup \begin{pmatrix} 1 & 1 & 0 & 0 \\ 1 & 0 & 1 & 0 \\ 1 & 0 & 0 & 1 \end{pmatrix}_{3 \times 4}.$$

The bicode words are
{(1 1 1 1 1), (0 0 0 0 0)} ∪ {(1 1 1 1), (0 0 0 0)}.

We now proceed on to define Parity-check bicode.

**DEFINITION 3.7.3:** *Parity check bicode is a (n, n-1) = ($n_1$, $n_1$-1) ∪ ($n_2$, $n_2$ – 1) bicode with parity check bimatrix H =* $\left(1\,1\,...\,1\right)_{1\times n_1} \cup \left(1\,1\,...\,1\right)_{1\times n_2}$.

*Each bicode word has one check symbol and all bicode words are given by all binary bivectors of length n = $n_1$ ∪ $n_2$ with an even number of ones. Thus if the sum of ones of a received bicode word is 1 at least an error must have occurred at the transmission.*

***Example 3.7.6:*** Let (n, k) = (4, 3) ∪ (5, 4) be a bicode with parity check bimatrix.

$$H = (1\ 1\ 1\ 1) \cup (1\ 1\ 1\ 1\ 1) = H_1 \cup H_2.$$

The bicodes related with the parity check bimatrix H is given by

$$\begin{Bmatrix} 0000 & 1001 \\ 0101 & 0011 \\ 1100 & 1010 \\ 0110 & 1111 \end{Bmatrix} \cup$$

$$\begin{Bmatrix} 00000 & 00011 & 11000 & 10100 \\ 10001 & 01100 & 00110 & 01010 \\ 01001 & 11011 & 10110 & 01111 \\ 00101 & 10010 & 11101 & 11110 \end{Bmatrix}.$$

Now the concept of linear bialgebra i.e. bimatrices are used when we define cyclic bicodes. So we define a cyclic bicode as the union of two cyclic codes for we know if C (n, k) is a cyclic code then if ($x_1 \ldots x_n$) ∈ C (n, k) it implies ($x_n$ $x_1 \ldots x_{n-1}$) belongs to C (n, k).

We just give an example of a cyclic code using both generator bimatrix and the parity check bimatrix or to be more precise using a generator bipolynomial and a parity check bipolynomial.

Suppose C (n, k) = $C_1$ ($n_1$, $k_1$) $\cup$ $C_2$ ($n_2$, $k_2$) be a cyclic bicode then we need generator polynomials $g_1(x) \mid x^{n_1} - 1$ and $g_2(x) \mid x^{n_2} - 1$ where $g(x) = g_1(x) \cup g_2(x)$ is a generator bipolynomial with degree of g(x) = ($m_1$, $m_2$) i.e. degree of $g_1(x) = m_1$ and degree of $g_2(x) = m_2$. The C(n, k) linear cyclic bicode is generated by the bimatrix;
$G = G_1 \cup G_2$ which is given below

$$\begin{pmatrix} g_0^1 & g_1^1 & g_{m_1}^1 & 0 & 0 \\ 0 & g_0^1 & g_{m_1-1}^1 & g_{m_1}^1 & 0 \\ \vdots & & & & \\ 0 & 0 & g_0^1\, g_1^1 & \cdots & g_{m_1}^1 \end{pmatrix}$$

$$\cup \begin{pmatrix} g_0^2 & g_1^2 & g_{m_2}^2 & 0 & 0 \\ 0 & g_0^2 & g_{m_2-1}^2 & g_{m_2}^2 & 0 \\ \vdots & & & & \\ 0 & 0 & g_0^2 & \cdots & g_{m_2}^2 \end{pmatrix}$$

$$= \begin{pmatrix} g_1 \\ xg_1 \\ \vdots \\ x^{k_1-1}g_1 \end{pmatrix} \cup \begin{pmatrix} g_2 \\ xg_2 \\ \vdots \\ x^{k_2-1}g_2 \end{pmatrix}$$

Then C (n, k) is cyclic bicode.

***Example 3.7.7:*** Suppose $g = g_1 \cup g_2 = 1 + x^3 \cup 1 + x^2 + x^3$ be the generator bipolynomial of a C (n, k) (= $C_1$ (6, 3) $\cup$ $C_2$ (7, 3)) bicode. The related generator bimatrix is given by

$$G = G_1 \cup G_2$$

i.e., $\begin{pmatrix} 1 & 0 & 0 & 1 & 0 & 0 \\ 0 & 1 & 0 & 0 & 1 & 0 \\ 0 & 0 & 1 & 0 & 0 & 1 \end{pmatrix} \cup \begin{pmatrix} 1 & 0 & 1 & 1 & 0 & 0 & 0 \\ 0 & 1 & 0 & 1 & 1 & 0 & 0 \\ 0 & 0 & 1 & 0 & 1 & 1 & 0 \\ 0 & 0 & 0 & 1 & 0 & 1 & 1 \end{pmatrix}$

Clearly the cyclic bicode is given by

$$\begin{Bmatrix} 000000 & 100100 \\ 001001 & 101101 \\ 010010 & 110110 \\ 011011 & 111111 \end{Bmatrix} \cup$$

$$\begin{Bmatrix} 0000000 & 0001011 & 0110001 & 1101001 \\ 1000101 & 1100010 & 0101100 & 1011000 \\ 0100111 & 1010000 & 0011101 & 0111010 \\ 0010110 & 1001110 & 1110100 & 1111111 \end{Bmatrix}$$

It is easily verified that the bicode is a cyclic bicode. Now we define how a cyclic bicode is generated by a generator bipolynomial. Now we proceed on to give how the check bipolynomial helps in getting the parity check bimatrix. Let $g = g_1 \cup g_2$ be the generator bipolynomial of a bicode C (n, k) = $C_1$ ($n_1$, $k_1$) $\cup$ $C_2$ ($n_2$, $k_2$).

Then $h_1 = \dfrac{x^{n_1-1}}{g_1}$ and $h_2 = \dfrac{x^{n_2-1}}{g_2}$; $h = h_1 \cup h_2$ is called the check bipolynomial of the bicode C(n, k). The parity check bimatrix $H = H_1 \cup H_2$

$$= \begin{pmatrix} 0 & 0 & 0 & h^1_{k_1} & \cdots & h^1_1 & h^1_0 \\ 0 & 0 & h^1_{k_1} & h^1_{k_1-1} & \ldots & h^1_0 & 0 \\ \vdots & & & & \vdots & & \vdots \\ h^1_{k_1} & & \ldots & & h^1_0 & \ldots & 0 \end{pmatrix} \cup$$

$$\begin{pmatrix} 0 & 0 & 0 & h^2_{k_2} & \cdots & h^2_1 & h^2_0 \\ 0 & 0 & h^2_{k_2} & h^2_{k_2-1} & \dots & h^2_0 & 0 \\ \vdots & & & & \vdots & & \vdots \\ h^2_{k_2} & & \dots & & h^2_0 & \dots & 0 \end{pmatrix}.$$

From the above example the parity check bipolynomial $h = h_1 \cup h_2$ is given by

$$(x^3 + 1) \cup (x^4 + x^3 + x^2 + 1).$$

The parity check bimatrix associated with this bipolynomial is given by $G = G_1 \cup G_2$ where,

$$G = \begin{pmatrix} 0 & 0 & 1 & 0 & 0 & 1 \\ 0 & 1 & 0 & 0 & 1 & 0 \\ 1 & 0 & 0 & 1 & 0 & 0 \end{pmatrix} \cup \begin{pmatrix} 0 & 0 & 1 & 0 & 1 & 1 & 1 \\ 0 & 1 & 0 & 1 & 1 & 1 & 0 \\ 1 & 0 & 1 & 1 & 1 & 0 & 0 \end{pmatrix}.$$

The linear bialgebra concept is used for defining biorthogonal bicodes or dual bicodes.

Let $C(n, k) = C_1(n_1, k_1) \cup C_2(n_2, k_2)$ be a linear bicode over $Z_2$. The dual (or orthogonal) bicode $C^{\perp}$ of C is defined by

$$\begin{aligned} C^{\perp} &= \{u \mid u.v = 0, \forall\, v \in C\} \\ &= \{u_1 \mid u_1.v_1 = 0 \forall v_1 \in C_1\} \cup \{u_2 \mid u_2.v_2 = 0 \forall v_2 \in C_2\}. \end{aligned}$$

Thus in case of dual bicodes we have if for the bicode C (n, k), where $G = G_1 \cup G_2$ is the generator bimatrix and if its parity check bimatrix is $H = H_1 \cup H_2$ then for the orthogonal bicode or dual bicode, the generator bimatrix is $H_1 \cup H_2$ and the parity check bimatrix is $G = G_1 \cup G_2$.

Thus we have $HG^T = GH^T$.

## 3.8 Best biapproximation and its application

In this section we define the notion of best biapproximation and pseudo inner biproduct. Here we give their applications to bicodes.

The applications of pseudo best approximation have been carried out in (2005) [50] to coding theory. We just recall the definition of pseudo-best approximation.

**DEFINITION 3.8.1:** *Let V be a vector space defined over a finite field $Z_p$ with pseudo inner product $\langle\ \rangle_p$. Let W be the subspace of V. Let $\beta \in V$, the pseudo best approximation to $\beta$ related to W is defined as follows. Let $\{\alpha_1, \alpha_2, \ldots, \alpha_k\}$ be the chosen basis for W, the pseudo best approximation to $\beta$ if it exists is given by*

$$\sum_{i=1}^{k} \langle \beta, \alpha_i \rangle_\rho \ \alpha_i \neq 0 \ if \sum_{i=1}^{k} \langle \beta, \alpha_i \rangle_\rho \ \alpha_i = 0$$

*then we say that the pseudo best approximation does not exist for this set of basis $\{\alpha_1, \ldots, \alpha_k\}$. In this case we chose another set of basis for W say $\{\gamma_1, \ldots, \gamma_k\}$ and calculate*

$$\sum_{i=1}^{k} \langle \beta, \gamma_i \rangle_p \gamma_i$$

*is the pseudo best approximation to $\beta$.*

Now we apply this pseudo best approximation, to get the most likely transmitted code word. Let C be code over $Z_q^n$. Clearly C is a subspace of $Z_q^n$. $Z_q^n$ is a vector space over $Z_p$ where $q = p^t$, (p – a prime, t > 1).

We take C = W in the definition and apply the best biapproximation to β where β is the received code word for some transmitted word from C but $\beta \notin W = C$, for if $\beta \in W = C$ we accept it as the correct message. If $\beta \notin C$ then we apply the notion of pseudo best approximation to β related to the subspace C in $Z_q^n$.

Let $\{c_1, \ldots, c_k\}$ be chosen as the basis of C then

$$\sum_{i=1}^{k} \langle \beta, c_i \rangle_p c_i$$

gives the best approximation to β clearly

$$\sum_{i=1}^{k} \langle \beta, c_i \rangle_p c_i$$

belongs to C provided

$$\sum_{i=1}^{k} \langle \beta, c_i \rangle_p c_i \neq 0.$$

It is easily seen that this is the most likely received message. If

$$\sum_{i=1}^{k} \langle \beta, c_i \rangle_p c_i = 0,$$

we can choose another basis for C so that

$$\sum_{i=1}^{k} \langle \beta, c_i \rangle_p c_i \neq 0 .$$

Now we just adopt this argument in case of bicodes.

**DEFINITION 3.8.2:** *Let $V = V_1 \cup V_2$ be a bivector space over the finite field $Z_p$, with some pseudo inner biproduct $\langle , \rangle_p$ defined on V. Let $W = W_1 \cup W_2$ be the subbispace of $V = V_1 \cup V_2$. Let $\beta \in V_1 \cup V_2$ i.e. $\beta = \beta_1 \cup \beta_2$ related to $W = W_1 \cup W_2$ is defined as follows: Let*

$$\{\alpha_1, \ldots, \alpha_k\} = \{\alpha_1^1, \ldots, \alpha_{k_1}^1\} \cup \{\alpha_1^2, \ldots, \alpha_{k_2}^2\}$$

*be the chosen basis of the bisubspace $W = W_1 \cup W_2$ the pseudo best biapproximation to $\beta = \beta_1 \cup \beta_2$ if it exists is given by*

$$\sum_{i=1}^{k}\langle \beta, \alpha_i\rangle_p \alpha_i = \sum_{i=1}^{k_1}\langle \beta_1, \alpha_i^1\rangle_{p_1} \alpha_i^1 \cup \sum_{i=1}^{k_2}\langle \beta_2, \alpha_i^2\rangle_{p_2} \alpha_i^2 .$$

*If* $\sum_{i=1}^{k}\langle \beta, \alpha_i\rangle_p \alpha_i = 0$ *then we say that the pseudo best biapproximation does not exist for the set of basis*

$$\{\alpha_1, \alpha_2, ..., \alpha_k\} = \{\alpha_1^1, \alpha_2^1, ..., \alpha_{k_1}^1\} \cup \{\alpha_1^2, \alpha_2^2, ..., \alpha_{k_2}^2\}.$$

*In this case we choose another set of basis* $\{\alpha'_1, ..., \alpha'_{k_1}\}$ *and find* $\sum_{i=1}^{k}\langle \beta, \alpha'_i\rangle_p \alpha'_1$ *which is taken as the best biapproximation to* $\beta = \beta_1 \cup \beta_2$.

Now we apply it in case of bicodes in the following way to find the most likely bicode word. Let $C = C_1 \cup C_2$ be a bicode in $Z_q^n$. Clearly C is a bivector subspace of $Z_q^n$. Take in the definition C = W and apply the pseudo best biapproximation. If some bicode word $x = x_1 \cup x_2$ in $C = C_1 \cup C_2$ is transmitted and $\beta = \beta_1 \cup \beta_2$ is received bicode word then if $\beta \in C = C_1 \cup C_2$ then it is accepted as the correct message if $\beta \notin C = C_1 \cup C_2$ then we apply pseudo best biapproximation to $\beta = \beta_1 \cup \beta_2$ related to the subbispace $C = C_1 \cup C_2$ in $Z_q^n$.

Here three cases occur if $\beta = \beta_1 \cup \beta_2 \notin C = C_1 \cup C_2$.

1) $\beta_1 \in C_1$ and $\beta_2 \notin C_2$ so that $\beta_1 \cup \beta_2 \notin C$
2) $\beta_1 \notin C_1$ and $\beta_2 \in C_2$ so that $\beta_1 \cup \beta_2 \notin C$
3) $\beta_1 \notin C_1$ and $\beta_2 \notin C_2$ so that $\beta_1 \cup \beta_2 \notin C$

We first deal with (3) then indicate the working method in case of (1) and (2).

Given $\beta = \beta_1 \cup \beta_2$ , with $\beta \notin C$ ; $\beta_1 \notin C_1$ and $\beta_2 \notin C_2$. choose a basis

$$(c_1, c_2, \ldots, c_k) = \{c_1^1, c_2^1, \ldots, c_{k_1}^1\} \cup \{c_1^2, c_2^2, \ldots, c_{k_2}^2\}$$

of the subbispace C. To find the best biapproximation to β in C find

$$\sum_{i=1}^{k} \left\langle \beta / c_i \right\rangle_p c_i = \sum_{i=1}^{k_1} \left\langle \beta_1 / c_i^1 \right\rangle_1 c_i^1 \cup \sum_{i=1}^{k_2} \left\langle \beta_2 / c_i^2 \right\rangle_2 c_i^2 .$$

If both

$$\sum_{i=1}^{k_1} \left\langle \beta^1 / c_i^1 \right\rangle_1 c_i^1 \neq 0$$

and

$$\sum_{i=1}^{k_2} \left\langle \beta_2 / c_i^2 \right\rangle_2 c_i^2 \neq 0$$

then

$$\sum_{i=1}^{k} \left\langle \beta / c_i \right\rangle_p c_i \neq 0$$

is taken as the best biapproximation to β.
If one of

$$\sum_{i=1}^{k_t} \left\langle \beta_t \, c_i^t \right\rangle_t c_i^t = 0 \;\; (t = 1, 2)$$

say

$$\sum_{i=1}^{n} \left\langle \beta_i \, c_i^1 \right\rangle c_i^1 = 0$$

then choose a new basis, say $\left\{c'_1 \ldots c'_{k_t}\right\}$ and calculate

$$\sum_{i=1}^{k_t} \left\langle \beta_t \, c'_i \right\rangle_t c'_i \; ;$$

that is

$$\sum_{i=1}^{k_1} \left\langle \beta_1 c'_i \right\rangle_1 c'_i \cup \sum_{i=1}^{k_2} \left\langle \beta_2 c_i^2 \right\rangle_2 c_i^2$$

will be the best biapproximation to β in C. If both

$$\sum_{i=1}^{k_1} \left\langle \beta_1 c_i^1 \right\rangle_1 c_i^1 = 0 \text{ and } \sum_{i=1}^{k_2} \left\langle \beta_2 c_i^2 \right\rangle_2 c_i^2 = 0$$

then choose a new basis for C.

$$\{b_1,\ldots,b_k\} = \{b_1^1,\ldots,b_{k_1}^1\} \cup \{b_1^2,\ldots,b_{k_2}^2\}$$

and find

$$\sum_{i=1}^{k} \langle \beta, b_i \rangle_p b_i = \sum_{i=1}^{k_1} \langle \beta_1, b_i^1 \rangle_1 b_i^1 \cup \sum_{i=1}^{k_2} \langle \beta_2, b_i^2 \rangle_2 b_i^2 .$$

If this is not zero it will be the best biapproximation to $\beta$ in C.

Now we proceed on to work for cases (1) and (2) if we work for one of (1) or (2) it is sufficient. Suppose we assume

$$\sum_{i=1}^{k} \langle \beta c_i \rangle_p c_i = \sum_{i=1}^{k_1} \langle \beta_1 c_i^1 \rangle_1 c_i^1 \cup \sum_{i=1}^{k_2} \langle \beta_2 c_i^2 \rangle_2 c_i^2$$

and say

$$\sum_{i=1}^{k_2} \langle \beta_2 c_i^2 \rangle_2 c_i^2 = 0$$

and

$$\sum_{i=1}^{k_1} \langle \beta_1 c_i^1 \rangle_1 c_i^1 \neq 0$$

then we choose only a new basis for $c_2$ and calculate the pseudo best approximation for $\beta_2$ relative to the new basis of $C_2$.

Now we illustrate this by the following example:

***Example 3.8.1:*** Let C = $C_1$ (6, 3) $\cup$ $C_2$ (8, 4) bicode over $Z_2$ generated by the parity check bimatrix H = $H_1 \cup H_2$ given by

$$\begin{bmatrix} 0 & 1 & 1 & 1 & 0 & 0 \\ 1 & 0 & 1 & 0 & 1 & 0 \\ 1 & 1 & 0 & 0 & 0 & 1 \end{bmatrix} \cup \begin{bmatrix} 1 & 1 & 0 & 1 & 1 & 0 & 0 & 0 \\ 0 & 0 & 1 & 1 & 0 & 1 & 0 & 0 \\ 1 & 0 & 1 & 0 & 0 & 0 & 1 & 0 \\ 1 & 1 & 1 & 1 & 0 & 0 & 0 & 1 \end{bmatrix}.$$

The bicode is given by $C = C_1 \cup C_2$

$$= \begin{Bmatrix} (000000) & (011011) & (110110) & (001110) \\ (100011) & (111000) & (100101) & (101101) \end{Bmatrix} \cup$$

{(0 0 0 0 0 0 0 0), (1 0 0 0 1 0 1 1), (0 1 0 0 1 1 0 0), (0 0 1 0 0 1 1 1), (0 0 0 1 1 1 0 1), (1 1 0 0 0 0 1 0 ), (0 1 1 0 1 1 1 0), (0 0 1 1 1 0 1 0), (0 1 0 1 0 1 0 0), (1 0 1 0 1 1 0 0), (1 0 0 1 0 1 1 0), (1 1 1 0 0 1 0 1), (0 1 1 1 0 0 1 1), (1 1 0 1 1 1 1 1), (1 0 1 1 0 0 0 1), (1 1 1 1 1 0 0 0 )}.

Choose a basis $B = B_1 \cup B_2$ = {(0 0 1 1 1 0), (1 1 1 0 0 0), (0 1 0 1 0 1)} $\cup$ {(1 1 1 0 0 1 0 1), (1 1 1 1 1 0 0 0 )}. Let $\beta$ = (1 1 1 1 1 1) $\cup$ (1 1 1 1 1 1 1 1) be the received bicode. Clearly $\beta \notin C_1 \cup C_2 = C$. Now we find the best biapproximation to $\beta$ in C relative to the basis $B = B_1 \cup B_2$ under the pseudo inner biproduct.

$\alpha_1 \cup \alpha_2$ = {⟨(1 1 1 1 1 1) / (0 0 1 1 1 0)⟩$_1$ (0 0 1 1 1 0) + ⟨(1 1 1 1 1 1) / (1 1 1 0 0 0)⟩$_1$ (1 1 1 0 0 0) + ⟨(1 1 1 1 1 1) / (0 1 0 1 0 1)⟩$_1$ (0 1 0 1 0 1) } $\cup$ {⟨(1 1 1 1 1 1 1 1 ) / (0 1 0 0 1 0 0 1)⟩$_2$ (0 1 0 0 1 0 0 1) + ⟨(1 1 1 1 1 1 1 1) / (1 1 1 0 0 1 0 1)⟩$_2$ (1 1 1 0 0 1 0 1) + ⟨(1 1 1 1 1 1 1 1) / (1 1 1 1 1 0 0 0)⟩$_2$ (1 1 1 1 1 0 0 0)}

= {(0 0 1 1 1 0) + (1 1 1 0 0 0) + (0 1 0 1 0 1)} $\cup$ {(0 1 0 0 1 0 0 1) + (1 1 0 0 0 0 1 0) + (1 1 1 0 0 1 0 1) + (1 1 1 1 1 0 0 0 )}

= (1 0 0 0 1 1 ) $\cup$ (1 0 0 1 0 1 1 0) $\in C_1 \cup C_2 = C$.

Thus this is the pseudo best biapproximation to the received bicode word {(1 1 1 1 1 1) $\cup$ (1 1 1 1 1 1 1 1)}.

Thus the method of pseudo best biapproximation to a received bicode word which is not a bicode word is always guaranteed. If one wants to get best of the best biapproximations one can vary the basis and find the resultant biapproximated bicodes, compare it with the received message, the bicode which gives the least Hamming bidistance from the received word is taken as the pseudo best biapproximated bicode.

*Note:* We just define the Hamming distance of two bicodes $x = x_1 \cup x_2$ and $y = y_1 \cup y_2$, which is given by $d^B (x, y) = d (x_1\ y_1) \cup d (x_2\ y_2)$ where $d (x_1\ y_1)$ and $d(x_2, y_2)$ are the Hamming distance. The least of Hamming bidistances is taken as minimum of the sum of $d (x_1\ y_1) + d (x_2\ y_2)$.

## 3.9 Markov bichains – biprocess

We just indicate the applications of vector bispaces and linear bialgebras in Markov biprocess; for this we have to first define the notion of Markov biprocess and its implications to linear bialgebra / bivector spaces. We may call it as Markov biprocess or Markov bichains.

Suppose a physical or mathematical system is such that at any moment it occupy two of the finite number of states (Incase of one of the finite number of states we apply Markov chains or the Markov process). For example say about a individuals emotional states like happy sad suppose a system more with time from two states or a pair of states to another pair of sates; let us constrict a schedule of observation times and a record of states of the system at these times. If we find the transition from one pair of state to another pair of state is not predetermined but rather can only be specified in terms of certain probabilities depending on the previous history of the system then the biprocess is called a stochastic biprocess. If in addition these transition probabilities depend only on the immediate history of the system that is if the state of the system at any observation is dependent only on its state at the immediately proceeding observations then the process is called Markov biprocess on Markov bichain.

The bitransition probability $p_{ij} = p^1_{i_1 j_1} \cup p^2_{i_2 j_2}$ (i, j = 1, 2,…, k) is the probabilities that if the system is in state j =

$(j_1, j_2)$ at any observation, it will be in state $i = (i_1, i_2)$ at the next observation. A transition matrix

$$P = [p_{ij}] = \left[p^1_{i_1 j_1}\right] \cup \left[p^2_{i_2 j_2}\right].$$

is any square bimatrix with non negative entries all of which bicolumn sum is $1 \cup 1$. A probability bivector is a column bivector with non negative entries whose sum is $1 \cup 1$.

The probability bivectors are said to be the state bivectors of the Markov biprocess. If $P = P_1 \cup P_2$ is the transition bimatrix of the Markov biprocess and $x^n = x_1^n \cup x_2^n$ is the state bivector at the $n^{th}$ observation then $x^{(n+1)} = P\, x^{(n)}$ and thus $x_1^{(n+1)} \cup x_2^{(n+1)} = P_1 x_1^{(n)} \cup P_2 x_2^{(n)}$. Thus Markov bichains find all its applications in bivector spaces and linear bialgebras.

## Chapter Four

# NEUTROSOPHIC LINEAR BIALGEBRA AND ITS APPLICATION

In this chapter for the first time we introduce the notion of neutrosophic linear bialgebra and examine some of its application.

To introduce the notion of neutrosophic bialgebra itself we need to introduce several of the neutrosophic algebraic structures. This chapter has five sections. In the first section of this chapter we introduce some basic neutrosophic algebraic structures including the neutrosophic bivector spaces and some Smarandache algebraic structure. In the second section we introduce the notion of Smarandache neutrosophic linear bialgebra and some of its properties. The third section gives Smarandache representation of finite Smarandache neutrosophic bisemigroup. Section four gives a brief introduction to Smarandache Markov bichains using S-neutrosophic bivector spaces. The final section gives Smarandache neutrosophic Leontief economic models.

## 4.1 Some basic Neutrosophic Algebraic Structures

In this section we just recall and in some cases introduce the basic neutrosophic algebraic structures essential for the introduction of neutrosophic linear bialgebra.

The notion of neutrosophic logic was created by Florentin Smarandache which is an extension / combination of the fuzzy logic in which the notion of indeterminacy is introduced. The concept of indeterminacy is not like the concept of imaginary for we can say it is an indeterminate. To be more logical the notion of indeterminacy in several places is better accepted than imaginary value the more in practical problems in industries, so a chosen value which is compatible can be accepted. For if we have polynomial equation in a single variable and we say the equation has imaginary roots say even for the simple equation $x^2 + 1 = 0$ in R [x]. One will not and cannot understand what one means by an imaginary solution. But on the other if some one says this equation has no solution or the solutions to this equation are indeterminate, immediately one understands that the solutions for this equation cannot be determined. So for a layman, the term cannot be determined is closer or better known to him than the concept of imaginary.

So throughout this book I will denote the indeterminate. We make a further assumption $I.I = I^2 = I$. It is to be noted an event may be true 'T' or false F or indeterminate I. So unlike our usual probability were T + F = 1 we have T + F + I and more over $T + F + I \gtreqless 0$. Now we are not interested in this book to study the neutrosophic logic i.e. a logic which solely deals with the concept of indeterminacy in the fuzzy logic or in plain logic.

But in this book we try to develop neutrosophic algebra. Already neutrosophic algebraic graphs, neutrosophic matrices etc have been very recently developed in [43]. Here we develop these neutrosophic algebraic structures that will help us in defining the notion of neutrosophic linear bialgebra.

**DEFINITION 4.1.1:** *Let K be the field of reals. We call the smallest field generated by $K \cup I$ to be the neutrosophic field for it involves the indeterminacy factor in it. We define $I^2 = I$, $I + I = 2I$ i.e., $I + \ldots + I = nI$, and if $k \in K$ then $k.I =$* *kI, 0I = 0. We denote the neutrosophic field by K(I) which is generated by K ∪ I that is K (I) = 〈K ∪ I〉.*

***Example 4.1.1:*** Let R be the field of reals. The neutrosophic field is generated by 〈R ∪ I〉 i.e. R(I) clearly R ⊂ 〈R ∪ I〉.

***Example 4.1.2:*** Let Q be the field of rationals. The neutrosophic field is generated by Q and I i.e. Q ∪ I denoted by Q(I).

**DEFINITION 4.1.2:** *Let K(I) be a neutrosophic field we say K(I) is a prime neutrosophic field, if K(I) has no proper subfield which is a neutrosophic field.*

***Example 4.1.3:*** Q(I) is a prime neutrosophic field where as R(I) is not a prime neutrosophic field for Q(I) ⊂ R (I).

It is very important to note that all neutrosophic fields are of characteristic zero. Likewise we can define neutrosophic subfield.

**DEFINITION 4.1.3:** *Let K(I) be a neutrosophic field, P ⊂ K(I) is a neutrosophic subfield of P if P itself is a neutrosophic field; K(I) will also be called as the extension neutrosophic field of the neutrosophic field P.*

Now we proceed on to define neutrosophic vector spaces, which are only defined over neutrosophic fields. We can define two types of neutrosophic vector spaces one when it is a neutrosophic vector space over ordinary field other being neutrosophic vector space over neutrosophic fields. To this end we have to define neutrosophic group under addition.

**DEFINITION 4.1.4:** *We know Z is the abelian group under addition. Z(I) denote the additive abelian group generated by the set Z and I, Z(I) is called the neutrosophic abelian group under '+'. Thus to define basically a neutrosophic group under addition we need a group under addition. So we proceed on to define neutrosophic abelian group under* *addition. Suppose G is an additive abelian group under '+'. G(I) = ⟨G ∪ I⟩, additive group generated by G and I, G(I) is called the neutrosophic abelian group under '+'.*

***Example 4.1.4:*** Let Q be the group under '+'; Q (I) = ⟨Q ∪ I⟩ is the neutrosophic abelian group under addition; '+'.

***Example 4.1.5:*** R be the additive group of reals, R(I) = ⟨R ∪ I⟩ is the neutrosophic group under addition.

***Example 4.1.6:*** $M_{n\times m}(I) = \{(a_{ij}) \mid a_{ij} \in Z(I)\}$ be the collection of all n × m matrices under '+' $M_{n\times m}(I)$ is a neutrosophic group under '+'.

Now we proceed on to define neutrosophic subgroup.

**DEFINITION 4.1.5:** *Let G(I) be the neutrosophic group under addition. P ⊂ G(I) be a proper subset of G(I). P is said to be neutrosophic subgroup of G(I) if P itself is a neutrosophic group i.e. P = ⟨$P_1$ ∪ I⟩ where $P_1$ is an additive subgroup of G.*

***Example 4.1.7:*** Let Z(I) = ⟨Z ∪ I⟩ be a neutrosophic group under '+'. ⟨2Z ∪ I⟩ = 2Z(I) is the neutrosophic subgroup of Z (I).

In fact Z(I) has several neutrosophic subgroups.

Now we proceed on to define the notion of neutrosophic quotient group.

**DEFINITION 4.1.6:** *Let G (I) = ⟨G ∪ I⟩ be a neutrosophic group under '+', suppose P (I) be a neutrosophic subgroup of G (I) then the neutrosophic quotient group*

$$\frac{G(I)}{P(I)} = \{a + P(I) \mid a \in G\,(I)\}.$$

***Example 4.1.8:*** Let Z (I) be a neutrosophic group under addition, Z the group of integers under addition P = 2Z(I) is a neutrosophic subgroup of Z(I), the neutrosophic quotient group using subgroup 2Z(I) is;

$$\frac{Z(I)}{P} = \{a + 2Z(I) \mid a \in Z(I)\} = \{(2n+1) + (2n+1)\, I \mid n \in Z\}.$$

Clearly $\frac{Z(I)}{P}$ is a group. For P = 2Z (I) serves as the additive identity. Take a, b $\in \frac{Z(I)}{P}$. If a, b $\in$ Z(I) \ P then two possibilities occur.

a + b is odd times I or a + b is odd or a + b is even times I or even if a + b is even or even times I then a + b $\in$ P. if a + b is odd or odd times I a + b $\in \frac{Z(I)}{P = 2Z(I)}$.

It is easily verified that P acts as the identity and every element in

$$a + 2Z\,(I) \in \frac{Z(I)}{2Z(I)}$$

has inverse. Hence the claim.

Now we proceed on to define the notion of neutrosophic vector spaces over fields and then we define neutrosophic vector spaces over neutrosophic fields.

**DEFINITION 4.1.7:** *Let G(I) by an additive abelian neutrosophic group. K any field. If G(I) is a vector space over K then we call G(I) a neutrosophic vector space over K.*

Now we give the notion of strong neutrosophic vector space.

**DEFINITION 4.1.8:** *Let G(I) be a neutrosophic abelian group. K(I) be a neutrosophic field. If G(I) is a vector space over K(I) then we call G(I) the strong neutrosophic vector space.*

**THEOREM 4.1.1:** *All strong neutrosophic vector spaces over K(I) are a neutrosophic vector spaces over K; as K $\subset$ K(I).*

*Proof:* Follows directly by the very definitions.

Thus when we speak of neutrosophic spaces we mean either a neutrosophic vector space over K or a strong neutrosophic vector space over the neutrosophic field K(I). By basis we mean a linearly independent set which spans the neutrosophic space.

Now we illustrate this by an example.

***Example 4.1.9:*** Let R(I) × R(I) = V be an additive abelian neutrosophic group over the neutrosophic field R(I). Clearly V is a strong neutrosophic vector space over R(I). The basis of V are {(0,1), (1,0)}.

***Example 4.1.10:*** Let V = R(I) × R(I) be a neutrosophic abelian group under addition. V is a neutrosophic vector space over R. The neutrosophic basis of V are {(1,0), (0,1), (I,0), (0,I)}, which is a basis of the vector space V over R.

A study of these basis and its relations happens to be an interesting form of research.

**DEFINITION 4.1.9:** *Let G(I) be a neutrosophic vector space over the field K. The number of elements in the neutrosophic basis is called the neutrosophic dimension of G(I).*

**DEFINITION 4.1.10:** *Let G(I) be a strong neutrosophic vector space over the neutrosophic field K(I). The number of elements in the strong neutrosophic basis is called the strong neutrosophic dimension of G(I).*

We denote the neutrosophic dimension of G(I) over K by $N_k$ (dim) of G (I) and that the strong neutrosophic dimension of G (I) by $SN_{K(I)}$ (dim) of G(I).

Now we define the notion of neutrosophic matrices.

**DEFINITION 4.1.11:** *Let $M_{nxm} = \{(a_{ij}) \,/\, a_{ij} \in K(I)\}$, where K(I), is a neutrosophic field. We call $M_{nxm}$ to be the neutrosophic matrix.*

***Example 4.1.11:*** Let $Q(I) = \langle Q \cup I \rangle$ be the neutrosophic field.

$$M_{4x3} = \begin{pmatrix} 0 & 1 & I \\ -2 & 4I & 0 \\ 1 & -I & 2 \\ 3I & 1 & 0 \end{pmatrix}$$

is the neutrosophic matrix, with entries from rationals and the indeterminacy I. We define product of two neutrosophic matrices whenever the product is defined as follows:
Let

$$A = \begin{pmatrix} -1 & 2 & -I \\ 3 & I & 0 \end{pmatrix}_{2\times 3}$$

and

$$B = \begin{pmatrix} -I & 1 & 2 & 4 \\ 1 & I & 0 & 2 \\ 5 & -2 & 3I & -I \end{pmatrix}_{3x4}$$

$$AB = \begin{bmatrix} -6I+2 & -1+4I & -2-3I & I \\ -2I & 3+I & 6 & 12+2I \end{bmatrix}_{2x4}$$

(we use the fact $I^2 = I$).

To define Neutrosophic Cognitive Maps (NCMs) we direly need the notion of Neutrosophic Matrices. We use square neutrosophic matrices for NCMs and use rectangular neutrosophic matrices for Neutrosophic Relational Maps (NRMs). We need the notion of neutrosophic graphs basically to obtain NCMs which will be nothing but directed neutrosophic graphs. Similarly neutrosophic relational maps will also be directed neutrosophic graphs.

It is no coincidence that graph theory has been independently discovered many times since it may quite properly be regarded as an area of applied mathematics. The subject finds its place in the work of Euler. Subsequent rediscoveries of graph theory were by Kirchhoff and Cayley. Euler (1707-1782) became the father of graph theory as well as topology when in 1936 he settled a famous unsolved problem in his day called the Konigsberg Bridge Problem.

Psychologist Lewin proposed in 1936 that the life space of an individual be represented by a planar map. His view point led the psychologists at the Research center for Group Dynamics to another psychological interpretation of a graph in which people are represented by points and interpersonal relations by lines. Such relations include love, hate, communication ability and power. In fact it was precisely this approach which led the author to a personal discovery of graph theory, aided and abetted by psychologists L. Festinger and D. Cartwright. Here it is pertinent to mention that the directed graphs of an FCMs or FRMs are nothing but the psychological inter-relations or feelings of different nodes; but it is unfortunate that in all these studies the concept of indeterminacy was never given any place, so in this section for the first time we will be having graphs in which the notion of indeterminacy i.e. when two vertex should be connected or not is never dealt with. If graphs are to display human feelings then this point is very vital for in several situations certain relations between concepts may certainly remain an indeterminate. So this section will purely cater to the properties of such graphs, the edges of certain vertices may not find its connection i.e., they are indeterminates, which we will be defining as neutrosophic graphs.

The world of theoretical physics discovered graph theory for its own purposes. In the study of statistical mechanics by Uhlenbeck the points stands for molecules and two adjacent points indicate nearest neighbor interaction of some physical kind, for example magnetic attraction or repulsion. But it is forgotten in all these situations we may have molecule structures which need not attract or repel but remain without action or not able to predict the action for such analysis we can certainly adopt the concept of neutrosophic graphs.

In a similar interpretation by Lee and Yang the points stand for small cubes in Euclidean space where each cube may or may not be occupied by a molecule. Then two points are adjacent whenever both spaces are occupied.

Feynmann proposed the diagram in which the points represent physical particles and the lines represent paths of the particles after collisions. Just at each stage of applying graph theory we may now feel the neutrosophic graph theory may be more suitable for application.

Now we proceed on to define the neutrosophic graph.

**DEFINITION 4.1.12:** *A neutrosophic graph is a graph in which at least one edge is an indeterminacy denoted by dotted lines.*

*Notation:* The indeterminacy of an edge between two vertices will always be denoted by dotted lines.

***Example 4.1.12:*** The following are neutrosophic graphs:

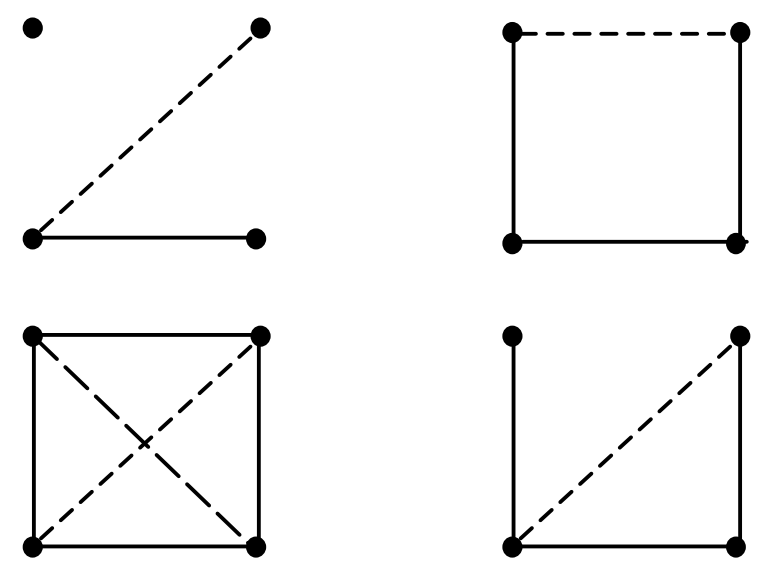

FIGURE:4.1.

All graphs in general are not neutrosophic graphs.

***Example 4.1.13:*** The following graphs are not neutrosophic graphs given in Figure 4.1.2:

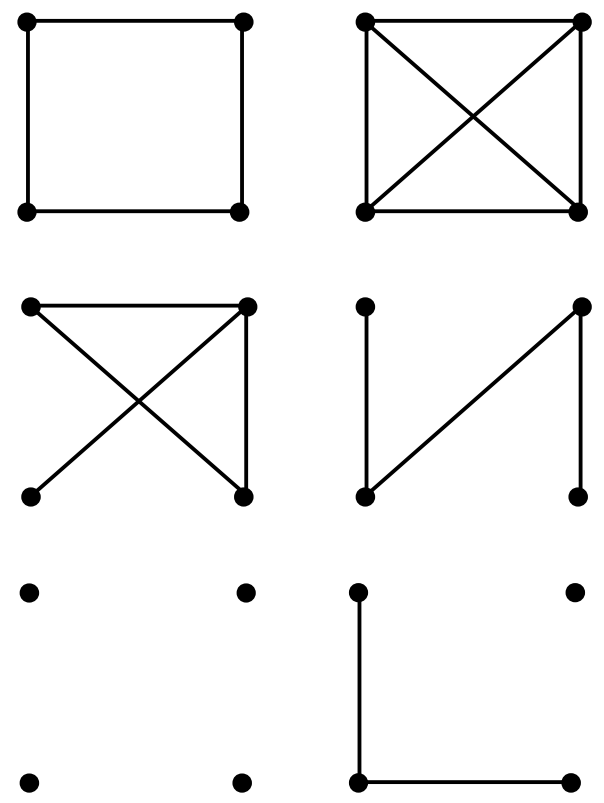

FIGURE: 4.1.2

**DEFINITION 4.1.13:** *A neutrosophic directed graph is a directed graph which has at least one edge to be an indeterminacy.*

**DEFINITION 4.1.14:** *A neutrosophic oriented graph is a neutrosophic directed graph having no symmetric pair of directed indeterminacy lines.*

**DEFINITION 4.1.15:** *A neutrosophic subgraph H of a neutrosophic graph G is a subgraph H which is itself a neutrosophic graph.*

**THEOREM 4.1.2:** *Let G be a neutrosophic graph. All subgraphs of G are not neutrosophic subgraphs of G.*

*Proof:* By an example. Consider the neutrosophic graph

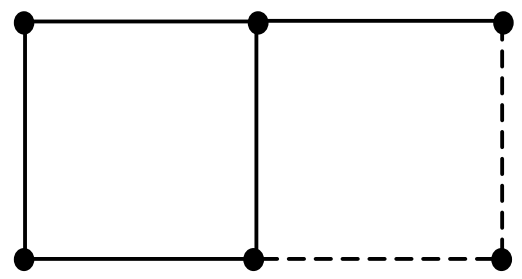

FIGURE: 4.1.3

given in Figure 4.1.3.

This has a subgraph given by Figure 4.1.4.

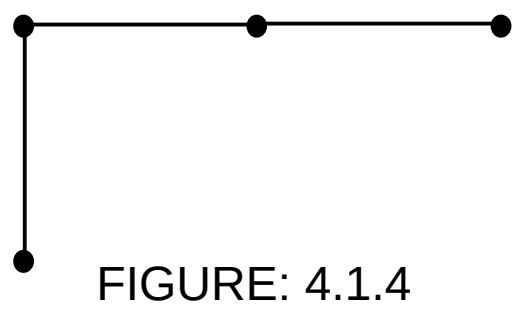

FIGURE: 4.1.4

which is not a neutrosophic subgraph of G.

**THEOREM 4.1.3:** *Let G be a neutrosophic graph. In general the removal of a point from G need not be a neutrosophic subgraph.*

*Proof:* Consider the graph G given in Figure 4.1.5.

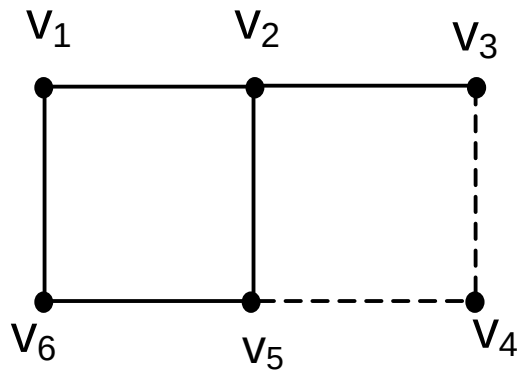

FIGURE: 4.1.5

G \ $v_4$ is only a subgraph of G but is not a neutrosophic subgraph of G.

Thus it is interesting to note that this is a main feature by which a graph differs from a neutrosophic graph.

**DEFINITION 4.1.16:** *Two graphs G and H are neutrosophically isomorphic if*

- *i.* *They are isomorphic.*
- *ii.* *If there exists a one to one correspondence between their point sets which preserve indeterminacy adjacency.*

**DEFINITION 4.1.17:** *A neutrosophic walk of a neutrosophic graph G is a walk of the graph G in which at least one of the lines is an indeterminacy line. The neutrosophic walk is* *neutrosophic closed if $v_0$ = $v_n$ and is neutrosophic open otherwise.*

*It is a neutrosophic trial if all the lines are distinct and at least one of the lines is a indeterminacy line and a neutrosophic path, if all points are distinct (i.e. this necessarily means all lines are distinct and at least one line is a line of indeterminacy). If the neutrosophic walk is neutrosophic closed then it is a neutrosophic cycle provided its n points are distinct and $n \geq 3$.*

*A neutrosophic graph is neutrosophic connected if it is connected and at least a pair of points are joined by a path. A neutrosophic maximal connected neutrosophic subgraph of G is called a neutrosophic connected component or simple neutrosophic component of G.*

*Thus a neutrosophic graph has at least two neutrosophic components then it is neutrosophic disconnected. Even if one is a component and another is a neutrosophic component still we do not say the graph is neutrosophic disconnected.*

***Example 4.1.14:*** Neutrosophic disconnected graphs are given in Figure 4.1.6.

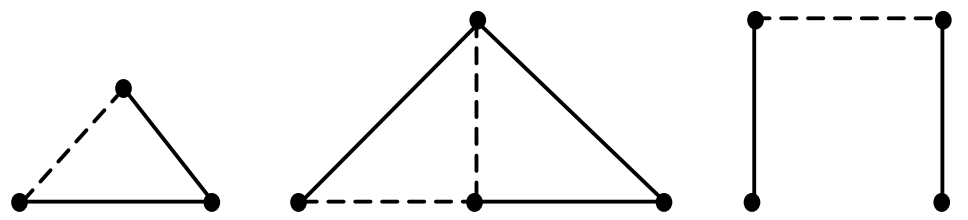

FIGURE: 4.1.6

***Example 4.1.15:*** Graph which is not neutrosophic disconnected is given by Figure 4.1.7.

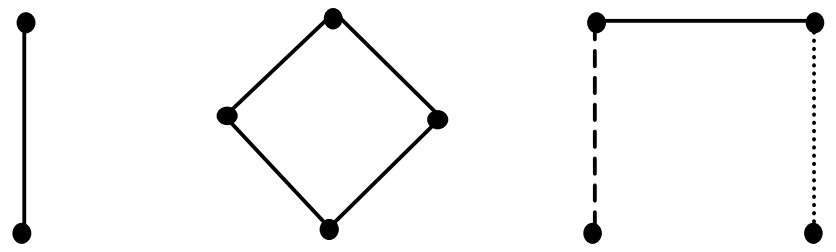

FIGURE: 4.1.7

Several results in this direction can be defined and analyzed.

**DEFINITION 4.1.18:** *A neutrosophic bigraph, G is a bigraph, G whose point set V can be partitioned into two*

*subsets $V_1$ and $V_2$ such that at least a line of G which joins $V_1$ with $V_2$ is a line of indeterminacy.*

This neutrosophic bigraphs will certainly play a role in the study of NRMs and in fact we give a method of conversion of data from NRMs to NCMs.

As both the models NRMs and NCMs work on the adjacency or the connection matrix we just define the neutrosophic adjacency matrix related to a neutrosophic

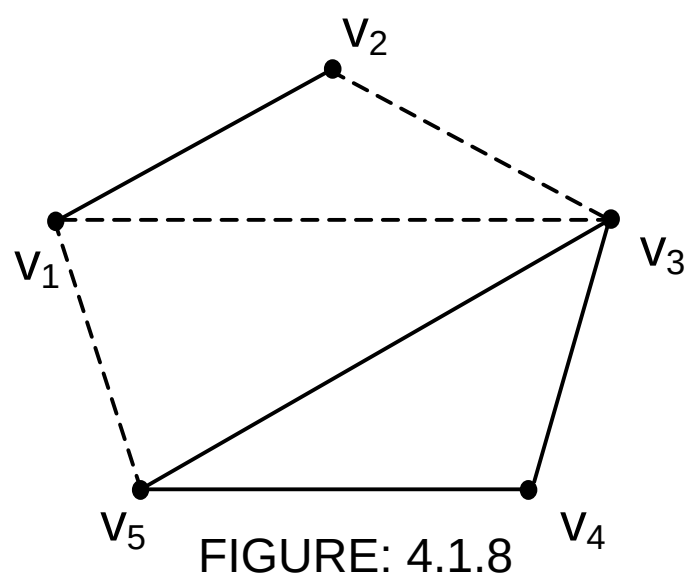

FIGURE: 4.1.8

graph G given by Figure 4.1.8.
The neutrosophic adjacency matrix is N(A)

$$N(A) = \begin{bmatrix} 0 & 1 & I & 0 & I \\ 1 & 0 & I & 0 & 0 \\ I & I & 0 & 1 & 1 \\ 0 & 0 & 1 & 0 & 1 \\ I & 0 & 1 & 1 & 0 \end{bmatrix}.$$

Its entries will not only be 0 and 1 but also the indeterminacy *I*.

**DEFINITION 4.1.19:** *Let G be a neutrosophic graph. The adjacency matrix of G with entries from the set (I, 0, 1) is called the neutrosophic adjacency matrix of the graph.*

Now as our main aim is the study of neutrosophic linear bialgebra we do not divulge into a very deep study of neutrosophic graphs or its properties but have given only the basic and the appropriate notions which are essential for studying of this book.

**DEFINITION 4.1.20:** *Let G be a neutrosophic group and K any field. If G is a vector space over K then we say G is a neutrosophic vector space over K.*

***Example 4.1.16:*** Let G = {The set of all 2 × 2 matrices with entries from Q (I) i.e. the neutrosophic field generated by Q and I}. G is a neutrosophic vector space over Q.

***Example 4.1.17:*** Let G = {the set of all m × n matrices with entries from Q (I)}. G is a neutrosophic vector space over Q.

Thus we saw the definition and a few examples of neutrosophic vector space. Now we proceed on to define. The notion of neutrosophic linear algebra.

**DEFINITION 4.1.21:** *Let G be a neutrosophic group under '+'. Let G is a neutrosophic vector space over the field K. If in G we have for every pair of elements, $g_1$, $g_2 \in G$, $g_1$ o $g_2 \in G$ 'o' a multiplicative operation on G, then we call G a neutrosophic linear algebra over K.*

***Example 4.1.18:*** Let G = {set of all 3 × 3 matrices with entries from Q (I)}. Clearly G is a neutrosophic linear algebra over Q, for G satisfies matrix multiplication.

We see that as in case of vector spaces and linear algebras every neutrosophic vector space in general is not a neutrosophic linear algebra. To this end we give an example of a neutrosophic vector space which is not a neutrosophic linear algebra.

***Example 4.1.19:*** Let G = {set of all 2 × 5 matrices with entries from Q (I)}. G is clearly a neutrosophic vector space over Q. But G is not a neutrosophic linear algebra over Q for we cannot define multiplication in the neutrosophic group G i.e., the set of $2 \times 5$ neutrosophic matrices.

It is clear by the very definition every neutrosophic linear algebra is a neutrosophic vector space.

Now the bases, dimension of a neutrosophic vector space are defined as in case of usual vector spaces. Now we proceed on to define the notion of strong neutrosophic vector space and strong neutrosophic linear algebra.

**DEFINITION 4.1.22:** *Let G be a neutrosophic group under '+'. If G is a neutrosophic vector space (linear algebra) over the neutrosophic field K(I) then we call G a strong neutrosophic vector space (linear algebra)..*

It is clear from the definition that every strong neutrosophic vector space (linear algebra) is a neutrosophic vector space (linear algebra).

But in general a neutrosophic vector space (linear algebra) is not a strong neutrosophic vector space (linear algebra).

***Example 4.1.20:*** Let G = {set of all $4 \times 2$ matrices with entries from Q (I)}. Clearly G is a neutrosophic group under addition. G is a neutrosophic vector space over Q and G is a strong neutrosophic vector space over Q (I).

**Remark:** Suppose V = {set of all $5 \times 5$ matrices with entries from Q}. V is a group. But V is not a neutrosophic vector space over Q (I).

The difference between the neutrosophic vector space and strong neutrosophic vector space can be also seen from the dimension or the number of base elements.

**DEFINITION 4.1.23:** *Let V (I) be a neutrosophic group and V (I) be a neutrosophic vector space over a field K. A non empty subset W of V (I) is said to be a neutrosophic subspace of V(I) if*

*(1) W is a neutrosophic subgroup of V (I).*
(2) *W is itself a neutrosophic vector space over K.*

***Example 4.1.21:*** Let V (I) = {the collection of all 5 × 5 matrices with entries from Q (I)}. V (I) is a group under '+'. V (I) is a neutrosophic vector space over Q. Consider W = {set of all 5 × 5 upper triangular matrices with entries from Q (I)}. Clearly W is a subset of V and W is a neutrosophic subspace of V. Suppose $W_1$ = {set of all 5 × 5 matrices with entries from Q}. Clearly $W_1$ is a subset of V (I), but $W_1$ is not a neutrosophic subgroup of V (I), hence $W_1$ is not a neutrosophic subspace of V (I).

So we can have in general a neutrosophic vector space which can have subspaces which are not neutrosophic vector space but only a vector space such subvector spaces are called as pseudo neutrosophic subspaces. This can happen only in neutrosophic vector spaces and never in strong neutrosophic vector spaces.

***Example 4.1.22:*** Let V = {R (I) [x] be the set of all neutrosophic polynomials with coefficient from R (I)}; V under '+' is a group. V is a neutrosophic vector space over Q, the field of rationals.

Now W = {R [x] the set of all polynomials in the variable x with coefficient in R}. W is a subset of V, W is also a subspace of V over Q. Clearly W is not a neutrosophic subspace of V but only a pseudo neutrosophic subspace of V over Q.

Now we define a linear transformation between two neutrosophic vector spaces which are defined over the same field. Let V and W be two neutrosophic vector spaces over the field K. A linear transformation $T_N : V \to W$ is said to be a neutrosophic linear transformation if T (I) = I and T : V $\to$ W is a linear transformation of vector space i.e. T (cI) = cT (I).

$$T(\alpha + \beta I) = T(\alpha) + T(\beta)I$$

i.e. the linear transformation leaves the indeterminacy invariant. However we adjoin the zero linear transformation

$T_o$ i.e. $T_o$ (υ) = 0 for all υ ∈ V and $T_o$ (I) = 0. Thus it is mater of routine to check that $Hom_K$ (V, W) is a neutrosophic vector space over K. In fact as in case of vector spaces, we have a one to one correspondence between the set of all m × n neutrosophic matrices with these neutrosophic linear transformations where dim V = n and dim W = m.

**DEFINITION 4.1.24:** *Let V be a finite dimensional strong neutrosophic vector space over Q (I). The set of all linear operators on V, denoted by $Hom_Q$ (V, V) is a strong neutrosophic vector space over Q (I).*

As in case of vector spaces we have in case of neutrosophic vector spaces, to every neutrosophic linear operator a neutrosophic square matrix is associated with it.

***Example 4.1.23:*** Let V = Q (I) × Q (I) be a neutrosophic vector space over Q.

Now we wish to calculate a basis for V. V is generated by {(0, 1), (1 0), (0 *I* ), (*I* 0)}.

Thus dimension of V is 4. Thus we cannot have V ≅ $Q^2$ Suppose V = Q (I) × Q (I) is strong neutrosophic vector space over Q (I) . Then a basis for V over Q(I) is {(0, 1), (1 0)}. Thus dimension of V is 2 and V ≅ $(Q (I))^2$. Thus it is interesting to note the same neutrosophic group as a neutrosophic vector space over a real field F and as a strong neutrosophic vector space over the neutrosophic field F (I) have varied dimensions.

***Example 4.1.24:*** Let V = {$(a_{ij})_{2 \times 2}$ i.e. the set of all 2 × 2 matrices with entries from Q (I)}. V is a neutrosophic group. V is a strong neutrosophic vector space over Q (I). The dimension of V is 4, the basis is

$$\left\{ \begin{pmatrix} 0 & 1 \\ 0 & 0 \end{pmatrix}, \begin{pmatrix} 1 & 0 \\ 0 & 0 \end{pmatrix}, \begin{pmatrix} 0 & 0 \\ 1 & 0 \end{pmatrix}, \begin{pmatrix} 0 & 0 \\ 0 & 1 \end{pmatrix} \right\}.$$

Suppose consider the same V i.e. V is the set of all $2 \times 2$ matrices with entries from Q (I). Now V is a neutrosophic vector space over Q. Now dimension V over Q is 8 i.e. a basis for V is

$$\left\{\begin{pmatrix}1 & 0\\ 0 & 0\end{pmatrix},\begin{pmatrix}0 & 1\\ 0 & 0\end{pmatrix},\begin{pmatrix}0 & 0\\ 1 & 0\end{pmatrix},\begin{pmatrix}0 & 0\\ 0 & 1\end{pmatrix},\right.$$

$$\left.\begin{pmatrix}I & 0\\ 0 & 0\end{pmatrix},\begin{pmatrix}0 & I\\ 0 & 0\end{pmatrix},\begin{pmatrix}0 & 0\\ I & 0\end{pmatrix},\begin{pmatrix}0 & 0\\ 0 & I\end{pmatrix}\right\}.$$

***Example 4.1.25:*** Consider V = {Q (I) [x] i.e. collection of all polynomials with coefficients from Q (I)}. Clearly V is a strong neutrosophic vector space over Q (I) and a neutrosophic vector space over Q. However dim of V is infinite in both cases, but of course different infinities.

***Example 4.1.26:*** Let V = {Q (I) [x] i.e. all polynomials of degree less than or equal to 3 with coefficients from Q (I)}. V is a neutrosophic group. V is a strong neutrosophic vector space over Q (I). A set of basis for V is $\{1, x, x^2, x^3\}$. Thus dimension of V as a strong neutrosophic vector space is 4. Now V can be realized as a neutrosophic vector space over Q.

Now a basis of V as a neutrosophic vector space over Q is $\{1, x_1, x^2, x^3, I, Ix, I x^2, I x^3\}$. Dimension of V as a neutrosophic vector space over Q is 8.

Thus we leave it as a problem for the reader. Suppose V is any neutrosophic group. Suppose V is a strong neutrosophic vector space over F (I) with dim V = n. Suppose V is considered just as a neutrosophic vector space over F. Is dim V = 2n ? Prove your claim.

Now we proceed on to define the notion of neutrosophic eigen values, neutrosophic eigen vectors and neutrosophic characteristic polynomial in case of square neutrosophic matrices as well as linear operators. If V is a neutrosophic group; V is a neutrosophic vector space over F, let dim V = 2n.

Then V as a strong neutrosophic vector space over F (I) is n i.e. $\dim_{F(I)} V = n$.

Thus even the same neutrosophic group will have different sets of eigen values and eigen vectors. For one the characteristic polynomial would be a degree n and for the other of degree 2n.

Now we first define the neutrosophic eigen values vectors of a neutrosophic matrix.

**DEFINITION 4.1.25:** *Let A be a $n \times n$ neutrosophic matrix with entries from Q(I); a neutrosophic characteristic value of A in Q(I) is a scalar c or an indeterminate dI + g such that the matrix ($A - c\,I_{n \times n}$) is singular (we see c can be dI + g and d and g are any scalar and I is the indeterminacy).*

***Example 4.1.27:*** Let A be a 3 × 3 neutrosophic matrix given below:

$$A = \begin{bmatrix} 3 & I & -1 \\ 2 & 2I & -1 \\ 2 & 2 & 0 \end{bmatrix}$$

The neutrosophic characteristic value is [$A - cI_{3 \times 3}$ ] such that | $A - c\,I_{3\times3}$| is singular.

Suppose f (c) = det |$A - c\,I_{n\times n}$| then f is called the characteristic polynomial of A. If is important to note that f is a monic polynomial which has exactly degree n. Thus the characteristic polynomial for the above 3 × 3 matrix given in this example is

$$\begin{vmatrix} 3-x & I & -1 \\ 2 & 2I-x & -1 \\ 2 & 2 & -x \end{vmatrix} = f(x)$$

i.e.,

$(3 - x) [(2I - x) (-x) + 2] - I [-2x + 2] -1 [4 - 2 (2I - x)].$

$= (3 - x) (-2Ix + x^2 + 2) + 2Ix - 2I - 4 + 4I - 2x$

$= (3 - x) (x^2 - 2Ix + 2) + 2Ix + 2I - 2x - 4$

$= 3x^2 - 6Ix + 6 - x^3 + 2Ix^2 - 2x + 2Ix + 2I - 2x - 4$

$= \quad -x^3 + (3 + 2I)\, x^2 - (4I + 4)\, x + 2I + 2 = 0$

is the neutrosophic characteristic polynomial which is clearly a neutrosophic polynomial of degree 3 in the variable x.

***Example 4.1.28:*** Let

$$A = \begin{bmatrix} 2 & I \\ -2 & 1 \end{bmatrix}$$

be a neutrosophic matrix with entries form Q(I). The characteristic neutrosophic polynomial is

$$\begin{vmatrix} x-2 & -I \\ 2 & x-1 \end{vmatrix} = 0$$

i.e.

$$\begin{aligned} (x-2)(x-1) + 2I &= 0 \\ x^2 - 3x + 2 + 2I &= 0 \\ x^2 - 3x + (2 + 2I) &= 0. \end{aligned}$$

The associated neutrosophic characteristic values are

$$\begin{aligned} x &= +3 \pm \frac{\sqrt{9-4(2+2I)}}{2} \\ &= \frac{3 \pm \sqrt{9-(2+2I)}}{2} \\ &= \frac{2 \pm \sqrt{7-2I}}{2} \notin Q(I) \end{aligned}$$

So this neutrosophic characteristic polynomial has no neutrosophic characteristic value in Q (I).

***Example 4.1.29:*** Let us consider the neutrosophic matrix

$$A = \begin{bmatrix} 2I+1 & 0 \\ -1 & 2 \end{bmatrix}.$$

The neutrosophic characteristic polynomial is

$$f(x) = \begin{vmatrix} x-\overline{2I+1} & 0 \\ 1 & x-2 \end{vmatrix}.$$

i.e. f (x) = (x – 2) (x – $\overline{2I+1}$) so the neutrosophic characteristic value related with the matrix are x = 2 and 2I + 1. Thus one of the value is real and the other is just a neutrosophic value. It is essential to note that neutrosophic roots in general need not occur in conjugate pairs.

Now we proceed on to define characteristic neutrosophic polynomial and neutrosophic characteristic values and related neutrosophic characteristic vector related with a linear operator of the neutrosophic vector spaces.

**DEFINITION 4.1.26:** *Let V be a strong neutrosophic vector space over the neutrosophic field F(I) and T be a linear neutrosophic operator on V. A neutrosophic characteristic value of T is a scalar C in F (I) such that for a non zero vector $\alpha$, $T\alpha = C\alpha$. If C is a neutrosophic characteristic value of T then*

*(a) any $\alpha$ such that $T\alpha = C\alpha$ is called a neutrosophic characteristic vector of T associated with the neutrosophic characteristic value C.*

*(b) The collection of all $\alpha$ such that $T\alpha = C\alpha$ is called the neutrosophic characteristic space associated with C.*

So if T : V $\to$ V where V is a strong neutrosophic vector space over F (I) and if dimension of V over F (I) is n to every neutrosophic linear operator T, there is associated a n $\times$ n neutrosophic matrix and we can find the neutrosophic characteristic polynomial, characteristic value and characteristic vector of T whenever they exist. However if we study the neutrosophic vector space over any field there are several changes to be taken in to consideration.

(1) The theory if V is a finite n-dimensional neutrosophic vector space over F then V $\ncong$ $F^n$.

(2) If T is a linear operator of V we cannot associate a $n \times n$ neutrosophic matrix with T.

As our main motivation is to introduce the notion of neutrosophic linear bialgebra we proceed to skip some of the results and just sketch only the very essentials.

The next important concept we wish to introduce is the inner product in neutrosophic vector spaces and strong neutrosophic vector spaces.

**DEFINITION 4.1.27:** *Let V be a neutrosophic vector space over the field of reals. A function $\langle\ \rangle$, from $V \times V$ to F is said to be a neutrosophic inner product if $\langle\ \rangle$ is the inner product on V and for $I \in V$; $\langle I\ I\rangle = I\ \langle 1, 1\rangle$. I is taken to be positive if we consider – I then the indeterminacy is negative.*

All properties of inner product on vector spaces can be carried out on strong neutrosophic vector space and neutrosophic vector space with suitable and appropriate modifications.

It is important to note that if V is a neutrosophic vector space over F then it is not possible to define linear functionals. For I will be in the vector space and $I \notin F$. This is evident if we take the simple neutrosophic vector space $V = Q(I) \times Q(I)$ over Q. This is occurring for Q (I) is a neutrosophic group under addition Q(I) is a strong neutrosophic vector space of dimension one over Q (I), but Q (I) as a neutrosophic vector space over Q is of dimension 2. This is the vast difference, which gives very often easy extension of results as in case of vector spaces.

Now we proceed on to define neutrosophic bivector spaces for we need this basically to define neutrosophic bigroups for the first time.

**DEFINITION 4.1.28:** *Let $G = G_1 \cup G_2$ if both $G_1$ and $G_2$ are proper subsets of G i.e. ($G_1 \not\subseteq G_2$ or $G_2 \not\subseteq G_1$) and if both $G_1$ and $G_2$ are neutrosophic groups then we say G is a neutrosophic bigroup.*

Bigroups are themselves very useful in the application of industrial problems, but as all industrial problems also have a concept of indeterminacy it is deem fit to introduce neutrosophic bigroups.

Now we illustrate these concepts by some examples so that reader can have a better understanding of the problem.

***Example 4.1.30:*** Consider the set G = {Q (I) × Q (I)} ∪ (Q (I) [x]) = $G_1 \cup G_2$. Clearly $G_1$ and $G_2$ are neutrosophic groups under '+' and '×' respectively. Thus G is a bigroup.

***Example 4.1.31:*** Let us consider G = {R (I) [x]} ∪ {Q (I) [x]}, R (I) [x] is the set of all neutrosophic polynomials in the indeterminate x with coefficients from R (I). Now Q (I) [x] is also the set of all neutrosophic polynomials with coefficients from the neutrosophic field Q (I). Clearly G = $G_1 \cup G_2$ is not a neutrosophic bigroup as Q (I) [x] ⊂ R (I) [x] i.e. $G_2 \underset{\neq}{\subset} G_1$.

***Example 4.1.32:*** Let G = $G_1 \cup G_2$ where $G_1$ is the set of all 3 × 2 neutrosophic matrices with entries from Q (I), $G_1$ is a neutrosophic group under matrix addition. $G_2$ be the set of all 2 × 3 matrices with entries from Q (I). $G_2$ is a neutrosophic group under matrix addition. Thus G = $G_1 \cup G_2$ is a neutrosophic bigroup.

Now we proceed on to define the neutrosophic subbigroup of the neutrosophic the bigroup G.

**DEFINITION 4.1.29:** *Let G = $G_1 \cup G_2$ be a neutrosophic bigroup. A proper subset H of G is said to be a neutrosophic subbigroup of G if (1) H = $H_1 \cup H_2$; $H_1$ is a neutrosophic subgroup of $G_1$ and $H_2$ is a neutrosophic subgroup of $G_2$ and $H_1$ and $H_2$ are proper subsets of H.*

We first illustrate this by an example.

***Example 4.1.33:*** Let $G = G_1 \cup G_2$ = {Q (I) × Q(I) × Q(I)} ∪ {set of all 3 × 3 matrices with entries form Q(I)}. Clearly both $G_1$ and $G_2$ are neutrosophic groups under addition. Further $G_1$ and $G_2$ are proper subsets of G. Thus $G = G_1 \cup G_2$ is the neutrosophic bigroup of G. Consider the proper subset $H = H_1 \cup H_2$ of G where $H_1$ = {(Q (I) × Q(I) × {0}} and H = {All 3 × 3 lower triangular neutrosophic matrices with entries from Q(I)}.

Clearly $H_1$ is neutrosophic subgroup of the neutrosophic group $G_1$ under component wise addition. Further $H_2$ is a neutrosophic subgroup of the neutrosophic group $G_2$ under matrix addition. Thus $H = H_1 \cup H_2$ is a neutrosophic subbigroup of the neutrosophic bigroup $G = G_1 \cup G_2$.

*Note:* A neutrosophic bigroup $G = G_1 \cup G_2$ is commutative if both the neutrosophic groups $G_1$ and $G_2$ are commutative if one of $G_1$ or $G_2$ alone is commutative we call G is a quasi commutative neutrosophic bigroup.

***Example 4.1.34:*** Let $G = G_1 \cup G_2$ where $G_1$ = Q(I) [x] and $G_2$ = {all 2 × 2 matrices with entries from Q(I) such that the matrix are non singular}. Clearly $G_1$ is a neutrosophic group under '+' and $G_2$ is a neutrosophic group under multiplication. Thus $G = G_1 \cup G_2$ is a neutrosophic bigroup. Clearly G is quasi commutative neutrosophic bigroup as $G_1$ is a commutative neutrosophic group whereas $G_2$ is a non commutative neutrosophic group.

But for us to define neutrosophic vector spaces we need only neutrosophic bigroups which are commutative so we do not indulge in defining several properties about neutrosophic bigroups in general.

Now we proceed on to define the notion of neutrosophic bifield.

**DEFINITION 4.1.30:** *Let $F = F_1 \cup F_2$ where $F_1$ and $F_2$ are proper subsets of F. If both $F_1$ and $F_2$ are neutrosophic fields then we say $F = F_1 \cup F_2$ is a neutrosophic bifield.*

***Example 4.1.35:*** Let F = Q (I) $\cup$ $Z_3$ (I). Clearly F is a neutrosophic field.

***Example 4.1.36:*** Let F = Q (I) $\cup$ R (I) = $F_1 \cup F_2$. We see both Q (I) and R (I) are neutrosophic fields but F is not a neutrosophic bifield, as Q(I) $\subset$ R(I) so Q(I) and R(I) are not proper subsets of F.

However we would be using only neutrosophic fields and not neutrosophic bifields in the definition of bivector spaces. Now we proceed on to define the notion of neutrosophic bivector spaces.

**DEFINITION 4.1.31:** *Let V = $V_1 \cup V_2$ be a neutrosophic bigroup under addition. Let F (I) be a neutrosophic field. If $V_1$ and $V_2$ are neutrosophic vector spaces over F (I) then we call V = $V_1 \cup V_2$ the strong neutrosophic bivector space, over the neutrosophic field F(I)..*

We now illustrate this by the following example.

***Example 4.1.37:*** Let V = $V_1 \cup V_2$ where $V_1$ = {Q (I) × Q (I)}, clearly $V_1$ is an additive neutrosophic group. Let $V_2$ = {The set of all 3 × 2 matrices with entries from Q(I) }; $V_2$ is a neutrosophic group under matrix addition. Thus V = $V_1 \cup V_2$ is a neutrosophic bigroup. Consider Q (I) the neutrosophic field. V = $V_1 \cup V_2$ is a strong neutrosophic bivector space over Q (I).

**Remark:** As in case of vector spaces we see the neutrosophic bigroup V = $V_1 \cup V_2$ given in the above example is not a neutrosophic bivector space over R (I). Thus the definition neutrosophic bivector space is also dependent on the neutrosophic field over which it is defined.

***Example 4.1.38:*** Let V = $V_1 \cup V_2$ where $V_1$ = {set of all 3 × 3 neutrosophic matrices with entries from R(I)}. $V_1$ is a neutrosophic group under matrix addition. $V_2$ = {R (I) × Q(I) × R (I) × Q (I) } = {$x_1$, $x_2$, $x_3$, $x_4$) / $x_1$, $x_3$ ∈ R (I) and $x_2$, $x_4$ ∈ Q (I)}. Clearly $V_2$ is a neutrosophic group under component wise addition.

Now $V = V_1 \cup V_2$ is a strong neutrosophic bivector space over Q (I). Clearly $V = V_1 \cup V_2$ is not a strong neutrosophic bivector space over R (I). Infact $V = V_1 \cup V_2$ is not even a bivector space over R (I).

Now we proceed on to define the notion of neutrosophic bivector space over a field.

**DEFINITION 4.1.32:** *Let $V = V_1 \cup V_2$ be a neutrosophic bigroup. Suppose both $V_1$ and $V_2$ are neutrosophic vector spaces over a field F (F not a neutrosophic field) then we call $V = V_1 \cup V_2$ a neutrosophic bivector space over the field F.*

It is important to note a neutrosophic bivector space is never a strong neutrosophic bivector space for in case of a strong neutrosophic bivector space we need the underlying field over which the neutrosophic bivector space is defined to be a neutrosophic field, but in case of neutrosophic bivector space we demand only the underlying field to be just a field and not a neutrosophic field. Now we proceed on to give examples of such structures.

***Example 4.1.39:*** Let $V = V_1 \cup V_2$ where $V_1$ = R (I) × Q (I) × R (I) be the neutrosophic group under addition and $V_2$ = {The set of all 3 × 2 neutrosophic matrices with entries from R (I)}, $V_2$ is a neutrosophic group under matrix addition. Clearly V is a neutrosophic strong bivector space over Q (I). The dimension of V is infinite dimensional for $V_1$ as a neutrosophic vector space over Q (I) is infinite dimensional and $V_2$ as a neutrosophic vector space over Q (I) is also infinite dimensional; hence $V = V_1 \cup V_2$ is an infinite dimensional strong neutrosophic bivector space over Q (I).

***Note:*** Clearly $V = V_1 \cup V_2$ is not a bivector space over R (I) for $V_1$ is not even a vector space over R (I).
Now we proceed on to give an example of a finite dimensional strong neutrosophic bivector space.

***Example 4.1.40:*** Let $V = V_1 \cup V_2$ with $V_1 = \{Q (I) \times Q (I)\}$ and $V_2$ = {set of all 3 × 2 matrices with entries from Q (I)}. Clearly $V = V_1 \cup V_2$ is a neutrosophic bigroup under addition. Further V is a strong neutrosophic bivector space over Q (I). Now a basis of $V = V_1 \cup V_2$ is given by

$$\{(0\ 1), (1\ 0)\} \cup$$
$$\left\{ \begin{pmatrix} 1 & 0 \\ 0 & 0 \\ 0 & 0 \end{pmatrix}, \begin{pmatrix} 0 & 1 \\ 0 & 0 \\ 0 & 0 \end{pmatrix}, \begin{pmatrix} 0 & 0 \\ 1 & 0 \\ 0 & 0 \end{pmatrix}, \begin{pmatrix} 0 & 0 \\ 0 & 1 \\ 0 & 0 \end{pmatrix}, \begin{pmatrix} 0 & 0 \\ 0 & 0 \\ 1 & 0 \end{pmatrix}, \begin{pmatrix} 0 & 0 \\ 0 & 0 \\ 0 & 1 \end{pmatrix} \right\}$$
$$= B_1 \cup B_2.$$

Thus dimension of $V = V_1 \cup V_2$ as a strong neutrosophic bivector space over Q (I) is (2, 6) i.e. 8.

***Note:*** If we consider $V = V_1 \cup V_2$ the neutrosophic bigroup as a neutrosophic vector space over Q. The dimension of is finite but not (2, 6) or 8 for a basis for $V = V_1 \cup V_2$ over Q is given by

B = {0, 1), (1 0), (0 I), (I, 0) } ∪

$$\left\{ \begin{pmatrix} 1 & 0 \\ 0 & 0 \\ 0 & 0 \end{pmatrix}, \begin{pmatrix} 0 & 1 \\ 0 & 0 \\ 0 & 0 \end{pmatrix}, \begin{pmatrix} 0 & 0 \\ 1 & 0 \\ 0 & 0 \end{pmatrix}, \begin{pmatrix} 0 & 0 \\ 0 & 1 \\ 0 & 0 \end{pmatrix}, \begin{pmatrix} 0 & 0 \\ 0 & 0 \\ 1 & 0 \end{pmatrix}, \right.$$
$$\left. \begin{pmatrix} I & 0 \\ 0 & 0 \\ 0 & 0 \end{pmatrix}, \begin{pmatrix} 0 & I \\ 0 & 0 \\ 0 & 0 \end{pmatrix}, \begin{pmatrix} 0 & 0 \\ I & 0 \\ 0 & 0 \end{pmatrix}, \begin{pmatrix} 0 & 0 \\ 0 & I \\ 0 & 0 \end{pmatrix}, \begin{pmatrix} 0 & 0 \\ 0 & 0 \\ I & 0 \end{pmatrix}, \begin{pmatrix} 0 & 0 \\ 0 & 0 \\ 0 & I \end{pmatrix} \right\}$$

$= B_1^1 \cup B_2^1$.

Thus dimension of $V = V_1 \cup V_2$ as a neutrosophic bivector space over Q is given by B and the dimension of $V = V_1 \cup V_2$ as a neutrosophic bivector space over Q is (4, 12) = 16. Thus we note the difference.

It fact it is left as an exercise for the interested and innovative reader to prove that in case of neutrosophic vector spaces V over a field F and the same bigroup V over the neutrosophic field F(I) the dimension is 2n and n respectively where, n = dim $V_1$ + dim $V_2$ as a strong neutrosophic bivector space and 2n = dim $V_1$ + dim $V_2$ as a neutrosophic bivector space.

Several interesting results about basis linearly independent elements etc can be derived and analyzed in the case of strong neutrosophic bivector space and neutrosophic bivector space.

We just give yet another example of a strong neutrosophic bivector space and just a neutrosophic bivector space and find the dimensions of them.

***Example 4.1.41:*** Let $V = V_1 \cup V_2$ be a neutrosophic bigroup, $V_1$ = {set of all 3 × 3 neutrosophic matrices with entries from Q(I)} and $V_2$ = {set of all 2 × 4 neutrosophic matrices with entries from Q (I)}. Clearly $V = V_1 \cup V_2$ is a neutrosophic bigroup under addition. $V = V_1 \cup V_2$ can be realized as a strong neutrosophic bivector space over Q(I). A basis for the strong neutrosophic bivector space $V = V_1 \cup V_2$ is given by $B = B_1 \cup B_2$;

$B_1 \cup B_2 =$

$$\left\{\begin{pmatrix}0&0&1\\0&0&0\\0&0&0\end{pmatrix},\begin{pmatrix}0&1&0\\0&0&0\\0&0&0\end{pmatrix},\begin{pmatrix}1&0&0\\0&0&0\\0&0&0\end{pmatrix},\begin{pmatrix}0&0&0\\1&0&0\\0&0&0\end{pmatrix},\begin{pmatrix}0&0&0\\0&1&0\\0&0&0\end{pmatrix}\right.$$

$$\left.\begin{pmatrix}0&0&0\\0&0&1\\0&0&0\end{pmatrix},\begin{pmatrix}0&0&0\\0&0&0\\1&0&0\end{pmatrix},\begin{pmatrix}0&0&0\\0&0&0\\0&1&0\end{pmatrix},\begin{pmatrix}0&0&0\\0&0&0\\0&0&1\end{pmatrix}\right\} \cup$$

$$\left\{\begin{pmatrix}1&0&0&0\\0&0&0&0\end{pmatrix},\begin{pmatrix}0&1&0&0\\0&0&0&0\end{pmatrix},\begin{pmatrix}0&0&1&0\\0&0&0&0\end{pmatrix},\begin{pmatrix}0&0&0&1\\0&0&0&0\end{pmatrix}\right.$$

$$\left.\begin{pmatrix}0&0&0&0\\1&0&0&0\end{pmatrix},\begin{pmatrix}0&0&0&0\\0&1&0&0\end{pmatrix},\begin{pmatrix}0&0&0&0\\0&0&1&0\end{pmatrix},\begin{pmatrix}0&0&0&0\\0&0&0&1\end{pmatrix}\right\}.$$

Clearly dimension of the strong neutrosophic bivector space is 9 + 8 = 17.

Now consider the same neutrosophic bigroup $V = V_1 \cup V_2$ described in the above example. Clearly $V = V_1 \cup V_2$ is a neutrosophic bivector space over Q. Now what is the dimension of $V = V_1 \cup V_2$ as a neutrosophic bivector space over Q. The basis of $V = V_1 \cup V_2$ over Q.

The basis of $V = V_1 \cup V_2$ is given by $B^1 = B_1^1 \cup B_2^1 =$

$$\left\{\begin{pmatrix}0&0&1\\0&0&0\\0&0&0\end{pmatrix},\begin{pmatrix}0&1&0\\0&0&0\\0&0&0\end{pmatrix},\begin{pmatrix}1&0&0\\0&0&0\\0&0&0\end{pmatrix},\begin{pmatrix}0&0&0\\1&0&0\\0&0&0\end{pmatrix},\begin{pmatrix}0&0&0\\0&1&0\\0&0&0\end{pmatrix}\right.$$

$$\begin{pmatrix}0&0&0\\0&0&1\\0&0&0\end{pmatrix},\begin{pmatrix}0&0&0\\0&0&0\\1&0&0\end{pmatrix},\begin{pmatrix}0&0&0\\0&0&0\\0&1&0\end{pmatrix},\begin{pmatrix}0&0&0\\0&0&0\\0&0&0\end{pmatrix},\begin{pmatrix}I&0&0\\0&0&0\\0&0&0\end{pmatrix}$$

$$\begin{pmatrix}0&I&0\\0&0&0\\0&0&0\end{pmatrix},\begin{pmatrix}0&0&I\\0&0&0\\0&0&0\end{pmatrix},\begin{pmatrix}0&0&0\\0&I&0\\0&0&0\end{pmatrix},\begin{pmatrix}0&0&0\\I&0&0\\0&0&0\end{pmatrix}$$

$$\left.\begin{pmatrix}0&0&0\\0&0&I\\0&0&0\end{pmatrix},\begin{pmatrix}0&0&0\\0&0&0\\I&0&0\end{pmatrix},\begin{pmatrix}0&0&0\\0&0&0\\0&I&0\end{pmatrix},\begin{pmatrix}0&0&0\\0&0&0\\0&0&I\end{pmatrix}\right\} \cup$$

$$\left\{\begin{pmatrix}1&0&0&0\\0&0&0&0\end{pmatrix},\begin{pmatrix}0&1&0&0\\0&0&0&0\end{pmatrix},\begin{pmatrix}0&0&1&0\\0&0&0&0\end{pmatrix},\begin{pmatrix}0&0&0&0\\0&0&0&1\end{pmatrix}\right.$$

$$\begin{pmatrix}0&0&0&1\\0&0&0&0\end{pmatrix},\begin{pmatrix}0&0&0&0\\1&0&0&0\end{pmatrix},\begin{pmatrix}0&0&0&0\\0&1&0&0\end{pmatrix},\begin{pmatrix}0&0&0&0\\0&0&1&0\end{pmatrix}$$

$$\begin{pmatrix}I&0&0&0\\0&0&0&0\end{pmatrix},\begin{pmatrix}0&I&0&0\\0&0&0&0\end{pmatrix},\begin{pmatrix}0&0&I&0\\0&0&0&0\end{pmatrix},\begin{pmatrix}0&0&0&I\\0&0&0&0\end{pmatrix}$$

$$\left.\begin{pmatrix}0&0&0&0\\I&0&0&0\end{pmatrix},\begin{pmatrix}0&0&0&0\\0&I&0&0\end{pmatrix},\begin{pmatrix}0&0&0&0\\0&0&I&0\end{pmatrix},\begin{pmatrix}0&0&0&0\\0&0&0&I\end{pmatrix}\right\}$$

Thus dimension of $V = V_1 \cup V_2$ as a neutrosophic bivector space over Q is given by 18 + 16 = 34. Thus the dimension of the neutrosophic bivector space over Q is twice the dimension of the strong neutrosophic bivector space over Q (I).

The above statement remains to be true whenever the spaces have the same underlying bigroup and the same field and its neutrosophic field. If these facts are not taken care of the above results cannot be true, this we explain by a very simple example.

***Example 4.1.42:*** Let $V = V_1 \cup V_2$ be a neutrosophic bigroup where $V_1 = \{R (I) \times R (I) \}$, clearly $V_1$ is a neutrosophic group under addition. Take $V_2$ = {set of all polynomials of degree less than or equal to 3 with coefficients from R (I)}; $V_2$ is a neutrosophic group under '+'. Now $V = V_1 \cup V_2$ is a strong neutrosophic bivector space of finite dimension over R (I). Clearly a basis for $V = V_1 \cup V_2$ over R (I) is given by {(1 0), (0 1)} $\cup$ $\{x, x_1 x^2, x^3\}$.

So dimension of the strong neutrosophic bivector space over R(I) is $\{2 \cup 4\}$= 6. Now $V = V_1 \cup V_2$ is a neutrosophic bivector space over Q and dimension of V as a neutrosophic bivector space of over Q is infinite. Thus instead of R if we use the field Q we get no relation between dimensions of strong neutrosophic bivector space and neutrosophic bivector space.

Now we proceed on to define the notion of linear bitransformation and linear bioperator on bivector spaces.

**DEFINITION 4.1.33:** *Let $V = V_1 \cup V_2$ and $W = W_1 \cup W_2$ be neutrosophic bigroups which are strong neutrosophic bivector spaces over the same neutrosophic field F (I). A neutrosophic linear bitransformation $T = T_1 \cup T_2$. (The symbol '$\cup$' is just only a notation and it has nothing to do with the usual union.) is defined from V to W by the neutrosophic linear transformation $T_1 : V_1 \to W_1$ and $T_2: V_2 \to W_2$ such that $T_i (I) = I$, i = 1, 2.*

On similar lines we define neutrosophic linear bitransformation of neutrosophic bivector spaces over any real field F.

*Note:* The set of all linear bitransformation of a strong neutrosophic bivectors spaces over F (I) form a strong neutrosophic bivector space over F (I). Suppose

$$T = T_1 \cup T_2: V \to W \ (V_1 \cup V_2 \to W_1 \cup W_2)$$

be the linear bitransformation of strong neutrosophic bivector spaces over F (I). Then the bikernel of T is defined to be

$$\text{BKer } T = \{\text{BKer } T_1\} \cup \{\text{BKer } T_2\}.$$

Clearly BKer T is a subset $V = V_1 \cup V_2$. Since $T_1 (I) = I$, $I \notin \text{BKer } T_1$ or $I \notin \text{BKer } T_2$. So

$$I \notin \text{BKer } (T_1 \cup T_2) = \text{BKer } T_1 \cup \text{BKer } T_2 = \text{BKer } T.$$

Further if $T = T_1 \cup T_2$ is a zero linear bitransformation, if $T(x) = (0)$ i.e. $T_1 (x_1) = 0$ and $T_2 (x_2)$ for all $x = x_1 \cup x_2 \in V = V_1 \cup V_2$ we denote this linear bitransformation by $T^o = T_1^o \cup T_2^o$.

In this case Bker $T^o = V = V_1 \cup V_2$.

Thus this study of neutrosophic spaces stands as a separate algebraic structure for indeterminacy is different from zero. If the indeterminate in general is mapped on to zero it has no meaning, the space itself has no meaning. Further we make an assumption only in case of the zero bitransformation as all elements are marked to zero so also I. Thus a linear bitransformation never maps I to zero but I to I only.

Thus as in case of linear transformation of vector spaces or bivector spaces we cannot say in case of neutrosophic vector spaces and neutrosophic bivectors spaces the kernel is a neutrosophic subspace or (bi subspace) of the domain space. Some amount of caution must be applied while studying bitransformations. We will define a pseudo subspace of a neutrosophic bivector space.

**DEFINITION 4.1.34:** *Let $V = V_1 \cup V_2$ be a strong neutrosophic bivector space over the neutrosophic field $\langle F \cup I \rangle$ = F(I). If $W = W_1 \cup W_2$ is a subbigroup of the neutrosophic bigroup $V_1 \cup V_2$ but not a neutrosophic subbigroup of G and if $W = W_1 \cup W_2$ happens to be the bivector space over $F \subset F \subset I$) then we call $W = W_1 \cup W_2$ the pseudo neutrosophic bivector space. If $V = V_1 \cup V_2$ happens to be a neutrosophic bivector space over F and if $W = W_1 \cup W_2$ is a subbigroup of V but not a neutrosophic subbigroup of G and if W happens to be a bivector space over F then we call W the pseudo neutrosophic bisubspace of the neutrosophic bivector space of $V = V_1 \cup V_2$.*

It is left as an exercise for the reader to prove for any linear bitransformation $T = T_1 \cup T_2$ from $V = V_1 \cup V_2$ to $W = W_1 \cup W_2$ the Bker T = Bker $T_1 \cup$ Bker $T_2$ is the pseudo neutrosophic bisubspace of $V = V_1 \cup V_2$.

***Example 4.1.43:*** Let $V = V_1 \cup V_2$ where $V_1$ = Q (I) × Q (I) × Q(I) be a neutrosophic group under addition. Let $V_2$ = {the collection of all 3 × 3 neutrosophic matrices with entries from Q (I)}, $V_2$ is a neutrosophic group under '+'. Clearly $V = V_1 \cup V_2$ is a neutrosophic bigroup. Further $V = V_1 \cup V_2$ is a strong neutrosophic bivector space over Q (I) or $V = V_1 \cup V_2$ is just the neutrosophic bivector space over Q. Now consider $W = W_1 \cup W_2$ where $W_1 = Q \times Q \times Q$ and $W_2$ = {the set of all $3 \times 3$ matrices with entries from Q}. Clearly $W = W_1 \cup W_2$ is a subbigroup of V. Infact W is pseudo neutrosophic bivector space over Q.

Now we have the following interesting result for which we have to define the notion of complete bivector pseudo subspace of a neutrosophic bivector space V.

**DEFINITION 4.1.35:** *Let $V = V_1 \cup V_2$ be a strong neutrosophic bivector (neutrosophic bivector space) over the neutrosophic field F (I). Suppose $W = W_1 \cup W_2$ is a subbigroup is $V_1 \cup V_2 = V$ and suppose $W = W_1 \cup W_2$ is the pseudo neutrosophic subspace of $V = V_1 \cup V_2$ and if $U = U_1 \cup U_2$ is any other pseudo neutrosophic subspace of V then $W \subset U$ implies and is implied by U = W i.e. $U_1 = W_1$ and $U_2 = W_2$ then we call W the complete pseudo neutrosophic bisubspace of V.*

**THEOREM 4.1.4:** *Let $V = V_1 \cup V_2$ be a strong neutrosophic bivector space over F (I) (or a neutrosophic bivector space over F). Every strong neutrosophic bivector space (or neutrosophic bivector space) has a complete pseudo neutrosophic bivector space.*

***Proof:*** The proof follows from the very definition.

Now we can state or visualize the notion of strong neutrosophic bivector space or (neutrosophic bivector space) as an extension of a bivector space over the field or in general a generalization of the bivector space.

***Example 4.1.44:*** Let $V = V_1 \cup V_2$ where $V_1$ = {set of all $2 \times 5$ matrices with entries from R (I)} and $V_2$ = {R (I) × R (I) × R (I)}. Now $V = V_1 \cup V_2$ is a strong neutrosophic bivector space over R (I). Consider the subspace $W = W_1 \cup W_2$ of $V = V_1 \cup V_2$, here $W_1$ = {set of all $2 \times 5$ matrices with entries from Q} and $W_2$ = {Q × Q × Q}. Clearly W =

$W_1 \cup W_2$ is a pseudo neutrosophic subbispace of V only over Q. Take $W^1 = W_1^1 \cup W_2^1$, $W_1^1$ = {set of all 2 × 5 matrices over R} and $W_2$ = {R × R × R}. $W_1^1 \cup W_2^1 = W^1$ is a pseudo neutrosophic subbispace of V over Q, still $W^1$ is not the complete pseudo neutrosophic bispace. Take U = $U_1 \cup U_2$ where $U_1$ = {set of all 2 × 5 with entries from R} and $U_2$ = {R × R × R}. U = $U_1 \cup U_2$ is a pseudo neutrosophic subbispace over R. U is the complete pseudo neutrosophic subbispace over R.

Now we proceed on to define quasi neutrosophic bivector space and quasi neutrosophic subbivector space over a field.

**DEFINITION 4.1.36:** *Let V = $V_1 \cup V_2$ be a bigroup such that $V_1$ is group and $V_2$ is a neutrosophic group. Then we call V = $V_1 \cup V_2$ to be the quasi neutrosophic bigroup.*

Now we just give an example before we proceed on to define the notion quasi neutrosophic bivector space and show we cannot define strong quasi neutrosophic bivector space.

***Example 4.1.45:*** Let V = $V_1 \cup V_2$ be a quasi neutrosophic bigroup, where $V_1$ = {Q (I) × Q (I) × Q (I)} a neutrosophic group under addition, $V_2$ = {set of all 3 × 2 matrices with entries from Q}, $V_2$ is a group under matrix addition. Clearly V = $V_1 \cup V_2$ is a quasi neutrosophic bivector space over Q.

*Note:* We cannot define the notion of strong quasi neutrosophic bivector space for we cannot make a quasi bigroup to be a strong neutrosophic bivector space over a neutrosophic field.

**DEFINITION 4.1.37:** *Let V = $V_1 \cup V_2$ be a neutrosophic bivector space over a field F. If a subbispace W = $W_1 \cup W_2$ where W is just a quasi subbigroup of V is a bivector space* *over F then we call $W = W_1 \cup W_2$ to be a quasi neutrosophic subbispace of V.*

Now we illustrate this situation by the following example:

***Example 4.1.46:*** Let $V = V_1 \cup V_2$ where $V_1 = \{Q (I) \times Q(I) \times Q(I) \}$ a neutrosophic group under addition and $V_2 = \{Q(I) (x) \times Q (I) (x)$ i.e. the collection of all polynomials with coefficient from $Q (I)\}$, $V_2$ is also a neutrosophic bigroup under component wise polynomial addition. Consider $(Q \times Q \times Q) = W_1$ and $W_2 = \{Q(I) \times Q(I)\}$. Clearly $W = W_1 \cup W_2$ is a quasi neutrosophic subbispace over Q.

Now we proceed on to define the notion of quasi neutrosophic bivector space linear bitransformation.

**DEFINITION 4.1.38:** *Let $V = V_1 \cup V_2$ and $W = W_1 \cup W_2$ be a two quasi neutrosophic bivector spaces over the field F. Define a linear bitransformation T from V to W i.e. $T = T_1 \cup T_2$ where $T_1 : V_1 \to W_1$ is a linear transformation of the neutrosophic vector space $V_1$ to the neutrosophic vector space $W_1$ and $T_2 : V_2 \to W_2$ is the linear transformation of vector spaces $V_2$ to $W_2$. The we call $T = T_1 \cup T_2$ to be a neutrosophic quasi linear bitransformation.*

We illustrate this by the following example:

***Example 4.1.47:*** Let $V = V_1 \cup V_2$ be a quasi neutrosophic bivector space over Q, where $V_1$ = {all polynomials of degree less than or equal to 3 with coefficients from Q} and $V_2 = \{Q(I) \times Q(I)\}$ and $W = W_1 \cup W_2$ be a quasi neutrosophic bivector space with $W_1 = \{Q \times Q \times Q \times Q\}$ group under component wise addition and

$$W_2 = \left\{ \begin{pmatrix} a & b \\ c & d \end{pmatrix} \Big| \, a, b, c, d \in Q(I) \right\}.$$

Define a quasi linear bitransformation T from V to W as follows :

$$T = T_1 \cup T_2 : V = V_1 \cup V_2 \to W = W_1 \cup W_2$$

$$T_1 (x, y, z, w) = (x + y, y + z, z + w, x + y + w)$$

and

$$T_2 (a, b) = \begin{pmatrix} a & b \\ 0 & 0 \end{pmatrix}; a, b \in V_2.$$

Clearly $T = T_1 \cup T_2$ is a quasi linear bitransformation of V to W.

Now we proceed on to define the notion of quasi linear bioperator.

**DEFINITION 4.1.39:** *Let $V = V_1 \cup V_2$ be a quasi neutrosophic bivector space defined over a field F. A linear bioperator T from V to V is called the quasi linear bioperator on V.*

We illustrate this by the following example:

***Example 4.1.48:*** Let $V = V_1 \cup V_2$ be a quasi neutrosophic bivector space over Q, where $V_1 = Q (I) \times Q (I) \times Q (I)$ and $V_2$ = {all $2 \times 2$ matrices with entries from Q}.

Let

$T : V \to V$

$T_1 \cup T_2 : V_1 \cup V_2 \to V_1 \cup V_2$

i.e. $(T_1 : V_1 \to V_1 ) \cup (T_2 : V_2 \to V_2)$

$T_1 (\alpha_1, \alpha_2, \alpha_3) = ( \alpha_1 + \alpha_3, \alpha_2 - \alpha_1, \alpha_3)$

$T_2 (x_1, x_2, x_3, x_4) = (x_1 + x_2, x_1 - x_3, x_2, x_4)$.

$T = T_1 \cup T_2$ is a quasi neutrosophic bioperator of the quasi neutrosophic bivector space $V = V_1 \cup V_2$.

Now one can as in case of bivector spaces define eigen bivalues and eigen bivectors for any related characteristic bipolynomial. Clearly the bivalues, bivector and bipolynomials are in general neutrosophic. All results derived in case of bivector spaces can also be derived in case of strong neutrosophic bivector spaces, neutrosophic bivector spaces, pseudo neutrosophic bivector spaces and quasi neutrosophic bivector spaces with minor and appropriate modifications in the respective definitions and derived all results.

Now we proceed on to define the notion of for the first time Smarandache neutrosophic vector spaces and Smarandache neutrosophic bivector spaces. To this end we need to define the notion of Smarandache neutrosophic semigroup, Smarandache neutrosophic bisemigroup and Smarandache neutrosophic fields.

**DEFINITION 4.1.40:** *Let S be a non empty set. S is said to be neutrosophic semigroup if S is a semi group and I $\in$ S.*

***Example 4.1.49:*** Let S = {set of all 2 × 2 matrices with entries from Q (I)}. S under matrix multiplication is a semigroup called the neutrosophic semigroup.

Now proceed on to define the notion of Smarandache neutrosophic semigroup and illustrate it with examples.

**DEFINITION 4.1.41:** *Let S be a neutrosophic semigroup. S is said to be a Smarandache neutrosophic semigroup (S-neutrosophic semigroup) if S contains a subset a such that S under the operations of S is a neutrosophic group.*

***Example 4.1.50:*** Let S = {set of all n × n matrices with entries from Q(I)}. S under matrix multiplication is a neutrosophic semigroup.

Consider the subset A of S where A = {set of all n × n matrices with entries from Q (I) such that the matrices are non singular i.e. their determinant is non zero}. Clearly S is a Smarandache neutrosophic semigroup. (S-neutrosophic semigroup).

Let us give one more example of a Smarandache neutrosophic semigroup.

***Example 4.1.51:*** Let S = {Q (I) × Q (I)}, S under component wise multiplication is a semigroup. Let A = {Q (I) \ {0}) × (Q (I) \ {0})} which is a proper subset of S. A under multiplication is a group. Thus S is a Smarandache neutrosophic semigroup.

Suppose we have a neutrosophic semigroup S which has a proper subset A such that A is only a group under the operations of S and not a neutrosophic group.

We define a notion of Smarandache quasi neutrosophic semigroup.

**DEFINITION 4.1.42:** *Let S be a neutrosophic semigroup. If in S we have a non empty subset A which is only a group but not a neutrosophic group then we call S a Smarandache quasi neutrosophic semigroup.*

We illustrate this with examples.

***Example 4.1.52:*** Let S be a neutrosophic semigroup given by S = {Z (I)}. S is a neutrosophic semigroup under multiplication. Take A = {-1, 1}. A is a group. So is a Smarandache quasi neutrosophic semigroup.

Now we proceed on to define the notion of Smarandache neutrosophic ring.

**DEFINITION 4.1.43:** *Let S(I) denote a neutrosophic set. If S(I) under '+' is a group and S(I) under '•' Is a semigroup so that (S(I), +, •) is a ring i.e., both the distributive laws are satisfied then we call S(I) a neutrosophic ring.*

Let Z (I) = {the ring generated by Z and I} = {a + bI / $I^2 = I$ such that a, b ∈ Z}. Thus for any ring R we can build a neutrosophic ring. Let $M^I_{n\times n}$ = {set of all n × n matrices with entries from Z (I)}. $M^{I}_{n\times n}$ is a neutrosophic ring under usual matrix addition and matrix multiplication.

Let R = Q (I) × Q (I), R is a neutrosophic ring under component wise addition and multiplication. Now we proceed on to define the notion of Smarandache Neutrosophic Ring.

**DEFINITION 4.1.44:** *Let S (I) be a neutrosophic ring if there exists proper subset A of S (I) which is a neutrosophic field then we call S (I) a Smarandache Neutrosophic ring (S-neutrosophic ring).*

***Example 4.1.53:*** Let S (I) = Q (I) × Q (I) be a neutrosophic ring. S (I) is a Smarandache Neutrosophic ring for A = Q(I) × {0} is a proper subset of S (I) such that A is a neutrosophic field. Thus S (I) is a Smarandache neutrosophic ring.

***Example 4.1.54:*** Let S(I) = {set of all n × n matrices with entries from Q (I)}.
Let

$$A = \left\{ \begin{pmatrix} a_{11} 0 & .... \\ .... & 0 \end{pmatrix} \middle| a_{11} \in Q(I) \right\}.$$

A is a neutrosophic field and is a proper subset of S (I). Thus S (I) is a Smarandache neutrosophic ring.

***Example 4.1.55:*** Let Q (I) [x] be the neutrosophic polynomial ring. Clearly A = Q (I) $\underset{\neq}{\subset}$ Q (I) [x]. A is a neutrosophic field so Q (I) [x] is a Smarandache neutrosophic ring.

Now we proceed on to define get another concept called the Smarandache quasi neutrosophic ring.

**DEFINITION 4.1.45:** *Let S (I) be a neutrosophic ring. S (I) is said to be a Smarandache quasi neutrosophic ring (S-quasi*

*neutrosophic ring) if S (I) has a proper nonempty subset, A such A is a field and not a neutrosophic field.*
*Thus a neutrosophic ring in general can be both S-neutrosophic as well as S-quasi neutrosophic ring.*

Now we illustrate this concept by the following examples:

***Example 4.1.56:*** Let Q (I) [x] be a neutrosophic ring Q (I) [x] is a Smarandache quasi neutrosophic ring as Q $\subset$ Q (I) [x] and Q is a ring.

***Example 4.1.57:*** Let R (I) be a neutrosophic ring. R(I) is a Smarandache quasi neutrosophic ring for Q $\subseteq$ R (I) and Q is a field. R (I) is also a Smarandache neutrosophic ring as Q(I) $\subseteq$ R (I) is a neutrosophic field.

Now we proceed on to define the notion of Smarandache neutrosophic vector spaces.

## 4.2 Smarandache neutrosophic linear bialgebra

In this section we just define the notion of neutrosophic linear bialgebra and illustrate it with examples.

**DEFINITION 4.2.1:** *Let R be a S-ring V is said to be a neutrosophic bimodule over R, if V is an abelian neutrosophic bigroup under operation + and for every r $\in$ R and v $\in$ V we have r v and v r are in V subject to*

1. $r (v_1 + v_2) = rv_1 + rv_2$
2. $r (s\ v) = (r\ s)\ v$
3. $(r + s)\ v = r\ v + s\ v$

*for all v, $v_1 \in$ V and r, s t $\in$ R.*

Now we illustrate this by the following example.

***Example 4.2.1:*** Let S = {set of all 3 × 2 matrices with entries from Z (I)}. S under '+' is a neutrosophic abelian group. Clearly S is a neutrosophic module over Z.

***Example 4.2.2:*** Let Z (I) [x] be the polynomials in x with coefficients from Z (I). Z (I) [x] is a neutrosophic abelian group and Z (I) [x] is a neutrosophic module over Z.

Now we proceed on to define Smarandache neutrosophic vector space.

**DEFINITION 4.2.2:** *Let S be a Smarandache ring. V be an abelian neutrosophic group. Let V be a neutrosophic module over S V is said to be a Smarandache neutrosophic vector space (S-neutrosophic vector space) if V is a neutrosophic vector space over a proper subset K of S where K is a field If V is also closed under an external binary operation multiplication which is associative, then we call V a Smarandache neutrosophic linear algebra (S-neutrosophic linear algebra).*

Now we illustrate this by the following examples:

***Example 4.2.3:*** Let Q [x] be a S-ring V = {R (I) [x] i.e. the set of all polynomials with coefficients from the neutrosophic field R (I)}. V under addition is an abelian group. V is a neutrosophic module over Q [x]. Now V is a Smarandache Neutrosophic Vector space over Q [x] as V is a vector space over Q.

***Example 4.2.4:*** Let V = {3 × 2 matrices with entries from R(I) } V is a neutrosophic abelian group under addition V is a neutrosophic module over R and V is also Smarandache neutrosophic vector space over R.

As our study pertains to only to linear bialgebra and neutrosophic linear bialgebra.

We proceed on to define Smarandache neutrosophic vector spaces.

***Example 4.2.5:*** Let V = Q (I) × Q (I) ∪ R (I) (x) $V_1 \cup V_2$ be a neutrosophic bisemigroup, i.e. $V_1$ is a neutrosophic semigroup under component wise multiplication and $V_2$ is a semigroup under polynomial multiplication. Clearly $V = V_1 \cup V_2$ is a neutrosophic bisemigroup. We call V to be a Smarandache neutrosophic bisemigroup if bother $V_1$ and $V_2$ are Smarandache neutrosophic semigroups. In the above example $V_1$ = Q (I) × { 1} is a group under component wise multiplication R (I) \ {o} ⊂ R (I) [x] is group under multiplication so R (I) [x] is a Smarandache neutrosophic semigroup.

**DEFINITION 4.2.3:** *Let $S = S_1 \cup S_2$ be a neutrosophic bisemigroup. S is said to be a Smarandache neutrosophic bisemigroup (S-neutrosophic bisemigroup) if S has a proper subset $A = A_1 \cup A_2$ such that A under the operations of S is a neutrosophic bigroup.*

We illustrate this get by an example.

***Example 4.2.6:*** Let Q (I) × Q [x] × Q [x] ∪ R (I) = S be a neutrosophic bisemigroup under multiplication. Clearly S is a Smarandache neutrosophic bisemigroup for A = Q(I) \ {0} ∪ {1} ∪ {1} ∪ R(I) \ {0} a neutrosophic bigroup under component wise multiplication. In fact this S has more than one neutrosophic bigroup.

We give a few types of Smarandache neutrosophic bivector spaces which we call as Type $I_B$ and Type $II_B$.

**DEFINITION 4.2.4:** *Let R be a S- ring $V = V_1 \cup V_2$ is a Smarandache neutrosophic bivector space of Type $I_B$ (S-neutrosophic bivector space of type $I_B$ ) if V is a bivector space over a proper subset K of R where K is a field. i.e. if $V_1$ is a neutrosophic module over a S-ring R and $V_2$ is a neutrosophic module over R and both $V_1$ and $V_2$ are* *neutrosophic vector spaces over a proper subset $K \underset{\neq}{\subset} R$ where K is a field.*

We define here yet another new type of Smarandache neutrosophic bivector space which we choose to call as Type $II_B$. However we wish to mention these two types type $I_B$ and type $II_B$ are unrelated.

**DEFINITION 4.2.5:** *Let F be a field. Suppose $B = B_1 \cup B_2$ be a Smarandache neutrosophic bisemigroup. If B has a neutrosophic bigroup. $V = V_1 \cup V_2$ i.e. $V_1 \subset B_1$ and $V_2 \subset B_2$ such that V is a bivector space over F then we call B a Smarandache neutrosophic bivector space of type $II_B$ over the field F.*

Now we define still a new version of Smarandache neutrosophic bivector space which we call as strong Smarandache neutrosophic bivector space.

**DEFINITION 4.2.6:** *Let R be a S-ring. $B = B_1 \cup B_2$ be a Smarandache neutrosophic bisemigroup. If V is a neutrosophic bigroup of B and V is a neutrosophic bivector space over the field K where $K \subseteq R$, then we say V is a Strong Smarandache neutrosophic bivector space (S-neutrosophic bivector space) over R..*

***Example 4.2.7:*** Let $B = B_1 \cup B_2$ where $B_1 = R (I) \times R$ and $B_2 = Q (I) \times Q (x)$; $B = B_1 \cup B_2$ is a neutrosophic bisemigroup under component wise multiplication $G = R (I) \setminus \{0\} \times \{1\} \cup \{\{I\} \cup Q \setminus \{0\}\}$ is a neutrosophic bigroup. So B is a Smarandache neutrosophic bisemigroup. Thus B is a strong Smarandache neutrosophic bivector space over $Q \subseteq R$ where R is the S-ring R the field of reals.

***Example 4.2.8:*** Let $B = B_1 \cup B_2 = R (I) [x] \cup$ {set of all $n \times n$ matrices with entries from B is a neutrosophic bisemigroup. In fact B is a Smarandache neutrosophic bisemigroup for $G = R (I) \setminus \{0\} \cup$ {set of all $n \times n$ matrix with determinant non zero with entries from R (I)} is a bigroup. Clearly G is a bivector space over Q ⊂ R (R the field of reals is a S-ring). So B is a strong Smarandache neutrosophic bivector space over R.

Interested readers can get several such examples. All notions of Smarandache basis etc can be defined as in case of Smarandache vector spaces. Further linear bitransformation and linear bioperator can also be derived as in case of neutrosophic bivector spaces. Thus these derivation and results can be treated as a routine with appropriate modification. Now we proceed on to give some of the important applications as we have elaborately done these work in case neutrosophic bivector spaces and Smarandache vector spaces in the following sections.

We now give some basic application of Smarandache neutrosophic bivector spaces to Smarandache representation of finite Smarandache bisemigroup, a smattering of neutrosophic logic using Smarandache bivector spaces of type $II_B$. Smarandache Markov chains using S-bivector spaces II Smarandache Leontief economic models and Smarandache anti linear bialgebra.

## 4.3 Smarandache representation of finite Smarandache bisemigroup

Here we for the first time define Smarandache representation of finite S-bisemigroup. We know every S-bisemigroup, $S = S_1 \cup S_2$ contains a bigroup $G = G_1 \cup G_2$. The Smarandache representation S-bisemigroups depends on the S-bigroup G which we choose. Thus this method widens the Smarandache representations. We first define the notion of Smarandache pseudo neutrosophic bisemigroup.

**DEFINITION 4.3.1:** *Let $S = S_1 \cup S_2$ be a neutrosophic bisemigroup. If S has only bigroup which is not a neutrosophic bigroup, then we all S a Smarandache pseudo neutrosophic bisemigroup (S - pseudo neutrosophic bisemigroup).*

***Example 4.3.1:*** Let $S = S_1 \cup S_2$ where $S_1 = Q(I) \times Q(I)$ and $S_2 = \{2 \times 2$ matrices with entries from $Q(I)\}$ both $S_1$ and $S_2$ under multiplication is a semigroup. Thus S is a neutrosophic bisemigroup. Take $G = G_1 \cup G_2$ where $G_1 = \{(Q \setminus \{0\}) \times (Q \setminus \{0\})\}$ and $G_2 = \{$set of all $2 \times 2$ matrices A with entries from Q such that $|A| \neq 0\}$. $G_1$ and $G_2$ are groups under multiplication. So S is a pseudo Smarandache Neutrosophic bisemigroup.

Now we give the Smarandache representation of finite pseudo Smarandache neutrosophic bisemigroups.

**DEFINITION 4.3.2:** *Let $G = G_1 \cup G_2$ be a Smarandache neutrosophic bisemigroup and $V = V_1 \cup V_2$ be a bivector space. A Smarandache birepresentation of G on V is a mapping $S_\rho = S_\rho^1 \cup S_\rho^2$ from $H_1 \cup H_2$ ($H_1 \cup H_2$ is a subbigroup of G which is not a neutrosophic bigroup) to invertible linear bitransformation on $V = V_1 \cup V_2$ such that*

$$S_{\rho_{xy}} = S^1_{\rho_{x_1 y_1}} \cup S^2_{\rho_{x_2 y_2}} = \left(S^1_{\rho_{x_1}} \circ S^1_{\rho_{y_1}}\right) \cup \left(S^2_{\rho_{x_2}} \circ S^2_{\rho_{y_2}}\right)$$

*for all $x_1, y_1 \in H_1$ and for all $x_2, y_2 \in H_2$. $H_1 \cup H_2 \subset G_1 \cup G_2$. Here $S_{\rho_x} = S^1_{\rho_{x_1}} \cup S^2_{\rho_{x_2}}$; to denote the invertible linear bitransformation on $V = V_1 \cup V_2$ associated to $x = x_1 \cup x_2$ on $H = H_1 \cup H_2$, so that we may write*

$$S_{\rho_x}(v) = S_{\rho_x}(v_1 \cup v_2) = S^1_{\rho_{x_1}}(v_1) \cup S^2_{\rho_{x_2}}(v_2)$$

*for the image of the vector $v = v_1 \cup v_2$ in $V = V_1 \cup V_2$ under $S_{\rho_x} = S^1_{\rho_{x_1}} \cup S^2_{\rho_{x_2}}$. As a result we have that $S_{\rho_e} = S^1_{\rho_{e_1}} \cup S^2_{\rho_{z_2}} = I^1 \cup I^2$ denotes the identity bitransformation on $V = V_1 \cup V_2$ and*

$$S_{\rho_{x^{-1}}} = S^1_{\rho_{x_1^{-1}}} \cup S^2_{\rho_{x_2^{-1}}} = \left(S^1_{\rho_{x_1}}\right)^{-1} \cup \left(S^2_{\rho_{x_2}}\right)^{-1}$$

*for all* $x = x_1 \cup x_2 \in H_1 \cup H_2 \subset G_1 \cup G_2 = G$.

In other words a birepresentation of $H = H_1 \cup H_2$ on $V = V_1 \cup V_2$ is a bihomomorphism from H into GL (V) i.e. ($H_1$ into GL ($V_1$)) $\cup$($H_2$ into GL ($V_2$)). The bidimension of $V = V_1 \cup V_2$ is called the bidegree of the representation.

Thus depending on the number of subbigroup of the S-neutrosophic bisemigroup we have several S-birepresentations of the finite S-neutrosophic bisemigroup.

Basic example of birepresentation would be Smarandache left regular birepresentation and Smarandache right regular birepresentation over a field K defined as follows.

We take $V_H = V_{H_1} \cup V_{H_2}$ to be a bivector space of bifunctions on $H_1 \cup H_2$ with values in K (where $H = H_1 \cup H_2$ is a subbigroup of the S-neutrosophic bisemigroup where H is not a neutrosophic bigroup). For Smarandache left regular birepresentation (S-left regular biregular representative) relative to $H = H_1 \cup H_2$ we define

$$SL_x = S^1 L^1_{x_1} \cup S^2 L^2_{x_2} = \left(S^1 \cup S^2\right)\left(L^1 \cup L_2\right)_{x_1 \cup x_2}$$

from $V_{H_1} \cup V_{H_2} \to V_{H_1} \cup V_{H_2}$ for each $x_1 \cup x_2 \in H = H_1 \cup H_2$ by for each $x = x_1 \cup x_2$ in $H = H_1 \cup H_2$ by $SL_x$ (f ) (z) = f ($x^{-1}$z) for each function f (z) in $V_H = V_{H_1} \cup V_{H_2}$

i.e. $S^1 L^1_{x_1} f_1(z_1) \cup S^2 L^2_{x_2} f_2(z_2)$

$= f_1(x_1^{-1} z_1) \cup f_2(x_2^{-1} z_2)$.

For the Smarandache right regular birepresentation (S-right regular birepresentation) we define $SR_x = SR_{x_1 \cup x_2}$ : $V_{H_1} \cup V_{H_2} \to V_{H_1} \cup V_{H_2}$ ; $H_1 \cup H_2 = H$ for each $x = x_1 \cup x_2 \in H_1 \cup H_2$ by $SR_x$ (f ) (z) = f (zx).

$S^1 R^1_{x_1} f_1(z_1) \cup S^2 R^2_{x_2}\left(f_2(z_2)\right)$

$= f_1(z_1 x_1) \cup f_2(z_2 x_2)$

for each function $f_1(z_1) \cup f_2(z_2) =$

$f(z)$ in $V_H = V_{H_1} \cup V_{H_2}$.

Thus if $x = x_1 \cup x_2$ and $y = y_1 \cup y_2$ are elements $H_1 \cup H_2 \subset G_1 \cup G_2$

Then $\left(SL_x \circ SL_y\right)\left(f(z)\right)$

$$
\begin{aligned}
&= SL_x (SL_y) (f) (z) \\
&= (SL_y (f)) x^{-1} z \\
&= f (y^{-1} x^{-1} z) \\
&= f_1 \left(y_1^{-1} x_1^{-1} z_1\right) \cup f_2 \left(y_2^{-1} x_2^{-1} z_2\right) \\
&= f_1 ((x_1 y_1)^{-1} z_1) \cup f_2 ((x_2 y_2)^{-1} z_2) \\
&= S^1 L^1_{x_1 y_1} (f_1)(z_1) \cup S^2 L^2_{x_2 y_2} f_2 (z_2) \\
&= \left[\left(S^1 L^1_{x_1} \cup S^2 L^2_{x_2}\right) \left(S^1 L^1_{y_1} \cup S^2 L^2_{y_2}\right)\right] f(z_1 \cup z_2) \\
&= \left[\left(S^1 L^1_{x_1} \cup S^2 L^2_{x_2}\right) \circ \left(S^1 L^1_{y_1} \cup S^2 L^2_{y_2}\right)\right]\left(f_1 (z_1) \cup f_2 (z_2)\right).
\end{aligned}
$$

and

$(SR_x \circ SR_y) (f) (z)$

$$
\begin{aligned}
&= \left(S^1 R^1_{x_1} \cup S^2 R^2_{x_2}\right) \circ \left(S^1 R^1_{y_1} \cup S^2 R^2_{y_2}\right) (f)(z) \\
&= SR_x (SR_y (f)) (z) \\
&= \left(S^1 R^1_{x_1} \cup S^2 R^2_{x_2}\right) \circ \left(S^1 R^1_{y_1} \cup S^2 R^2_{y_2}\right)(f)(z) \\
&= \left(S^1 R^1_{y_1} (f_1) \cup S^2 R^2_{y_2} (f_2)\right)\left(z_1 x_1 \cup z_2 x_2\right) \\
&= f_1 (z_1 x_1 y_1) \cup f_2 (z_2 x_2 y_2) \\
&= S^1 R^1_{x_1 y_1} f_1 (z_1) \cup S^2 R^2_{x_2 y_2} f_2 (z_2) \\
&= S R_{xy} (f) (z).
\end{aligned}
$$

Thus for a given S-neutrosophic bisemigroup we will have several V's associated with them i.e. bivector space functions on each $H_1 \cup H_2 \subset G_1 \cup G_2$, H a subbigroup of the S-neutrosophic bisemigroup with values from K. This study in this direction is innovative.

We have yet another Smarandache birepresentation which can be convenient is the following. For each $w = w_1 \cup w_2$ in $H = H_1 \cup H_2$, H bisubgroups of the S-neutrosophic bisemigroup $G = G_1 \cup G_2$.

Define a bifunction

$$\phi_w(z) = \phi^1_{w_1}(z_1) \cup \phi^2_{w_2}(z_2)$$

on $H_1 \cup H_2 = H$ by

$$\phi^1_{w_1}(z_1) \cup \phi^2_{w_2}(z_2) = 1 \cup 1$$

where $w = w_1 \cup w_2 = z = z_1 \cup z_2$

$$\phi^1_{w_1}(z_1) \cup \phi^2_{w_2}(z_2) = 0 \cup 0$$

when $z \neq w$.

Thus the functions $\phi_w = \phi^1_{w_1} \cup \phi^2_{w_2}$ for $w = w_1 \cup w_2$ in $H = H_1 \cup H_2$ ($H \subset G$) form a basis for the space of bifunctions on each $H = H_1 \cup H_2$ contained in $G = G_1 \cup G_2$.

One can check that

$$SL_x\left(\phi_w\right) = \left(\phi_{xw}\right)$$

$$\text{i.e. } S^1 L^1_{x_1}\left(\phi_{w_1}\right) \cup S^2 L^2_{x_2}\left(\phi_{w_2}\right) = \phi^1_{x_1 w_1} \cup \phi^2_{x_{12} w_2}$$

$$SR_x(\phi_w) = \phi_{xw}$$

$$\text{i.e. } S^1 R^1_{x_1}\left(\phi^1_{w_1}\right) \cup S^2 R^2_{x_2}\left(\phi^2_{w_2}\right) = \phi^1_{x_1 w_1} \cup \phi^2_{x_2 w_2}$$

for all $x \in H_1 \cup H_2 \subset G/$.
Observe that

$$SL_x \circ SR_y = SR_y \circ SL_x$$

$$\text{i.e. } \left(S^1L^1_{x_1} \cup S^2L^2_{x_2}\right) \circ \left(S^1L^1_{y_1} \cup S^2L^2_{y_2}\right)$$

$$\left(S^1L^1_{y_1} \cup S^2L^2_{y_2}\right) \circ \left(S^1L^1_{x_1} \cup S^2L^2_{x_2}\right)$$

for all $x = x_1 \cup x_2$ and $y = y_1 \cup y_2$ in $G = G_1 \cup G_2$.

More generally suppose we have a bihomomorphism from the bigroups $H = H_1 \cup H_2 \subset G = G_1 \cup G_2$ (G a S-neutrosophic bisemigroup) to the bigroup of permutations on a non empty finite biset. $E^1 \cup E^2$. That is suppose for each $x_1$ in $H_1 \subset G_1$ and $x_2$ in $H_2$, $H_2 \subset G_2$, x in $H_1 \cup H_2 \subset G_1 \cup G_2$ we have a bipermutation $\pi^1_{x_1} \cup \pi^1_{x_2}$ on $E_1 \cup E_2$ i.e. one to one mapping of $E_1$ on to $E_1$ and $E_2$ onto $E_2$ such that

$$\pi_x \circ \pi_y = \pi^1_{x_1} \circ \pi^1_{y_1} \cup \pi^2_{x_2} \circ \pi^2_{y_2}\ ,\ \pi e = \pi^1_{e_1} \cup \pi^2_{e_2}$$

is the biidentity bimapping of $E_1 \cup E_2$ and $\pi_{x^{-1}} = \pi^1_{x_1^{-1}} \cup \pi^1_{x_2^{-1}}$ is the inverse mapping of $\pi_x = \pi^1_{x_1} \cup \pi^2_{x_2}$ on $E_1 \cup E_2$. Let $V_H = V^1_{H_1} \cup V^2_{H_2}$ be the bivector space of K-valued bifunctions on $E_1 \cup E_2$.

Then we get the Smarandache birepresentation of $H_1 \cup H_2$ on $V_{H_1} \cup V_{H_2}$ by associating to each $x = x_1 \cup x_2$ in $H_1 \cup H_2$ the linear bimapping

$$\pi_x = \pi^1_{x_1} \cup \pi^2_{x_2}\ :\ V_{H_1} \cup V_{H_2} \to V_{H_1} \cup V_{H_2}$$

defined by

$$\pi_x (f) (a) = f (\pi_x (a))$$

$$\text{i.e. } \left(\pi_{x_1} \cup \pi_{x_2}\right)\left(f^1 \cup f^2\right)\left(a_1 \cup a_2\right)$$

$$= f^1 \left(\pi_{x_1} (a_1)\right) \cup f^2 \left(\pi_{x_2} (a_2)\right)$$

for every $f^1 (a_1) \cup f^2 (a_2) = f (a)$ in $V_{H_1} \cup V_{H_2}$.

This is called the Smarandache bipermutation birepresentation corresponding to the bihomomorphism $x \mapsto \pi x$ i.e. $x_1 \mapsto \pi_{x_1} \cup x_2 \mapsto \pi_{x_2}$ from $H = H_1 \cup H_2$ to permutations on $E = E_1 \cup E_2$.

It is indeed a Smarandache birepresentation for we have several E's and $V_H = V^1_{H_1} \cup V^2_{H_2}$'s depending on the number of proper subsets $H = H_1 \cup H_2$ in $G_1 \cup G_2$ (G the S-bisemigroup) which are bigroups under the operations of $G = G_1 \cup G_2$ because for each $x = x_1 \cup x_2$ and $y = y_1 \cup y_2$ in $H = H_1 \cup H_2$ and each function $f (a) = f_1 (a_1) \cup f_2 (a_2)$ in $V_H = V^1_{H_1} \cup V^2_{H_2}$ we have

$$\left(\pi_x \circ \pi_y\right) (f) (a)$$

$$= \left(\pi^1_{x_1} \cup \pi^2_{x_2}\right) \circ \left(\pi^1_{y_1} \cup \pi^2_{y_2}\right)\left(f_1 \cup f_2\right)\left(a_1 \cup a_2\right)$$

$$= \left(\pi^1_{x_1} \circ \pi^1_{y_1}\right) (f_1) (a_1) \cup \left(\pi^2_{x_2} \circ \pi^2_{y_2}\right)\left(f_2\right)\left(a_2\right)$$

$$= \pi^1_{x_1}\left(\pi^1_{y_1}(f_1)(a_1)\right) \cup \pi^2_{x_2}\left(\pi^2_{y_2}(f_2)(a_2)\right)$$
$$= \pi^1_{y_1}(f_1)\left(\pi^1_{x_1}(I^1(a_1))\right) \cup \pi^2_{y_2}(f_2)\left(\pi^2_{x_2}(I^2(a_2))\right)$$
$$= f_1\left(\pi^1_{y_1}1(\pi^1_{x_1}(1(a_1))\right) \cup f_2\left(\pi^2_{y_2}1(\pi^2_{x_2}(1(a_2))\right)$$
$$= f_1\left(\pi^1_{(x_1y_1)}1(a_1)\right) \cup f_2\left(\pi^2_{(x_2y_2)}1(a_2)\right).$$

Alternatively for each $b = b_1 \cup b_2 \in E_1 \cup E_2$ defined by

$$\psi_b(a) = \psi^1_{b_1}(a_1) \cup \psi^2_{b_2}(a_2)$$

be the function on $E_1 \cup E_2$ defined by $\psi_b(a) = 1$ i.e.,

$$\psi^1_{b_1}(a_1) \cup \psi^2_{b_2}(a_2) = 1 \cup 1.$$

when $a = b$ i.e. $a_1 \cup b_1 = a_2 \cup b_2$. $\psi_b(a) = 0$ when $a \neq b$ i.e. $\psi^1_{b_1}(a_1) \cup \psi^2_{b_2}(a_2) = 0 \cup 0$ when $a_1 \cup b_1 \neq a_2 \cup b_2$. Then the collection of functions $\psi_b$ for $b \in E_1 \cup E_2$ is a basis for $V_H = V^1_{H_1} \cup V^2_{H_2}$ and $\pi_x(\psi) = \psi_{\pi_x(b)}$ $\forall x$ in H and b in E i.e.

$$\pi_{x_1}(\psi^1) \cup \pi_{x_2}(\psi^2) = \psi^1_{\pi_{x_1(b_1)}} \cup \psi^2_{\pi_{x_2(b_2)}}$$

for $x = x_1 \cup x_2$ in $H = H_1 \cup H_2$ and $b_1 \cup b_2$ in $E_1 \cup E_2$. This is true for each proper subset $H = H_1 \cup H_2$ in the S-neutrosophic semigroup $G = G_1 \cup G_2$ and the bigroup $H = H_1 \cup H_2$ associated with the bipermutations of the non empty finite set $E = E_1 \cup E_2$.

Next we shall discuss about Smarandache isomorphic bigroup representation. To this end we consider two bivector spaces $V = V_1 \cup V_2$ and $W = W_1 \cup W_2$ defined over the same field K and that T is a linear biisomorphism from V on to W.

Assume

$$\rho H = \rho^1 H_1 \cup \rho^2 H_2$$

and

$$\rho'_H = \rho'^1 H_1 \cup \rho'^2 H_2$$

are Smarandache birepresentations of the subbigroup $H = H_1 \cup H_2$ in $G = G_1 \cup G_2$ (G a pseudo S-neutrosophic bisemigroup) on V and W respectively.

$$T \text{ o } (\rho H)_x = (\rho' H)_x \circ T$$

for all $x = x_1 \cup x_2 \in H = H_1 \cup H_2$
i.e.

$$(T_1 \cup T_2) \circ \left(\rho^1 H_1 \cup \rho^2 H_2\right)_{x_1 \cup x_2}$$

$$= T_1 \left(\rho^1 H_1\right)_{x_1} \cup T_2 \left(\rho^2 H_2\right)_{x_2}$$

$$= \left(\rho' H_1\right)^1_{x_1} \circ T_1 \cup \left(\rho' H_2\right)^2_{x_2} \circ T_2$$

then we say $T = T_1 \cup T_2$ determines a Smarandache bi-isomorphism between the birepresentation $\rho H$ and $\rho' H$. We may also say that $\rho H$ and $\rho' H$ are Smarandache biisomorphic S-bisemgroup birepresentations.

However it can be verified that Smarandache biisomorphic birepresentation have equal degree but the converse is not true in general.

Suppose V = W be the bivector space of K-valued functions on $H = H_1 \cup H_2 \subset G_1 \cup G_2$ and define T on V = W by
$T (f) (a) = f (a^{-1})$
i.e. $T_1 (f_1) (a_1) \cup T_2 (f_2) (a_2) = f_1 \left(a_1^{-1}\right) \cup f_2 \left(a_2^{-1}\right)$.

This is one to one linear bimapping from the space of K-valued bifunctions $H_1$ on to itself and To $SR_x = SL_x \circ T$
i.e. $\left(T_1 \circ S^1 R^1_{x_1}\right) \cup \left(T_2 \circ S^2 R^2_{x_2}\right) = \left(S^1 L^1_{x_1} \circ T_1\right) \cup \left(S^2 L^2_{x_2} \circ T_2\right)$

for all $x = x_1 \cup x_2$ in $H = H_1 \cup H_2$. For if f (a) is a bifunction on $G = G_1 \cup G_2$ then

$$\begin{aligned}
(T \text{ o } SR_x)(f)(a) &= T(SR_x(f))(a) \\
&= SR_x(f)(a^{-1}) \\
&= f(a^{-1} x) \\
&= T(f)(x^{-1} a) \\
&= SL_x(T(f))(a) \\
&= (SL_x \text{ o } T)(f)(a).
\end{aligned}$$

Therefore we see that S-left and S-right birepresentations of $H = H_1 \cup H_2$ are biisomorphic to each other.

Suppose now that $H = H_1 \cup H_2$ is a subbigroup of the S-bisemigroup G and $\rho H = \rho^1 H_1 \cup \rho^2 H_2$ is a birepresentation of $H = H_1 \cup H_2$ on the bivector space $V_H = V^1_{H_1} \cup V^2_{H_2}$ over the field K and let $\nu_1,\ldots,\nu_n$ be a basis of $V_H = V^1_{H_1} \cup V^2_{H_2}$. For each $x = x_1 \cup x_2$ in $H = H_1 \cup H_2$ we can associate to $(\rho H)_x = (\rho^1 H_1)_{x_1} \cup (\rho^2 H_2)_{x_2}$ an invertible $n \times n$ bimatrix with entries in K using this basis we denote this bimatrix by $(M^1 H_1)_{x_1} \cup (M^2 H_2)_{x_2} = (MH)_x$ where $M = M_1 \cup M_2$.

The composition rule can be rewritten as

$(MH)_{xy} = (MH)_x (MH)_y$

$(M^1 H_1)_{x_1 y_1} \cup (M^2 H_2)_{x_2 y_2}$

$\left[(M^1 H_1)_{x_1} \cup (M^2 H_2)_{x_2}\right]\left[(M^1 H_1)_{y_1} \cup (M^2 H_2)_{y_2}\right]$

$(M^1 H_1)_{x_1} (M^1 H_1)_{y_1} \cup (M^2 H_2)_{x_2} (M^2 H_2)_{y_2}$

where the bimatrix product is used on the right side of the equation. We see depending on each $H = H_1 \cup H_2$ we can have different bimatrices $MH = M^1 H_1 \cup M_2 H_2$, and it need not in general be always a $n \times n$ bimatrices it can also be a $m \times m$ bimatrix $m \neq n$. A different choice of basis for $V = V_1 \cup V_2$ will lead to a different mapping

$$x \mapsto Nx \text{ i.e. } x_1 \cup x_2 \mapsto N^1_{x_1} \cup N^2_{x_2}$$

from H to invertible $n \times n$ bimatrices.

However the two mappings

$$x \mapsto M_x = M^1_{x_1} \cup M^2_{x_2}$$

$$x \mapsto N_x = N^1_{x_1} \cup N^2_{x_2}$$

will be called as Smarandache similar relative to the subbigroup $H = H_1 \cup H_2 \subset G = G_1 \cup G_2$ in the sense that there is an invertible $n \times n$ bimatrix $S = S^1 \cup S^2$ with entries in K such that $N_x = SM_x S^{-1}$ i.e.

$$N^1_{x_1} \cup N^2_{x_2} = S^1 M^1_{x_1} (S^1)^{-1} \cup S^2 M^2_{x_2} (S^2)^{-1}$$

for all $x = x_1 \cup x_2 \subset G = G_1 \cup G_2$. It is pertinent to mention that when a different H' is taken $H \neq H'$ i.e. $H^1 \cup H^2 \neq \left(H^{'}\right)^1 \cup \left(H^{'}\right)^2$ then we may have a different $m \times m$ bimatrix. Thus using a single S-neutrosophic bisemigroup we have very many such bimappings depending on each $H \subset G$. On the other hand one can begin with a bimapping $x \mapsto M_x$ from H into invertible $n \times n$ matrices with entries in K i.e. $x_1 \mapsto M^1_{x_1} \cup x_2 \mapsto M^2_{x_2}$ from $H = H_1 \cup H_2$ into invertible $n \times n$ matrices. Thus now one can reformulate the condition for two Smarandache birepresentations to be biisomorphic.

If one has two birepresentation of a fixed subbigroup $H = H_1 \cup H_2$, H a subbigroup of the S-neutrosophic bisemigroup G on two bivector spaces V and W ($V = V^1 \cup V^2$ and $W = W^1 \cup W^2$) with the same scalar field K then these two Smarandache birepresentations are Smarandache biisomorphic if and only if the associated bimappings from $H = H_1 \cup H_2$ to invertible bimatrices as above, for any choice of basis on $V = V^1 \cup V^2$ and $W = W^1 \cup W^2$ are bisimilar with the bisimilarity bimatrix S having entries in K.

Now we proceed on to give a brief description of Smarandache biirreducible birepresentation, Smarandache biirreducible birepresentation and Smarandache bistable representation and so on. Now we proceed on to define Smarandache bireducibility of finite S-neutrosophic bisemigroups.

Let G be a finite neutrosophic S-bisemigroup when we say G is a S-finite bisemigroup or finite S-bisemigroup we only mean all proper subset in G which are subbigroups in $G = G_1 \cup G_2$ are of finite order $V_H$ be a bivector space over a field K and $\rho H$ a birepresentation of H on $V_H$.

Suppose that there is a bivector space $W_H$ of $V_H$ such that $(\rho H)_x W_H \subseteq W_H$ here $W_H = W^1_{H_1} \cup W^2_{H_2}$ where $H = H_1 \cup H_2$ and $V_H = V^1_{H_1} \cup V^2_{H_2}$, $H = H_1 \cup H_2$, $\rho H = \rho^1 H_1 \cup \rho^2 H_2$ and $x = x_1 \cup x_2 \in H$ i.e. $x_1 \in H_1$ and $x_2 \in H_2$.
This is equivalent to saying that $(\rho H)_x (W_H) = W_H$ i.e.

$$\left[(\rho^1 H_1)_{x_1} \cup (\rho^2 H_2)_{x_2}\right]\left[W^1_{H_1} \cup W^2_{H_2}\right] = W^1_{H_1} \cup W^2_{H_2}$$

for all $x = x_1 \cup x_2 \in H_1 \cup H_2$ as $(\rho H)_{x^{-1}} = \left[(\rho H)_x\right]^{-1}$

$$\text{i.e. } (\rho^1 H_1 \cup \rho^2 H_2)_{(x_1 \cup x_2)^{-1}} = \left[(\rho^1 H_1 \cup \rho^2 H_2)_{(x_1 \cup x_2)}\right]^{-1}$$

$$(\rho^1 H_1)_{x_1^{-1}} \cup (\rho^2 H_2)_{x_2^{-1}} = \left[(\rho^1 H_1)_{x_1}\right]^{-1} \cup \left[(\rho^2 H_2)_{x_2}\right]^{-1}$$

We say $W_H = W^1_{H_1} \cup W^2_{H_2}$ is Smarandache biinvariant or Smarandache bistable under the birepresentation $\rho H = \rho^1 H_1 \cup \rho^2 H_2$.

We say the bisubspace $Z_H = Z^1_{H_1} \cup Z^2_{H_2}$ of $V_H = V^1_{H_1} \cup V^2_{H_2}$ to be a Smarandache bicomplement of a subbispace

$$W_H = W^1_{H_1} \cup W^2_{H_2} \text{ if } W_H \cap Z_H = \{0\}$$

and

$$W_H + Z_H = V_H \text{ i.e. } \left(W^1_{H_1} \cap Z^1_{H_1}\right) \cup \left(W^2_{H_2} \cap Z^2_{H_2}\right) = \{0\} \cup \{0\}$$

and

$$\left(W^1_{H_1} + Z^1_{H_1}\right) \cup \left(W^2_{H_2} + Z^2_{H_2}\right) = V^1_{H_1} + V^2_{H_2}$$

here $W^i_{H_i} + Z^i_{H_i}$, i = 1, 2 denotes the bispan of $W_H$ and $Z_H$ which is a subbispace of $V_H$ consisting of bivectors of the form $w + z = (w_1 + z_1) \cup (w_2 + z_2)$ where $w \in W_H$ and $z \in Z_H$. These conditions are equivalent to saying that every bivector $\nu = \nu_1 \cup \nu_2 \in V^1_{H_1} \cup V^2_{H_2}$ can be written in a unique way as w + z = ($w_1 + z_1$) ∪ ($w_2 + z_2$), $w_i \in W^i_{H_i}$ and $z_i \in Z^i_{H_i}$, i = 1, 2.

Complementary bispaces always exists because of basis for a bivector subspace of a bivector space can be enlarged to a basis of a whole bivector space.

If $Z_H$ = $Z^1_{H_1} \cup Z^2_{H_2}$ and $W_H$ = $W^1_{H_1} \cup W^2_{H_2}$ are complementary subbispaces (bisubspaces) of a bivector space $V_H$ = $V^1_{H_1} \cup V^2_{H_2}$ then we get a linear bimapping $P_H$ = $P^1_{H_1} \cup P^2_{H_2}$ on $V_H$ = $V^1_{H_1} \cup V^2_{H_2}$ on to $W_H$ = $W^1_{H_1} \cup W^2_{H_2}$ along $Z_H$ = $Z^1_{H_1} \cup Z^2_{H_2}$ and is defined by

$P_H$ (w + z) w for all w ∈ $W_H$ and z ∈ $Z_H$. Thus $I_H - P_H$ is the biprojection of $V_H$ on to $Z_H$ along $W_H$ where $I_H$ denotes the identity bitransformation on $V_H$ = $V^1_{H_1} \cup V^2_{H_2}$.

**Note:**

$$
\begin{aligned}
(P_H)^2 &= \left(P^1_{H_1} \cup P^2_{H_2}\right)^2 \\
&= \left(P^1_{H_1}\right)^2 \cup \left(P^2_{H_2}\right)^2 \\
&= P^1_{H_1} \cup P^2_{H_2}
\end{aligned}
$$

when $P_H$ is a biprojection.

Conversely if $P_H$ is a linear bioperator on $V_H$ such that $(P_H)^2 = P_H$ then $P_H$ is the biprojection of $V_H$ on to the bisubspace of $V_H$ which is the biimage of $P_H$ = $P^1_{H_1} \cup P^2_{H_2}$ along the subspace of $V_H$ which is the bikernel of $\rho H = \rho^1 H_1 \cup \rho^2 H_2$.

It is important to mention here unlike usual complements using a finite bigroup we see when we used pseudo S-neutrosophic bisemigroups. The situation is very varied. For each proper subset H of G ($H_1 \cup H_2 \subset G_1 \cup G_2$) where H is a subbigroup of G we get several important S-bicomplements and several S-biinvariant or S-bistable or S-birepresentative of $\rho H = \rho^1 H_1 \cup \rho^2 H_2$.

Now we proceed on to define Smarandache biirreducible birepresentation. Let G be a S-finite neutrosophic bisemigroup, $V_H = V^1_{H_1} \cup V^2_{H_2}$ be a bivector space over a field K, $\rho H = \rho^1 H_1 \cup \rho^2 H_2$ be a birepresentation of H on $V_H$ and $W_H$ is a subbispace of $V_H = V^1_{H_1} \cup V^2_{H_2}$ which is invariant under $\rho H = \rho^1 H_1 \cup \rho^2 H_2$. Here we make an assumption that the field K has characteristic 0 or K has positive characteristic and the number of elements in each H = $H^1 \cup H^2$ is not divisible by the characteristic K, $H_1 \cup H_2 \subset G_1 \cup G_2$ is a S-bisemigroup.

Let us show that there is a bisubspace $Z_H = Z^1_{H_1} \cup Z^2_{H_2}$ of $V^1_{H_1} \cup V^2_{H_2} = V_H$ such that $Z_H$ is a bicomplement of $W_H = W^1_{H_1} \cup W^2_{H_2}$ and $Z_H$ is also biinvariant under the birepresentation $\rho H$ of H i.e. $\rho^1 H_1 \cup \rho^2 H_2$ of $H_1 \cup H_2$ on $V_H = V^1_{H_1} \cup V^2_{H_2}$. To do this we start with any bicomplements $(Z_H)_o = \left(Z^1_{H_1}\right)_o \cup \left(Z^2_{H_2}\right)_o$ of $W_H = W^1_{H_1} \cup W^2_{H_2}$ of $V_H = V^1_{H_1} \cup V^2_{H_2}$ and let $(P_H)_o = \left(P^1_{H_1} \cup P^2_{H_2}\right)_o$ : $V_H = V^1_{H_1} \cup V^2_{H_2} \to V^1_{H_1} \cup V^2_{H_2}$ be the biprojection of $V_H = V^1_{H_1} \cup V^2_{H_2}$ on to $W^1_{H_1} \cup W^2_{H_2} = W_H$ along $(Z_H)_o$. Thus $(P_H)_o = \left(P^1_{H_1} \cup P^2_{H_2}\right)_o$ maps V to W and $(P_H)_o$ ow = w for all w $\in$ W.

Let $m = m_1 \cup m_2$ denote the number of elements in $H = H_1 \cup H_2 \subset G_1 \cup G_2$ i.e. $|H_i| = m_i$, i = 1, 2. Define a linear bimapping

$$P_H : V_H \to V_H$$

i.e. $P^1_{H_1} \cup P^2_{H_2} : V^1_{H_1} \cup V^2_{H_2} \to V^1_{H_1} \cup V^2_{H_2}$

by $P_H = P^1_{H_1} \cup P^2_{H_2}$

$$= \frac{1}{m_1} \sum_{x_1 \in H_1} \left(\rho^1 H_1\right)_{x_1} \circ \left(P^1_{H_1}\right) \circ \left(\rho^1 H_1\right)^{-1}_{x_1}$$

$$\cup \frac{1}{m_2} \sum_{x_2 \in H_2} \left(\rho^2 H_2\right)_{x_2} \circ \left(P^2_{H_2}\right) \circ \left(\rho^2 H_2\right)^{-1}_{x_2}$$

assumption on K implies that $\frac{1}{m_i}$ (i = 1, 2) makes sense as an element of K i.e. as the multiplicative inverse of a sum of m 1's in K where 1 refers to the multiplicative identity element of K. This expression defines a linear bimapping on $V_H = V^1_{H_1} \cup V^2_{H_2}$ because $\left(\rho H\right)_x$'s and $(P_H)_o$ are linear bimapping. We actually have that $P_H = P^1_{H_1} \cup P^2_{H_2}$ bimaps $V_H$ to $W_H$ i.e. $V^1_{H_1} \cup V^2_{H_2}$ to $W^1_{H_1} \cup W^2_{H_2}$ and because the

$$\left(P_H\right)_o = \left(P^1_{H_1} \cup P^2_{H_2}\right)_o$$

maps $V_H = V^1_{H_1} \cup V^2_{H_2}$ to $W_H = W^1_{H_1} \cup W^2_{H_2}$, and because the

$$\left(\rho_H\right)_x\text{'s} \left( = \left(\rho^1_{H_1}\right)_{x_1} \cup \left(\rho^2_{H_2}\right)_{x_2} \right)$$

maps $W_H = W^1_{H_1} \cup W^2_{H_2}$ to $W^1_{H_1} \cup W^2_{H_2}$ .

If $w \in W_H$ then

$$\left[(\rho H)_x\right]^{-1} w = \left[\left(\rho^1 H_1\right)_{x_1} \cup \left(\rho^2 H_2\right)_{x_2}\right]^{-1} (w_1 \cup w_2)$$

$$= \left(\rho^1 H_1\right)^{-1}_{x_1} (w_1) \cup \left(\rho^2 H_2\right)_{x_2} (w_2) \in W^1_{H_1} \cup W^2_{H_2}$$

for all $x = x_1 \cup x_2$ in $H = H_1 \cup H_2 \subset G = G_1 \cup G_2$ and then

$$\left(P_H\right)_o \left(\left(\rho H\right)_x\right)^{-1} \omega =$$

$$= \left(P_H\right)_o \left(\left(\rho^1 H_1\right)_{x_1}\right)^{-1} (w_1) \cup \left(P^2_{H_2}\right)_o \left(\left(\rho^2 H_2\right)_{x_2}\right)^{-1} (w_2)$$

$$= \left(\left(\rho^1 H_1\right)_{x_1}\right)^{-1} (w_1) \cup \left(\left(\rho^2 H_2\right)_{x_2}\right)^{-1} (w_2).$$

Thus we conclude that $(P_H)$ (w) = w i.e. $\left(P^1_{H_1}\right)(w_1) = w_1$ and $\left(P^2_{H_2}\right)(w_2) = w_2$ i.e. $P_H = P^1_{H_1} \cup P^2_{H_2}$ for all w = $(w_1 \cup w_2)$ in $W_H$ = $W^1_{H_1} \cup W^2_{H_2}$ by the very definition of $P_H$.
The definition of $P_H$ also implies that

$$(\rho H)_y \circ P_H \circ \left[(\rho H)_y\right]^{-1} = P_H$$

i.e.

$$\left(\rho^1 H_1\right)_{y_1} \circ P^1_{H_1} \circ \left(\left(\rho^1 H_1\right)_{y_1}\right)^{-1} \cup \left(\rho^2 H_2\right)_{y_2} \circ P^1_{H_2} \circ \left(\left(\rho^2 H_2\right)_{y_2}\right)^{-1}$$
$$= P^1_{H_1} \cup P^2_{H_2}$$

for all $y \in H = H_1 \cup H_2$.
The only case this does not occur is when $W_H$ = {0} i.e. $W^1_{H_1} \cup W^2_{H_2}$ = {0} ∪ {0}. Because $P_H$ $(V_H) \subset W_H$ and $P_H$ (w) = w for all $w \in W_H$ = $W^1_{H_1} \cup W^2_{H_2}$. $P_H$ = $P^1_{H_1} \cup P^2_{H_2}$ is a biprojection of $V_H$ onto $W_H$ i.e. $P^i_{H_i}$ is a projection of $V^i_{H_i}$ onto $W^i_{H_i}$, i = 1, 2 along some bisubspace $Z_H$ = $Z^1_{H_1} \cup Z^2_{H_2}$ of $V_H$ = $V^1_{H_1} \cup V^2_{H_2}$. Specifically one should take $Z_H$ = $Z^1_{H_1} \cup Z^2_{H_2}$ to be the bikernel of $P_H$ = $P^1_{H_1} \cup P^2_{H_2}$ . It is easy to see that $W_H \cap Z_H$ = {0} i.e. $W^1_{H_1} \cap Z^1_{H_1} = \{0\}$ and $W^2_{H_2} \cap Z^2_{H_2} = \{0\}$ since $P^i_{H_i}(w_i) = w_i$ for all $w_i \in W^i_{H_i}$, i = 1, 2.

On the other hand if $\nu = \nu_1 \cup \nu_2$ is any element of $V_H$ = $V^1_{H_1} \cup V^2_{H_2}$ then we can write $\nu = \nu_1 \cup \nu_2$ as

$$P_H(\nu) = P^1_{H_1}(\nu_1) \cup P^2_{H_2}(\nu_2) \text{ so } P_H(\nu) + (V - P_H(\nu)).$$

Thus $\nu - P_H(\nu)$ lies in $Z_H$, the bikernel of $P_H$. This shows that $W_H$ and $Z_H$ satisfies the essential bicomplement of $W_H$ in $V_H$. The biinvariance of $Z_H$ under the birepresentation ρH is evident.

Thus the Smarandache birepresentation ρH of H on $V_H$ is biisomorphic to the direct sum of H on $W_H$ and $Z_H$, that are the birestrictions of ρH to $W_H$ and $Z_H$.

There can be smaller biinvariant bisubspaces within these biinvariant subbispaces so that one can repeat the process for each H, H ⊂ G. We say that the subbispaces $(W_H)_1$ , $(W_H)_2$, …, $(W_H)_t$ of $V_H$ i.e.

$$\left(W^1_{H_1} \cup W^2_{H_2}\right)_1, \left(W^1_{H_1} \cup W^2_{H_2}\right)_2, \ldots, \left(W^1_{H_1} \cup W^2_{H_2}\right)_t$$

of $V^1_{H_1} \cup V^2_{H_2}$ form an Smarandache biindependent system related to each subbigroup $H = H_1 \cup H_2 \subset G = G_1 \cup G_2$. If $(W_H)_j \neq (0)$ for each j and if $w_j \in (W_H)_j$, $1 \le j \le t$ and

$$\sum_{j=1}^{t} w_j = \sum_{j=1}^{t} w^1_j \cup \sum_{j=1}^{t} w^2_j = 0 \cup 0$$

where $w_j = w^1_j \cup w^2_j$, $w^1_j \in W^1_{H_1}$ and $w^2_j \in W^2_{H_2}$ imply $w^i_j$ = 0, i = 1, 2; j = 1, 2, …, t. If in addition it spans $(W_H)_1, (W_H)_2, \ldots, (W_H)_t$ = $V^1_{H_1} \cup V^2_{H_2}$ = $V_H$ then every bivector $\nu = \nu^1 \cup \nu^2$ on $V^1_{H_1} \cup V^2_{H_2}$ can be written in a unique way as $\sum_{j=1}^{t} u_j$ with $u_j = u^1_j \cup u^2_j \in$ ( $W^1_{H_1} \cup W^2_{H_2}$ ) for each j.

Next we proceed on to give two applications to Smarandache Markov bichains and Smarandache Leontief economic bimodels.

### 4.4 Smarandache Markov bichains using S-neutrosophic bivector spaces

Suppose a physical or a mathematical system is such that at any movement it can occupy one of a finite number of states when we view them as stochastic bioprocess or Markov bichains we make an assumption that the system moves with time from one state to another so that a schedule of observation times keep the states of the system at these times. But when we tackle real world problems say even for simplicity the emotions of a persons it need not fall under the category of sad, cold, happy, angry, affectionate, disinterested, disgusting, many times the emotions of a person may be very unpredictable depending largely on the situation, and the mood of the person and its relation with another, so such study cannot fall under Markov chains, for at a time more than one emotion may be in a person and also such states cannot be included and given as next pair of observation, these changes and several feelings atleast two at a time will largely affect the very transition bimatrix

$$P = P_1 \cup P_2 = \left[p_{ij}^{1}\right] \cup \left[p_{ij}^{2}\right]$$

with non negative entries for which each of the column sums are one and all of whose entries are positive. This has relevance as even the policy makers are humans and their view is ultimate and this rules the situation. Here it is still pertinent to note that all decisions are not always possible at times certain of the views may be indeterminate at that period of time and may change after a period of time but all our present theory have no place for the indeterminacy only the neutrosophy gives the place for the concept of indeterminacy, based on which we have built neutrosophic vector spaces, neutrosophic bivector spaces, then now the notion of Smarandache -neutrosophic bivector spaces and so on.

So to over come the problem we have indecisive situations we give negative values and indeterminate situations we give negative values so that our transition neutrosophic bimatrices individual columns sums do not add to one and all entries may not be positive.

Thus we call the new transition neutrosophic bimatrix which is a square bimatrix which can have negative entries and I the indeterminate also falling in the set $[-1\ 1] \cup \{I\}$ and whose column sums can also be less than 1 and I as the Smarandache neutrosophic transition bimatrix.

Further the Smarandache neutrosophic probability bivector will be a bicolumn vector which can take entries from [–1, 1] ∪ [–I, I] whose sum can lie in the biinterval [–1, 1] ∪ [–I I]. The Smarandache neutrosophic probability bivectors $x^{(n)}$ for n = 0, 1, 2, … are said to be the Smarandache state neutrosophic bivectors of a Smarandache neutrosophic Markov bioprocess. Clearly if P is a S-transition bimatrix of a Smarandache Markov bioprocess and $x^{(n)} = x_1^{(n_1)} \cup x_2^{(n_2)}$ is the Smarandache state neutrosophic bivectors at the $n^{th}$ pair of observation then

$$x^{(n+1)} \neq px^{(n)}$$

$$\text{i.e. } x_1^{(n+1)} \cup x_2^{(n_2+1)} \neq p_1 x_1^{(n_1)} \cup p_2 x_2^{(n_2)}.$$

Further research in this direction is innovative and interesting.

## 4.5 Smarandache neutrosophic Leontief economic bimodels

Matrix theory has been very successful in describing the inter relation between prices outputs and demands in an economic model. Here we just discuss some simple bimodels based on the ideals of the Nobel laureate Massily Leontief. We have used not only bimodel structure based on bimatrices also we have used the factor indeterminacy. So our matrices would be only Neutrosophic bimatrices. Two types of models which we wish to discuss are the closed or input-output model and the open or production model each of which assumes some economics parameter which describe the inter relations between the industries in the economy under considerations. Using neutrosophic bimatrix theory we can combine and study the effect of price bivector. Before the basic equations of the input-output model are built we just recall the definition of fuzzy neutrosophic bimatrix. For we need this type of matrix in our bimodel.

**DEFINITION 4.5.1:** *Let $M_{nxm} = \{(a_{ij}) \,/\, a_{ij} \in K(I)\}$, where K(I), is a neutrosophic field. We call $M_{nxm}$ to be the neutrosophic rectangular matrix.*

***Example 4.5.1:*** Let Q($I$ ) = ⟨Q ∪ $I$ ⟩ be the neutrosophic field.

$$M_{4\times 3} = \begin{pmatrix} 0 & 1 & I \\ -2 & 4I & 0 \\ 1 & -I & 2 \\ 3I & 1 & 0 \end{pmatrix}$$

is the neutrosophic matrix, with entries from rationals and the indeterminacy $I$.

We define product of two neutrosophic matrices and the product is defined as follows:

Let

$$A = \begin{pmatrix} -1 & 2 & -I \\ 3 & I & 0 \end{pmatrix}_{2\times 3}$$

and

$$B = \begin{pmatrix} I & 1 & 2 & 4 \\ 1 & I & 0 & 2 \\ 5 & -2 & 3I & -I \end{pmatrix}_{3x4}$$

$$AB = \begin{bmatrix} -6I+2 & -1+4I & -2-3I & I \\ -4I & 3+I & 6 & 12+2I \end{bmatrix}_{2x4}.$$

(we use the fact $I^2 = I$).

Let $M_{n\times n}$ = {($a_{ij}$) | ($a_{ij}$) ∈ Q($I$ )}, $M_{n\times n}$ is a neutrosophic vector space over Q and a strong neutrosophic vector space over Q($I$ ).

Now we proceed onto define the notion of fuzzy integral neutrosophic matrices and operations on them, for more about these refer [43].

**DEFINITION 4.5.2:** *Let N = [0, 1] ∪ I where I is the indeterminacy. The m × n matrices $M_{m \times n}$ = {($a_{ij}$) / $a_{ij}$ ∈* [0, 1] *∪ I} is called the fuzzy integral neutrosophic matrices. Clearly the class of m × n matrices is contained in the class of fuzzy integral neutrosophic matrices.*

***Example 4.5.2:*** Let

$$A = \begin{pmatrix} I & 0.1 & 0 \\ 0.9 & 1 & I \end{pmatrix},$$

A is a 2 × 3 integral fuzzy neutrosophic matrix.

We define operation on these matrices. An integral fuzzy neutrosophic row vector is 1 × n integral fuzzy neutrosophic matrix. Similarly an integral fuzzy neutrosophic column vector is a m × 1 integral fuzzy neutrosophic matrix.

***Example 4.5.3:*** A = (0.1, 0.3, 1, 0, 0, 0.7, *I*, 0.002, 0.01, *I*, 0.12) is a integral row vector or a 1 × 11, integral fuzzy neutrosophic matrix.

***Example 4.5.4:*** $B = (1, 0.2, 0.111, I, 0.32, 0.001, I, 0, 1)^T$ is an integral neutrosophic column vector or B is a 9 × 1 integral fuzzy neutrosophic matrix.

We would be using the concept of fuzzy neutrosophic column or row vector in our study.

**DEFINITION 4.5.3:** *Let P = ($p_{ij}$) be a m × n integral fuzzy neutrosophic matrix and Q = ($q_{ij}$) be a n × p integral fuzzy neutrosophic matrix. The composition map P • Q is defined by R = ($r_{ij}$) which is a m × p matrix where $r_{ij}$ = $\max_k \min(p_{ik} q_{kj})$ with the assumption max($p_{ij}$, I) = I and min($p_{ij}$, I) = I where $p_{ij}$ ∈* [0, 1]. *min (0, I) = 0 and max(1, I) = 1.*

***Example 4.5.5:*** Let

$$P = \begin{bmatrix} 0.3 & I & 1 \\ 0 & 0.9 & 0.2 \\ 0.7 & 0 & 0.4 \end{bmatrix}, \ Q = (0.1, I, 0)^T$$

be two integral fuzzy neutrosophic matrices.

$$P \bullet Q = \begin{bmatrix} 0.3 & I & 1 \\ 0 & 0.9 & 0.2 \\ 0.7 & 0 & 0.4 \end{bmatrix} \bullet \begin{bmatrix} 0.1 \\ I \\ 0 \end{bmatrix} = (I, I, 0.1).$$

***Example 4.5.6:*** Let

$$P = \begin{bmatrix} 0 & I \\ 0.3 & 1 \\ 0.8 & 0.4 \end{bmatrix}$$

and

$$Q = \begin{bmatrix} 0.1 & 0.2 & 1 & 0 & I \\ 0 & 0.9 & 0.2 & 1 & 0 \end{bmatrix}.$$

One can define the max-min operation for any pair of integral fuzzy neutrosophic matrices with compatible operation.

Now we proceed onto define the notion of fuzzy neutrosophic matrices. Let $N_s = [0, 1] \cup nI / n \in (0, 1]\}$; we call the set $N_s$ to be the fuzzy neutrosophic set.

**DEFINITION 4.5.4:** *Let $N_s$ be the fuzzy neutrosophic set.*

$$M_{n \times m} = \{(a_{ij}) / a_{ij} \in N_s\}$$

*we call the matrices with entries from $N_s$ to be the fuzzy neutrosophic matrices.*

***Example 4.5.7:*** Let $N_s = [0,1] \cup \{nI/ n \in (0,1]\}$ be the set

$$P = \begin{bmatrix} 0 & 0.2I & 0.31 & I \\ I & 0.01 & 0.7I & 0 \\ 0.31I & 0.53I & 1 & 0.1 \end{bmatrix}$$

P is a 3 × 4 fuzzy neutrosophic matrix.

***Example 4.5.8:*** Let $N_s = [0, 1] \cup \{nI / n \in (0, 1]\}$ be the fuzzy neutrosophic matrix. A = [0, 0.12*I*, *I*, 1, 0.31] is the fuzzy neutrosophic row vector:

$$B = \begin{bmatrix} 0.5I \\ 0.11 \\ I \\ 0 \\ -1 \end{bmatrix}$$

is the fuzzy neutrosophic column vector.

Now we proceed on to define operations on these fuzzy neutrosophic matrices.

Let $M = (m_{ij})$ and $N = (n_{ij})$ be two m × n and n × p fuzzy neutrosophic matrices. $M \bullet N = R = (r_{ij})$ where the entries in the fuzzy neutrosophic matrices are fuzzy indeterminates i.e. the indeterminates have degrees from 0 to 1 i.e. even if some factor is an indeterminate we try to give it a degree to which it is indeterminate for instance 0.9*I* denotes the indeterminacy rate; it is high where as 0.01*I* denotes the low indeterminacy rate. Thus neutrosophic matrices have only the notion of degrees of indeterminacy. Any other type of operations can be defined on the neutrosophic matrices and fuzzy neutrosophic matrices. The notion of these matrices have been used to define neutrosophic relational equations and fuzzy neutrosophic relational equations.

Here we give define the notion of neutrosophic bimatrix and illustrate them with examples. Also we define fuzzy neutrosophic matrices.

**DEFINITION 4.5.5:** *Let $A = A_1 \cup A_2$ where $A_1$ and $A_2$ are two distinct neutrosophic matrices with entries from a neutrosophic field. Then $A = A_1 \cup A_2$ is called the neutrosophic bimatrix.*

*It is important to note the following:*

*(1) If both $A_1$and $A_2$ are neutrosophic matrices we call A a neutrosophic bimatrix.*

*(2) If only one of $A_1$ or $A_2$ is a neutrosophic matrix and other is not a neutrosophic matrix then we all $A = A_1 \cup A_2$ as the semi neutrosophic bimatrix. (It is clear all neutrosophic bimatrices are trivially semi neutrosophic bimatrices).*

*It both $A_1$ and $A_2$ are $m \times n$ neutrosophic matrices then we call $A = A_1 \cup A_2$ a $m \times n$ neutrosophic bimatrix or a rectangular neutrosophic bimatrix.*

*If $A = A_1 \cup A_2$ be such that $A_1$ and $A_2$ are both $n \times n$ neutrosophic matrices then we call $A = A_1 \cup A_2$ a square or a $n \times n$ neutrosophic bimatrix. If in the neutrosophic bimatrix $A = A_1 \cup A_2$ both $A_1$ and $A_2$ are square matrices but of different order say $A_1$ is a $n \times n$ matrix and $A_2$ a $s \times s$ matrix then we call $A = A_1 \cup A_2$ a mixed neutrosophic square bimatrix. (Similarly one can define mixed square semi neutrosophic bimatrix).*

*Likewise in $A = A_1 \cup A_2$ if both $A_1$ and $A_2$ are rectangular matrices say $A_1$ is a $m \times n$ matrix and $A_2$ is a $p \times q$ matrix then we call $A = A_1 \cup A_2$ a mixed neutrosophic rectangular bimatrix. (If $A = A_1 \cup A_2$ is a semi neutrosophic bimatrix then we call A the mixed rectangular semi neutrosophic bimatrix).*

Just for the sake of clarity we give some illustration.

**Notation:** We denote a neutrosophic bimatrix by $A_N = A_1 \cup A_2$.

***Example 4.5.9:*** Let

$$A_N = \begin{bmatrix} 0 & I & 0 \\ 1 & 2 & -1 \\ 3 & 2 & I \end{bmatrix} \cup \begin{bmatrix} 2 & I & 1 \\ I & 0 & I \\ 1 & 1 & 2 \end{bmatrix}$$

$A_N$ is the $3 \times 3$ square neutrosophic bimatrix.

At times one may be interested to study the problem at hand (i.e. the present situation) and a situation at the $r^{th}$ time period the predicted model.

All notion and concept at all times is not determinable. For at time a situation may exist for a industry that it cannot say the monetary value of the output of the $i^{th}$ industry needed to satisfy the outside demand at one time, this notion may become an indeterminate (For instance with the advent of globalization the electronic goods manufacturing industries are facing a problem for in the Indian serenio when an exported goods is sold at a cheaper rate than manufactured Indian goods will not be sold for every one will prefer only an exported good, so in situation like this the industry faces only a indeterminacy for it cannot fully say anything about the movements of the manufactured goods in turn this will affect the $\sigma_{ij}$. $\sigma_{ij}$ may also tend to become an indeterminate. So to study such situation simultaneously the neutrosophic bimatrix would be ideal we may have the newly redefined production vector which we choose to call as Smarandache neutrosophic production bivector which will have its values taken from +ve value or –ve value or an indeterminacy.

So Smarandache neutrosophic Leontief open model is got by permitting.

$$x \geq 0, d \geq 0, c \geq 0$$
$$x \leq 0, d \leq 0, c \leq 0$$

and

x can be I

d can take any value and c can be a neutrosophic bimatrix. We can say $(1 - c)^{-1} \geq 0$ productive $(1 - c)^{-1} < 0$ non productive or not up to satisfaction and $(1 - c^{-1}) = nI$, I the indeterminacy i.e. the productivity cannot be determined i.e. one cannot say productive or non productive but cannot be determined. $c = c_1 \cup c_2$ is the consumption neutrosophic bimatrix.

$c_1$ at time of study and $c_2$ after a stipulated time period.

x, d, c can be greater than or equal to zero less than zero or can be an indeterminate.

$$x = \begin{bmatrix} x_1^1 \\ \vdots \\ x_k^1 \end{bmatrix} \cup \begin{bmatrix} x_1^2 \\ \vdots \\ x_k^2 \end{bmatrix}$$

production neutrosophic bivector at the times $t_1$ and $t_2$ the demand neutrosophic bivector $d = d^1 \cup d^2$

$$d = \begin{bmatrix} d_1^1 \\ \vdots \\ d_k^1 \end{bmatrix} \cup \begin{bmatrix} d_1^2 \\ \vdots \\ d_k^2 \end{bmatrix}$$

at time $t_1$ and $t_2$ respectively. Consumption neutrosophic bimatrix

$c = c_1 \cup c_2$

$$c_1 = \begin{bmatrix} \sigma_{11}^1 & \cdots & \sigma_{1k}^1 \\ \sigma_{21}^1 & \cdots & \sigma_{2k}^1 \\ \vdots & & \\ \sigma_{k1}^1 & \cdots & \sigma_{kk}^1 \end{bmatrix}$$

$$c_2 = \begin{bmatrix} \sigma_{11}^2 & \cdots & \sigma_{1k}^2 \\ \sigma_{21}^2 & \cdots & \sigma_{2k}^2 \\ \vdots & & \\ \sigma_{k1}^2 & \cdots & \sigma_{kk}^2 \end{bmatrix}$$

at times $t_1$ and $t_2$ respectively.

$$\sigma_{i1} x_1 + \sigma_{12} x_2 + \ldots + \sigma_{ik} x_k$$

$$= \left(\sigma_{i1}^1 x_1^1 + \sigma_{i2}^1 x_2^1 + \ldots + \sigma_{ik}^1 x_k^1\right) \cup \left(\sigma_{i1}^2 x_1^2 + \sigma_{i2}^2 x_2^2 + \ldots + \sigma_{ik}^2 x_k^2\right)$$

is the value of the output of the $i^{th}$ industry needed by all k industries at the time periods $t_1$ and $t_2$ to produce a total output specified by the production neutrosophic bivector x = $x^1 \cup x^2$. Consumption neutrosophic bimatrix c is such that; production if $(1 – c)^{-1}$ exists and

$$(1 - c)^{-1} \geq 0$$

i.e. c = $c_1 \cup c_2$ and $(1 - c_1)^{-1} \cup (1 - c_2)^{-1}$ exists and each of $(1 - c_1)^{-1}$ and $(1 - c_2)^{-1}$ is greater than or equal to zero. A consumption neutrosophic bimatrix c is productive if and only if there is some production bivector $x \geq 0$ such that

$$x > cx$$

$$\text{i.e. } x^1 \cup x^2 > c^1 x^1 \cup c^2 x^2.$$

A consumption bimatrix c is productive if each of its birow sum is less than one. A consumption bimatrix c is productive if each of its bicolumn sums is less the one. Non productive if bivector $x < 0$ such that $x < c\,x$.

Now quasi productive if one of $x^1 \geq 0$ and $x^1 > c^1 x^1$ or $x^2 \geq 0$ and $x^1 > c^1 x^1$.

Now production is indeterminate if x is indeterminate x and cx are indeterminates or x is indeterminate and c x is determinate. Production is quasi indeterminate if at $t_1$ or $t_2$, $x^i \geq 0$ and $x^i > c^i x^i$ are indeterminates quasi non productive and indeterminate if one of $x^i < 0$, $c^i x^i < 0$ and one of $x^i$ and $I^i x^i$ are indeterminate. Quasi production if one of $c^i x^i > 0$ and $x^i > 0$ and $x^i < 0$ and $I^i x^i < 0$. Thus 6 possibilities can occur at anytime of study say at times $t_1$ and $t_2$ for it is but very natural as in any industrial problem the occurrences of any factor like demand or production is very much dependent on the people and the government policy and other external factors.

Chapter Five

# SUGGESTED PROBLEMS

In this chapter we introduce problems, which would certainly help any reader to become more familiar with this new algebraic structure. Most of the problems are easy, only a few of them are difficult. Certainly these problems will throw more light on this new subject.

1. Prove any interesting condition on a linear bioperator T to be bidiagonalizable.

2. Under what conditions on a, b, c and $a^1$, $b^1$, $c^1$, $d^1$ the square bimatrix

$$A = \begin{bmatrix} 0 & 0 & 0 & 0 \\ a & 0 & 0 & 0 \\ 0 & b & 0 & 0 \\ 0 & 0 & c & 0 \end{bmatrix} \cup \begin{bmatrix} a^1 & 0 & 0 & 0 \\ 0 & b^1 & 0 & 0 \\ 0 & 0 & 0 & c^1 \\ 0 & 0 & d^1 & 0 \end{bmatrix}$$

is bidiagonalizable?

3. Prove if T is a linear bioperator on $V = V_1 \cup V_2$ which is a (m, m) dimensional bivector space and suppose $T = T_1 \cup T_2$ has distinct characteristic bivalues; prove T is bidiagonalizable.

4. Define the notion of annihilating bipolynomials in case of a linear bioperators.

5. Find the annihilating bipolynomial for the following bimatrix $A = A_1 \cup A_2$

$$A = \begin{bmatrix} 0 & 1 & 0 & 1 \\ 1 & 0 & 1 & 0 \\ 0 & 1 & 0 & 1 \\ 1 & 0 & 1 & 0 \end{bmatrix} \cup \begin{bmatrix} 1 & 1 & 0 & 0 \\ -1 & -1 & 0 & 0 \\ -2 & -2 & 2 & 1 \\ 1 & 1 & -1 & 0 \end{bmatrix}.$$

6. Is the characteristic bipolynomial of $A = A_1 \cup A_2$ given in problem 5 same as the annihilating bipolynomial?

7. Define minimal bipolynomial for any linear bioperator T.

8. Find the minimal bipolynomial for the bimatrix $A = A_1 \cup A_2$ given in problem 5.

9. Find a $3 \times 3$ bimatrix and the related linear bioperator whose minimal bipolynomial is $\{x^2\} \cup \{x^3\}$.

10. Define for a bivector space $V = V_1 \cup V_2$ a subbispace W to be biinvariant where $T = T_1 \cup T_2$. (T a linear bioperator on V).

11. Using the definition of a linear bioperator, $T = T_1 \cup T_2$ prove T is bitriangulable if and only if both $T_1$ and $T_2$ are triangulable:

"In a finite dimensional bivector space $V = V_1 \cup V_2$ over the field F and $T = T_1 \cup T_2$ a linear bioperator on V, prove $T = T_1 \cup T_2$ is bitriangulable if and only if minimal bipolynomial for T is the product of linear bipolynomials over F."

*Note:* It may happen that only one of $T_1$ or $T_2$ may be triangulable in $T = T_1 \cup T_2$ in which case we say the linear bioperator T is semi bitriangulable.

12. Give an example of a linear bioperator $T = T_1 \cup T_2$, which is only semibitriangulable.

13. Prove every bimatrix $A = A_1 \cup A_2$ such that in which $A^2 = A_1^2 \cup A_2^2$ where $A^2 = A$ i.e. $A_1^2 = A_1$ and $A_2^2 = A_2$ is similar to a bidiagonal bimatrix.

14. If V and W be bivector spaces over the same field F and T be a linear bitransformation from V into W. Suppose V is finite dimensional prove; Rank T + Nullity T = dim V i.e. rank $(T_1 \cup T_2)$ + nullity $(T_1 \cup T_2)$ = dim $V_1 \cup$ dim $V_2$ i.e. (rank $T_1$ + nullity $T_1$) $\cup$ (rank $T_2$ + millity $T_2$) = dim $V_1 \cup$ dim $V_2$.

15. Describe the range bispace (rank) and null bispace (nullity) for the linear bitransformation. $T = T_1 \cup T_2 : V \to W$ ($V = V_1 \cup V_2$ and $W = W_1 \cup W_2$), given by $T (x^1 \cup x^2) = T_1 (x_1, y_1, z_1) \cup T (x_2, y_2, z_2)$ , $w_2 = (x_1 - z_1, y_1 + z_1) \cup (x_2 + y_2, z_2, z_2 - w_2, w_2 + z_2, x_2 + w_2 + y_2)$.

16. For the above linear bitransformation T, given in problem 15 is it possible to find a U: $W \to V$

such that $T \circ U$ is identity U also a linear bitransformation from W to V where $U = U_1 \cup U_2$?

17. Let $T: V \to V$ be a linear bioperator defined by (when $T = T_1 \cup T_2 : V_1 \cup V_2 \to V_1 \cup V_2$) $T_1 (x_1, x_2, x_3) \cup T_2 (a, b, c, d) = (3x_1, x_1 - x_2, 2x_1 + x_2 + x_3) \cup (a + b + c + d, b + c + d , c + d, d)$

   a. Find $(T^2 - I) (T - 2I)$ i.e. $\left(T_1^2 \cup T_2^2 - I^1 \cup I^2\right)\left(T_1 \cup T_2 - 2(I^1 \cup I^2)\right) = \left(T_1^2 - I^1\right)\left(T_1 - 2I^1\right) \cup (T_2^2 - I^1)\left(T_2 - 2I^2\right)$.
   b. Is T non singular?
   c. Does the eigen bivalues exist over Q?
   d. Does the eigen bivalues exist over R?
   e. If eigen bivalues exists find the eigen bivectors?
   f. Is T a bidiagonalizable bioperator?

18. Suppose a linear bitransformation $f = f_1 \cup f_2$ from $V = V_1 \cup V_2$ into the scalar field $F \cup F$ is defined such that $f (c\,\alpha + \beta) = c\, f (\alpha) + f (\beta)$; where

$$c = c_1 \cup c_2$$
$$\alpha = \alpha_1 \cup \alpha_2$$
$$\beta = \beta_1 \cup \beta_2.$$

Suppose $f (c\,\alpha + \beta)$

$$= (f_1 \cup f_2) ((c_1 \cup c_2) (\alpha_1 \cup \alpha_2) + (\beta_1 \cup \beta_2))$$
$$= f_1 (c_1 \alpha_1 + \beta_1) \cup f_2 (c_2 \alpha_2 + \beta_2)$$
$$= [c_1 f_1 (\alpha_1 + f_1 (\beta_1)] \cup [c_2 (f_2 (\alpha_2)) + f_2 (\beta_2)]$$
$$= (c_1 \cup c_2) [f_1(\alpha_1) \cup f_2(\alpha_2)] + [f_1(\beta_1) \cup f_2(\beta_2)]$$
$$= c\, f (\alpha) + f(\beta).$$

(for $c \in F$, $\in V$)

Then $f = f_1 \cup f_2$ is called the linear bifunctional on the bivector space $V = V_1 \cup V_2$. Give an example of a linear bifunctional.

19. Define dual bivector space $V^* = V_1^* \cup V_2^*$. (Hint $V^* = L^B(V_1, F) \cup L^B(V_2, F)$).

20. Prove the existence of a unique dual basis for given basis $B = \{\alpha_1, \ldots, \alpha_n\}$ where $\alpha_i = \alpha_i^1 \cup \alpha_i^2$ i = 1, 2, …, n for the dual bivector space $V^*$ of the bivector space V.

21. Prove the following are equivalent, where T is a linear bioperator on a finite dimensional bivector space $V = V_1 \cup V_2$. If $c_1, \ldots, c_k$ be the distinct characteristic bivalues of ($c_i = c_i^1 \cup c_i^2$ i = 1, 2, …, k) and let $W_i = (T - c_iI) = (T_1 - c_i^1 I^1) \cup (T_2 - c_i^2 I^2)$; $1 \le i \le k$.

1. T is bidiagonalizable
2. The characteristic bipolynomial of $T = T_1 \cup T_2$ is $f = f_1 \cup f_2$

$$= \{(x - c_1^1)^{d_1^1} \ldots (x - c_k^1)^{d_k^1}\} \cup \left\{\left(x - x_1^2\right)^{d_1^2} \ldots \left(x - c_k^2\right)^{d_k^2}\right\}$$

and $\dim W_i = \dim W_i^1 \cup \dim W_i^2 = d_i^1 \cup d_i^2$ ; i = 1, 2,…, k

3. $\dim W_1 + \ldots + \dim W_k = \dim V = \dim V_1 \cup \dim V_2$

$= (\dim W_1^1 \cup \dim W_1^2) + \ldots + (\dim W_k^1 \cup \dim W_k^2)$

$= (\dim W_1^1 + \ldots + \dim W_k^1) \cup (\dim W_1^2 + \ldots + \dim W_k^2)$.

22. Is $T : V \to V$ defined by the bimatrix.

$$\begin{bmatrix} -9 & 4 & 4 \\ -8 & 3 & 4 \\ -16 & 8 & 7 \end{bmatrix} \cup \begin{bmatrix} 6 & -3 & -2 \\ 4 & -1 & -2 \\ 10 & -5 & -3 \end{bmatrix}$$

bidiagonalizable over $\mathbb{C}$ ?

23. Let $V = \{V_1 \cup V_2\}$ where $V_1 = Q \times Q \times Q \times Q$ and $V_2$ = {set of all polynomials of degree less than 6}, be a bivector inner biproduct space over Q, under the inner biproduct $( / ) = ( / )_1 \cup ( / )_2$ where $( / )_1$ is the standard inner product $( / )_2 =$

$$\int_{-1}^{1} f(x)\, g(x)\, dx .$$

Let $S = S_1 \cup S_2$

$= \{(5\ 3\ -1,\ 2),\ (6\ 7\ 0\ 3),\ (-\ 4\ 2\ 0\ 0)\} \cup \{1, x^3 - 2, x^4 + 1\}$

find $S^{\perp} = S_1^{\perp} \cup S_2^{\perp}$ .

24. Prove $\{0\}^{\perp} \cup \{0\}^{\perp} = V_1^{\perp} \cup V_2^{\perp}$ for any bivector space $V = V_1 \cup V_2$.

25. Prove if $V = V_1 \cup V_2$ is a finite dimensional inner biproduct space and if a linear bifunctional on $V = V_1 \cup V_2$ then there exists a unique vector $\beta = \beta_1 \cup \beta_2 \in V(\beta_1 \in V_1, \beta_2 \in V_2)$ such that $f(\alpha) = (\alpha \mid \beta)$ for all $\alpha \in V$ i.e. $f_1(\alpha_1) = (\alpha_1 \mid \beta_1)$ and $f_2(\alpha_2) = (\alpha_2 \mid \beta_2)$.

26. Illustrate the result in problem (25) by an explicit example.

27. Prove for any linear bioperator $T = T_1 \cup T_2$ on a finite dimensional inner biproduct bivector space $V = V_1 \cup V_2$, there exists a unique linear bioperator $T^* = T_1^* \cup T_2^*$ on $V = V_1 \cup V_2$ such that $(T\alpha \mid \beta) = (\alpha \mid T^* \beta)$ for all $\alpha, \beta$ in $V = V_1 \cup V_2$ i.e. $(T_1 \alpha_1 \mid \beta_1) \cup (T_2 \alpha_2 \mid \beta_2) = (\alpha_1 \mid T_1^* \beta_1) \cup (\alpha_2 \mid T_2^* \beta_2)$ where $\alpha = \alpha_1 \cup \alpha_2$ and $\beta = \beta_1 \cup \beta_2$.

28. Prove the related result in case of biorthonormal basis for $V_1 \cup V_2 = V$.

29. Prove the biadjoint of a linear bioperator depends also on the inner biproduct on V.

30. Prove for any bioperators $U = U_1 \cup U_2$ and $T = T_1 \cup T_2$ .

a. $(T + U)^* = T^* + U^*$.
b. $(cT)^* = c\,T^*$ .
c. $(TU)^* = U^*\,T^*$.
d. $(T^*)^* = T$.

31. Prove the product of two self biadjoint bioperators is self biadjoint if and only if the two bioperators commute.

32. Prove if $V = V_1 \cup V_2$ is an inner biproduct space and $T = T_1 \cup T_2$ is a self biadjoint linear bioperator on $V = V_1 \cup V_2$, then each characteristic bivalue of $T = T_1 \cup T_2$ is real and characteristic bivectors of T associated with distinct characteristic bivalues are biorthogonal.

33. Prove in case of complex $n \times n$ bimatrix $A = A_1 \cup A_2$ there is a biunitary bimatrix U such that $U^{-1} A U$ upper triangular.

34. Prove analogous of spectral theorem for a binormal bioperator on a finite dimensional complex inner biproduct space $V = V_1 \cup V_2$ or a self biadjoint bioperator on a finite dimensional real inner biproduct space.

35. If $T = T_1 \cup T_2$ is a bidiagonalizable normal bioperator on a finite dimensional inner biproduct space and

$$T = \sum_{j=1}^{K} C_j E_j \quad \left(C_j = C_j^1 \cup C_j^2;\; E_j = E_j^1 \cup E_j^2\right)$$

$$\text{i.e. } T_1 \cup T_2 = \left(\sum_{j=1}^{K} C_j^1 E_j^1\right) \cup \left(\sum_{j=1}^{K} C_j^2 E_j^2\right)$$

be its spectral biresolution illustrate the spectral biresolution for a finite dimensional inner biproduct space.

36. Prove if U and T are binormal bioperators which commute then U + T and UT are binormal bioperators.

37. Derive the modified form of Cayley Hamilton theorem in case of linear bioperator $T = T_1 \cup T_2$ on a finite dimensional bivector space.

38. State and prove the modified form of the primary decomposition, theorem in case of a linear bioperator $T = T_1 \cup T_2$ on a finite dimensional bivector space $V = V_1 \cup V_2$. Define the notion of T-ad-missible subbispace of V.

40. Give the modified form of Cyclic Decomposition Theorem in case of a linear bioperator T on a finite dimensional bivector space $V = V_1 \cup V_2$.

41. State prove a theorem analogous of Generalized Cayley Hamilton Theorem for a linear bioperator T on a finite dimensional bivector space $V = V_1 \cup V_2$.

42. Obtain some interesting results on pseudo inner biproducts.

43. Obtain the modified form of spectral theorem using pseudo inner biproducts.

44. Find the Jordan canonical biform for the following bimatrices over $\mathbb{C}$.

a. $$A = A_1 \cup A_2 = \begin{pmatrix} 1 & 3 \\ 0 & 1 \end{pmatrix} \cup \begin{pmatrix} 1 & 1 \\ 0 & 1 \end{pmatrix}$$

b. $X = \begin{pmatrix} -1 & 0 & 0 & 0 \\ 0 & 1 & 0 & 0 \\ 0 & 0 & 1 & 0 \\ 0 & 0 & 0 & 2 \end{pmatrix} \cup \begin{pmatrix} i & 0 \\ 0 & -i \end{pmatrix}$

c. $P = \begin{pmatrix} 1 & 0 & 0 & 0 \\ 0 & -1 & 0 & 0 \\ 0 & 0 & 2 & 0 \\ 0 & 0 & 0 & 1 \end{pmatrix} \cup \begin{pmatrix} 0 & 1 \\ -1 & 0 \end{pmatrix}$

d. $D = \begin{pmatrix} -1 & 0 & 0 & 0 \\ 0 & 1 & 0 & 0 \\ 0 & 0 & 1 & 0 \\ 0 & 0 & 0 & 2 \end{pmatrix} \cup \begin{pmatrix} 1 & 0 & 0 & 0 \\ 0 & -1 & 0 & 0 \\ 0 & 0 & 2 & 0 \\ 0 & 0 & 0 & 1 \end{pmatrix}$

45. Find the Jordan canonical forms over the field of complex numbers for the bimatrix

$$A = \begin{pmatrix} 2 & -1 \\ 1 & -1 \end{pmatrix} \cup \begin{pmatrix} 0 & 0 & 1 \\ 1 & 0 & 0 \\ 0 & 1 & 0 \end{pmatrix}.$$

46. Calculate the characteristic and minimum bipolynomial for the transformation bimatrix;

$$M = \begin{pmatrix} -5 & -9 & 0 & 0 \\ 4 & 7 & 0 & 0 \\ 0 & 0 & -1 & -1 \\ 0 & 0 & 4 & 3 \end{pmatrix} \cup \begin{pmatrix} 3 & 1 & 0 & -1 \\ -1 & 2 & 1 & 1 \\ 0 & -1 & 1 & 0 \\ 0 & 1 & 1 & 2 \end{pmatrix}$$

Find the Jordan biform of M.

47. Find the Jordan biform of the bimatrix.

$$\begin{bmatrix} -12 & -7 & -4 & -2 & -1 \\ 16 & 8 & 4 & 1 & 0 \\ 52 & 36 & 23 & 15 & 9 \\ -48 & -28 & -16 & -9 & -4 \\ -64 & -48 & -32 & -20 & -13 \end{bmatrix} \cup \begin{bmatrix} -18 & -31 & 12 & -15 & -19 \\ 36 & 50 & -19 & 27 & 29 \\ -16 & -20 & 10 & -12 & -12 \\ 16 & 32 & -12 & 14 & 20 \\ -64 & -80 & 32 & -48 & -46 \end{bmatrix}.$$

48. Does bimatrix

$$\begin{pmatrix} 2 & 0 & 0 & 2 \\ 0 & 1 & -2 & 0 \\ 0 & 1 & -1 & 0 \\ -4 & 0 & 0 & 2 \end{pmatrix} \cup \begin{pmatrix} 3 & 1 & 0 & -1 \\ -1 & 2 & 1 & 1 \\ 0 & -1 & 1 & 0 \\ 0 & 1 & 1 & 2 \end{pmatrix}$$

have a Jordan biform over R? why?

49. For the (5, 3) $\cup$ (5, 2) bicode given by the generator bimatrix, $G = G_1 \cup G_2$ where

$$G = \begin{pmatrix} 0 & 1 & 1 & 1 & 1 \\ 1 & 0 & 0 & 1 & 0 \end{pmatrix} \cup \begin{pmatrix} 1 & 1 & 1 & 0 & 0 \\ 0 & 0 & 1 & 1 & 0 \\ 1 & 1 & 1 & 1 & 1 \end{pmatrix}$$

Find its dual bicode.

50. Find the parity check bimatrix and dual bicode for the bicodes given by the generator by matrix.

$$G = \begin{pmatrix} 1 & 1 & 1 & 0 \\ 0 & 1 & 1 & 0 \\ 0 & 0 & 1 & 1 \end{pmatrix} \cup \begin{pmatrix} 1 & 0 & 1 & 1 \\ 0 & 1 & 1 & 1 \\ 1 & 0 & 0 & 1 \end{pmatrix}.$$

51. For the parity check bimatrix, $H = H_1 \cup H_2$ where

$$H = \begin{pmatrix} 0 & 0 & 0 & 1 & 1 & 1 & 1 \\ 0 & 1 & 1 & 0 & 0 & 1 & 1 \\ 1 & 0 & 1 & 0 & 1 & 0 & 1 \end{pmatrix} \cup \begin{pmatrix} 1 & 0 & 0 & 1 & 1 & 0 & 1 \\ 0 & 1 & 0 & 1 & 0 & 1 & 1 \\ 0 & 0 & 1 & 0 & 1 & 1 & 1 \end{pmatrix}$$

Find the dual bicode.

52. Determine the bicode words with the generator bipolynomial $(1 + x^2 + x^3 + x^4) \cup (x^4 + x^2 + x + 1) = g_1 \cup g_2$. $g_1 / x^2 + 1$ and $g_2 / x^8 + 1$ Find their dual bicode.

53. Find the adjacency bimatrix of the following bigraphs :

a. $G = G_1 \cup G_2$.

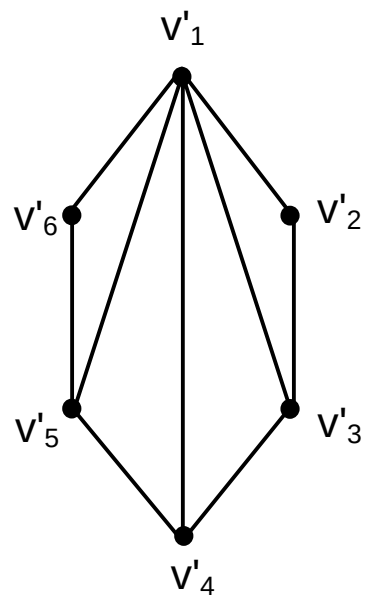

$G_1$

FIGURE: 5.1

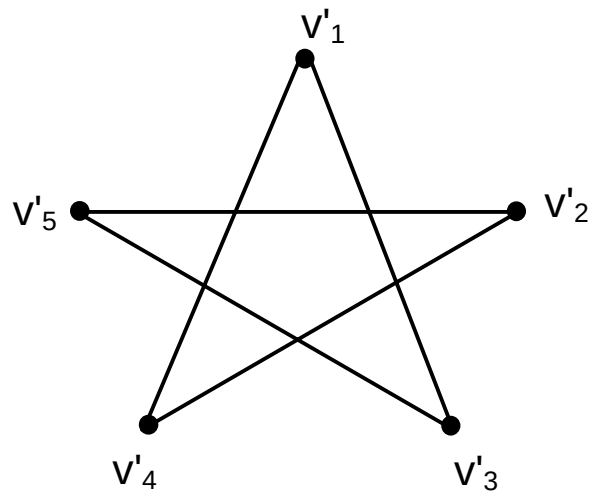

$G_2$

FIGURE: 5.2

b.

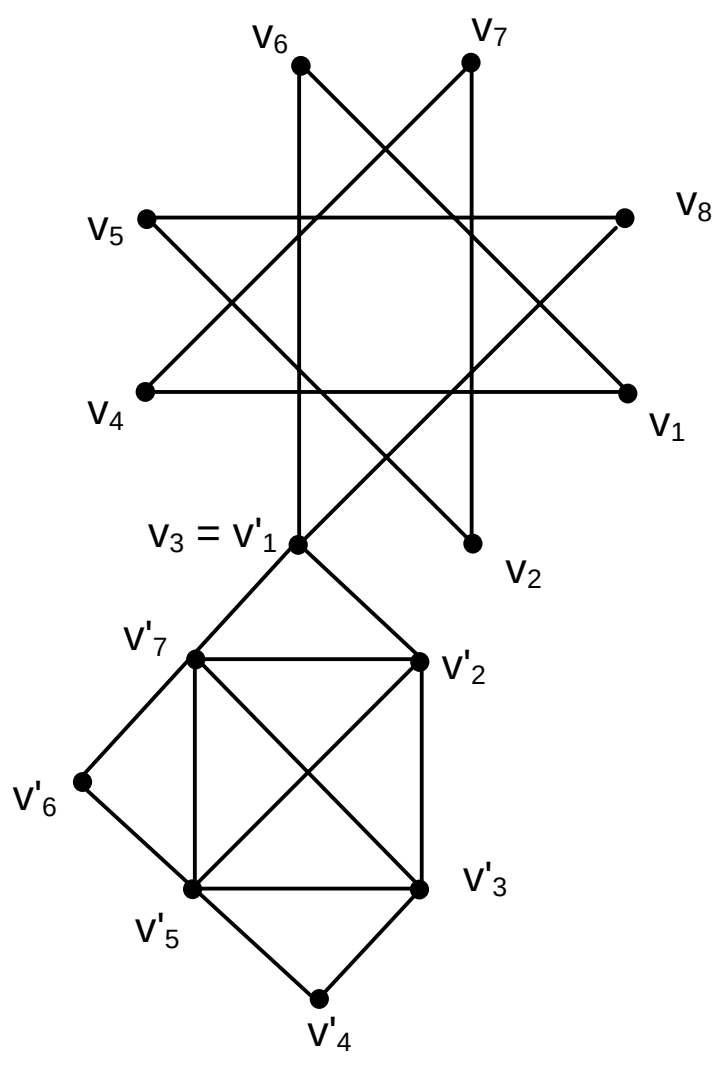

FIGURE: 5.3

c.

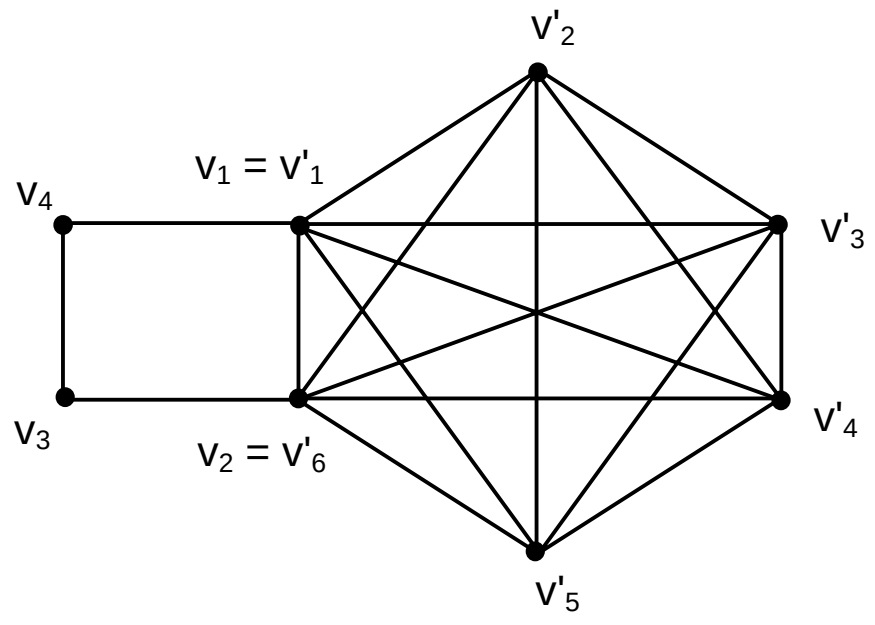

FIGURE: 5.4

54. For the bigraph given in problem 53 (iii) draw the cycle bimatrix.

55. Find the cycle bimatrix of the following bigraph:

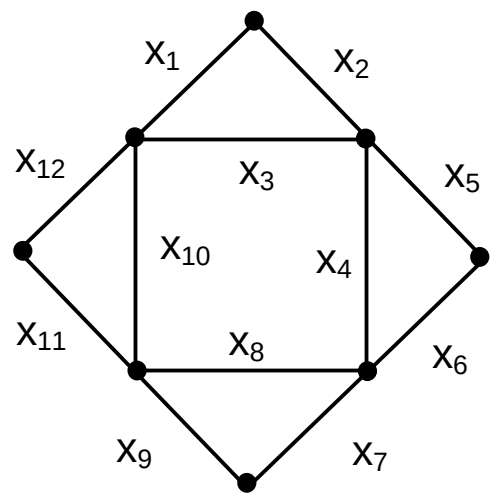

FIGURE: 5.5

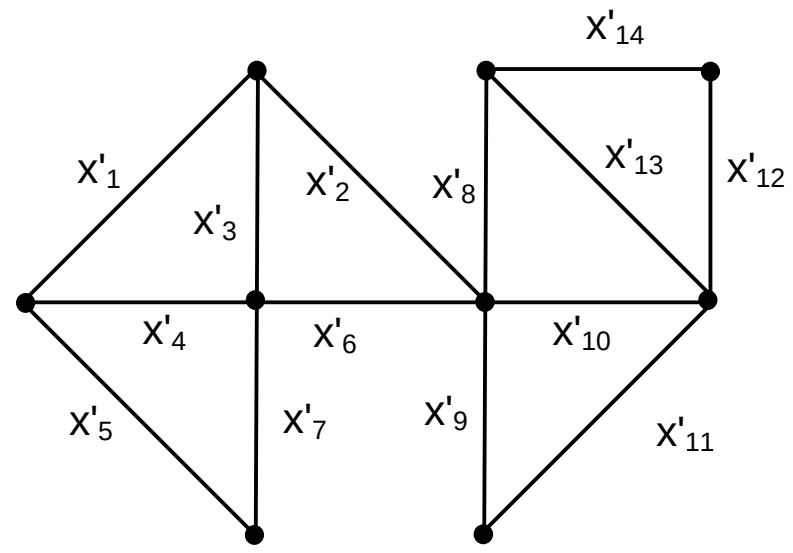

FIGURE: 5.6

56. Let C = C (7, 4) be a code generated by the parity check matrix.

$$H = \begin{pmatrix} 0 & 0 & 0 & 1 & 1 & 1 & 1 \\ 0 & 1 & 1 & 0 & 0 & 1 & 1 \\ 1 & 0 & 1 & 0 & 1 & 0 & 1 \end{pmatrix}$$

over $Z_2$. If the received words are (1 0 0 1 0 1 0) and (1 1 0 1 0 1 1) find the sent words using pseudo best approximation.

57. Let C be a binary (7, 3) generated by the generator matrix

$$G = \begin{pmatrix} 0 & 0 & 0 & 1 & 1 & 1 & 1 \\ 0 & 1 & 1 & 0 & 0 & 1 & 1 \\ 1 & 0 & 1 & 0 & 1 & 0 & 1 \end{pmatrix}$$

over $Z_2$, answer the following:

a. If (1 1 1 1 1 1 1) is the received word find the correct code word using method of pseudo best approximation.
b. Find the dual of C (7, 3). If (1 1 1 1 1 1 1) is the received code some sent code word in $C^{\perp}$ (7, 3) find the correct code using pseudo best approximation.
c. Is the received words different or same in C (7, 3) and $C^{\perp}$ (7, 3)?

58. Using pseudo best biapproximation find the sent bicode if the received bicode word is $x = x_1 \cup x_2 = (1\ 1\ 1\ 0\ 1\ 1\ 1) \cup (1\ 1\ 1\ 1\ 1\ 1\ 1\ 1)$.

The bicode is calculated using the parity check bimatrix $H = H_1 \cup H_2 =$

$$\begin{bmatrix} 0 & 0 & 1 & 0 & 1 & 1 & 1 \\ 0 & 1 & 0 & 1 & 1 & 1 & 0 \\ 1 & 0 & 1 & 1 & 1 & 0 & 0 \end{bmatrix} \cup \begin{bmatrix} 1 & 1 & 0 & 0 & 1 & 0 & 0 & 0 \\ 0 & 1 & 1 & 0 & 0 & 1 & 0 & 0 \\ 1 & 0 & 0 & 1 & 0 & 0 & 1 & 0 \\ 0 & 1 & 0 & 1 & 0 & 0 & 0 & 1 \end{bmatrix}.$$

59. Using the generator matrix

$$G = \begin{pmatrix} 1 & 1 & 1 & 0 & 0 \\ 0 & 0 & 1 & 1 & 0 \\ 1 & 1 & 1 & 1 & 1 \end{pmatrix}$$

and the parity check matrix

$$H = \begin{pmatrix} 0 & 0 & 0 & 1 & 1 & 1 & 1 \\ 0 & 1 & 1 & 0 & 0 & 1 & 1 \\ 1 & 0 & 1 & 0 & 1 & 0 & 1 \end{pmatrix}$$

the bicode $C = C_1 (5, 3) \cup C_2 (7, 4)$ is obtained. If the received bicode is $(1\ 0\ 1\ 0\ 1) \cup (0\ 1\ 0\ 1\ 0\ 1\ 0)$ find using the pseudo best biapproximations the sent message.

60. Give an example of a bisemilinear algebra / semi linear bialgebra (Both the terms are used synonymously).

61. Does there exist a bisemivector space which is not a bisemilinear algebra? Justify your answer with proof or counter example.

62 Obtain some interesting results about bisemi linear algebras which are not true in bisemivector spaces.

63. Let G = {The set of all 3 × 3 matrices with entries from Q (I)}, G is a group under '+'. Clearly G is a neutrosophic vector space over Q. Find the dimension of G over Q. Also find a set of basis for this neutrosophic vector space G.

64. Find the dimension and a basis for the neutrosophic vector space G = {3 × 3 matrices with entries from Q(I) } over Q (I).

65. Find a basis and the dimension of the strong neutrosophic vector space V = {set of all 7 × 3 matrices with entries from R(I) } over R (I). Find the basis and dimension of V when V is a strong neutrosophic vector space over Q (I).

66. Let V = {set of all 3 × 2 matrices with entries from Q (I)} and W = {Q (I) × Q (I) × Q (I)} be neutrosophic group under addition. V and W are strong neutrosophic vector spaces over Q (I). Find dimension of $Hom_{Q(I)} (V, W)$. Is $Hom_{Q (I)} (V, W)$ a strong neutrosophic vector spaces over Q (I) ? Justify your claim. Give a pair of basis for $Hom_{Q(I)}$ ($V_1$ W).

67. Find pseudo neutrosophic subspaces of V and W given in problem (66).

68. Let T : V → W, T (x, y, z) = (x + y, x – y, x + z, z + y) where V and W are strong neutrosophic vector spaces over Q (I). Give a basis for V and a basis for W.

69. Find KerT for the linear transformation given in problem (68).

70. Let Q (I) be the neutrosophic field over which a n dimensional neutrosophic vector space V is defined. Is $V \cong [Q (I)]^n$? Justify your claim!

71. Suppose V is a neutrosophic vector space over Q of dimension n. Is $V \cong Q^m$? Give examples to justify your claim.

72. If V is a neutrosophic vector space over a field F and if it is a strong neutrosophic vector space over the field F (I). If dim V as a strong neutrosophic vector space over F(I) is n then will dim V as a neutrosophic vector space over F be 2n?

73. Find a means to define linear operator of a neutrosophic vector space over any field F.

74. Suppose V = Q (I) × Q (I) × Q (I) be a neutrosophic group under component wise addition and V is a neutrosophic vector space over Q. Find a linear operator on V. What is the order of the matrix associated with it.

75. Define a inner product on V = {set of all polynomials of degree less than or equal to 10 with coefficients from Q (I)}, V is a neutrosophic group under addition. V is a neutrosophic vector space over Q, which is not a standard inner product.

76. Write a result analogous to Gram Schmidt orthogonalization process for a neutrosophic vector space.

77. Derive some interesting notions on orthogonal neutrosophic vectors.

78. Define adjoint for a linear operator on a neutrosophic vector space and strong neutrosophic vector space.

79. If $T : V \to V$ is a linear operator on a neutrosophic vector space.

Is $(T^*)^* = T$ ?
Is $(T + U)^* = T^* + U^*$ ?
Is $(T U)^* = U^* T^*$?
where $U : V \to V$ .

80. Can dual neutrosophic vector space be defined ? Why?

81. Define dual space for the strong neutrosophic vector spaces V over F (I). Is $V^* \cong V$ ?

82. Obtain interesting results about dual spaces in case of strong neutrosophic vector spaces.

83. Give an example of a neutrosophic bigroup which is quasi commutative.

84. Find proper neutrosophic subbigroups of the neutrosophic bigroup $G = G_1 \cup G_2$ where $G_1$ =

{Set of all (7 × 2) neutrosophic matrices with entries from Q (I)}, $G_1$ is a neutrosophic group under matrix addition, $G_2$ = {Q(I)(x), all neutrosophic polynomials with coefficients from Q (I)}; under polynomial addition $G = G_1 \cup G_2$ is a neutrosophic bigroup. Find at least three neutrosophic subbigroups of G?

85. Let $G = G_1 \cup G_2$ where $G_1$ = {R (I) × Q (I) × Q (I)} a neutrosophic group under addition and $G_2$ = {Set of all 4 × 4 matrices with entries from Q (I) which are non singular}. $G_2$ is a neutrosophic group under multiplication. Find neutrosophic subgroups of G.

86. Prove kernel of a linear bitransformation T is a pseudo neutrosophic bisubspace of a domain space, where $T : V \to W$ here V and W are strong neutrosophic bivector space over F (I).

87. Give an example of a quasi neutrosophic bigroup.

88. Is it always possible to construct a quasi neutrosophic subbigroup from the neutrosophic bigroup?

89. Does the neutrosophic bigroup $V = V_1 \cup V_2$ where $V_1$ = Q (I) × Q (I) × Q × R, a neutrosophic group under addition and $V_2$ = {set of all 7 × 2 matrices with entries from Q (I)} have a quasi neutrosophic subbigroup. If so find it.

90. For the neutrosophic bivector space $V = V_1 \cup V_2$ over Q where $V_1$ = {All polynomials with coefficients from Q (I)} under '+' and $V_2$ = {Set of all 2 × 2 matrices with entries from Q(I)}. Define a linear bioperator T with a non trivial bikernel.

91. Find a quasi neutrosophic subbivectorspace. Find a pseudo neutrosophic subbivector space V. Does V have a complete pseudo neutrosophic subbivector space? (V given in problem 90)

92. Find the complete pseudo neutrosophic subbivector space of $V = V_1 \cup V_2$ where $V_1$ = {The set of all $m \times n$ matrices with entries from Q (I)}, $V_1$ is a neutrosophic group under matrix addition and $V_2$ = {The set of all $m \times m$ matrices with entries from Q (I)}, $V_2$ also neutrosophic group under matrix addition. $V = V_1 \cup V_2$ is just the neutrosophic bivector space over Q.

93. Does $V = V_1 \cup V_2$ defined in problem 92 have any quasi neutrosophic subbivector space?

94. Is the set of all quasi linear bitransformation of two quasi neutrosophic bivector spaces V and W over F a quasi neutrosophic bivector space over F. Justify your answer.

95. Let $V = V_1 \cup V_2$ , $V_1 = ( Q \times Q \times Q \times Q)$ group under component wise addition and $V_2$ = {The set of all $2 \times 3$ matrices with entries from Q (I)}. $V = V_1 \cup V_2$ is clearly a quasi neutrosophic bivector space over Q. Let $W = W_1 \cup W_2$ where $W_1$ = {Set of all $2 \times 2$ matrices with entries from Q}, under the operation matrix addition is a group, $W_2$ = {all polynomials of degree less than or equal to 3 with coefficient from Q (I)} is a group under polynomial addition. W is a quasi neutrosophic bivector space over Q. Define $T : V \to W$ a quasi neutrosophic linear bitransformation;

   a. such that $T = B \text{ Ker } T_1 \cup B \text{ Ker } T_2 = 0 \cup 0$.
   b. such that $B \text{ Ker } T \neq (0) \cup (0)$

96. Let $V = V_1 \cup V_2$ be a (m, n) dimension quasi neutrosophic bivector space over Q, where $V_1 = \{Q \times Q \times Q \times Q\}$, a group under component wise addition and $V_2$ = {The set of all $3 \times 2$ neutrosophic matrices with entries from Q (I) }; $V_2$ is a neutrosophic group under matrix addition:

   a. What is dimension of V ?
   b. Find a quasi linear bioperator T on V such that B Ker T = {0}.

97. Let $V = V_1 \cup V_2$ be a quasi neutrosophic bivector space over R, where $V_1 = \{R \times R \times R\}$ and $V_2$ = {set of all $3 \times 2$ neutrosophic matrices with entries from R (I)}.

   a. Find dimension of $V = V_1 \cup V_2$ .
   b. Find quasi linear bioperator $T = T_1 \cup T_2$ on $V = V_1 \cup V_2$ such that B Ker T = $\{0\} \cup \{0\}$.
   c. Find a quasi linear bioperator $T = T_1 \cup T_2$ on $V = V_1 \cup V_2$ such that B Ker T $\neq \{0\} \cup \{0\}$.

98. Does there exists a neutrosophic ring which is only Smarandache quasi neutrosophic ring and not a Smarandache neutrosophic ring?

99. Does there exist a neutrosophic ring which is only Smarandache neutrosophic ring and not a Smarandache quasi neutrosophic ring? Justify your claim!

100. Does there exists a neutrosophic semigroup which is a Smarandache neutrosophic semigroup and not a Smarandache quasi neutrosophic semigroup.

101. Does there exists a neutrosophic semigroup which is a Smarandache quasi neutrosophic semigroup and not a Smarandache neutrosophic semigroup. Justify your answer.

102. Give an example of a neutrosophic semigroup which is not a Smarandache neutrosophic semigroup.

103. Is the set of all $3 \times 3$ matrices with entries from Z(I) a Smarandache neutrosophic semigroup?

104. Is S = {Z (I) × Z (I) × Z (I), ×} a Smarandache quasi neutrosophic semigroup under component wise multiplication?

105. Is S = Z × Z (I) a Smarandache quasi neutrosophic semigroup under component wise multiplication?

106. Give an example of a Smarandache neutrosophic vector space.

107. Is V = R(I) × R a Smarandache neutrosophic vector space over R × R?

108. Is V = Q (I) × Q (I) × Q (I) a Smarandache neutrosophic vector space over the S-ring, S = Q × Q × Q?

109. Give an example of a quasi neutrosophic bivector space.

110. Illustrate by an example a Smarandache neutrosophic bisemigroup.

111. Is. $S = S_1 \cup S_2$ where $S_1$ = R × Q (I) × Q (I) [x] and $S_2$ = {R (I) } × Q (I) a Smarandache

neutrosophic bisemigroup under the operation component wise multiplication.

a. Illustrate by an example a Smarandache neutrosophic bisemigroup.
b. Strong Smarandache neutrosophic bivector space of finite dimension.

112. Give an example of a pseudo Smarandache neutrosophic bisemigroup.

113. Does there exist any relation between Smarandache neutrosophic bisemigroup and a pseudo Smarandache neutrosophic bisemigroup?

114. Give an industry real model and apply it in Smarandache neutrosophic Leontief bimodel using Smarandache neutrosophic bimatrices.

# BIBLIOGRAPHY

It is worth mentioning here that we are only citing the texts that apply directly to linear algebra, and the books which have been referred for the purpose of writing this book. To supply a complete bibliography on linear algebra is not only inappropriate owing to the diversity of handling, but also a complex task in itself, for, the subject has books pertaining from the flippant undergraduate level to serious research. We have limited ourselves to only listing those research-level books on linear algebra which ingrain an original approach in them. Longer references/ bibliographies and lists of suggested reading can be found in many of the reference works listed here.

# INDEX

## A

## B

### C

### D

### E

### F

### G

### H

### I

### J

### K

## L

## M

## N

**O**

## P

## Q

## R

## S

**T**

**U**

**V**

**W**

# About the Authors

**Dr.W.B.Vasantha Kandasamy** is an Associate Professor in the Department of Mathematics, Indian Institute of Technology Madras, Chennai, where she lives with her husband Dr.K.Kandasamy and daughters Meena and Kama. Her current interests include Smarandache algebraic structures, fuzzy theory, coding/ communication theory. In the past decade she has guided nine Ph.D. scholars in the different fields of non-associative algebras, algebraic coding theory, transportation theory, fuzzy groups, and applications of fuzzy theory of the problems faced in chemical industries and cement industries. Currently, six Ph.D. scholars are working under her guidance. She has to her credit 287 research papers of which 209 are individually authored. Apart from this, she and her students have presented around 329 papers in national and international conferences. She teaches both undergraduate and post-graduate students and has guided over 45 M.Sc. and M.Tech. projects. She has worked in collaboration projects with the Indian Space Research Organization and with the Tamil Nadu State AIDS Control Society. This is her 19th book.

She can be contacted at vasantha@iitm.ac.in
You can visit her work on the web at: http://mat.iitm.ac.in/~wbv

---

**Dr.Florentin Smarandache** is an Associate Professor of Mathematics at the University of New Mexico, Gallup Campus, USA. He published over 60 books and 80 papers and notes in mathematics, philosophy, literature, rebus. In mathematics his research papers are in number theory, non-Euclidean geometry, synthetic geometry, algebraic structures, statistics, and multiple valued logic (fuzzy logic and fuzzy set, neutrosophic logic and neutrosophic set, neutrosophic probability). He contributed with proposed problems and solutions to the Student Mathematical Competitions. His latest interest is in information fusion were he works with Dr.Jean Dezert from ONERA (French National Establishment for Aerospace Research in Paris) in creasing a new theory of plausible and paradoxical reasoning (DSmT).

He can be contacted at smarand@unm.edu

---

**K. Ilanthenral** is the editor of The Maths Tiger, Quarterly Journal of Maths. She can be contacted at ilanthenral@gmail.com